\DocumentMetadata{lang=en-US,testphase=phase-I} % Tag the PDF links for TeXShop's preview.
 \documentclass[11pt, letterpaper]{amsart}
\calclayout
\usepackage[thinlines]{easytable}
\usepackage{graphicx}
\usepackage{amsmath,amssymb,amsfonts,amsthm,mathtools,mathrsfs}
\usepackage{enumitem}
\usepackage[margin=1in]{geometry}
\usepackage{microtype}

\usepackage[colorlinks=true]{hyperref}

\numberwithin{equation}{section}
\allowdisplaybreaks
\newtheorem{theorem}{Theorem}[section]
\newtheorem{proposition}[theorem]{Proposition}
\newtheorem{lemma}[theorem]{Lemma}
\newtheorem{corollary}[theorem]{Corollary}
\theoremstyle{definition}
\newtheorem{definition}[theorem]{Definition}

\theoremstyle{remark}
\newtheorem{remark}[theorem]{Remark}
\newtheorem*{remark*}{Remark}

\newcommand{\SL}{\operatorname{SL}}
\newcommand{\GL}{\operatorname{GL}}
\newcommand{\PGL}{\operatorname{PGL}}

\newcommand{\SO}{\operatorname{SO}}
\renewcommand{\O}{\operatorname{O}}
\newcommand{\CT}{\operatorname{CT}}

\newcommand{\Lie}{\operatorname{Lie}}

\newcommand{\St}{\operatorname{St}}
\newcommand{\A}{\mathbb A}
\newcommand{\C}{\mathbb C}
\newcommand{\Q}{\mathbb Q}
\newcommand{\R}{\mathbb R}
\newcommand{\Z}{\mathbb Z}
\newcommand{\D}{\mathcal D}

\title[Vanishing near the virtual cohomological dimension]
{Vanishing near the virtual cohomological dimension of $\SL_n(\Z)$}

\author{Taiwang Deng}
\address{Beijing Institute of Mathematical Sciences and Applications (BIMSA), Huairou District, 100084, Beijing, China}
\email{dengtaiw@bimsa.cn}

\date{}

\begin{document}

\begin{abstract}
We prove the Church--Farb--Putman conjecture that
\[
 H^{N_n-i}(\SL_n(\Z);\Q)=0,
 \qquad N_n=\binom n2,\quad n\ge2,\quad 0\le i<n-1.
\]
Franke's comparison theorem reduces the assertion to a
complementary-degree bound for automorphic cohomology.
The trivial cuspidal support is treated using Franke's
complex of compact duals.  For nontrivial cuspidal support, we group the
inducing factors into blocks of an intermediate Levi subgroup and use
constant terms to reduce the vanishing bound to degree estimates on its
blocks.  A weighted count of Harish--Chandra parameters, together with
the rank-dependent motivic-weight bounds of Chenevier--Lannes, gives the
required estimate except possibly for mixed blocks with one cuspidal factor
of rank two or four, or two cuspidal factors of rank two.  In these cases,
we bound the cohomological degree using cuts between consecutive inducing
factors and the corresponding cuspidal constant terms.
\end{abstract}

\maketitle

\tableofcontents

\section{Introduction}\label{sec:introduction}

\subsection{Background}

Let
\[
 \Gamma_n=\SL_n(\Z),\qquad N_n=\binom n2.
\]
Borel--Serre \cite{BorelSerre1973} prove that $N_n$ is the virtual
cohomological dimension of $\Gamma_n$.  Its rational cohomology
vanishes above $N_n$.  We ask how far this vanishing extends below
$N_n$.

The spherical Tits building $\mathcal B_n$ of $\Q^n$ is the simplicial
complex of flags of nonzero proper $\Q$-subspaces of $\Q^n$
\cite[Section~2]{CFP2014}.  Its top reduced homology is the rational
Steinberg module
\[
 \St_n=\widetilde H_{n-2}(\mathcal B_n;\Q).
\]
Borel--Serre virtual duality
\cite[Theorem~11.4.2]{BorelSerre1973} gives
\[
 H^{N_n-i}(\Gamma_n;\Q)
 \cong H_i(\Gamma_n;\St_n).
\]
Church, Farb, and Putman \cite[Conjecture~2]{CFP2014} conjectured
that these groups vanish for $0\le i<n-1$.  

For the automorphic argument we express this question in terms of
compact supports.  Set
\[
 X_n=\SL_n(\R)/\SO_n,\qquad Y_n=\Gamma_n\backslash X_n,\qquad
 d_n=\dim X_n=N_n+n-1.
\]
Poincar\'e--Lefschetz duality identifies the
dual of $H^{N_n-i}(\Gamma_n;\C)$ with compactly supported cohomology
in degree $n-1+i$.  The conjecture is therefore equivalent to
\[
 H_c^j(Y_n;\C)=0\qquad(j<2n-2).
\]
Here the part with $j<n-1$ corresponds to the already known
vanishing above $N_n$.  

\subsection{Main result and previous work}

\begin{theorem}\label{thm:main}
For every $n\ge2$ and every $0\le i<n-1$,
\begin{equation}\label{eq:main}
 H^{N_n-i}(\Gamma_n; \Q)=0.
\end{equation}
\end{theorem}

Theorem~\ref{prop:general-mixed-support-vanishing} gives
\[
 H^q(Y_n;\C)_{\mathfrak c}=0\qquad(q>d_n-2n+2)
\]
for every full-level cuspidal support $\mathfrak c$ compatible with the
infinitesimal character of the trivial representation on $\mathfrak{sl}_n$.
Its proof combines a weighted degree estimate with a calculation of
cuspidal constant terms for groups containing one or two cuspidal factors.

Lee--Szczarba \cite[Theorem~4.1]{LeeSzczarba1976} proved the top degree in the
conjecture, Church--Putman \cite{ChurchPutman2017} proved the
codimension-one case, and Br\"uck--Miller--Patzt--Sroka--Wilson
\cite[Theorem~B]{BMPSW2024} proved the codimension-two case.  These arguments
use the Steinberg-module formulation of Borel--Serre duality.  
Our approach uses Franke's automorphic description of cohomology
\cite{Franke1998}.  The vanishing for trivial cuspidal support
follows from Franke's explicit computation
\cite[Theorems~7.10--7.11]{Franke2008}. We give the degree estimate
in Theorem~\ref{prop:constant}.
For nontrivial cuspidal support,
Theorem~\ref{prop:general-mixed-support-vanishing}
gives the equivalent ordinary-degree bound; the passage to group cohomology is
spelled out in Corollary~\ref{cor:main-theorem}.

Grbac and Grobner
\cite[Theorem~A]{GrbacGrobner2025} construct
degenerate Eisenstein classes outside the image of the Borel map in degree
$n-1$ for certain even $n$.  In particular, they show that $H^3(\SL_4(\Z);\C)\ne0$, which implies that our range of vanishing $i<n-1$ is optimal as a
uniform bound in $n$.

Our results should be compared with those of Li and Sun
\cite[Theorem~1.8 and Section~5.4]{LiSun2019}, who prove that the Borel map
\[
 H^q(\mathfrak{sl}_n(\R),\SO_n;\C)
 \longrightarrow H^q(\SL_n(\Z);\C)
\]
is an isomorphism for $0\le q\le n-2$.  This stable range is sharp
as a uniform bound in $n$.

\subsection{The automorphic framework}

Franke's comparison theorem computes ordinary arithmetic cohomology
from a relative Lie algebra complex of automorphic forms
\cite[Theorem~18]{Franke1998}.  For $\Gamma_n=\SL_n(\Z)$, Shapiro's lemma
\eqref{eq:shapiro-decomposition} expresses $H^q(\Gamma_n;\C)$ as the
direct sum of the two groups $H^q(\GL_n(\Z);E_\varepsilon)$ for the
trivial representation $E_1$ and the determinant representation $E_{\det}$.
Both have the infinitesimal character of the trivial representation on
$\mathfrak{sl}_n$.  For fixed cuspidal support and growth bound, the
associated graded pieces of Franke's filtration are described by Laurent
coefficients of Eisenstein series \cite[Theorem~14 and Corollary~1]{Franke1998}.
By Corollary~\ref{cor:franke-degree-reduction}, if a cuspidal-support summand has
nonzero cohomology in degree $q$ for $E_\varepsilon$, then so does
$M(t)\otimes E_\varepsilon$ for some inducing datum $t$, where the module
$M(t)$ (cf. \eqref{eq:Mt} ) contributes to one graded piece of the
Franke filtration (Theorem \ref{thm:franke-graded}).

By M\oe glin--Waldspurger's classification of the discrete spectrum
of $\GL_n$ \cite{MoeglinWaldspurger1989}, the relevant square-integrable Levi representations
are products of Speh representations $J(\ell,\pi)$, where $\pi$ is
cuspidal and $\ell\ge1$.  Up to a common translation, the
Harish--Chandra parameter of $J(\ell,\pi)$ consists of a string of
$\ell$ consecutive entries for each entry of the parameter of $\pi$.
After including the inducing twists, we normalize the full parameter
to have coordinate sum zero and shift it by $(n-1)/2$, as in
Subsection~\ref{subsec:speh-weights}.  For the $\nu$-th Speh
factor, let $B_\nu$ be the set of entries it contributes to this
normalized and shifted parameter.  The infinitesimal-character condition
requires the sets $B_\nu$ to be pairwise disjoint with union
$U_n=\{0,1,\ldots,n-1\}$.

The same tiling enters the formula for the complementary degree $j=d_n-q$
in Proposition~\ref{prop:degree}.  For the level-one cuspidal factors of
rank at least two, we use the bounds of Chenevier--Lannes on the motivic
weight recalled in Lemma~\ref{lem:w11}.  Their use in restricting the
cuspidal supports of $H^*(\SL_n(\Z);\C)$ already appears in
Grbac--Grobner \cite[Theorem~2.2]{GrbacGrobner2025}.

\subsection{The degree bound for nontrivial support}

We proceed in three steps.
\begin{enumerate}[label=(\arabic*),leftmargin=2.2em]
\item \emph{Separate the support into components.}
At the beginning of Section~\ref{sec:numerical},
we place the character weights and the centres of the cuspidal
factors on $\R$. If the $\nu$-th Speh factor is $J(\ell,\pi)|\det|^z$, its
cuspidal support consists of the twists $\pi|\det|^{z-(\ell-1)/2+j}$ ($0\le j<\ell$),
whose shifted weight sets have union $B_\nu$ in \eqref{eq:packet}.
Starting from the tiling \eqref{eqn: tiling-inf-char-7.1}, we join
character weights differing by $1$, and we
join centres differing by $1$ when the underlying cuspidal
representations are isomorphic.  Let $K$ range over the connected
components of the union of these points and intervals.  Each pair
of twists of the same cuspidal representation whose exponents differ
by $1$ has both centres in the same component.  List the connected components
from right to left and let $n_1,\ldots,n_s$ be their total ranks.  These
ranks determine the standard parabolic $Q=P_{n_1,\ldots,n_s}$.
Propositions~\ref{lem:block-coefficient-map} and
\ref{lem:separated-cuspidal-parameters} therefore apply with the identity
permutation.  At full level they give
\[
 \mathcal A_{\mathfrak c}^{K_f}\simeq\mathcal I_Q(0)^{K_f}\otimes\D_y,
\]
where $\mathcal A_{\mathfrak c}^{K_f}$ is the full-level summand for the
cuspidal support $\mathfrak c$ (Subsection~\ref{subsec:topological-trivial-datum})
and $\mathcal I_Q(0)$ is normalized induction from the modules on
the Levi blocks of $Q$ with their determinant twists. The relative spectral variables
$y=(y_1,\ldots,y_s)$ vary the twist on the $b$-th block by $|\det|^{y_b}$,
subject to $\sum_{b=1}^s n_b y_b=0$. The space $\D_y$ consists of finite
linear combinations of derivative evaluations at the origin of this hyperplane.

\item \emph{Estimate each component.}
Components consisting only of characters are treated by
Theorem~\ref{prop:constant} and duality~\eqref{eq:duality} for both determinant parities.
For every other component $K$, let $n_K$ be its total rank.
If the corresponding $\GL_{n_K}$-block contributes in cohomological degree
$q_K$, Lemma~\ref{lem:general-component-degree} gives $d_{n_K}-q_K\ge J_K^0$.
If the cuspidal representations with centres in $K$ have ranks $k_\nu$
and centres $c_\nu$, set
$H_K=\sum_{\nu:\,c_\nu\in K}\lfloor k_\nu/2\rfloor$.
For components without characters or with $H_K\ge4$, this bound,
Lemma~\ref{lem:weighted-cuspidal-depth}, and the weight bounds of
Lemma~\ref{lem:w11} suffice in the combined degree estimate below.
For a mixed component with $H_K=3$, the proof of
Proposition~\ref{lem:small-rank-two-component} shows that the bound $J_K^0$
also suffices in the combined degree estimate of Proposition~\ref{lem:global-degree-from-components}.
The remaining mixed components have cuspidal rank patterns
\[
 (2),\qquad (4),\qquad (2,2).
\]
For $(2)$ and $(4)$ with $m$ character factors,
Lemma~\ref{lem:component-depth-at-most-two} applies Theorem~\ref{lem:one-copy-internal-bound}
with $h=1$ and $h=2$, giving the respective bounds $3m-2$ and $4m+1$.
For $(2,2)$ with $m$ character factors, the same lemma uses the bound $4m-4$
of Theorem~\ref{prop:two-cuspidal-bound} when $\Phi_K<0$, with $\Phi_K$
as in Proposition~\ref{lem:global-degree-from-components}.

\item \emph{Combine the bounds.}
For the $b$-th component, take $J_b=2n_b-2$ if it consists only of characters.
For a component $K$ containing cuspidal factors, the bounds in (2)
give $J_b=J_K$ with $J_K\ge J_K^0-\Phi_K$.
For a contributing Kostant summand indexed by $v$, equation
\eqref{eq:separated-support-degree} in the proof of
Proposition~\ref{lem:global-degree-from-components} gives
\[
 d_n-q\ge\sum_{b=1}^s J_b+
       \dim\mathfrak n_Q-\ell_W(v)+(s-1).
\]
Here $\mathfrak n_Q$ is the nilradical of $Q$ and $\ell_W(v)$
is the Weyl-group length of $v$.
The term $\dim\mathfrak n_Q-\ell_W(v)$ accounts for roots
joining distinct Levi blocks of $Q$, and $s-1$ is the dimension of the relative
split centre.  For these choices of $J_K$,
Proposition~\ref{lem:global-degree-from-components} gives
\[
 d_n-q\ge 2n-2+\sum_K(\Phi_K+J_K-J_K^0)\ge 2n-2,
\]
where $K$ ranges over the components containing cuspidal factors.
\end{enumerate}

\subsection{Organization of the paper}

Section~\ref{sec:automorphic-setup} fixes the automorphic and
cohomological conventions and treats trivial cuspidal support.
Section~\ref{sec:packets} proves the degree formula.
Section~\ref{sec:eisenstein-coefficients} recalls the constant-term
and residue formulas and separates the cuspidal parameters.
Section~\ref{sec:numerical} combines the degree bounds and proves
Theorem~\ref{thm:main}, assuming Theorems~\ref{lem:one-copy-internal-bound}
and~\ref{prop:two-cuspidal-bound} for the remaining small groups.
Section~\ref{sec:residue-quotients} treats the coefficient modules
for a consecutive interval of characters and for two copies of a
cuspidal representation. Section~\ref{sec:internal-cohomology} uses
cuts in these modules to prove the two bounds.
Subsection~\ref{subsec:two-cuspidal-example} illustrates the proof of
Theorem~\ref{prop:two-cuspidal-bound} on $\GL_{14}$ with ten character
factors and the two cuspidal representations $\Delta_{17}$ and $\Delta_{15}$.

\subsection{Scope and integral questions}

The theorem concerns full level and trivial rational coefficients.
Full level is used both in the spherical intertwiner calculation and
in \eqref{eq:rank-weight-lower-bound}.  Trivial coefficients give the
consecutive weight set $U_n$. The compact-dual calculation in
Theorem~\ref{prop:constant} treats the full trivial cuspidal-support summand.
The top-degree
rational cohomology of principal congruence subgroups of prime level
has nonzero finite-field Steinberg quotients
\cite[Section~1]{MillerPatztPutman2021}.
For an irreducible algebraic coefficient $F$ nontrivial on $\SL_n$,
the shifted Harish--Chandra parameter of $F^\vee$ has gaps.
The estimates that use the consecutive weights in $U_n$ would therefore
have to be reconsidered.
For $\SL_4(\Z)$ and $\GL_4(\Z)$, Bajpai--Deng
\cite{BajpaiDeng2026} compute homological Euler characteristics with
irreducible rational algebraic coefficients and derive vanishing
results and lower bounds for the dimensions of cohomology groups.

Let
\[
 \St_n^{\Z}=\widetilde H_{n-2}(\mathcal B_n;\Z).
\]
The integral strengthening
\[
 H_i(\SL_n(\Z);\St_n^{\Z})=0\qquad(0\le i<n-1)
\]
does not follow from our rational argument, which does not control
torsion in integral Steinberg homology.
For integral results and conjectures, see
\cite[Sections~1.6--1.9]{MillerPatztPutman2025}.

Integral Steinberg homology also enters the rank filtration in
algebraic $K$-theory.  Rognes's component-filtration spectral sequence
computes the homology of the rank-filtration layers from group homology
with coefficients built from Steinberg modules
\cite[Theorems~7.8--7.9]{Rognes2000}.
For related iterated bar constructions, see
\cite[Sections~1.1--1.3]{MillerPatztWilson2025}.
The rank spectral sequence converges to spectrum homology
$H_*(K(\Z);\Z)$ \cite[(8.3)]{Rognes2000}.
Extracting $K$-groups by this approach requires resolving differentials
and extensions and further input relating spectrum homology to homotopy.

\subsection*{Acknowledgements}
\addtocontents{toc}{\protect\setcounter{tocdepth}{1}}
While a postdoctoral fellow at the Max Planck Institute for
Mathematics, the author learned of the Church--Farb--Putman
conjecture from Professor G. Harder, when he was preparing a chapter on
Franke's theorem for his book \cite{Harder2022}.
The author is supported by the National Natural Science Foundation
of China (Grant No.~12401013) and the Beijing Natural Science
Foundation (Grant No.~1244042).

The author used ChatGPT (OpenAI) for 
editorial assistance in improving the language, exposition, and notational inconsistencies of the manuscript. All mathematical content and references were verified by the author, who take full responsibility for the manuscript.

\section{Cohomological and automorphic setup}\label{sec:automorphic-setup}

\subsection{Duality and Shapiro's lemma}
We first relate the cohomology groups of $\SL_n$ and $\GL_n$
by Shapiro's lemma.

Let
\[
 N_n=\binom n2,\qquad
 d_n=\dim\bigl(\SL_n(\R)/\SO_n\bigr)=N_n+n-1.
\]
Let
\[
 Y_n:=\SL_n(\Z)\backslash\SL_n(\R)/\SO_n.
\]

Choose a torsion-free normal subgroup
$\Gamma'\triangleleft\SL_n(\Z)$ of finite index and put
\[
 M'=\Gamma'\backslash\SL_n(\mathbb R)/\SO_n,\qquad
 F=\SL_n(\Z)/\Gamma'.
\]
Normality of $\Gamma'$ makes Poincar\'e--Lefschetz duality on
$M'$ equivariant for the action of $F$:
\[
                     H^q(M';\mathbb C)^\vee
                        \cong H_c^{d_n-q}(M';\mathbb C).
\]
Over $\mathbb C$, taking invariants under the finite group $F$ yields
\begin{equation}\label{eq:duality}
                 H^q(Y_n;\mathbb C)^\vee
                    \cong H_c^{d_n-q}(Y_n;\mathbb C).
\end{equation}
We use $\GL_n$ for the automorphic calculation.  Put
\[
 \widetilde Y_n=
 \GL_n(\Z)\backslash\GL_n(\R)/(\R_{>0}\O_n).
\]
Let ${\det}_r$ denote the one-dimensional algebraic representation
of $\GL_r$ given by $g\mapsto{\det}g$ and omit the subscript when
the rank is clear.
The normal inclusion $\SL_n(\Z)\subset\GL_n(\Z)$ has quotient of
order two.  Let $\mathscr L_1$ and $\mathscr L_{\det}$ be the
rank-one local systems defined by the trivial and determinant
characters of $\GL_n(\Z)/\SL_n(\Z)$.  Then
$\operatorname{Ind}_{\SL_n(\Z)}^{\GL_n(\Z)}\mathbf1
\cong\mathbf1\oplus{\det}$.  Shapiro's lemma
\cite[Chapter~IX, Section~2]{BorelWallach2000} therefore gives
\begin{equation}\label{eq:shapiro-decomposition}
 H^*(Y_n;\C)
 \cong H^*(\widetilde Y_n;\mathscr L_1)
 \oplus H^*(\widetilde Y_n;\mathscr L_{\det}),
\end{equation}
and the same formula holds with $H^*$ replaced by $H_c^*$.
The automorphic cochain constructions preserve the summands indexed
by $\varepsilon=1$ and $\varepsilon={\det}$, which we treat
separately.

For $0\le i<n-1$ and $q=N_n-i$ we have $q-(d_n-2n+2)=n-1-i>0$.
It therefore suffices to prove, for both determinant parities,
\begin{equation}\label{eq:target}
 H^q(\widetilde Y_n;\mathscr L_\varepsilon)=0
 \qquad(q>d_n-2n+2,\ \varepsilon\in\{1,{\det}\}).
\end{equation}
We use duality when an input is stated in compactly supported
cohomology.  It may interchange the determinant parities, so those
inputs will also be used for both parities.

\subsection{Global notation and conventions}
\label{subsec:global-notation}
Let
\[
 G_n=\GL_n,\qquad K_{n,\infty}=\O_n,\qquad \A=\A_\Q,
\]
where $\A_\Q$ is the ring of adeles of $\Q$.  When the rank is
clear, write $G=\GL_n$, $K_\infty=\O_n$,
$\mathfrak g=\Lie(G(\R))$, and
$\mathfrak k=\Lie(K_\infty)$.
Let $\mathfrak g_\C=\mathfrak g\otimes_\R\C$ be the
complexified Lie algebra.  We
denote by $U(\mathfrak g_\C)$ its universal enveloping algebra.

Let $P_0$ be the
Borel subgroup of upper-triangular matrices.  Unless stated otherwise,
every parabolic subgroup is standard, i.e., it contains $P_0$.  If
$P=M_PN_P$ is a parabolic of $G$, then $M_P$ is its Levi factor
and $N_P$ its unipotent radical.  We denote the normalized
constant-term map along $P$ by $\CT_P$.
Throughout, for representations $\sigma_i$ of $\GL_{n_i}$
over $\A$ or $\Q_v$, we use
\[
 \sigma_1\times\cdots\times\sigma_r
 \quad\text{or}\quad \mathop{\times}_{i=1}^r\sigma_i
\]
to denote the normalized parabolic induction of
$\sigma_1\boxtimes\cdots\boxtimes\sigma_r$ to $\GL_{n_1+\cdots+n_r}$.
For the nonarchimedean normalization and product notation, see
\cite[Sections~1.8 and 4.1]{BernsteinZelevinsky1977}; for the
archimedean and adelic conventions, see
\cite[Sections~I.1 and II.1]{MoeglinWaldspurger1989}.
For an inducing section $f$ and a spectral parameter $\lambda$,
\[
 E_P^G(f,\lambda)
\]
denotes the Eisenstein series formed using normalized parabolic
induction.  Let $\rho_P$ be
the half-sum
of the roots in $N_P$, following
\cite[Sections~3.3 and 6.1--6.2]{Shahidi2010}.

For integers $p\le q$, the notation $[p,q]$ means the integer set
$\{p,p+1,\ldots,q\}$.  We write $U_r=[0,r-1]$ for $r\ge1$.
The closed real interval with endpoints
$x\le y$ is written $[x,y]_{\R}$.  For an automorphic
representation $\pi$, $\pi^\vee$ is its contragredient; for a
finite-dimensional vector space or finite-length local algebra $V$,
$V^\vee$ is its algebraic dual.  We use $\boxtimes$ for an external
tensor product.  Every standard $L$-function denoted by
$\Lambda(s,\pi)$ is completed, including its archimedean factors;
see \cite{GodementJacquet1972}.

An \emph{associate class of cuspidal data} is a Weyl-associate class
of pairs $(M,\sigma)$, with $M$ a Levi subgroup and $\sigma$
cuspidal on $M(\A)$ \cite[Section~2]{GrbacGrobner2024}.
We call such a class a \emph{cuspidal support}.  Fix one class
$\mathfrak c$ and work in the corresponding summand.
The contragredient class $\mathfrak c^\vee$ is represented by
$(M,\sigma^\vee)$.  A
$K_\infty$-type is an irreducible finite-dimensional representation
of the maximal compact group $K_\infty$.

\subsection{Automorphic forms and the comparison theorem}
\label{subsec:automorphic-comparison}

We now specify the automorphic complex for the two coefficient
systems $\mathscr L_1$ and $\mathscr L_{\det}$ at full level, following \cite[Sections~3--7]{Franke1998}.

For an algebraic group $H$, denote by $H^\circ$ its identity
component.  For a real Lie group, the same notation denotes the
identity component in the real topology, and
$\pi_0(H)=H/H^\circ$ denotes its group of connected components.

Let $A_G$ be the maximal $\Q$-split torus in the centre of $G$,
and put
\[
 K_f=\prod_p\GL_n(\Z_p),\qquad
 A_G(\R)^\circ=\R_{>0}I_n.
\]

For a standard parabolic $P=M_PN_P$, let $A_P$ denote the maximal
$\Q$-split torus in the centre of $M_P$. Let $X_*(A_P)$ denote the cocharacter lattice of $A_P$, and put
\[
 \mathfrak a_P=X_*(A_P)\otimes_\Z\R,\qquad
 \mathfrak a_P^G=\mathfrak a_P/\mathfrak a_G.
\]
For the fixed minimal parabolic $P_0$, set
\[
 A_0=A_{P_0},\qquad
 \mathfrak a_0=\mathfrak a_{P_0},\qquad
 \mathfrak a_0^G=\mathfrak a_{P_0}^G.
\]
Using the standard invariant inner product, we identify the quotient
$\mathfrak a_P^G$ with the orthogonal complement of
$\mathfrak a_G$ in $\mathfrak a_P$.

For a standard parabolic $P$, including $P=G$, the
\emph{Harish--Chandra height map}
$H_P:G(\A)\to\mathfrak a_P$ is defined as follows.
Choose an Iwasawa decomposition $g=pk$, with $p\in P(\A)$
and $k\in K_\infty K_f$.  Then $H_P(g)$ is characterized by
\[
 \langle\xi,H_P(g)\rangle
 =\sum_v\log|\xi(p_v)|_v
\]
for every $\Q$-rational character $\xi$ of $P$.  Here the
sum runs over all places of $\Q$, and $|\cdot|_v$ is the
normalized local absolute value.  The value is independent of the
Iwasawa decomposition.

For $\lambda\in(\mathfrak a_P)^*_\C$, the \emph{spectral
character} on $M_P(\A)$ is
\[
 m\longmapsto e^{\langle\lambda,H_P(m)\rangle}.
\]
In normalized parabolic induction, this character multiplies the
modulus factor $\delta_P^{1/2}$.  The function
$g\mapsto e^{\langle\lambda,H_P(g)\rangle}$ is its extension
to $G(\A)$ via Iwasawa decomposition.

Let
\[
 [G]:=G(\Q)A_G(\R)^\circ\backslash G(\A).
\]
An \emph{automorphic form} on $[G]$ is smooth at infinity,
right invariant under a compact open subgroup of $G(\A_f)$,
and finite under both $K_\infty$ and the centre
$\mathfrak Z_G=Z(U(\mathfrak g_\C))$ of the universal enveloping
algebra \cite{BorelJacquet1979}; see also
\cite[Section~3.1]{Borel2007}.  We require uniform moderate growth:
for an adelic norm $\lVert\cdot\rVert$ as in
\cite[Page~919]{Arthur1978}, there is an integer $N$, independent
of $D$, such that on every Siegel set $\mathfrak S$,
\[
  |R(D)\phi(g)|
  \le C_{D,\mathfrak S}\lVert g\rVert^N,
  \qquad
  D\in U(\mathfrak g_\C),\quad g\in\mathfrak S.
\]
Here $R(D)$ denotes right differentiation by $D$.
Under the preceding $K_\infty$- and $\mathfrak Z_G$-finiteness
conditions, this is equivalent to the usual moderate-growth condition;
see \cite[Section~4.3]{Borel2007}.  Write $\mathcal A(G)$ for this space.

For a $(\mathfrak g,K_\infty)$-module $V$, define its relative
cochain complex, following
\cite[Chapter~I, Sections~1 and 5]{BorelWallach2000}, by
\[
 C^q(\mathfrak g,K_\infty;V)
 =\operatorname{Hom}_{K_\infty}
   \left(\textstyle\bigwedge^q
   \bigl((\mathfrak g/\mathfrak k)\otimes_\R\C\bigr),V\right),
\]
with the relative Chevalley--Eilenberg differential.  Its cohomology is
denoted by
\[
 H^\bullet(\mathfrak g,K_\infty;V)
 =H^\bullet\!\left(C^\bullet(\mathfrak g,K_\infty;V)\right).
\]

To remove the scalar central direction, let
\[
 G^1:=
 \left(\bigcap_{\chi\in X_\Q^*(G)}\ker\chi\right)^\circ
 =\SL_n,
 \qquad
 \mathfrak g^1=\Lie(G^1(\R)),
\]
and set
\[
 \mathfrak Z_G^1=Z\bigl(U(\mathfrak g^1_\C)\bigr).
\]
Thus $\mathfrak g^1=\mathfrak{sl}_n(\R)$.  If
$J\subset\mathfrak Z_G^1$ is an ideal of finite codimension and $V$
is a $\mathfrak Z_G^1$-module, set, following
\cite[Section~4, Equation~(1)]{Franke1998},
\[
 \operatorname{Fin}_J V
 =\{\,v\in V\mid J^Nv=0\quad\text{for some }N\ge1\,\}.
\]
Thus $\operatorname{Fin}_J V$ is the $J$-power-torsion submodule.
Its infinitesimal-character support is contained in the finite set of
characters $\chi:\mathfrak Z_G^1\to\C$ satisfying
$J\subset\ker\chi$.

Let
\[
 \mathfrak k^1=\mathfrak k\cap\mathfrak g^1=\mathfrak{so}(n)
\]
and put
\[
 \mathfrak p_G
 =(\mathfrak g^1/\mathfrak k^1)\otimes_\R\C
 \cong
 \bigl(\mathfrak g/(\mathfrak k\oplus\mathfrak a_G)\bigr)
 \otimes_\R\C.
\]
For $\varepsilon\in\{1,{\det}\}$, let $E_\varepsilon$ denote,
respectively, the trivial representation or the algebraic determinant
representation of $G(\R)$.  Its restriction to $\mathfrak g^1$ is
trivial, while its restriction to $K_\infty=\O_n$ is the
one-dimensional representation $\varepsilon$.  Accordingly, set
\[
 \mathcal A_\varepsilon(G)^{K_f}
 =\mathcal A(G)^{K_f}\otimes\varepsilon ,
\]
viewed as a $(\mathfrak g^1,K_\infty)$-module.  This is the
automorphic coefficient module corresponding to the local system
$\mathscr L_\varepsilon$ of
\eqref{eq:shapiro-decomposition}.

Let
\[
 J_0=\operatorname{Ann}_{\mathfrak Z_G^1}(\mathbf1)
\]
be the maximal ideal giving the trivial
$\mathfrak g^1_\C$-module.  Since $\mathfrak g^1$ acts trivially on
both $E_1$ and $E_{\det}$, the same ideal $J_0$ occurs for the two
coefficient systems.

Define the automorphic cochain complex by
\begin{equation}\label{eq:aut-complex}
 C_{\rm aut}^q(\varepsilon,J_0)= C^q(\mathfrak g^1,K_\infty;\operatorname{Fin}_{J_0}
   \mathcal A_\varepsilon(G)^{K_f})
 =\operatorname{Hom}_{K_\infty}
 \left(
   \bigwedge^q\mathfrak p_G,\,
   \operatorname{Fin}_{J_0}
   \mathcal A_\varepsilon(G)^{K_f}
 \right),
\end{equation}
with relative Chevalley--Eilenberg differential $d_{\rm aut}$.
The argument space $\mathfrak p_G$ already omits the split
central direction $\mathfrak a_G$.

A cochain $c\in C_{\rm aut}^q(\varepsilon,J_0)$ is determined
by its values on a basis of $\bigwedge^q\mathfrak p_G$.
Equivariance makes $\operatorname{im}(c)$ a finite-dimensional
$K_\infty$-submodule.  We will choose the growth bound for a
cocycle from these finitely many values.

To relate this complex to geometric cohomology, replace its
automorphic coefficient module by
$C^\infty([G])^{K_f}\otimes E_\varepsilon$.
The resulting relative cochain complex identifies with the de Rham
complex of $\widetilde Y_n$ with coefficients in
$\mathscr L_\varepsilon$
\cite[Section~2.3, Equation~(16), pp.~197--198; Section~7.4,
p.~255]{Franke1998}.
Thus inclusion of automorphic forms into smooth functions gives a
map of complexes, with $d_{\rm aut}$ corresponding to the
exterior derivative.  Franke's comparison theorem gives the
following cohomology isomorphism.

\begin{theorem}\label{thm:franke-comparison}
For $\varepsilon\in\{1,{\det}\}$, the inclusion
\[
 \operatorname{Fin}_{J_0}\mathcal A_\varepsilon(G)^{K_f}
 \hookrightarrow
 C^\infty([G])^{K_f}\otimes E_\varepsilon
\]
induces an isomorphism on
$(\mathfrak g^1,K_\infty)$-cohomology:
\begin{equation}\label{eq:franke-comparison}
 H^q\bigl(C_{\rm aut}^\bullet(\varepsilon,J_0)\bigr)
 \cong H^q(\widetilde Y_n;\mathscr L_\varepsilon).
\end{equation}
\end{theorem}

\begin{proof}
Apply \cite[Theorem~18, p.~255]{Franke1998} with $E_1=\mathbf1$
or $E_{\det}={\det}$.  Both have infinitesimal character $J_0$
on $\mathfrak g^1$, and their restrictions to $K_\infty$
are the prescribed characters.

Taking $K_f$-fixed vectors is exact because $K_f$ is compact.
The finite-adelic twist $(\xi_{E_\varepsilon})$ in Franke's
comparison \cite[Section~2.3, Theorem~4, p.~198]{Franke1998}
restricts trivially to $K_f$, since
$H_G(k_f)=0$ for $k_f\in K_f$.  At full level,
\[
 G(\Q)A_G(\R)^\circ\backslash G(\A)/(K_fK_\infty)
 \cong
 \GL_n(\Z)\backslash\GL_n(\R)/(\R_{>0}\O_n)
 =\widetilde Y_n,
\]
and the restrictions of $E_1$ and $E_{\det}$ to
$\GL_n(\Z)$ give precisely the local systems
$\mathscr L_1$ and $\mathscr L_{\det}$.
Thus Franke's comparison is exactly
\eqref{eq:franke-comparison}.
\end{proof}

\subsection{Constant-term exponents and the Franke filtration}

The Franke filtration records how far the real parts of the constant-term
exponents lie in the positive chamber.  We first fix the normalization
of these exponents, then use their positive parts to define a finite
filtration of the weighted complex.
For a recent account of the Franke filtration for general linear groups,
including its decomposition by cuspidal support, see
\cite[Section~2]{GrbacGrobner2024}.

Let $P=M_PN_P$ be standard.  Normalize the Haar measure $dn$ on
$N_P(\Q)\backslash N_P(\A)$ to have total volume $1$.  For
$\phi\in\mathcal A(G)$, define its normalized constant term by
\[
 \CT_P\phi(g)=
 e^{-\langle\rho_P,H_P(g)\rangle}
 \int_{N_P(\Q)\backslash N_P(\A)}\phi(ng)\,dn.
\]
There is a finite set
\[
 \mathcal E_P(\phi)\subset(\mathfrak a_P^G)^*_\C
\]
such that, writing $g=m\exp(X)k$ with
$X\in\mathfrak a_P^G$,
\begin{equation}\label{eq:constant-exponents}
 \CT_P\phi(m\exp(X)k)
 =
 \sum_{\lambda\in\mathcal E_P(\phi)}
 e^{\langle\lambda,X\rangle}
 \Phi_{P,\lambda}(X;m,k).
\end{equation}
Here $\Phi_{P,\lambda}(X;m,k)$ is polynomial in $X$, with
coefficients in automorphic forms on the Levi $M_P$:
\[
 \Phi_{P,\lambda}(X;m,k)
 =
 \sum_{j=1}^{r_\lambda}
 p_{P,\lambda,j}(X)\phi_{P,\lambda,j}(m,k)
\]
for suitable polynomials $p_{P,\lambda,j}$ and Levi automorphic
forms $\phi_{P,\lambda,j}$.  An element
$\lambda\in\mathcal E_P(\phi)$ for which
$\Phi_{P,\lambda}\ne0$ is called a
\emph{constant-term exponent} of $\phi$ along $P$.
This is the normalized form of
\cite[Section~6, Equations~(2)--(4)]{Franke1998}.

\begin{definition}\label{def:cuspidal-constant-term-component}
For a standard parabolic $P$, let
$(\CT_P\phi)^{\mathrm{cusp}}$ denote the projection of $\CT_P\phi$ to the space of cuspidal automorphic forms on $M_P$.  It
is characterized by having the same Petersson pairing with every cusp
form on $M_P(\Q)\backslash M_P(\A)^1$ as $\CT_P\phi$
\cite[Section~6.2, Equation~(6.8)]{BernsteinLapid2024}.  Here
$M_P(\A)^1=\ker(H_P|_{M_P(\A)})$.  Decompose
$(\CT_P\phi)^{\mathrm{cusp}}$ by its exponents as in
\eqref{eq:constant-exponents}.  A nonzero summand belonging to one
exponent is called a \emph{cuspidal constant-term component} of $\phi$
along $P$.  For a meromorphic automorphic family $F(\lambda)$, the
same term denotes the corresponding meromorphic family of summands.
\end{definition}

Franke's filtration uses the positive parts of these exponents.
Let $\Delta_P^G$ be the simple roots of $A_P$ in $N_P$, modulo
$A_G$, and let $\omega_\alpha^G$ be the corresponding fundamental
weights.  For $\mu\in(\mathfrak a_P^G)^*$, Langlands' combinatorial lemma
gives a unique subset $I(\mu)\subseteq\Delta_P^G$ and decomposition
\begin{equation}\label{eq:langlands-positive-part}
 \mu=\sum_{\alpha\notin I(\mu)}a_\alpha\omega_\alpha^G
       -\sum_{\beta\in I(\mu)}b_\beta\beta,
 \qquad a_\alpha>0,\quad b_\beta\ge0.
\end{equation}
Set $\mu_+=\sum_{\alpha\notin I(\mu)}a_\alpha\omega_\alpha^G$.
See \cite[Section~6, Lemma~1 and (5)--(7)]{Franke1998} and
\cite[Section~3.4, Lemma~3.1]{LiSun2019}.
Let
\[
 \mathscr C_P^G=\left\{\sum_{\alpha\in\Delta_P^G}c_\alpha\alpha:
                                     c_\alpha\ge0\right\},
 \qquad \mathscr C_0^G=\mathscr C_{P_0}^G.
\]
We view $(\mathfrak a_P^G)^*$ as a subspace of
$(\mathfrak a_0^G)^*$ by the canonical inclusion.

Fix $J$, a cuspidal support $\mathfrak c$, and a dominant
$\tau\in(\mathfrak a_0^G)^*$.  Let
$\{P\}_{\mathfrak c}$ be the associate class of parabolics with
Levi conjugate to a representative of $\mathfrak c$.
Write $\mathcal A(G)_{\mathfrak c}$ for its automorphic summand.
Define
\[
 \mathcal A_{\mathfrak c,\varepsilon}
 =\left\{\phi\in
     (\operatorname{Fin}_J\mathcal A(G)_{\mathfrak c})^{K_f}:
 \operatorname{Re}\lambda\in\tau-\mathscr C_0^G
 \text{ for all constant-term exponents }\lambda\right\}
 \otimes\varepsilon.
\]
By \cite[Theorem~15]{Franke1998}, the displayed exponent
condition defines Franke's bounded-growth module.  Increasing $\tau$ exhausts
$(\operatorname{Fin}_J\mathcal A(G)_{\mathfrak c})^{K_f}\otimes\varepsilon$.
For fixed $J$, only finitely many exponents occur
\cite[Section~6, paragraph preceding (8)]{Franke1998}.
Choose an integer-valued function $T$ on their positive parts such that
\begin{equation}\label{eq:T-order}
 T(\mu)<T(\nu)\quad\text{if}\quad
 \mu\ne\nu,\qquad \nu\in\mu-\mathscr C_0^G.
\end{equation}
Franke's finite decreasing filtration is
\begin{equation}\label{eq:franke-filtration}
 F^p\mathcal A_{\mathfrak c,\varepsilon}
 =\left\{\phi\in\mathcal A_{\mathfrak c,\varepsilon}:T((\operatorname{Re}\lambda)_+)\ge p
 \text{ for every exponent along }P\in\{P\}_{\mathfrak c}\right\}.
\end{equation}
Its terms are $(\mathfrak g^1,K_\infty)$-submodules.
On the relative cochain complex
$C^\bullet=C^\bullet(\mathfrak g^1,K_\infty;
\mathcal A_{\mathfrak c,\varepsilon})$, put
\begin{equation}\label{eq:filtered-aut-complex}
 F^pC^\bullet
 =\operatorname{Hom}_{K_\infty}
   \left(\bigwedge^\bullet\mathfrak p_G,
                           F^p\mathcal A_{\mathfrak c,\varepsilon}\right).
\end{equation}

Let $x=(x_1,\ldots,x_d)$, let
$R_x=\C\{x_1,\ldots,x_d\}$ be the ring of convergent holomorphic
germs at $x=0$, and let
$\mathfrak m_x=(x_1,\ldots,x_d)$.  The duals of the natural quotient
maps form an increasing system, and we set
\[
 \D_x=\bigcup_{N\ge0}(R_x/\mathfrak m_x^{N+1})^\vee.
\]
Here $(-)^\vee=\operatorname{Hom}_\C(-,\C)$.  Thus $\D_x$ is the
space of finite-order distributions supported at $x=0$, acting on
holomorphic germs by finite sums of derivative evaluations, as in
\cite[Section~6, p.~234, before (11)]{Franke1998}.
The $R_x$-module action on $\D_x$ is
\[
 (r\lambda)(f)=\lambda(rf),
 \qquad r,f\in R_x,\quad \lambda\in\D_x.
\]

\begin{definition}
\label{def:distribution-evaluation}\label{def:coefficient-map}
Let $W$ be a complete locally convex complex vector space, and
let $\mathscr O_x(W)$ denote the space of $W$-valued holomorphic
germs at $x=0$.
If $\lambda\in\D_x$ factors through
$R_x/\mathfrak m_x^{N+1}$ and $F\in\mathscr O_x(W)$, define
\begin{equation}
\label{eqn: extract-coefficient}
 \lambda_W(F)
   =(\lambda\otimes\operatorname{id}_W)(j_0^N F)\in W,
\end{equation}
where
$j_0^N F\in(R_x/\mathfrak m_x^{N+1})\otimes_\C W$
is the Taylor expansion of $F$ at $x=0$ modulo terms of total
degree greater than $N$.  The value $\lambda_W(F)$ is independent
of $N$, provided that $\lambda$ factors through
$R_x/\mathfrak m_x^{N+1}$.

Let $V$ be a complex vector space and let $H(x):V\to W$ be a
family of linear maps such that $x\mapsto H(x)v$ is holomorphic
near $0$ for every $v\in V$.  Its \emph{coefficient map} is
\[
 \Xi_H:V\otimes_\C\D_x\longrightarrow W,\qquad
 \Xi_H(v\otimes\lambda)=\lambda_W\bigl(H(\,\cdot\,)v\bigr).
\]
For $\alpha\in\Z_{\ge0}^d$, put
$\delta_\alpha(f)=\partial_x^\alpha f(0)/\alpha!$.  Then
\[
 \Xi_H(v\otimes\delta_\alpha)
 =\left.\frac1{\alpha!}\partial_x^\alpha\bigl(H(x)v\bigr)\right|_{x=0}.
\]
Thus $\operatorname{im}\Xi_H$ is the span of these Taylor coefficients.
We also use $\Xi_H$ with codomain any specified subspace of $W$
containing them.  For $\eta\in V\otimes\D_x$, write
$\eta(H):=\Xi_H(\eta)$.
\end{definition}

\subsection{Laurent-coefficient modules and the associated graded}
\label{subsec:laurent-coefficient-modules}

The description of a filtration quotient has three ingredients: a
compact realization of an induced representation
(Definition~\ref{def:compact-realization}),
finite derivatives at a spectral point, and the identifications made
by standard intertwining operators.  We introduce them in this order.

After decomposing the square-integrable Levi spectrum, an indexing datum
$t=(R,\sigma_t,\lambda_t)$ consists of
\begin{enumerate}[label=(\roman*),leftmargin=2em]
\item a standard parabolic $R=M_RN_R$ containing a member of
$\{P\}_{\mathfrak c}$;
\item a square-integrable automorphic representation $\sigma_t$ of
$M_R(\A)$.  Write $\sigma_{t,\infty}$ for its component at the real
place and $\sigma_{t,f}$ for its finite-adelic part.  We require
$\sigma_{t,f}$ to have nonzero
$K_f\cap M_R(\A_f)$-fixed vectors;
\item a spectral point
$\lambda_t\in(\mathfrak a_R^G)^*_\C$, with dominant real part, such
that, after the canonical inclusion in
$(\mathfrak a_0^G)^*_\C$,
$\operatorname{Re}\lambda_t\in
\tau-\mathscr C_0^G$, and the cuspidal support of
$\sigma_t\otimes e^{\langle\lambda_t,H_R(\cdot)\rangle}$
belongs to $\mathfrak c$ \cite[Section~2]{GrbacGrobner2024}.
\end{enumerate}
We call $\sigma_t$ the \emph{square-integrable Levi representation}
of the inducing datum $t$.
Fix a datum $t=(R,\sigma_t,\lambda_t)$, and write
\[
 M_R=\prod_{\nu=1}^s\GL_{r_\nu},\qquad
 \sigma_t=\boxtimes_{\nu=1}^s\Sigma_\nu.
\]
For $\lambda\in(\mathfrak a_R^G)^*_\C$, define its block
coordinates $\lambda_\nu$ by
\[
 e^{\langle\lambda,H_R(m)\rangle}
 =\prod_{\nu=1}^s|{\det} m_\nu|_\A^{\lambda_\nu},
 \qquad m=(m_1,\ldots,m_s)\in M_R(\A),
\]
where $|\cdot|_\A$ is the adelic norm.  The induced family is
\[
 I_t(\lambda):=\mathop{\times}_{\nu=1}^s
       \bigl(\Sigma_\nu|{\det}|_\A^{\lambda_\nu}\bigr).
\]
For each place $v$, let $\Sigma_{\nu,v}$ be the local
component of $\Sigma_\nu$, and put
\[
 I_t(\lambda)_v:=\mathop{\times}_{\nu=1}^s
       \bigl(\Sigma_{\nu,v}|{\det}|_v^{\lambda_\nu}\bigr),
\]
with normalized induction to $G(\Q_v)$.
In particular, $I_t(\lambda)_\infty$ is the induced
representation of $G(\R)$, and
$I_t(\lambda)_f:=\bigotimes_{v<\infty}'I_t(\lambda)_v$
is its finite-adelic counterpart.
In compact realizations and relative Lie algebra cohomology,
we use the $K_\infty$-finite vectors of
$I_t(\lambda)_\infty$.
Define
\[
 \chi_{\sigma_t,\lambda_t}:\mathfrak Z_G^1\longrightarrow\C
\]
to be the $G^1$-infinitesimal character of $I_t(\lambda)$ at
$\lambda=\lambda_t$.
We require the inducing datum $t$ to be \emph{$J$-compatible}, meaning that
\[
 J\subseteq\ker\chi_{\sigma_t,\lambda_t}.
\]
This condition ensures that every finite-order derivative at
$\lambda_t$ is killed by some power of $J$.
For $J=J_0$, compatibility means that $I_t(\lambda_t)_\infty$ has the same
infinitesimal character on $\mathfrak{sl}_n$ as the trivial representation.

\begin{definition}\label{def:compact-realization}
Fix an automorphic realization of $\sigma_t$, and let
$V_{\sigma_{t,\infty}}$ and $V_{\sigma_{t,f}}$ be the
representation spaces of its real and finite-adelic parts.
The \emph{compact realization} of an induced representation is
the model obtained as follows.
\begin{enumerate}[label=(\roman*),leftmargin=2em]
\item \emph{Local representations.}  At the real place, put
$K_{R,\infty}=K_\infty\cap M_R(\R)$.  Define
\begin{equation}\label{eq:compact-picture}
 \mathcal I_{R,K}(\sigma_{t,\infty})=
 \left\{
 F\in C^\infty(K_\infty,V_{\sigma_{t,\infty}}):
 F(k_Rk)=\sigma_{t,\infty}(k_R)F(k)
 \right\}_{K_\infty\text{-finite}},
\end{equation}
where $k_R\in K_{R,\infty}$, $k\in K_\infty$, and
$K_\infty$ acts by right translation.

At a finite place $v$, let $K_v=\GL_n(\Z_v)$ and let
$V_{\sigma_{t,v}}$ be the space of the local component
$\sigma_{t,v}$.  Define $\mathcal I_{R,K_v}(\sigma_{t,v})$
to be the space of locally constant functions
$F:K_v\to V_{\sigma_{t,v}}$ satisfying
\[
 F(m_Rn_Rk)=\sigma_{t,v}(m_R)F(k),
\]
for $m_R\in M_R(\Q_v)\cap K_v$,
$n_R\in N_R(\Q_v)\cap K_v$, and $k\in K_v$.
The group $K_v$ acts by right translation.

\item \emph{Adelic representations.}  Put
$K_{R,f}=M_R(\A_f)\cap K_f$ and
$V_{\sigma_t}=V_{\sigma_{t,\infty}}\otimes_\C
V_{\sigma_{t,f}}$.  Its compact realization is
\begin{equation}\label{eq:adelic-compact-realization}
 \mathcal I_{R,K}(\sigma_t)=
 \left\{
 F\in C^\infty(K_\infty K_f,V_{\sigma_t}):
 F(m_Rn_Rk)=\sigma_t(m_R)F(k)
 \right\}_{K_\infty\text{-finite}},
\end{equation}
where $m_R\in K_{R,\infty}K_{R,f}$,
$n_R\in N_R(\A_f)\cap K_f$, and $k\in K_\infty K_f$.
The group $K_\infty K_f$ acts by right translation.

\item \emph{Full level.}  Write $u_t=(R,\sigma_t)$ for the pair obtained
by omitting the spectral point from $t$, and define
\[
 W(u_t)=\mathcal I_{R,K}(\sigma_{t,\infty})
          \otimes_\C V_{\sigma_{t,f}}^{K_{R,f}}.
\]
Every $f\in W(u_t)$ is a finite sum $\sum_iF_i\otimes v_i$
in this algebraic tensor product.
The compact realization of a
$K_f$-fixed adelic section means this space $W(u_t)$.
\end{enumerate}
\end{definition}

\begin{lemma}\label{lem:compact-realization}
For every $\lambda$, restriction $f\mapsto f|_K$ gives
the isomorphisms
\[
\begin{aligned}
 \operatorname{res}_{K_\infty}:
 \bigl(I_t(\lambda)_\infty\bigr)_{K_\infty\text{-finite}}
 &\xrightarrow{\ \sim\ }\mathcal I_{R,K}(\sigma_{t,\infty}),\\
 \operatorname{res}_{K_\infty K_f}:
 I_t(\lambda)_{K_\infty\text{-finite}}
 &\xrightarrow{\ \sim\ }\mathcal I_{R,K}(\sigma_t).
\end{aligned}
\]
These maps intertwine right translation by $K_\infty$ and
$K_\infty K_f$, respectively.  At each finite place $v$,
restriction gives an isomorphism of $K_v$-modules
\[
\begin{aligned}
 \operatorname{res}_{K_v}:I_t(\lambda)_v
 &\xrightarrow{\ \sim\ }\mathcal I_{R,K_v}(\sigma_{t,v}),\\
 f&\longmapsto f|_{K_v},
\end{aligned}
\]
where the target is the space of locally constant functions with
the covariance condition of
Definition~\ref{def:compact-realization}\textup{(i)}.
The spaces in Definition~\ref{def:compact-realization} and
their compact-group actions are independent of $\lambda$.
At full level, restriction further gives isomorphisms of
$K_\infty$-modules
\[
 \bigl(I_t(\lambda)^{K_f}\bigr)_{K_\infty\text{-finite}}
 \cong\mathcal I_{R,K}(\sigma_t)^{K_f}\cong W(u_t).
\]
\end{lemma}

\begin{proof}
The real, finite-place, and adelic assertions use the same
extension construction.  In $G(\R)$, $G(\Q_v)$, and
$G(\A)$, Iwasawa decomposition writes $g=nmk$, with
$n\in N_R$ and $m\in M_R$ over the same field or adeles,
and with $k\in K_\infty$, $K_v$, or $K_\infty K_f$,
respectively.  At the real place and over $\A$, set
\[
\begin{aligned}
 \widetilde F_{\lambda,\infty}(nmk)
 &=e^{\langle\lambda+\rho_R,H_R(m)\rangle}
       \sigma_{t,\infty}(m)F(k),\\
 \widetilde F_\lambda(nmk)
 &=e^{\langle\lambda+\rho_R,H_R(m)\rangle}
       \sigma_t(m)F(k).
\end{aligned}
\]
At a finite place $v$, the same formula is
\[
 \widetilde F_{\lambda,v}(nmk)
 =\delta_{R,v}(m)^{1/2}
   \left(\prod_{\nu=1}^s|{\det} m_\nu|_v^{\lambda_\nu}\right)
   \sigma_{t,v}(m)F(k),
\]
where $m=(m_\nu)_\nu$ and $\delta_{R,v}$ is the modulus
character of $R(\Q_v)$.

In each case, two Iwasawa decompositions differ by an element of
$R\cap K$.  The spectral and modulus characters are trivial
on this intersection, so the covariance of $F$ in
Definition~\ref{def:compact-realization} makes the extension
well defined.  It satisfies the normalized inducing covariance
and restricts to $F$; conversely, that covariance recovers any
induced section from its restriction.  The formulas commute with
right translation by $K$, whose action on $F$ is independent
of $\lambda$.  They preserve smoothness and
$K_\infty$-finiteness at infinity and local constancy at finite
places.  This proves all three restriction isomorphisms; see
\cite[Section~II.1, p.~643]{MoeglinWaldspurger1989}.

At full level, evaluation at the identity identifies
$I_t(\lambda)_f^{K_f}$ with
$V_{\sigma_{t,f}}^{K_{R,f}}$, by
$G(\A_f)=R(\A_f)K_f$.  Tensoring with the real compact
model gives the last isomorphism.  This is the full-level
space in Franke's construction
\cite[Section~5.2, p.~218, and Section~6, (10)--(11), p.~234]{Franke1998}.
\end{proof}

For $f=\sum_iF_i\otimes v_i\in W(u_t)$, 
$f(k_\infty)=\sum_iF_i(k_\infty)\otimes v_i$ is viewed as an
automorphic form on $M_R$ in the chosen realization of
$\sigma_t$.  The scalar inducing function used in the Eisenstein
sum is
\[
 f_\lambda(n_Rm k_\infty k_f)
   :=e^{\langle\lambda+\rho_R,H_R(m)\rangle}
        \bigl[f(k_\infty)\bigr](m),
\]
where $n_R\in N_R(\A)$, $m\in M_R(\A)$,
$k_\infty\in K_\infty$, and $k_f\in K_f$.  Here
\[
 \bigl[f(k_\infty)\bigr](m)
 =\sum_i\bigl[F_i(k_\infty)\otimes v_i\bigr](m).
\]
Each tensor is viewed as an automorphic form through the chosen
realization of $\sigma_t$, then evaluated at $m$.

Covariance in \eqref{eq:compact-picture} and $K_{R,f}$-invariance make the
formula independent of the Iwasawa decomposition.  For fixed $f$,
only the exponential varies with $\lambda$.  The family
$\lambda\mapsto f_\lambda$ is called the \emph{flat section}
associated with $f$.

Let $\mathcal O_{\lambda_t}$ be the local ring of holomorphic germs
at $\lambda=\lambda_t$, with maximal ideal
$\mathfrak m_{\lambda_t}$.  The space
\begin{equation}\label{eq:Dt}
D_t=\bigcup_{N\ge0}
 \left(\mathcal O_{\lambda_t}/
       \mathfrak m_{\lambda_t}^{N+1}\right)^\vee
\end{equation}
is the space of finite-order holomorphic functionals supported at
$\lambda_t$.  Such a functional has order at most $N$ if it
factors through
$\mathcal O_{\lambda_t}/\mathfrak m_{\lambda_t}^{N+1}$; equivalently,
it is a linear combination of derivatives of total degree at most
$N$.  Define
\begin{equation}\label{eq:Mt}
 M(t):=W(u_t)\otimes D_t.
\end{equation}
The following Leibniz
rule makes $M(t)$ a $(\mathfrak g,K_\infty)$-module.
For each $\lambda$, the identification $f\mapsto f_\lambda$
transports the $\mathfrak g$-action on $I_t(\lambda)$ to
the fixed space $W(u_t)$.  For $X\in\mathfrak g$, the
resulting operator $\pi_\lambda(X):W(u_t)\to W(u_t)$ is
characterized by
\[
 \bigl(\pi_\lambda(X)f\bigr)_\lambda=R(X)f_\lambda,
 \qquad f\in W(u_t),
\]
where $R(X)$ is right differentiation in the real group variable.
The $K_\infty$-action is independent of
$\lambda$, so on $M(t)$ it is
\[
 k\cdot(f\otimes\delta)=(k\cdot f)\otimes\delta.
\]
Choose linear coordinates at $\lambda_t$, and let
$\partial^\alpha_{\lambda_t}\in D_t$ denote the functional
\[
 h\longmapsto
 \left.\partial^\alpha_\lambda h(\lambda)\right|_{\lambda=\lambda_t}.
\]
The $\mathfrak g$-action is defined by the Leibniz rule
\begin{equation}\label{eq:g-action-jets}
 X\cdot\bigl(f\otimes\partial^\alpha_{\lambda_t}\bigr)
 =
 \sum_{\beta+\gamma=\alpha}
 \binom{\alpha}{\beta}
 \left.
   \partial^\gamma_\lambda\bigl(\pi_\lambda(X)f\bigr)
 \right|_{\lambda=\lambda_t}
 \otimes\partial^\beta_{\lambda_t}.
\end{equation}
For the jet action \eqref{eq:g-action-jets}, the condition
$J\subseteq\ker\chi_{\sigma_t,\lambda_t}$ makes $M(t)$
$J$-power-torsion.  

For $f\in W(u_t)$, the flat section $f_\lambda$ gives the
Eisenstein series, initially for $\operatorname{Re}\lambda$
sufficiently dominant,
\[
 E_{R,t}^G(f,\lambda)(g)
 =\sum_{\gamma\in R(\Q)\backslash G(\Q)}f_\lambda(\gamma g),
 \qquad g\in G(\A).
\]
The sum is linear in $f$, converges absolutely in this region,
and continues meromorphically to $(\mathfrak a_R^G)^*_{\C}$ as a
family of smooth automorphic functions
\cite[Sections~4, 6, and 7]{Langlands1976}.  We write
$E_R^G(f,\lambda)$ when $t$, and hence $\sigma_t$, is fixed.

\begin{definition}\label{def:franke-indexing-groupoid}
For $k\ge0$ and $p\in\Z$, let
$\mathscr T_{\mathfrak c,k,\tau,p}$ have the indexing data
$t=(R,\sigma_t,\lambda_t)$ above as objects, subject to
\[
 \dim\mathfrak a_R^G=k,\qquad
 \operatorname{Re}\lambda_t\in\tau-\mathscr C_0^G,\qquad
 T((\operatorname{Re}\lambda_t)_+)=p,\qquad
 J\subseteq\ker\chi_{\sigma_t,\lambda_t}.
\]
A morphism $w:t\to t'$ is an element
$w\in N_G(A_0)(\Q)/Z_G(A_0)(\Q)$ satisfying
\[
 wM_Rw^{-1}=M_{R'},\qquad
 [{}^w\sigma_t]=[\sigma_{t'}],\qquad
 w\lambda_t=\lambda_{t'}.
\]
Here $\dot w\in N_G(A_0)(\Q)$ represents $w$,
${}^w\sigma_t(m)=\sigma_t(\dot w^{-1}m\dot w)$, and
brackets denote isomorphism classes.  Composition is multiplication
of Weyl elements.  This is the groupoid of
\cite[Section~6, before (12)]{Franke1998}, restricted to the
fixed support and filtration degree.
\end{definition}

For these morphisms, the standard intertwiner
$M(w,\lambda):W(u_t)\to W(u_{t'})$ is holomorphic at
$\lambda_t$ \cite[Section~6, Lemma~2]{Franke1998}.
Transporting spectral coordinates by $w$, its action on distributions is
\begin{equation}\label{eq:intertwiner-jets}
 w_*(f\otimes\partial^\alpha_{\lambda_t})
 =\sum_{\beta+\gamma=\alpha}\binom\alpha\beta
   \left.\partial_\lambda^\gamma M(w,\lambda)f
   \right|_{\lambda=\lambda_t}\otimes\partial^\beta_{w\lambda_t}.
\end{equation}
These maps make $t\mapsto M(t)$ a functor.  Its colimit identifies
$v\in M(t)$ with $w_*v\in M(t')$ for every morphism $w:t\to t'$.

\begin{theorem}\label{thm:franke-graded}
There is an isomorphism of $(\mathfrak g^1,K_\infty)$-modules
\begin{equation}\label{eq:franke-graded}
 \bigoplus_k\left(\varinjlim_{t\in\mathscr T_{\mathfrak c,k,\tau,p}}
 M(t)\right)\otimes E_\varepsilon
 \xrightarrow{\ \sim\ }\operatorname{gr}_F^p
                                \mathcal A_{\mathfrak c,\varepsilon}.
\end{equation}
Here $0\le k\le\dim\mathfrak a_P^G$, for
$P\in\{P\}_{\mathfrak c}$.
\end{theorem}

\begin{proof}
This is \cite[Theorem~14 and Section~6, (9)--(14)]{Franke1998}.
For a holomorphic Eisenstein family, the map sends
$f\otimes\partial^\alpha_{\lambda_t}$ to its derivative
$\left.\partial_\lambda^\alpha E_R^G(f,\lambda)
\right|_{\lambda=\lambda_t}$, modulo the next filtration term.
At a pole, Franke replaces evaluation by the constant Laurent
coefficient on a line transverse to all polar hyperplanes.
His theorem asserts independence of that line and equivariance
\emph{after} passage to the filtration quotient.
Section~6, Proposition~1 of the same paper gives compatibility
with cuspidal support.  Taking $K_f$-fixed vectors and tensoring
with $E_\varepsilon$ preserve the isomorphism.
\end{proof}

\begin{lemma}\label{lem:groupoid-colimit-cohomology}
For each groupoid $\mathscr T_{\mathfrak c,k,\tau,p}$ of
Definition~\ref{def:franke-indexing-groupoid},
the
canonical map
\[
 \varinjlim_t H^\bullet(\mathfrak g^1,K_\infty;M(t))
 \longrightarrow
 H^\bullet\!\left(\mathfrak g^1,K_\infty;
                   \varinjlim_tM(t)\right)
\]
is an isomorphism.  In particular, every class on the right is a
finite sum of images of closed cochains belonging to individual
objects $M(t)$.
\end{lemma}

\begin{proof}
This is the passage to cohomology used in
\cite[Section~7.4, Equation~(1), p.~256]{Franke1998}.
The colimit is a direct sum of coinvariants under finite
automorphism groups.  Taking coinvariants under a finite group
is exact over $\C$, and relative cochains commute with these
coinvariants and direct sums.  The last assertion follows from
finite support in the direct sum.
\end{proof}

\begin{corollary}\label{cor:franke-degree-reduction}
Let $F$ be a finite-dimensional algebraic representation of $G$ whose dual
has infinitesimal character $\chi$ on $\mathfrak g^1$, and put $J=\ker\chi$.
If
\[
 H^q\!\left(\mathfrak g^1,K_\infty;
   (\operatorname{Fin}_J\mathcal A(G)_{\mathfrak c})^{K_f}\otimes F\right)\ne0,
\]
then $H^q(\mathfrak g^1,K_\infty;M(t)\otimes F)\ne0$ for some
$J$-compatible inducing datum $t$ of cuspidal support $\mathfrak c$.
\end{corollary}

\begin{proof}
A cocycle representing a nonzero class is determined by its values on a
basis of $\bigwedge^q\mathfrak p_G$.  These finitely many values lie in
$\mathcal A_{\mathfrak c,1}\otimes F$ for some growth bound $\tau$.
Its class in the corresponding relative cochain complex $C^\bullet$ is
still nonzero.  The finite filtration of \eqref{eq:filtered-aut-complex}
with coefficients in $F$ gives
\begin{equation}\label{eq:franke-spectral-sequence}
 E_1^{p,r}=H^{p+r}(\operatorname{gr}_F^pC^\bullet)
 \quad\Longrightarrow\quad H^{p+r}(C^\bullet).
\end{equation}
Thus $H^q(\operatorname{gr}_F^pC^\bullet)\ne0$ for some $p$.
Apply Theorem~\ref{thm:franke-graded} with $\varepsilon=1$ and tensor
with $F$.  The proof of Lemma~\ref{lem:groupoid-colimit-cohomology}
applies to $M(t)\otimes F$, since the intertwining maps are
$w_*\otimes\mathrm{id}_F$ and tensoring with $F$ commutes with direct
sums and coinvariants.  Hence some $M(t)\otimes F$ has nonzero cohomology
in degree $q$.
\end{proof}

\subsection{Cuspidal-support summands and the trivial support}
\label{subsec:topological-trivial-datum}

For the cuspidal support $\mathfrak c$ of
Subsection~\ref{subsec:global-notation}, let
$\mathcal A(G)_{\mathfrak c}$ be its associate-class summand.
The decomposition
\[
 \mathcal A(G)=\bigoplus_{\mathfrak c}\mathcal A(G)_{\mathfrak c}
\]
is recalled in \cite[Section~2]{GrbacGrobner2024}.
When $G$ is clear, we write $\mathcal A_{\mathfrak c}=\mathcal A(G)_{\mathfrak c}$
and denote its subspace fixed by $K_f=\prod_p\GL_n(\Z_p)$ by $\mathcal A_{\mathfrak c}^{K_f}$.
Theorem~\ref{thm:franke-comparison} and
\eqref{eq:shapiro-decomposition} give
\[
 H^q(Y_n;\C)_{\mathfrak c}
 :=\bigoplus_{\varepsilon\in\{1,{\det}\}}
 H^q\!\left(\mathfrak g^1,K_\infty;
 (\operatorname{Fin}_{J_0}\mathcal A(G)_{\mathfrak c})^{K_f}
                                     \otimes\varepsilon\right),
 \qquad
 H^q(Y_n;\C)=\bigoplus_{\mathfrak c}H^q(Y_n;\C)_{\mathfrak c}.
\]
Let $\mathrm{PD}:H_c^j(Y_n;\C)\to H^{d_n-j}(Y_n;\C)^\vee$
be the duality isomorphism in \eqref{eq:duality}.  Define
\[
 H_c^j(Y_n;\C)_{\mathfrak c}
 =\mathrm{PD}^{-1}\!\left(
       H^{d_n-j}(Y_n;\C)_{\mathfrak c^\vee}^{\vee}\right).
\]
The dual of one summand is viewed inside the dual of the full
cohomology by extending functionals by zero on the other summands.
The contragredient occurs because the adjoint of a spherical
Hecke operator $h$ is $h^\vee(g)=h(g^{-1})$.

Let ${\rm triv}$ be the support of the trivial representation:
\[
 \left(M_{P_0},\mathop{\boxtimes}_{i=1}^n
                         |\cdot|_\A^{(n+1-2i)/2}\right).
\]
It is self-dual, since negating the exponents reverses their order.
Franke identifies this support summand with the generalized trivial
Hecke eigenspace \cite[Section~5]{Franke2008}.

\begin{theorem}\label{prop:constant}
For $m\ge2$,
\begin{equation}\label{eq:constant}
 H_c^j(Y_m;\C)_{\rm triv}=0\qquad(j<2m-2).
\end{equation}
For $m=1$ its lowest degree is $0$.
\end{theorem}

\begin{proof}
For $m=1$, $Y_m$ is a point.  Assume $m\ge2$.
Franke's compact-dual calculation
\cite[Theorem~3.1, Corollary~3.3, and Theorems~7.10--7.11]{Franke2008}
gives two families of homogeneous basis classes for the full-level
trivial-support summand.  We recall their indexing and degree
formulas, then bound each degree.
Put
\[
 \mathcal O_m=\{q:3\le q\le m,\ q\text{ odd}\}.
\]
Each occurrence of $q$ records a primitive compact-dual generator
of degree $2q-1$.  For a multiplicity function
$l:\mathcal O_m\to\Z_{\ge0}$, put
\[
 d(l)=\sum_{q\in\mathcal O_m}(2q-1)l(q).
\]

The first family \cite[Theorem~7.10 and (7.19)--(7.20)]{Franke2008}
is indexed by $e\in\{0,1\}$, $N\ge0$, and nonempty subsets
$X_0,\ldots,X_N\subset\mathcal O_m$, with
\[
 l(q)=\#\{i:q\in X_i\},\qquad
 \sum_{i=0}^N\max X_i=m-e,\qquad
 j=d(l)+N+e.
\]
The term $N+e$ is the incidence degree: the associated parabolic
has $N+1+e$ Levi blocks.  For each $i$,
\[
 \sum_{q\in X_i}(2q-1)+1\ge2\max X_i.
\]
Summing gives
\[
 d(l)+N+1\ge2(m-e),\qquad
 j\ge2m-e-1\ge2m-2.
\]

The second family occurs only for even $m$ and contains the
product of the Levi Euler classes, of total degree $m$
\cite[Theorem~7.11]{Franke2008}.
For $m>2$, its indices satisfy
\[
 \sum_{i=0}^N(1+\max X_i)=m-2e,\qquad
 j=m+d(l)+N+e,
\]
Here $e\in\{0,1\}$, $N\ge0$, each $X_i$ is a nonempty subset of
$\mathcal O_m$, and $l(q)=\#\{i:q\in X_i\}$.
Since $2q-1\ge q$,
\[
 d(l)\ge\sum_{i=0}^N\max X_i=m-2e-(N+1),
 \qquad j\ge2m-e-1\ge2m-2.
\]
For $m=2$, the exceptional Euler class has degree $2=2m-2$.
These families exhaust the summand, proving \eqref{eq:constant}.
Finally, \eqref{eq:shapiro-decomposition} makes each
$H_c^j(\widetilde Y_m;\mathscr L_\varepsilon)_{\rm triv}$ a direct
summand of $H_c^j(Y_m;\C)_{\rm triv}$, so the bound holds for
both determinant parities.
\end{proof}

\begin{remark}
\label{rem: trivial-support-proof}
While revising this paper, we became aware of the related
observation of Church--Farb--Putman \cite[Remark~1.5]{CFP19}.
For an imaginary quadratic field $K$, they note that Franke's
results imply vanishing of the Hecke-trivial subspace of
$H^q(\SL_n\mathcal O_K;\Q)$ for $q>\nu_n-n+1$, where
$\nu_n=\operatorname{vcd}(\SL_n\mathcal O_K)$.
\end{remark}

\section{Speh representations and cohomological degrees}\label{sec:packets}

\subsection{Harish--Chandra parameters of Speh representations}\label{subsec:speh-weights}

We now use the inducing data to estimate the degrees on the first
page of \eqref{eq:franke-spectral-sequence}.  We express these bounds
in the complementary degree $j=d_n-q$.  The infinitesimal-character
condition partitions $U_n$ into subsets $B_\nu$, one for each
inducing factor.  We obtain these subsets from the Harish--Chandra
parameters as follows.
\begin{definition}\label{def:weights}
Let $\sigma$ be a representation of $\GL_s(\R)$ with infinitesimal
character.  Via the Harish--Chandra isomorphism
\cite[Chapter~VI, Section~0.2]{BorelWallach2000}, represent this
infinitesimal character by a semisimple conjugacy class in
$\mathfrak{gl}_s(\C)$.
Write $\{\lambda_1,\ldots,\lambda_s\}$ for the unordered multiset of
eigenvalues of this class, counted with multiplicity; this is the
\emph{Harish--Chandra parameter} of $\sigma$.  The Harish--Chandra
parameter of the trivial representation is
\[
 \left\{\frac{1-s}{2},\frac{3-s}{2},\ldots,\frac{s-1}{2}\right\}.
\]
For an automorphic representation, we use the parameter of its
real-place component.  Twisting by $|{\det}|^u$ adds $u$ to
each entry of the parameter.
\end{definition}

M\oe glin--Waldspurger classify the irreducible discrete automorphic
representations of general linear groups as Speh representations
\cite[p.~606 and Sections~II.2--II.3]{MoeglinWaldspurger1989}.
Consequently, an irreducible square-integrable automorphic
representation of a Levi subgroup is an external tensor product
of Speh representations.

In the cohomological summands below, the infinitesimal character on
$\mathfrak g^1$ is that of the trivial coefficient and is therefore regular
algebraic.  Let $\pi$ be a unitary cuspidal automorphic
representation of $\PGL_k(\A)$ whose pullback to $\GL_k(\A)$
becomes regular algebraic after a twist by $|{\det}|^s$,
in the sense of
\cite[Definitions~1.8 and~3.12, pp.~89--90 and~111]{Clozel1990}.
We use the unitary algebraicity convention of
\cite[Definition~8.2.7]{ChenevierLannes2019}: the
Harish--Chandra entries are distinct elements of
$\tfrac12\Z$ with integral pairwise differences.
References below to a regular algebraic unitary
cuspidal representation use this convention; the unitary
representation itself need not be algebraic in Clozel's
normalization.  All $L$-functions retain the unitary
normalization of $\pi$.
For $\ell\ge1$, let $J(\ell,\pi)$ denote the Speh
representation obtained as the Langlands quotient of the normalized
induction
\[
 \pi|{\det}|^{(\ell-1)/2}\times
 \pi|{\det}|^{(\ell-3)/2}\times\cdots\times
 \pi|{\det}|^{-(\ell-1)/2}.
\]
Here $\operatorname{rank}(\pi)=k$, whereas $J(\ell,\pi)$ is a
representation of $\GL_{k\ell}(\A)$.  Its cuspidal support is
represented on $\GL_k^\ell$ by
$\boxtimes_{j=0}^{\ell-1}\pi|{\det}|^{(\ell-1)/2-j}$
\cite[Section~3, Equations~(3.1)--(3.2)]{GrbacGrobner2024}.
We distinguish factors with
$k\ge2$ from the case $k=1$, in which $\pi=1$ and we write
$J(m,1)$.  With our Langlands-quotient convention,
\[
                         J(m,1)\cong\mathbf1_{\GL_m(\A)}.
\]
The parameter of $J(m,1)|{\det}|^s$ is therefore obtained by
adding $s$ to the trivial parameter in Definition~\ref{def:weights},
with rank $m$.  We keep the determinant twist separate from
$J(m,1)$.

Let $\lambda_1<\cdots<\lambda_k$ be the Harish--Chandra parameter
of $\pi$.  The entries have integral pairwise differences.
Purity applied after the twist that makes $\pi$ algebraic, followed
by unitary normalization, makes them real and symmetric
about zero \cite[Lemma~4.9]{Clozel1990}.  Put
\[
 w=w(\pi):=-2\lambda_1,\qquad a_i:=\lambda_i-\lambda_1,
 \qquad 0=a_1<a_2<\cdots<a_k=w.
\]
Here $w(\pi)$ is the motivic weight of
\cite[Definition~8.2.7]{ChenevierLannes2019}.

Every cuspidal factor of a $J_0$-compatible inducing datum
(as in Subsection~\ref{subsec:laurent-coefficient-modules}) has distinct
parameter entries with integral differences,
because its shifted entries form a subset of $U_n$.
A twist by $|{\det}|^s$ therefore makes it regular
algebraic in Clozel's normalization.  Returning to the
unitary normalization places it under the convention above.
Thus both integral and strictly half-integral unitary
parameters occur in the arguments below.

We apply this convention to all factors of one inducing datum.
Fix $p$, $k_0$, and
\[
 t=(R,\sigma_t,\lambda_t)
 \in\mathscr T_{\mathfrak c,k_0,\tau,p}.
\]
Let
\[
 M_R=\prod_{\nu=1}^s\GL_{r_\nu},\qquad
 \sigma_t=\Sigma_1\boxtimes\cdots\boxtimes\Sigma_s,
\]
where each $\Sigma_\nu$ is a Speh representation of rank
$r_\nu$.  Put $z_\nu=(\lambda_t)_\nu$ in the block
coordinates of Subsection~\ref{subsec:laurent-coefficient-modules},
and $x_\nu=\operatorname{Re}z_\nu$.  Thus the $\nu$-th factor
in the normalized induction is
$\Sigma_\nu\otimes|{\det}|_{\A}^{z_\nu}$.
The representation $I_t(\lambda_t)$ has the
infinitesimal character of the trivial $G^1$-representation.
Its Harish--Chandra parameter therefore differs from that of the
trivial $\GL_n(\R)$-representation by a common scalar $z_c$.
Put $x_c=\operatorname{Re}z_c$, and let
$\{\lambda_{\nu,1},\ldots,\lambda_{\nu,r_\nu}\}$ be the
Harish--Chandra parameter of $\Sigma_{\nu,\infty}$, as in
Definition~\ref{def:weights}.  The twisted factor
$\Sigma_\nu\otimes|{\det}|_{\A}^{z_\nu}$ has parameter
$\{\lambda_{\nu,i}+z_\nu:1\le i\le r_\nu\}$.
After normalizing the full parameter to have coordinate sum zero
and shifting it by $(n-1)/2$, the $\nu$-th factor contributes
the set
\[
 B_\nu
 =\left\{\lambda_{\nu,i}+z_\nu-z_c+\frac{n-1}{2}:
            1\le i\le r_\nu\right\}.
\]
The infinitesimal-character condition gives
\[
 U_n=B_1\sqcup\cdots\sqcup B_s,
 \qquad |B_\nu|=r_\nu.
\]
For $\Sigma_\nu=J(\ell,\pi)$, its parameter consists of
$\lambda_i+(1-\ell)/2+j$, with $1\le i\le k$ and
$0\le j<\ell$.  Hence
\begin{equation}\label{eq:packet}
\begin{aligned}
 B_\nu&=\bigsqcup_{i=1}^k[s_0+a_i,s_0+a_i+\ell-1]\subset U_n,\\
 s_0&=z_\nu-z_c+\frac{n-\ell-w(\pi)}2\in\Z.
\end{aligned}
\end{equation}
The entries are distinct, so the segments in \eqref{eq:packet}
are disjoint.  Consequently,
\begin{equation}\label{eq:regularity}
                         a_{i+1}-a_i\ge\ell\qquad(1\le i<k).
\end{equation}

To compute the centre, write $\Sigma_\nu=J(\ell,\pi)$.
The Speh shifts $(\ell-1)/2,(\ell-3)/2,\ldots,-(\ell-1)/2$
preserve the symmetry of the cuspidal parameter about zero.
Thus the parameter of $\Sigma_\nu$ is real and symmetric about zero.

Since $B_\nu\subset U_n$ and the parameter of $\Sigma_\nu$
is real, the translation
$z_\nu-z_c$ is real and equals $x_\nu-x_c$.
Translating the preceding symmetry therefore gives
\begin{equation}\label{eq:block-exponent}
 \operatorname{c}(B_\nu)
 :=\frac{\min B_\nu+\max B_\nu}{2}
 =\frac{1}{r_\nu}\sum_{b\in B_\nu}b
 =x_\nu-x_c+\frac{n-1}{2}.
\end{equation}
If $\Sigma_\nu=J(m,1)$, this specializes to
\[
 B_\nu=p_\nu+U_m=[p_\nu,p_\nu+m-1],\qquad
 p_\nu=x_\nu-x_c+\frac{n-m}{2}\in\Z,
 \qquad 0\le p_\nu\le n-m.
\]
Thus the subset $B_\nu$ contributed by the $J(m,1)$-factor
at $\lambda_t$ is an integer segment.
Later hypotheses concern its position and complement in $U_n$;
nearby complex parameters need not give integer entries.
In particular,
$x_\nu-x_\mu=\operatorname{c}(B_\nu)-\operatorname{c}(B_\mu)$.
For the upper-triangular Borel, the dominant standard module has weakly
decreasing real exponents.  We therefore number the factors so that
\[
 \operatorname{c}(B_1)\ge\operatorname{c}(B_2)\ge\cdots
 \ge\operatorname{c}(B_s).
\]
Factors with equal centres may be placed in either order.
Proposition~\ref{prop:degree} shows that \eqref{eq:degree} is
unchanged by permuting such factors.

\begin{lemma}\label{lem:w11}
Every regular algebraic unitary cuspidal representation of
$\PGL_k(\A_\Q)$, with $k\ge2$, which is unramified at every finite
place has motivic weight, in the normalization of
Subsection~\ref{subsec:speh-weights},
\begin{equation}\label{eq:w11}
                              w\ge11.
\end{equation}
More precisely,
\begin{equation}\label{eq:rank-weight-lower-bound}
 w\ge
 \begin{cases}
 11,&k=2,\\
 22,&k=3,\\
 19,&k=4,\\
 23,&k\ge5.
 \end{cases}
\end{equation}
\end{lemma}

\begin{proof}
Chenevier--Lannes \cite[Theorem~F]{ChenevierLannes2019} classify the level-one
algebraic cuspidal representations of motivic weight at most $22$,
without a bound on the rank.  Their list contains
\[
 1,\quad \Delta_{11},\Delta_{15},\Delta_{17},\Delta_{19},\Delta_{21},
 \quad \operatorname{Sym}^2\Delta_{11}.
\]
Here $\Delta_w$ is the unitary cuspidal representation associated
with the normalized level-one eigenform of weight $w+1$.  The list
also contains four rank-four representations of motivic weights
$19,21,21,21$.  The symmetric square has rank three and
motivic weight $22$.  This proves both bounds.
\end{proof}

\begin{remark}\label{rem:rank-two-motivic-parity}
For $k=2$, put $w=w(\pi)$. The discrete series $\pi_\infty=D_w$
has lowest positive $\SO_2$-weight $w+1$. On this weight space
the rotation $-I_2$ acts by $(-1)^{w+1}$. Since $\pi$ has trivial
central character, $(-1)^{w+1}=1$, so $w$ is odd.
Thus for a shifted weight set $B=\{b,b+w\}\subset U_n$, its centre is
\[
 \operatorname{c}(B)=\frac{b+(b+w)}2=b+\frac w2\in\tfrac12+\Z.
\]
\end{remark}

\subsection{The exact degree formula}

Continue with $t=(R,\sigma_t,\lambda_t)$ and
$\sigma_t=\boxtimes_{\nu=1}^s\Sigma_\nu$ from the preceding
subsection.  We first record the cohomological reductions used
to compute the complementary degree.  Put
\[
\begin{aligned}
 K_{R,\infty}&=K_\infty\cap M_R(\R),
 &\mathfrak n_R&=\Lie(N_R)_\C,\\
 \mathfrak m_R^1&=
   \bigl(\Lie(M_R(\R))\cap\mathfrak g^1\bigr)_\C,
 &\mathfrak r_R^1&=\mathfrak m_R^1\ltimes\mathfrak n_R.
\end{aligned}
\]
For $M(t)$ from \eqref{eq:Mt}, the inducing module is
\[
 V=(V_{\sigma_{t,\infty}})_{K_{R,\infty}\text{-finite}}
       \otimes V_{\sigma_{t,f}}^{K_{R,f}}\otimes D_t.
\]
The group $K_{R,\infty}$ acts on the first factor, and
$\mathfrak n_R$ acts trivially.  For $X\in\mathfrak m_R^1$, set
$h_X(\lambda)=\langle\lambda+\rho_R,H_R(X)\rangle$, where
$H_R(X)$ is the complex-linear extension of the differential of
$H_R$ at the identity.  The action is
\[
\begin{aligned}
 X\cdot(v_\infty\otimes v_f\otimes\delta)
 =d\sigma_{t,\infty}(X)v_\infty\otimes v_f\otimes\delta
+v_\infty\otimes v_f\otimes h_X\delta, 
\end{aligned}
\]
where for $f\in\mathcal O_{\lambda_t}$, 
\[
(h_X\delta)(f)=\delta(h_Xf).
\]
The shift $\rho_R$ accounts for normalized induction; the action
on $D_t$ is the jet action \eqref{eq:g-action-jets}.

\begin{lemma}\label{lem:parabolic-cohomology}
We have 
\[
 H^q(\mathfrak g^1,K_\infty;M(t)\otimes E_\varepsilon)
 \cong
 H^q(\mathfrak r_R^1,K_{R,\infty};V\otimes E_\varepsilon).
\]
The latter group is computed by the relative
Hochschild--Serre spectral sequence
\begin{equation}\label{eq:levi-hochschild-serre}
\begin{aligned}
 E_2^{a,b}
 =H^a\!\left(\mathfrak m_R^1,K_{R,\infty};
       V\otimes H^b(\mathfrak n_R;E_\varepsilon)\right)\Longrightarrow
 H^{a+b}\!\left(\mathfrak r_R^1,K_{R,\infty};
       V\otimes E_\varepsilon\right).
\end{aligned}
\end{equation}
\end{lemma}

\begin{proof}
Apply Shapiro's lemma and the tensor identity for the
finite-dimensional coefficient $E_\varepsilon$
\cite[Chapter~I, Section~6]{BorelWallach2000}.
Since $\mathfrak n_R$ acts trivially on $V$, we have an
identity of complexes
\[
 C^\bullet(\mathfrak n_R;V\otimes E_\varepsilon)
 =V\otimes C^\bullet(\mathfrak n_R;E_\varepsilon).
\]
The relative Hochschild--Serre spectral sequence for
$\mathfrak r_R^1=\mathfrak m_R^1\ltimes\mathfrak n_R$
then gives \eqref{eq:levi-hochschild-serre}
\cite[Chapter~III, Sections~2--3]{BorelWallach2000}.
\end{proof}

Let $W^R$ be the minimal-length representatives of $W_G/W_{M_R}$,
where $W_G$ and $W_{M_R}$ are the Weyl groups of $G$ and $M_R$.
Write $\ell_W$ for Coxeter length, $\rho$ for the half-sum of
the positive roots of $G$, and $\Lambda_\varepsilon$ for the
highest weight of $E_\varepsilon$.

\begin{lemma}\label{lem:kostant-cohomology}
As an $M_R$-module,
\begin{equation}\label{eq:kostant-decomposition}
 H^b(\mathfrak n_R;E_\varepsilon)
 \cong
 \bigoplus_{\substack{v\in W^R\\ \ell_W(v)=b}}
 V^{M_R}_{\mu_v},
 \qquad
 \mu_v=v^{-1}(\Lambda_\varepsilon+\rho)-\rho,
\end{equation}
where $V^{M_R}_{\mu_v}$ has highest weight $\mu_v$.
Consequently, the $v$-summand of \eqref{eq:levi-hochschild-serre} is
\[
 H^a\!\left(\mathfrak m_R^1,K_{R,\infty};
       V\otimes V^{M_R}_{\mu_v}\right)
\]
in total degree $a+\ell_W(v)$.
\end{lemma}

\begin{proof}
The decomposition is Kostant's theorem
\cite[Chapter~III, Sections~2--3]{BorelWallach2000}.
Substitution into the $E_2$-term of
\eqref{eq:levi-hochschild-serre} gives the last assertion.
\end{proof}

For $\eta\in(\mathfrak a_R)^*_{\C}$, let $\C_\eta$ be the
one-dimensional representation of $M_R(\R)$ on which $g$ acts
by $e^{\langle\eta,H_R(g)\rangle}$.

\begin{lemma}\label{lem:finite-derivative-vanishing}
For the inducing module $V$ above and $v\in W^R$, put
\[
 V_0=\sigma_{t,\infty}\otimes V_{\sigma_{t,f}}^{K_{R,f}}
       \otimes\C_{\lambda_t+\rho_R}\otimes V^{M_R}_{\mu_v}.
\]
If $H^a(\mathfrak m_R^1,K_{R,\infty};V_0)=0$ for all $a$,
then
\[
 H^a(\mathfrak m_R^1,K_{R,\infty};V\otimes V^{M_R}_{\mu_v})=0
 \qquad\text{for all }a.
\]
In particular, nonzero cohomology requires
$\sigma_{t,\infty}\otimes\C_{\lambda_t+\rho_R}$ and
$(V^{M_R}_{\mu_v})^\vee$ to have the same infinitesimal
character on $\mathfrak m_R^1$.
\end{lemma}

\begin{proof}
Filter $D_t$ by distribution order,
\[
 D_t^{\le N}
 =\left(\mathcal O_{\lambda_t}/\mathfrak m_{\lambda_t}^{N+1}\right)^\vee,
 \qquad D_t^{\le-1}=0.
\]
The jet action \eqref{eq:g-action-jets} preserves this filtration.
Write $\cdot_{\lambda_t}$ for the specialized coefficient action.
For $|\alpha|=N$, the leading term is
\[
 X(f\otimes\partial^\alpha_{\lambda_t})
 \equiv(X\cdot_{\lambda_t}f)\otimes\partial^\alpha_{\lambda_t}
 \pmod{V_0\otimes D_t^{\le N-1}}.
\]
Thus every graded quotient is a finite direct sum of copies of $V_0$.
The long exact sequences give vanishing at each finite order.
Every relative cochain has finite distribution order, so the same
vanishing holds for their union.  If the two infinitesimal characters
differ, the required vanishing for $V_0$ follows from
\cite[Chapter~I, Sections~4--5]{BorelWallach2000}.
\end{proof}

Let $F$ be an irreducible algebraic representation of $\GL_n$,
with highest weight $\Lambda_F$, and let
\[
 \gamma=\Lambda_F+\rho_n,\qquad
 \lambda=-\frac1n\sum_{j=1}^n\gamma_j.
\]
In the notation of Definition~\ref{def:weights}, the Harish--Chandra
parameter of $F$ is $\{\gamma_1,\ldots,\gamma_n\}$, with
$\gamma_1>\cdots>\gamma_n$. For $F^\vee$ the parameter is
$\{-\gamma_j:1\le j\le n\}$. Its trace-zero representative on
$\mathfrak{sl}_n$ is $\{-\gamma_j-\lambda:1\le j\le n\}$.

For a standard parabolic $R\subseteq\GL_n$ with Levi $M_R$ and
nilradical $\mathfrak n_R$, the Kostant decomposition
\eqref{eq:kostant-decomposition} for $F$ gives
\[
 H^b(\mathfrak n_R;F)
 \cong\bigoplus_{\substack{v\in W^R\\ \ell_W(v)=b}}V_{\mu_v}^{M_R},
 \qquad
 \mu_v=v^{-1}(\Lambda_F+\rho_n)-\rho_n=v^{-1}\gamma-\rho_n.
\]

\begin{lemma}\label{lem:consecutive-kostant-character}
Let $R=P_{n_1,\ldots,n_s}\subseteq\GL_n$, and put
$m_b=\sum_{a<b}n_a$.  Let $U$ be an
$(\mathfrak m_R^1,K_{R,\infty})$-module whose generalized
infinitesimal character is
represented, on the successive Levi blocks, by
\[
 \eta_b=(\eta_{b,1}<\cdots<\eta_{b,n_b})\quad(1\le b\le s),
 \qquad \sum_{b,j}\eta_{b,j}=0.
\]
Assume that all entries are real and pairwise distinct.

For $v\in W^R$, let $\mu_v=v^{-1}\gamma-\rho_n$.
If, for some $a$, 
\[
 H^a\!\left(\mathfrak m_R^1,K_{R,\infty};
       U\otimes\C_{\rho_R}\otimes V_{\mu_v}^{M_R}\right)\ne0,
\]
then
\begin{equation}\label{eq:ordered-kostant-parameter-comparison}
 -\gamma_{v(m_b+j)}=\eta_{b,j}+\lambda
 \qquad(1\le b\le s,\ 1\le j\le n_b).
\end{equation}
These equalities determine at most one such $v$.
If $n_1=m$ and $\eta_{1,j}=\eta_{1,1}+j-1$ for
$1\le j\le m$, then the $\GL_m$-factor of this Kostant
summand is a determinant character.
\end{lemma}

\begin{proof}
On block $b$, let
\[
 \rho_{R,b}=\frac{n-n_b-2m_b}{2},\qquad
 \rho_{n,m_b+j}=\rho_{n_b,j}+\rho_{R,b}.
\]
Put $g=\operatorname{diag}(g_1,\ldots,g_s)$. For $b<c$,
$\operatorname{Ad}(g)$ acts on the $(b,c)$-block of $\mathfrak n_R$
by $X\mapsto g_bXg_c^{-1}$. Hence
\[
 \begin{aligned}
 \delta_R^{1/2}(g)
 &=\prod_{b<c}|\det g_b|^{n_c/2}|\det g_c|^{-n_b/2}\\
 &=\prod_{b=1}^s|\det g_b|^{(\sum_{c>b}n_c-\sum_{c<b}n_c)/2}
 =\prod_{b=1}^s|\det g_b|^{\rho_{R,b}}.
 \end{aligned}
\]
Thus the inducing parameter on block $b$ is
$(\eta_{b,j}+\rho_{R,b})_j$, while the $\GL_{n_b}$-factor of
$(V_{\mu_v}^{M_R})^\vee$ has parameter
\[
 -\mu_{v,m_b+j}-\rho_{n_b,j}
 =-\gamma_{v(m_b+j)}+\rho_{R,b}.
\]
Note that  $U\otimes\C_{\rho_R}$ has parameter
$(\eta_{b,j}+\rho_{R,b})_{b,j}$ and $(V_{\mu_v}^{M_R})^\vee$
has parameter $(-\gamma_{v(m_b+j)}+\rho_{R,b})_{b,j}$.
By the infinitesimal-character criterion used in the proof of
Lemma~\ref{lem:finite-derivative-vanishing}, nonzero
cohomology requires these two modules to have
the same infinitesimal character on
\[
 \mathfrak m_R^1
 =\{(X_1,\ldots,X_s):\textstyle\sum_b\operatorname{tr}X_b=0\}.
\]
After a permutation within each block, the two parameters differ by a
multiple of the total trace, since their restrictions to the trace-zero
diagonal Cartan agree. Cancelling $\rho_{R,b}$ therefore gives
one number $\lambda'$
such that
\[
 \{-\gamma_{v(m_b+j)}:1\le j\le n_b\}
 =\{\eta_{b,j}+\lambda':1\le j\le n_b\}
 \qquad(1\le b\le s).
\]
Summing all coordinates gives
$-\sum_j\gamma_j=n\lambda'$, hence $\lambda'=\lambda$.
Minimality of $v\in W^R$ makes
$v(m_b+1)<\cdots<v(m_b+n_b)$.
For fixed $b$, both $(-\gamma_{v(m_b+j)})_j$ and
$(\eta_{b,j}+\lambda)_j$ are strictly increasing. The preceding
equality of multisets therefore gives
\eqref{eq:ordered-kostant-parameter-comparison} coordinatewise.
Distinctness also makes $v$ unique.

If the first block is consecutive, then
$\gamma_{v(j)}=\gamma_{v(1)}-(j-1)$.  Consequently
\[
 \mu_{v,j}
 =\gamma_{v(1)}-(j-1)-\frac{n+1-2j}{2}
 =\gamma_{v(1)}-\frac{n-1}{2}
 \qquad(1\le j\le m).
\]
All highest-weight coordinates on this block are the same
integer $d$, so its coefficient is ${\det}^d$, also
when $m=1$.
\end{proof}

For a fixed $v$, write
$V^{M_R}_{\mu_v}=\boxtimes_{\nu=1}^s E_{\nu,v}$ and put
\[
\begin{gathered}
 \mathfrak g_\nu^1=\mathfrak{sl}_{r_\nu}(\R),\qquad
 \mathfrak p_\nu=
   (\mathfrak g_\nu^1/\mathfrak{so}(r_\nu))\otimes_\R\C,\\
 d_{r_\nu}=\dim_\C\mathfrak p_\nu
          =\frac{r_\nu(r_\nu+1)}2-1.
\end{gathered}
\]

\begin{lemma}\label{lem:levi-kunneth-duality}
The cohomology of the semisimple Levi factors is
\[
\begin{aligned}
 H^h\!\left(\bigoplus_{\nu=1}^s\mathfrak g_\nu^1,
       \prod_{\nu=1}^s\O_{r_\nu};
       \bigotimes_{\nu=1}^s
         (\Sigma_{\nu,\infty}\otimes E_{\nu,v})\right)\cong
 \bigoplus_{q_1+\cdots+q_s=h}
 \bigotimes_{\nu=1}^s
 H^{q_\nu}(\mathfrak g_\nu^1,\O_{r_\nu};
       \Sigma_{\nu,\infty}\otimes E_{\nu,v}).
\end{aligned}
\]
\end{lemma}

\begin{proof}
Splitting the exterior arguments among the Lie-algebra factors
identifies the relative complex with the tensor product of the
factor complexes.  Apply the K\"unneth formula
\cite[Chapter~I, Section~1.3, Equation~(2)]{BorelWallach2000}.
The twists $|{\det}|^{x_\nu}$ are trivial on
$(\mathfrak g_\nu^1,\O_{r_\nu})$, so they do not alter this calculation.
\end{proof}

Fix $v\in W^R$.  Since $\sigma_t$ has trivial central character,
the Lie algebra $(\mathfrak a_R^G)_\C$ of the relative split centre
acts on $V\otimes V^{M_R}_{\mu_v}$ through the following action
on $D_t$:
\[
 H\cdot\delta
 =\bigl(\langle\lambda+\rho_R,H\rangle+\mu_v(H)\bigr)\delta.
\]
The decompositions
\[
 \mathfrak m_R^1=(\mathfrak a_R^G)_\C
       \oplus\bigoplus_{\nu=1}^s(\mathfrak g_\nu^1)_\C,
 \qquad K_{R,\infty}=\prod_{\nu=1}^s\O_{r_\nu}
\]
give the full K\"unneth formula
\begin{equation}\label{eq:levi-kunneth-centre}
\begin{aligned}
 &H^h\!\left(\mathfrak m_R^1,K_{R,\infty};
       V\otimes V^{M_R}_{\mu_v}\right)\\
 &\quad\cong V_{\sigma_{t,f}}^{K_{R,f}}\otimes_\C
 \bigoplus_{\substack{a,q_1,\ldots,q_s\ge0\\
                      a+q_1+\cdots+q_s=h}}
 \Bigl(H^a\bigl((\mathfrak a_R^G)_\C;D_t\bigr)\otimes_\C
       \bigotimes_{\nu=1}^s
 H^{q_\nu}\!\left(\mathfrak g_\nu^1,\O_{r_\nu};
       \Sigma_{\nu,\infty}\otimes E_{\nu,v}\right)\Bigr).
\end{aligned}
\end{equation}
The finite-place space $V_{\sigma_{t,f}}^{K_{R,f}}$ has trivial
$(\mathfrak m_R^1,K_{R,\infty})$-action.  The next lemma computes
the relative split-centre factor in \eqref{eq:levi-kunneth-centre}.

\begin{lemma}\label{lem:split-centre-cohomology}
The Lie algebra cohomology
$H^a\bigl((\mathfrak a_R^G)_\C;D_t\bigr)$, computed with
this action on the coefficient module $D_t$, vanishes for
every $a\ge0$ unless
\[
 \langle\lambda_t+\rho_R,H\rangle+\mu_v(H)=0
 \qquad\text{for every }H\in(\mathfrak a_R^G)_\C.
\]
When this condition holds,
\[
 H^a\bigl((\mathfrak a_R^G)_\C;D_t\bigr)
 =
 \begin{cases}
 \C\,\operatorname{ev}_{\lambda_t},&a=0,\\
 0,&a>0,
 \end{cases}
 \qquad \operatorname{ev}_{\lambda_t}(f)=f(\lambda_t).
\]
Consequently, every nonzero summand in
\eqref{eq:levi-kunneth-centre} satisfies
\[
 a=0,\qquad h=q_1+\cdots+q_s.
\]
\end{lemma}

\begin{proof}
Choose a basis $H_1,\ldots,H_{s-1}$ of $(\mathfrak a_R^G)_\C$
with dual basis $\theta_1,\ldots,\theta_{s-1}$, and set
\[
 x_i(\lambda)=\langle\lambda+\rho_R,H_i\rangle+\mu_v(H_i).
\]
Since $(\mathfrak a_R^G)_\C$ is abelian, the
Chevalley--Eilenberg complex for the stated action and its
differential are \cite[Corollary~7.7.3, p.~240]{Weibel1994}
\[
\begin{gathered}
 C^\bullet\bigl((\mathfrak a_R^G)_\C;D_t\bigr)
 =
 D_t\otimes_\C\bigwedge^\bullet
          \bigl((\mathfrak a_R^G)_\C\bigr)^*,\\
 d(\delta\otimes\omega)
   =\sum_{j=1}^{s-1}x_j\delta\otimes(\theta_j\wedge\omega).
\end{gathered}
\]
If $x_i(\lambda_t)\ne0$, then $x_i^{-1}$ is a holomorphic germ.
Multiplication by $x_i^{-1}$ on $D_t$, defined by
$(x_i^{-1}\delta)(f)=\delta(x_i^{-1}f)$, is inverse to
multiplication by $x_i$.  Define a degree-$-1$ map by
\[
 h(\delta\otimes\omega)
   =(x_i^{-1}\delta)\otimes\iota_{H_i}\omega,
\]
where $\iota_{H_i}$ is contraction of alternating forms with $H_i$.
The graded Leibniz rule for contraction
\cite[Lemma~14.13(b), Equation~(14.11), p.~358]{Lee2012} gives
\[
 \iota_{H_i}(\theta_j\wedge\omega)
   +\theta_j\wedge\iota_{H_i}\omega=\theta_j(H_i)\omega
\]
and hence
\[
 (dh+hd)(\delta\otimes\omega)
 =\sum_{j=1}^{s-1}(x_jx_i^{-1}\delta)
       \otimes\theta_j(H_i)\omega
 =\delta\otimes\omega.
\]
Hence the complex is contractible, as is its tensor product
with the semisimple Levi complex.

Suppose now that every $x_i(\lambda_t)=0$.
Since the $H_i$ form a basis, the affine-linear functions $x_i$
are coordinates centred at $\lambda_t$.
Thus $D_t=\bigotimes_{i=1}^{s-1}\D_{x_i}$, and the
complex $C^\bullet\bigl((\mathfrak a_R^G)_\C;D_t\bigr)$
is the tensor product of
\[
 \D_{x_i}\xrightarrow{x_i}\D_{x_i}
 \qquad\text{in degrees }0,1.
\]
For $e_r(f)=f^{(r)}(0)/r!$, multiplication satisfies
$x_i e_r=e_{r-1}$ for $r\ge1$ and $x_i e_0=0$.
Each factor therefore has cohomology $\C e_0$ in degree zero
and zero in degree one.  Their tensor product has cohomology
$\C\,\operatorname{ev}_{\lambda_t}$ in degree zero.
For $s=1$, the centre is zero and the empty tensor product is
$\C$ in degree zero.
\end{proof}

\begin{proposition}\label{prop:speh-cohomology-bound}
Let $\Sigma_\nu=J(\ell,\pi_\nu)$ be a factor of $\sigma_t$, with
$\ell\ge1$ and $k=\operatorname{rank}(\pi_\nu)\ge2$.
For $v\in W^R$, let $E_{\nu,v}$ be the algebraic coefficient
on this factor, as in Lemma~\ref{lem:levi-kunneth-duality}.
If
\[
 H^{q_\nu}\!\left(\mathfrak{sl}_{k\ell}(\R),\O_{k\ell};
       \Sigma_{\nu,\infty}\otimes E_{\nu,v}\right)\ne0,
 \qquad \Sigma_{\nu,\infty}=J(\ell,\pi_{\nu,\infty}),
\]
then
\begin{equation}\label{eq:b}
 d_{k\ell}-q_\nu\ge b(k,\ell):=
 \begin{cases}
 \dfrac{k\ell(k\ell-\ell+1)}4,&k\ {\rm even},\\[5pt]
 \dfrac{k\ell(k\ell-\ell+1)-\ell}4,&k\ {\rm odd}.
 \end{cases}
\end{equation}
\end{proposition}

\begin{proof}
Put $h=\lfloor k/2\rfloor$.
The absolute determinant twists $|{\det}|^s$ used in
Subsection~\ref{subsec:speh-weights} act trivially on
$(\mathfrak{sl}_{k\ell},\O_{k\ell})$ and therefore do
not change the relative cohomology considered here.
In the unitary normalization,
\[
 \pi_{\nu,\infty}\cong
 \begin{cases}
 D_{m_1}\times\cdots\times D_{m_h},&k=2h,\\
 D_{m_1}\times\cdots\times D_{m_h}\times\operatorname{sgn}^{e},
     &k=2h+1,
 \end{cases}
\]
where $m_1>\cdots>m_h>0$, $D_m$ is the unitary discrete-series
representation of $\GL_2(\R)$ with Harish--Chandra parameter
$(m/2,-m/2)$, and $e\in\{0,1\}$
\cite[Proposition~8.2.13(i), p.~198]{ChenevierLannes2019}.
Hence
\[
 \Sigma_{\nu,\infty}\cong
 \begin{cases}
 \displaystyle\mathop{\times}_{i=1}^h J(\ell,D_{m_i}),&k=2h,\\[3pt]
 \displaystyle\left(\mathop{\times}_{i=1}^h J(\ell,D_{m_i})\right)
       \times(\operatorname{sgn}^{e}\circ{\det}),&k=2h+1.
 \end{cases}
\]
The Speh normalization is that of
\cite[Section~2, pp.~1521--1522]{NairPrasad2021};
the final character for odd $k$ is on $\GL_\ell(\R)$.
By \eqref{eq:regularity}, the length-$\ell$ strings in the
infinitesimal character are disjoint.
Put $N=k\ell$.  Apply
\cite[Section~1, p.~1516; Section~8, Theorem~6, p.~1541,
and Proposition~13, p.~1547]{NairPrasad2021} with
\[
 G=\SL_N,\qquad K=\SO_N,\qquad
 E=(E_{\nu,v}|_{\SL_N})^\vee.
\]
The irreducible constituents of
$\Sigma_{\nu,\infty}|_{\SL_N(\R)}$ with nonzero cohomology
with coefficients in $E_{\nu,v}|_{\SL_N}$ are modules
$A_{\mathfrak q^1}(\chi^1)$, where
\[
 \mathfrak q^1=\mathfrak q\cap\mathfrak{sl}_N(\C),
 \qquad L^1=\widetilde L\cap\SL_N(\R),
\]
$\mathfrak q=\mathfrak l\oplus\mathfrak u$ is stable
under the Cartan involution, $\widetilde L$ is its real
Levi subgroup with complexified Lie algebra $\mathfrak l$,
and $\chi^1$ is a
one-dimensional unitary representation of $L^1$.
Cohomological induction is taken in degree
$\dim(\mathfrak u\cap\mathfrak{so}_N(\C))$.
For the Speh factors above, the partition in
\cite[Section~9, Example~6, p.~1551]{NairPrasad2021} has
\[
 n_1=\cdots=n_h=\ell,\qquad
 n_0=
 \begin{cases}
 0,&k=2h,\\
 \ell,&k=2h+1.
 \end{cases}
\]
Thus we may take
\[
 \widetilde L\cong
 \begin{cases}
 \GL_\ell(\C)^h,&k=2h,\\
 \GL_\ell(\C)^h\times\GL_\ell(\R),&k=2h+1.
 \end{cases}
\]

Set
\[
 \mathfrak p_{\widetilde L}
 :=\{X\in\mathfrak l:X^{\mathsf T}=X\},
 \qquad
 r_{\mathfrak q}=\dim(\mathfrak u\cap\mathfrak p_\nu).
\]
On $(\mathfrak{sl}_N,\SO_N)$ with coefficient
$E_{\nu,v}|_{\SL_N}$, the Vogan--Zuckerman formula
\cite[Theorem~5.5]{VoganZuckerman1984}, in the form of
\cite[Section~8, Proposition~13(1), p.~1547]{NairPrasad2021}, gives
\[
 q_\nu=r_{\mathfrak q}+t,\qquad
 0\le t\le\dim\mathfrak p_{\widetilde L}-1.
\]
Here the subtraction by $1$ removes the scalar line
$\C I_{k\ell}\subset\mathfrak p_{\widetilde L}$, leaving
the trace-zero subspace.
Taking invariants under $\O_N/\SO_N$ gives the factor in
Lemma~\ref{lem:levi-kunneth-duality} and preserves this
degree bound for either determinant parity.

Put $\mathfrak p=\{X\in\mathfrak{gl}_{k\ell}(\C):X^{\mathsf T}=X\}$,
whose trace-zero subspace identifies with $\mathfrak p_\nu$.
Taking the $(-1)$-eigenspaces of the Cartan involution in
$\mathfrak{gl}_{k\ell}(\C)
=\overline{\mathfrak u}\oplus\mathfrak l\oplus\mathfrak u$ gives
\[
 \mathfrak p=(\overline{\mathfrak u}\cap\mathfrak p)
       \oplus\mathfrak p_{\widetilde L}
       \oplus(\mathfrak u\cap\mathfrak p).
\]
Here the bar denotes complex conjugation with respect to
$\mathfrak{gl}_{k\ell}(\R)$.
Conjugation exchanges the two outer summands, each of dimension
$r_{\mathfrak q}$, since the nilradicals consist of trace-zero matrices.
As $\dim\mathfrak p=d_{k\ell}+1$, we obtain
\begin{equation}
\label{eqn: dim-formula}
 d_{k\ell}=2r_{\mathfrak q}+\dim\mathfrak p_{\widetilde L}-1.
\end{equation}
Therefore
\[
 d_{k\ell}-q_\nu
 =r_{\mathfrak q}+\dim\mathfrak p_{\widetilde L}-1-t
 \ge r_{\mathfrak q}.
\]
Finally, by \eqref{eqn: dim-formula}
\[
 r_{\mathfrak q}
 =\frac12\left(\frac{k\ell(k\ell+1)}2
                     -\dim\mathfrak p_{\widetilde L}\right).
\]
Here
\[
 \dim\mathfrak p_{\widetilde L}=
 \begin{cases}
 h\ell^2,&k=2h,\\
 h\ell^2+\ell(\ell+1)/2,&k=2h+1.
 \end{cases}
\]
Substituting $h=k/2$ or $(k-1)/2$ gives
\[
 b(k,\ell):=r_{\mathfrak q}=
 \begin{cases}
 \dfrac{k\ell(k\ell-\ell+1)}4,&k\ {\rm even},\\[4pt]
 \dfrac{k\ell(k\ell-\ell+1)-\ell}4,&k\ {\rm odd}.
 \end{cases}
\]
This proves \eqref{eq:b}.
\end{proof}

Fix a nonzero class
$\alpha\in H^q(\mathfrak g^1,K_\infty;M(t)\otimes E_\varepsilon)$,
and put $j=d_n-q$.
Via Shapiro's isomorphism in Lemma~\ref{lem:parabolic-cohomology},
regard $\alpha$ as a class in the abutment of
\eqref{eq:levi-hochschild-serre}.  Its nonzero associated graded
image lies in $E_\infty^{a,b}$ for some $a+b=q$.
Choose a representative of this image in $E_2^{a,b}$.
By \eqref{eq:kostant-decomposition}, choose $v\in W^R$ with
$\ell_W(v)=b$ such that this representative has a nonzero component in
\[
 H^a\!\left(\mathfrak m_R^1,K_{R,\infty};
       V\otimes V^{M_R}_{\mu_v}\right).
\]
Fix a nonzero component of this class in the K\"unneth decomposition
\eqref{eq:levi-kunneth-centre}, with semisimple degrees
$(q_1,\ldots,q_s)$.
Lemma~\ref{lem:split-centre-cohomology} gives $a=\sum_\nu q_\nu$.

For subsets $A,B$ contributed by distinct inducing factors,
ordered so that $\operatorname{c}(A)\ge\operatorname{c}(B)$, set
\[
 \kappa(A,B)=|A||B|-\#\{(a,b)\in A\times B:a<b\}.
\]

\begin{proposition}\label{prop:degree}
The complementary degree $j=d_n-q$ of $\alpha$ satisfies
\begin{equation}\label{eq:degree}
 j=\sum_\nu(d_{r_\nu}-q_\nu)+
       \sum_{\mu<\nu}
       \bigl(r_\mu r_\nu-\kappa(B_\mu,B_\nu)\bigr)+(s-1).
\end{equation}
\end{proposition}

\begin{proof}
Since $a+b=q$, $b=\ell_W(v)$, and $a=\sum_\nu q_\nu$, we have
\[
 q=\ell_W(v)+\sum_\nu q_\nu.
\]

Let $h$ be the Harish--Chandra parameter of
$\sigma_{t,\infty}\otimes\C_{\lambda_t}$, normalized to have
coordinate sum zero.  Nonzero Levi cohomology requires the
infinitesimal character of
$\sigma_{t,\infty}\otimes\C_{\lambda_t+\rho_R}$ to equal that
of $(V^{M_R}_{\mu_v})^\vee$; see
Lemma~\ref{lem:finite-derivative-vanishing}.
The highest-weight formula
\cite[Chapter~V, Section~5, Equation~(5.43), and
Theorem~5.62]{Knapp2002} and
\eqref{eq:kostant-decomposition} give
\[
 \mu_v+\rho_{M_R}
 =\Lambda_\varepsilon+v^{-1}\rho-\rho_R.
\]
After normalizing coordinate sums to zero, the dual parameter
is $-v^{-1}\rho+\rho_R$.  The character comparison therefore
gives $h+\rho_R=-v^{-1}\rho+\rho_R$, up to permutations within
each Levi block.  Choose the representative
\begin{equation}\label{eq:parameter-kostant-sign}
 h_i=-\rho_{v(i)},\qquad
 h_i+\frac{n-1}{2}=v(i)-1.
\end{equation}
Thus the entries $v(i)-1$ on block $\nu$ are exactly $B_\nu$.

The length of a permutation is its number of inversions
\cite[Chapter~II, Section~6, Proposition~2.70,
pp.~168--169]{Knapp2002}.
Since $v\in W^R$ is increasing within each block, its inversions
are precisely the pairs
$(a,b)\in B_\mu\times B_\nu$, $\mu<\nu$, with $a>b$.
Consequently,
\[
\begin{aligned}
 \dim\mathfrak n_R&=\sum_{\mu<\nu}r_\mu r_\nu,\\
 \ell_W(v)&=\sum_{\mu<\nu}\kappa(B_\mu,B_\nu),\\
 \dim\mathfrak n_R-\ell_W(v)
 &=\sum_{\mu<\nu}
       \bigl(r_\mu r_\nu-\kappa(B_\mu,B_\nu)\bigr).
\end{aligned}
\]
Using
$d_n=\sum_\nu d_{r_\nu}+\dim\mathfrak n_R+(s-1)$ now gives
\[
\begin{aligned}
 j=d_n-q
 &=\sum_\nu(d_{r_\nu}-q_\nu)
       +\dim\mathfrak n_R-\ell_W(v)+(s-1)\\
 &=\sum_\nu(d_{r_\nu}-q_\nu)
       +\sum_{\mu<\nu}
          \bigl(r_\mu r_\nu-\kappa(B_\mu,B_\nu)\bigr)+(s-1),
\end{aligned}
\]
as required.

\end{proof}

\section{Constant terms and Eisenstein coefficient modules}
\label{sec:eisenstein-coefficients}

We use full constant terms to separate groups of cuspidal parameters.
We first fix the intertwining and residue conventions needed to
identify the resulting coefficient maps.

\subsection{Constant terms and normalized intertwining operators}

Write $P_{n_1,\ldots,n_s}$ for the standard upper block-triangular
parabolic with Levi $\prod_i\GL_{n_i}$.

For standard parabolics $R\subset P$ and an automorphic form
$\varphi$ on $N_P(\A)M_P(\Q)\backslash G(\A)$, write
\[
 (C_{P,R}\varphi)(g)
 :=\int_{N_R(\Q)\backslash N_R(\A)}\varphi(ng)\,dn.
\]
Let $C_{P,R}^{\mathrm{cusp}}\varphi$ denote its cuspidal
projection on $M_R$, as in
Definition~\ref{def:cuspidal-constant-term-component}; these are
the unnormalized conventions of
\cite[Section~6.2]{BernsteinLapid2024}.

For standard parabolics $P,Q$, put
\[
 W^{\supset Q}(P)
 =\{w\in{}_{Q}W_{P}:M_Q\subset wM_Pw^{-1}\}.
\]
Here ${}_{Q}W_{P}$ consists of the minimal-length representatives
of $W_Q\backslash W/W_P$.  For $w\in W^{\supset Q}(P)$,
let $P_w\subset P$ be the standard parabolic with Levi
$M_{P_w}=w^{-1}M_Qw$.  The global unnormalized standard
intertwining operator $M(w,\lambda)$ maps normalized induction
from $P_w$ to normalized induction from $Q$.  For a function $h$ on
$N_Q(\A)M_Q(\Q)\backslash G(\A)$, write
$[h]_\nu(g)=e^{\langle\nu,H_Q(g)\rangle}h(g)$.

\begin{proposition}
\label{prop:bernstein-lapid-constant-term}
Let $E(\varphi,\lambda)$ be the Eisenstein series induced from
$P$.  For every standard parabolic $Q$,
\[
 C_{G,Q}^{\mathrm{cusp}}E(\varphi,\lambda)
 =\sum_{w\in W^{\supset Q}(P)}
  \left[M(w,\lambda)
        (C_{P,P_w}^{\mathrm{cusp}}\varphi)\right]_{w\lambda}.
\]
The equality holds meromorphically in $\lambda$.
\end{proposition}

\begin{proof}
This is \cite[Equation~(6.18), p.~208]{BernsteinLapid2024}.
\end{proof}

We use the standard and normalized intertwining operators defined in
\cite[Section~I.1, p.~607, and Section~II.1, pp.~643--644]{MoeglinWaldspurger1989}.
For $\sigma=\boxtimes_{j=1}^r\sigma_j$, a representation of a
Levi subgroup over $\Q_v$ or $\A$, $t\in\C^r$, and
$u\in S_r$, write
\begin{equation}\label{eq:normalized-intertwiner-notation}
 N(u,\sigma,t):
 \mathop{\times}_{j=1}^r\sigma_j|{\det}|^{t_j}
 \longrightarrow
 \mathop{\times}_{j=1}^r
       \sigma_{u^{-1}(j)}|{\det}|^{t_{u^{-1}(j)}}.
\end{equation}
Denote by $M(u,\sigma,t)$  the corresponding unnormalized
operator and by $r(u,\sigma,t)$ its normalizing scalar, so
\[
 M(u,\sigma,t)=r(u,\sigma,t)N(u,\sigma,t).
\]
We call $t=(t_1,\ldots,t_r)$ the \emph{exponent vector}
relative to the ordered factors $(\sigma_1,\ldots,\sigma_r)$:
the $j$-th factor is twisted by the determinant character
$|{\det}|^{t_j}$.
For $u\in S_r$, define $(ut)_j=t_{u^{-1}(j)}$; thus the
target has exponent vector $ut$ relative to the permuted
factors $(\sigma_{u^{-1}(1)},\ldots,\sigma_{u^{-1}(r)})$.
For adelic $\sigma$, the normalized global operator is
\[
 N(u,\sigma,t)=\bigotimes_v' N(u,\sigma_v,t).
\]
Thus the inducing representation distinguishes the local and
global operators.

For general $w\in{}_{Q}W_{P}$, define the standard parabolics
$P_w\subset P$ and $Q_w\subset Q$ by
\[
 M_{P_w}=M_P\cap w^{-1}M_Qw,\qquad
 M_{Q_w}=M_Q\cap wM_Pw^{-1}.
\]
The operator $M(w,\lambda)$ now intertwines from $P_w$ to $Q_w$,
with $\lambda$ restricted to the Levi of $P_w$.  Write
\[
 E_{Q_w}^{Q}(\psi,\nu)(g)
 :=\sum_{\gamma\in Q_w(\Q)\backslash Q(\Q)}
     e^{\langle\nu,H_{Q_w}(\gamma g)\rangle}\psi(\gamma g).
\]
For $m\in M_Q(\A)$ and $k\in K_\infty K_f$, the full normalized
constant-term formula is
\begin{equation}\label{eq:full-normalized-constant-term}
 \CT_Q E(\varphi,\lambda)(mk)
 =\delta_Q(m)^{-1/2}
   \sum_{w\in{}_{Q}W_{P}}
   E_{Q_w}^{Q}\!\left(M(w,\lambda)(C_{P,P_w}\varphi),w\lambda\right)(mk).
\end{equation}
This is \cite[Section~6.10, Equation~(6.17)]{BernsteinLapid2024}
with our normalization; the identity holds meromorphically in $\lambda$.
We project onto the prescribed cuspidal support and generalized
infinitesimal character on each target Levi block, retaining the
polynomial dependence on its split centre.
Proposition~\ref{prop:bernstein-lapid-constant-term} is the cuspidal
projection of the unnormalized formula.

At every finite place, the normalized operator sends the normalized
spherical vector to the normalized spherical vector.  All finite-place
operators below are restricted to these lines.  At infinity, we use
the $K_\infty$-finite compact realization of
Definition~\ref{def:compact-realization}.

\begin{lemma}\label{lem:cuspidal-L-nonvanishing}
Let $\pi$ be a unitary cuspidal automorphic representation of
$\GL_k(\A)$, with $k>1$.  The completed $L$-function
$\Lambda(t,\pi^\vee)$ is entire and nonzero for
$\Re t\ge1$ or $\Re t\le0$.  In particular,
\begin{equation}
\label{eqn: non-vanishing-results-L}
 t\in\tfrac12\Z\setminus\{\tfrac12\}
 \quad\Longrightarrow\quad \Lambda(t,\pi^\vee)\ne0.
\end{equation}
For $a,b\in\tfrac12\Z\setminus\{\tfrac12\}$, the ratio
\[
                 \frac{\Lambda(a+h,\pi^\vee)}
                      {\Lambda(b+h,\pi^\vee)}
\]
and its reciprocal are holomorphic and nonzero near $h=0$.
Moreover, $\Lambda(t,\pi)$ and $\Lambda(t,\pi^\vee)$ have the
same finite order of vanishing at $t=1/2$.
\end{lemma}

\begin{proof}
Entireness follows from Godement--Jacquet
\cite[Theorem~13.8]{GodementJacquet1972}.
For $\Re t>1$, absolute convergence of the finite Euler product
and the nonvanishing archimedean factors give
$\Lambda(t,\pi^\vee)\ne0$.
On the line $\Re t=1$, nonvanishing follows from
Jacquet--Shalika \cite[Theorem~(1.3)]{JacquetShalika1976}.
Applying the same results to $\pi$ and using the functional
equation gives $\Lambda(t,\pi^\vee)\ne0$ for $\Re t\le0$.
This proves \eqref{eqn: non-vanishing-results-L} and the assertion
about ratios.  The functional equation also gives
\[
 \Lambda(\tfrac12+h,\pi^\vee)
 =\epsilon(\tfrac12+h,\pi^\vee)\Lambda(\tfrac12-h,\pi),
\]
where the epsilon factor is holomorphic and nonzero.
Thus the two orders of vanishing coincide; they are finite
because these entire functions are not identically zero.
\end{proof}

\subsection{Residues of Eisenstein series and character quotients}
\label{subsec:speh-residue-quotients}

Let $d\ge1$, and let $\tau$ be an irreducible unitary
cuspidal automorphic representation of $\GL_r(\A)$.
Fix a realization of $\tau$ by cuspidal automorphic forms.
Let $S=M_SN_S=P_{r,\ldots,r}\subset\GL_{dr}$, with
$M_S=\GL_r^d$.  For $t=(t_1,\ldots,t_d)$ with $\sum_i t_i=0$,
choose a holomorphic family of sections
\[
 F_t\in\tau|{\det}|^{t_1}\times\cdots\times\tau|{\det}|^{t_d}
\]
of normalized induction.  Here holomorphy means that
$t\mapsto F_t|_K$, where $K=K_\infty K_f$, is holomorphic with
values in $\mathcal I_{S,K}(\tau^{\boxtimes d})$ from
Definition~\ref{def:compact-realization}\textup{(ii)}.
For an Iwasawa decomposition $g=umk$, the section is recovered from
its restriction to $K$ by
\[
 F_t(umk)=\delta_S(m)^{1/2}
       \prod_{i=1}^d|{\det} m_i|_\A^{t_i}
       \tau^{\boxtimes d}(m)F_t(k),
\]
for $u\in N_S(\A)$, $m=(m_1,\ldots,m_d)\in M_S(\A)$,
and $k\in K$; see Lemma~\ref{lem:compact-realization}.
Each value $F_t(g)$ is an automorphic form on $M_S$ in the chosen
realization of $\tau^{\boxtimes d}$, with the Levi subgroup acting
by right translation.
For $\Re(t_i-t_{i+1})$ sufficiently large, define
\begin{equation}\label{eq:speh-eisenstein-family}
 \mathcal E(t)(g)
 :=\sum_{\gamma\in S(\Q)\backslash\GL_{dr}(\Q)}
          [F_t(\gamma g)](1_{M_S}).
\end{equation}
Here $1_{M_S}$ is the identity of $M_S$, so each summand is
the value of the automorphic form $F_t(\gamma g)$ at that identity.
The sum converges absolutely in this region; its meromorphic
continuation defines $\mathcal E(t)$ elsewhere
\cite[Sections~2.2--2.3]{BernsteinLapid2024}.
Define $\Phi$ by taking successive residues of
$\mathcal E(t)$ along $t_i-t_{i+1}=1$, in the order
$i=d-1,\ldots,1$.  For $d=1$, the sum has one term and
$\mathcal E(0)(g)=[F_0(g)](1)=F_0(1)(g)$; set
$\Phi=\mathcal E(0)\in\tau$.
Denote by $\mathcal R(d,\tau)$ the complex linear span
of the forms $\Phi$ obtained as $t\mapsto F_t$ ranges over all
holomorphic families of sections above.
Thus $\mathcal R(d,\tau)$ is a space of automorphic
forms on $\GL_{dr}(\Q)\backslash\GL_{dr}(\A)$.
In particular, $\mathcal R(1,\tau)=\tau$ in the chosen realization.

For 
$h=C_{\GL_{dr},S}^{\mathrm{cusp}}\Phi$ with 
$\Phi \in \mathcal R(d,\tau)$
 and each
$g\in\GL_{dr}(\A)$, define a function $f_h(g)$ on
the Levi subgroup $M_S(\A)$ by
\begin{equation}\label{eq:constant-term-induced-section}
 f_h(g)(m):=\delta_S(m)^{-1/2}h(mg),
 \qquad g\in\GL_{dr}(\A),\quad m\in M_S(\A).
\end{equation}
Since $h$ is left $N_S(\A)$-invariant, for
$m_0\in M_S(\A)$ and $n\in N_S(\A)$ one has
\[
 f_h(m_0ng)(m)
 =\delta_S(m_0)^{1/2}f_h(g)(mm_0).
\]
This is the covariance relation for normalized induction,
with right translation in $m$ as the Levi action
\cite[Sections~I.1 and II.1]{MoeglinWaldspurger1989}.
Since $f_h(g)(1)=h(g)$, this realization retains the full
scalar section and commutes with right translation.

\begin{lemma}\label{lem:speh-residue-constant-term}
For every $\Phi\in\mathcal R(d,\tau)$, put
$h=C_{\GL_{dr},S}^{\mathrm{cusp}}\Phi$.  Then
\begin{equation}\label{eq:speh-levi-coefficient-condition}
 f_h(g)\in\boxtimes_{i=1}^d\tau|{\det}|^{(2i-d-1)/2}
 \qquad\text{for every }g\in\GL_{dr}(\A).
\end{equation}
Under the realization \eqref{eq:constant-term-induced-section},
the constant-term map is a nonzero
$\GL_{dr}(\A)$-equivariant injection
\[
C_{\GL_{dr},S}^{\mathrm{cusp}}:  \mathcal R(d,\tau)
 \hookrightarrow
 \tau|{\det}|^{(1-d)/2}\times\cdots\times
 \tau|{\det}|^{(d-1)/2}.
\]
The image is the Speh subrepresentation under right translation.
Consequently, $\mathcal R(d,\tau)\simeq J(d,\tau)$
as $\GL_{dr}(\A)$-representations.
\end{lemma}

\begin{proof}
For the inducing family $F_t$ above, put
\[
 t^0=\left(\frac{d+1-2i}{2}\right)_{i=1}^d,
 \qquad u_0(i)=d+1-i\quad(1\le i\le d).
\]
In the notation \eqref{eq:normalized-intertwiner-notation},
the normalized global operator for the $d$ copies of $\tau$
is $N(u_0,\tau^{\boxtimes d},t)$.
Taking successive residues along $t_i-t_{i+1}=1$ in the constant-term
formula of Proposition~\ref{prop:bernstein-lapid-constant-term}, as in
\cite[Section~II.2, p.~646, and Section~II.3, p.~649]{MoeglinWaldspurger1989},
gives
\begin{equation}\label{eq:speh-residue-constant-term}
 (C_{\GL_{dr},S}^{\mathrm{cusp}}\Phi)(g)
 =c_{d,\tau}\bigl[N(u_0,\tau^{\boxtimes d},t^0)F_{t^0}\bigr](g)(1_{M_S}),
 \qquad c_{d,\tau}\ne0.
\end{equation}
For $d=1$, no residue is taken and $c_{1,\tau}=1$.
By \cite[Section~II.3, pp.~647--649]{MoeglinWaldspurger1989},
$\mathcal R(d,\tau)\simeq J(d,\tau)$, and the image of
$N(u_0,\tau^{\boxtimes d},t^0)$ is the Speh subrepresentation of the induction
with exponent vector $u_0t^0$.
Since $(u_0t^0)_i=(2i-d-1)/2$, the residue identity gives
\eqref{eq:speh-levi-coefficient-condition}; the realization
\eqref{eq:constant-term-induced-section} then gives the
stated target.  As $F_{t^0}$ ranges over the inducing space,
the image is the Speh subrepresentation.
The constant-term map commutes with right translation and
is nonzero.  Its domain $\mathcal R(d,\tau)\simeq J(d,\tau)$
is irreducible, so the map is injective.
\end{proof}

\begin{lemma}\label{lem:adjacent-character-quotient}
The normalized operator
\[
 N\bigl((12),1\boxtimes1,(\tfrac12,-\tfrac12)\bigr)
\]
has image isomorphic to $J(2,1)=\mathbf1_{\GL_2(\A)}$.
After identifying its image with this trivial representation,
denote the resulting quotient by
\begin{equation}\label{eq:global-residue-quotient}
 R_{1,1}:|\cdot|_\A^{1/2}\times|\cdot|_\A^{-1/2}
             \twoheadrightarrow\mathbf1_{\GL_2(\A)}.
\end{equation}
Its kernel equals the kernel of the automorphic residue map.
A common twist of the two
characters by $|\cdot|_\A^s$ gives the quotient
$R_{1,1}\otimes|{\det}|_\A^s$; its kernel in the compact realization
is independent of $s$.
\end{lemma}

\begin{proof}
Apply Lemma~\ref{lem:speh-residue-constant-term} with $d=2$ and
$\tau=1$.  Its proof identifies the normalized operator's image
with $J(2,1)$.  Formula~\eqref{eq:speh-residue-constant-term}
and injectivity of the cuspidal constant-term map on
$\mathcal R(2,1)$ identify the two kernels.
A common determinant twist is trivial on the compact subgroup,
so the twisted quotients have the same compact realization and kernel.
\end{proof}

\subsection{Holomorphy and simultaneous residues}
\label{subsec:simultaneous-residues}

We locate the possible polar hyperplanes of Eisenstein series induced
from cuspidal automorphic representations and prove that residues
along these hyperplanes commute.  Let
\begin{equation}
\label{eqn: inducing-data-general-Levi}
 P=P_{r_1,\ldots,r_h},\qquad
 \tau=\boxtimes_{i=1}^h\tau_i,\qquad
 I_P(t)=\mathop{\times}_{i=1}^h\tau_i|{\det}|^{t_i},
 \qquad \sum_i r_it_i=0,
\end{equation}
where the $\tau_i$ are unitary cuspidal representations with
trivial central character and nonzero hyperspecial fixed
vectors at every finite place.  This includes the character
$1$ on $\GL_1$.  Fix a real point $t^0$ with
$t_1^0\ge\cdots\ge t_h^0$, and assume that the
Harish--Chandra parameter of $I_P(t^0)_\infty$ has distinct
entries.

Fix automorphic realizations of the $\tau_i$, and put
$G=\GL_n$, where $n=\sum_i r_i$.  A section $f\in I_P(t)$
takes values in the resulting space of automorphic forms
on $M_P=\prod_i\GL_{r_i}$.  For $\Re(t_i-t_{i+1})$
sufficiently large, define the linear operator
\[
\begin{aligned}
 E_P(t):I_P(t)&\longrightarrow C^\infty(G(\Q)\backslash G(\A)),\\
 (E_P(t)f)(g)
 &:=\sum_{\gamma\in P(\Q)\backslash G(\Q)}
                         [f(\gamma g)](1_{M_P}).
\end{aligned}
\]
Here $1_{M_P}$ is the identity of $M_P$ and each summand evaluates
the automorphic form $f(\gamma g)$ at that identity.
The sum converges absolutely in this region and defines $E_P(t)$
elsewhere by meromorphic continuation, as in
\eqref{eq:speh-eisenstein-family} of
Subsection~\ref{subsec:speh-residue-quotients}.
For the data of that construction, $\mathcal E(t)=E_S(t)F_t$.

All holomorphic sections below are
$K_f$-fixed and locally contained in a fixed finite sum of
$K_\infty$-types in the compact realization
of Definition~\ref{def:compact-realization}.

Consider
\begin{equation}\label{eq:cuspidal-polar-pairs}
 \mathcal H=\{(i,j):i<j,\ \tau_i\simeq\tau_j,\
                              t_i^0-t_j^0=1\},\qquad
 u_{ij}=t_i-t_j-1.
\end{equation}

\begin{proposition}\label{prop:cuspidal-normal-crossings}
The functions $u_e$, $e\in\mathcal H$, extend to coordinates
$(u,y)$, linear in $t-t^0$, on
$\sum_i r_it_i=0$.  For every holomorphic section $f(t)$,
\[
       \left(\prod_{e\in\mathcal H}u_e\right)E_P(t)f(t)
\]
is holomorphic near $t^0$.
Consequently, for every $S\subseteq\mathcal H$, the residue
\[
 \operatorname{Res}_S(E_Pf)
 :=\left.\left(\prod_{e\in S}u_e\right)E_P(t)f(t)
                                      \right|_{u_e=0,\ e\in S}
\]
is a well-defined meromorphic germ in the remaining
coordinates.  Taking these residues in any order gives the
same germ.  Here each residue means the coefficient of
$u_e^{-1}$, without an orientation sign.
\end{proposition}

\begin{proof}
Regularity of the total Harish--Chandra parameter implies
\[
 \tau_i\simeq\tau_j,\quad i\ne j
       \quad\Longrightarrow\quad t_i^0\ne t_j^0.
\]
For a fixed index $i$, there is at most one $j$ with
$(i,j)\in\mathcal H$: two distinct such indices $j,k$ would give
$\tau_j\simeq\tau_k$ and $t_j^0=t_k^0=t_i^0-1$, contradicting
the preceding implication.  Similarly, there is at most one $j$
with $(j,i)\in\mathcal H$.  Since for each pair $(i,j)\in \mathcal H$ we have $i<j$,
the pairs form disjoint paths of indices.
To extend the $u_{ij}$ to coordinates on the hyperplane
$\sum_i r_it_i=0$, it suffices to prove that their differentials
$du_{ij}=dt_i-dt_j$ are linearly independent on its tangent space.
A linear form vanishing on this tangent space is a scalar multiple
of $\sum_i r_i\,dt_i$.  Thus a possible relation among the
restricted differentials has the form
\begin{equation}
\label{eqn: linear-dependent-differ}
 \sum_{(i,j)\in\mathcal H}c_{ij}(dt_i-dt_j)
                       =c\sum_i r_i\,dt_i.
\end{equation}
Let $\partial t_1,\ldots,\partial t_h$ be the basis of the
ambient space $\C^h$ dual to $dt_1,\ldots,dt_h$.
Evaluating \eqref{eqn: linear-dependent-differ} on
$\sum_{i=1}^h \partial t_i$ gives
$0=c\sum_{i=1}^h r_i=cn$, hence $c=0$.
At an endpoint of each path the remaining identity gives the
coefficient of its unique edge as zero.  Removing that edge and
continuing gives $c_{ij}=0$ for every pair.
Since $u_{ij}(t^0)=0$, the $u_{ij}$ extend to the asserted coordinates.

Let $\mathcal O_{t^0}$ be the ring of holomorphic germs on
$\sum_i r_it_i=0$ at $t^0$.  For a block permutation $w$, let
\[
 \operatorname{Inv}(w):=\{(i,j):i<j,\ w(i)>w(j)\}.
\]
The standard operator has the normalization
\cite[Sections~I.1 and II.1]{MoeglinWaldspurger1989}
\[
\begin{aligned}
 M(w,\tau,t)&=r(w,\tau,t)N(w,\tau,t),\\
 r(w,\tau,t)&=a_w(t)
 \prod_{(i,j)\in\operatorname{Inv}(w)}
 \frac{\Lambda(t_i-t_j,\tau_i\times\tau_j^\vee)}
      {\Lambda(1+t_i-t_j,\tau_i\times\tau_j^\vee)},
\end{aligned}
\]
where $a_w\in\mathcal O_{t^0}^{\times}$ and $\Lambda$ denotes the
completed Rankin--Selberg $L$-function, including the local factors
at all archimedean places; see
\cite{JacquetPiatetskiShapiroShalika1983,JacquetShalika1990}.
The trivial central characters exclude nontrivial determinant twists:
\[
 \tau_i\simeq\tau_j|{\det}|^z
 \quad\Longrightarrow\quad
 x^{r_i z}=1\ (x>0)
 \quad\Longrightarrow\quad z=0.
\]
For $i<j$, we have $t_i^0-t_j^0\ge0$.
The analytic properties in
\cite[Section~II.2, properties~(1)--(3), p.~646]{MoeglinWaldspurger1989}
therefore give
\[
 \Lambda(1+t_i-t_j,\tau_i\times\tau_j^\vee)
                                      \in\mathcal O_{t^0}^{\times}, \qquad \begin{cases}
 u_{ij}\Lambda(t_i-t_j,\tau_i\times\tau_j^\vee)
                  \in\mathcal O_{t^0}^{\times},&(i,j)\in\mathcal H,\\
 \Lambda(t_i-t_j,\tau_i\times\tau_j^\vee)
                  \in\mathcal O_{t^0},&(i,j)\notin\mathcal H.
 \end{cases}
\]
In fact, a pole at  $t_i^0-t_j^0=0$ or $1$ requires $\tau_i\simeq\tau_j$.
For isomorphic pairs, the case $t_i^0=t_j^0$ is excluded at the beginning of the proof and
the remaining possible pole is simple and occurs at $t_i^0-t_j^0=1$.
Consequently,
\[
 \left(\prod_{e\in\mathcal H\cap\operatorname{Inv}(w)}u_e\right)
                         r(w,\tau,t)\in\mathcal O_{t^0}.
\]

Since the representations $\tau_{i,\infty}$ are unitary and
$t_i^0-t_j^0\ge0$ for $i<j$, 
\cite[Proposition~I.10 ]{MoeglinWaldspurger1989} implies that 
$N(w,\tau_\infty,t)$ is holomorphic near $t^0$.
At each finite place $v$, let $e_v(t)$ and $e_{w,v}(wt)$ be the
normalized spherical vectors in the source and target of $N(w,\tau_\infty,t)$,  then
\[
 N(w,\tau_v,t)e_v(t)=e_{w,v}(wt).
\]
Both vectors are constant in their compact realizations.
It follows that $N(w,\tau,t)f(t)$ is holomorphic for the sections
under consideration.  Hence
\[
\begin{aligned}
 \left(\prod_{e\in\mathcal H}u_e\right)M(w,\tau,t)f(t)=\left(\prod_{e\in\mathcal H\setminus\operatorname{Inv}(w)}u_e\right)
 \left[\left(\prod_{e\in\mathcal H\cap\operatorname{Inv}(w)}u_e\right)
                    r(w,\tau,t)\right]N(w,\tau,t)f(t)
\end{aligned}
\]
is holomorphic.  By Proposition~\ref{prop:bernstein-lapid-constant-term},
this proves holomorphy of every cuspidal constant term after multiplying
by $\prod_{e\in\mathcal H}u_e$; the factors $[\,\cdot\,]_{wt}$ are
entire in $t$.  The pole-detection criterion
\cite[Sections~I.3.4 and I.4.10]{MoeglinWaldspurger1995} gives
the holomorphicity near $(0,0)$ of 
\begin{equation}
\label{eqn: def-A(u,y)}
 A(u,y):=\left(\prod_{e\in\mathcal H}u_e\right)E_P(t)f(t).
\end{equation}
Here $u=(u_e)_{e\in\mathcal H}$.

For $S\subseteq\mathcal H$, it follows that
\[
 \operatorname{Res}_S(E_Pf)
 =\frac{\left.A(u,y)\right|_{u_e=0,\ e\in S}}
        {\prod_{e\in\mathcal H\setminus S}u_e}.
\]
Setting the coordinates $u_e$, $e\in S$, to zero in $A$ commutes,
which proves independence of the order of the residues.
\end{proof}

Under the hypotheses of Proposition~\ref{prop:cuspidal-normal-crossings},
fix a holomorphic section $f(t)\in I_P(t)$ and a standard parabolic $Q$.
For $w\in W^{\supset Q}(P)$ with $P_w\subsetneq P$, the
operator $C_{P,P_w}^{\mathrm{cusp}}$ in
Proposition~\ref{prop:bernstein-lapid-constant-term} takes a
constant term along the proper parabolic $P_w\cap M_P$ of $M_P$.
This term vanishes by cuspidality of $\tau$.
Thus only $P_w=P$, equivalently $wM_Pw^{-1}=M_Q$, can contribute,
and the formula reduces to
\begin{equation}
\label{eqn: constant-term-formula-special}
 \bigl(C_{G,Q}^{\mathrm{cusp}}(E_P(t)f(t))\bigr)(g)
 =\sum_{\substack{w\in W^{\supset Q}(P)\\P_w=P}}
       [M(w,\tau,t)f(t)](g)(1_{M_Q}).
\end{equation}

If $F$ has expansion
$\sum_{\gamma,\beta}F_{\gamma,\beta}u^\gamma y^\beta$ at $(u,y)=(0,0)$,
write $[u^\gamma y^\beta]F=F_{\gamma,\beta}$ for the coefficient of
$u^\gamma y^\beta$.  In a Laurent expansion the entries of $\gamma$
may be negative.

\begin{corollary}\label{cor:residues-and-coefficients}
Put $p(u)=\prod_{e\in\mathcal H}u_e$ and
$A(u,y)=p(u)E_P(t)f(t)$.  For nonnegative multiindices
$\alpha,\beta$ in the $u$- and $y$-coordinates,
\[
 [u^{\alpha-\mathbf1_{\mathcal H}}y^\beta]E_P(t)f(t)
       =[u^\alpha y^\beta]A(u,y),
\]
where $\mathbf1_{\mathcal H}$ has every entry equal to $1$.
Thus the Laurent coefficients of $E_P(t)f(t)$ and the Taylor
coefficients of $A$ span the same space, including all finite
derivatives normal to the residue hyperplanes.
\end{corollary}

\begin{proof}
By Proposition~\ref{prop:cuspidal-normal-crossings}, $A$ is holomorphic.
Write $A=\sum_{\alpha,\beta\ge0}A_{\alpha,\beta}u^\alpha y^\beta$.
Since $p(u)=u^{\mathbf1_{\mathcal H}}$, division by $p$ gives
\[
 E_P(t)f(t)=\sum_{\alpha,\beta\ge0}
 A_{\alpha,\beta}u^{\alpha-\mathbf1_{\mathcal H}}y^\beta.
\]
Thus both coefficient spans are generated by the same $A_{\alpha,\beta}$.
Since $(\partial_u^\alpha\partial_y^\beta A)(0,0)
=\alpha!\beta!\,A_{\alpha,\beta}$ for every $\alpha,\beta$, this span
contains all finite derivatives in the normal $u$-directions.
\end{proof}

\subsection{Separation by full constant terms}
\label{subsec:separation-full-constant-terms}

Use the data $(P,\tau,t^0,\mathcal H)$ of Proposition~\ref{prop:cuspidal-normal-crossings}.
Partition the $h$ inducing factors in \eqref{eqn: inducing-data-general-Levi}
by their indices $\{1,\ldots,h\}=B_1\sqcup\cdots\sqcup B_s$, keeping
$i$ and $j$ in the same block for every $(i,j)\in\mathcal H$.  Put
\[
 n_b=\sum_{i\in B_b}r_i,\qquad Q=P_{n_1,\ldots,n_s}.
\]
Set $h_0=0$ and $h_b=\sum_{a=1}^b|B_a|$.  Define
\[
 \omega(i)=h_{b-1}+\#\{j\in B_b:j\le i\}
 \qquad(i\in B_b).
\]
Thus $\omega(B_b)=\{h_{b-1}+1,\ldots,h_b\}$, and $\omega$
preserves the original order within each $B_b$.

Consider the parabolic  $P^\omega=P_{r_{\omega^{-1}(1)},\ldots,r_{\omega^{-1}(h)}}$,
and the induced representation
\[
 I_{P^\omega}(\omega t)=\mathop{\times}_{j=1}^h
       \tau_{\omega^{-1}(j)}|{\det}|^{(\omega t)_j},
 \qquad (\omega t)_j=t_{\omega^{-1}(j)}.
\]
Denote by $E^\omega(t):=E_{P^\omega}(\omega t)$ the corresponding Eisenstein operator.
By construction, we have 
\[
P^\omega\subseteq Q.
\]

\begin{remark}
\label{rem: Summation through Q}
Suppose the blocks $B_b$ are consecutive in the prescribed order, so
$P\subseteq Q$.  Let $P_b\subseteq\GL_{n_b}$ have block sizes
$(r_i)_{i\in B_b}$.  Since $N_Q\subseteq P$ and
$P\cap M_Q=\prod_b P_b$, the block diagonal embedding identifies
$\prod_b(P_b(\Q)\backslash\GL_{n_b}(\Q))$ with $P(\Q)\backslash Q(\Q)$.
For a flat section $v_t\in I_P(t)^{K_f}$, summation first over these
cosets gives, in the region of absolute convergence,
\begin{equation}
\label{eqn: Eis-decomposition-iteration}
 (E_P(t)v_t)(g)=
 \sum_{\gamma\in Q(\Q)\backslash G(\Q)}
 \sum_{\delta\in P(\Q)\backslash Q(\Q)}
       [v_t(\delta\gamma g)](1_{M_P}).
\end{equation}
Thus the inner sum forms the Eisenstein series on the Levi blocks,
and the outer sum forms the Eisenstein series for $Q$.
This identity of Eisenstein series extends meromorphically.

Put $K=K_\infty K_f$, $K_b=K\cap\GL_{n_b}(\A)$ and
$K_{n_b,f}=\prod_\ell\GL_{n_b}(\Z_\ell)$.
Let $V_b$ consist of the restrictions to $K_b$ of the $\O_{n_b}$-finite,
$K_{n_b,f}$-fixed sections of $\mathop{\times}_{i\in B_b}\tau_i|{\det}|^{t_i}$.
By Lemma~\ref{lem:compact-realization}, $V_b$ is independent of $t$.
Write $v=v_t|_K$.  In the realization of $I_P(t)$ induced from $Q$,
the value of $v_t$ at $k\in K$ is $v^\sharp(k)\in\bigotimes_b V_b$, where
\[
 v^\sharp(k)(k_1,\ldots,k_s)
     =v\bigl(\operatorname{diag}(k_1,\ldots,k_s)k\bigr)
     \qquad(k_b\in K_b).
\]
\end{remark}

\vspace{.5cm}

Let $c_b=n_b^{-1}\sum_{i\in B_b}r_it_i^0$.
For $i\in B_b$, put
\[
 y_b=\frac1{n_b}\sum_{j\in B_b}r_j(t_j-t_j^0),\qquad
 x_i=t_i-t_i^0-y_b.
\]
Then $\sum_{i\in B_b}r_ix_i=0$ and $\sum_b n_by_b=0$.
For $j_b=\max B_b$ take $x^{(b)}=(x_i)_{i\in B_b\setminus\{j_b\}}$
as coordinates on $\sum_{i\in B_b}r_ix_i=0$, and put
$x=(x^{(1)},\ldots,x^{(s)})$.
For $(i,j)\in\mathcal H$, both indices belong to the same $B_b$, so
$u_{ij}=t_i-t_j-1=x_i-x_j$.  Put
\[
 p_b=\prod_{\substack{(i,j)\in\mathcal H\\i,j\in B_b}}(x_i-x_j),
 \qquad H_b(x^{(b)})=p_bE_{P_b}^{\GL_{n_b}}
       \bigl((t_i^0-c_b+x_i)_{i\in B_b}\bigr).
\]
Thus $p(t):=\prod_{(i,j)\in\mathcal H}u_{ij}=\prod_b p_b$.  Since $\omega$ preserves the order in
each $B_b$, Proposition~\ref{prop:cuspidal-normal-crossings} applied
to the block $b$ gives holomorphy of $H_b$.

\begin{lemma}\label{lem:regularized-sum-through-Q}
Let $v_t\in I_{P^\omega}(\omega t)^{K_f}$ be a flat section and $v=v_t|_K$.
For $y=(y_1,\ldots,y_s)$, let $F_{x,y,v}$ be the section induced from $Q$
with twist $|\det|^{c_b+y_b}$ on block $b$.  For fixed $x$ and $v$, its
restriction to $K$ is independent of $y$ and is given by
\[
 F_{x,y,v}(k)=F_{x,0,v}(k)
 =\left(\bigotimes_b H_b(x^{(b)})\right)v^\sharp(k)
 \qquad(k\in K),
\]
where $v^\sharp$ is defined in Remark~\ref{rem: Summation through Q}.
Then for $g\in G(\A)$ we have the meromorphic identity
\[
 p(t)\sum_{\delta\in P^\omega(\Q)\backslash Q(\Q)}
       [v_t(\delta g)](1_{M_{P^\omega}})
   =[F_{x,y,v}(g)](1_{M_Q}).
\]
\end{lemma}

\begin{proof}
Consider $g=nmk$ with $n\in N_Q(\A)$, $m=(m_b)_b\in M_Q(\A)$ and $k\in K$.
By linearity it suffices to take $v^\sharp(k)=\bigotimes_b u_b$, $u_b\in V_b$.
Let $u_{b,t}$ be the flat section with compact restriction $u_b$ and
exponents $(t_i)_{i\in B_b}$.
Since $t_i-(t_i^0-c_b+x_i)=c_b+y_b$ for every $i\in B_b$, the common
twist on block $b$ is $|\det|^{c_b+y_b}$.  By the product formula
$|\det(\delta_b m_b)|_\A=|\det m_b|_\A$, so in the region of absolute convergence
\[
 p_b\sum_{\delta_b\in P_b(\Q)\backslash\GL_{n_b}(\Q)}
       [u_{b,t}(\delta_b m_b)](1_{M_{P_b}})
 =|\det m_b|_\A^{c_b+y_b}\,[H_b(x^{(b)})u_b](m_b).
\]
By Remark~\ref{rem: Summation through Q}, the cosets in
$P^\omega(\Q)\backslash Q(\Q)$ are represented by
$\delta=\operatorname{diag}(\delta_1,\ldots,\delta_s)$ with the $\delta_b$
ranging independently over the block cosets above.  Since $p(t)=\prod_b p_b$,
normalized induction from $Q$ gives
\[
 \begin{aligned}
 p(t)\sum_\delta [v_t(\delta g)](1_{M_{P^\omega}})
 &=\delta_Q(m)^{1/2}\prod_b
     \left(|\det m_b|_\A^{c_b+y_b}[H_b(x^{(b)})u_b](m_b)\right)\\
 &=[F_{x,y,v}(g)](1_{M_Q}).
 \end{aligned}
\]
The identity extends meromorphically in $t$.
\end{proof}

Let $\mathfrak c$ be the cuspidal support of $I_P(t^0)$, and let $\mathcal A_b$ be the corresponding
module on $[\GL_{n_b}]$ for the support
$\boxtimes_{i\in B_b}\tau_i|{\det}|^{t_i^0-c_b}$.

For $y$ in $Y=\{y\in\C^s:\sum_b n_by_b=0\}$, set
\[
 \mathcal I_Q(y)=\mathop{\times}_{b=1}^s
                    \bigl(\mathcal A_b|{\det}|^{c_b+y_b}\bigr).
\]
Denote by $\D_y$ the finite-order distributions at $0\in Y$ and $f_y$
 the flat section through $f\in\mathcal I_Q(0)^{K_f}$.

\begin{proposition}\label{lem:block-coefficient-map}
Assume that, for the permutation $\omega$ defined above, the global
standard intertwiner
\[
 M(\omega,\tau,t):I_P(t)^{K_f}\longrightarrow I_{P^\omega}(\omega t)^{K_f}
\]
and its inverse are holomorphic near $t^0$ in the compact realizations.
For every $f\in\mathcal I_Q(0)^{K_f}$, the family $E_Q(y)f_y$ is
holomorphic near $y=0$.  Coefficient evaluation as in
Definition~\ref{def:coefficient-map} defines a surjective
$(\mathfrak g^1,K_\infty)$-module map
\begin{equation}\label{eq:separated-support-induction}
 \Phi:\mathcal I_Q(0)^{K_f}\otimes\D_y
      \twoheadrightarrow\mathcal A_{\mathfrak c}^{K_f},
 \qquad f\otimes\eta\longmapsto\eta\bigl(E_Q(y)f_y\bigr),
\end{equation}
where the source has the action \eqref{eq:g-action-jets}.
\end{proposition}

\begin{proof}
Put $H(t)=p(t)E_P(t)$.
Proposition~\ref{prop:cuspidal-normal-crossings} gives holomorphy of $H$.
The Eisenstein functional equation
\cite[Theorem~6.1.7, p.~110]{Shahidi2010} gives
\[
 H^\omega(t):=pE^\omega(t)=H(t)M(\omega,\tau,t)^{-1}.
\]
Thus $H^\omega$ is holomorphic near $t^0$.  Take linear coordinates $z$
on $\sum_i r_it_i=0$ centred at $t^0$.  For compact sections
$v$ of $I_P(t^0)^{K_f}$ and $w$ of $I_{P^\omega}(\omega t^0)^{K_f}$,
write
\[
 M(\omega,\tau,t)^{-1}w=\sum_\mu v_\mu z^\mu,\qquad
 M(\omega,\tau,t)v=\sum_\mu w_\mu z^\mu.
\]
Since $M(\omega,\tau,t)$ and its inverse intertwine the compact group actions,
the $v_\mu$ are full-level compact sections for the parabolic $P$ and the $w_\mu$
are full-level compact sections for $P^\omega$.  Taking the coefficient
of $z^\alpha$ in $H^\omega(t)=H(t)M(\omega,\tau,t)^{-1}$ and
$H(t)=H^\omega(t)M(\omega,\tau,t)$ gives
\[
 \begin{aligned}
 [z^\alpha]\bigl(H^\omega(t)w\bigr)
   &=\sum_{\mu\le\alpha}[z^{\alpha-\mu}]\bigl(H(t)v_\mu\bigr),\\
 [z^\alpha]\bigl(H(t)v\bigr)
   &=\sum_{\mu\le\alpha}[z^{\alpha-\mu}]\bigl(H^\omega(t)w_\mu\bigr).
 \end{aligned}
\]
As $w$ varies, the first equality shows that every Taylor coefficient of
$H^\omega$ lies in the span of those of $H$.  The second
gives the reverse inclusion.
By \cite[Section~2]{GrbacGrobner2024}, the Taylor coefficients of the
regularized cuspidal Eisenstein family span $\mathcal A_{\mathfrak c}$.
Averaging the inducing sections over $K_f$ gives its $K_f$-fixed part.
For these sections $p$ clears the poles, so the coefficients of $H$,
and hence of $H^\omega$, span $\mathcal A_{\mathfrak c}^{K_f}$.

Fix a flat section $v_t\in I_{P^\omega}(\omega t)^{K_f}$ and write
$v=v_t|_K$.  By Lemma~\ref{lem:regularized-sum-through-Q}, summing over
$Q(\Q)\backslash G(\Q)$ in \eqref{eqn: Eis-decomposition-iteration} gives
\[
 H^\omega(x,y)v=E_Q(y)F_{x,y,v}.
\]
The compact value $F_{x,y,v}(k)$ is independent of $y$.  Define
$f_{\alpha,v}\in\mathcal I_Q(0)^{K_f}$ by its restriction to $K$,
\begin{equation}\label{eq:block-coefficient-section}
 f_{\alpha,v}(k)=[x^\alpha]F_{x,y,v}(k)\qquad(k\in K).
\end{equation}
For $y\in Y$, $(f_{\alpha,v})_y$ denotes the unique section in
$\mathcal I_Q(y)^{K_f}$ with this same restriction to $K$
(Lemma~\ref{lem:compact-realization}).

We show that these sections span $\mathcal I_Q(0)^{K_f}$.
For $U_b=V_b\otimes\D_{x^{(b)}}$, the block coefficient maps
\[
 C_b:U_b\twoheadrightarrow\mathcal A_b^{K_{n_b,f}},\qquad
 u\otimes\eta\longmapsto\eta(H_bu),
\]
are surjective by \cite[Section~2]{GrbacGrobner2024}.
Thus $C=\bigotimes_b C_b$ is an $L$-equivariant surjection for
$L=\prod_b\O_{n_b}$.  Fix $f\in\mathcal I_Q(0)^{K_f}$.
Its compact values $f(k)$ span a finite-dimensional $L$-module $W_f$.
Choose a linear lift $r_f:W_f\to\bigotimes_b U_b$ under $C$.
For Haar measure $d\ell$ of total mass $1$ on $L$, put
\[
 s_f(w)=\int_L\ell^{-1}r_f(\ell w)\,d\ell.
\]
The integrand lies in the finite-dimensional $L$-span of $r_f(W_f)$.
Then $Cs_f(w)=w$ and $s_f(\ell w)=\ell s_f(w)$.
Consider the finite expansion
\[
 s_f(f(k))=\sum_\alpha v_\alpha^\sharp(k)\otimes
                   \frac{\partial_x^\alpha|_0}{\alpha!}.
\]
Since $L$ acts trivially on the distributions, each $v_\alpha^\sharp$
satisfies $v_\alpha^\sharp(\ell k)=\ell v_\alpha^\sharp(k)$ for
$\ell\in L$.  It is $K_\infty$-finite and $K_f$-fixed, so by induction
in stages it corresponds to a compact section $v_\alpha$ for $P^\omega$.
Applying $C$ gives $f=\sum_\alpha f_{\alpha,v_\alpha}$.

Since $v$ is $K_\infty$-finite and $K_f$-fixed, the values $v^\sharp(k)$
span a finite-dimensional subspace of $\bigotimes_b V_b$.  Thus we can write
\[
 \begin{aligned}
 v^\sharp(k)&=\sum_{\nu=1}^N a_\nu(k)\bigotimes_b v_{b,\nu},
       \qquad v_{b,\nu}\in V_b,\\
 F_{x,y,v}(k)&=\sum_{\nu=1}^N a_\nu(k)
       \bigotimes_b H_b(x^{(b)})v_{b,\nu},
 \end{aligned}
\]
with continuous functions $a_\nu$ on $K$.  On $[\GL_{n_b}]$, each
$H_b(x^{(b)})v_{b,\nu}$ satisfies a moderate-growth bound uniform for
$x$ in a sufficiently small closed polydisc about $0$, by
\cite[Section~8.6, Corollary]{BernsteinLapid2024}.  The sum is finite
and the $a_\nu$ are bounded on $K$, so the resulting bound for
$F_{x,y,v}(k)$ is uniform in $k$ and $x$.
For $\operatorname{Re}y$ sufficiently far in the positive chamber for
$Q$, summation over $Q(\Q)\backslash G(\Q)$ as in
Remark~\ref{rem: Summation through Q}, applied to $P^\omega\subseteq Q$, gives
\[
 \bigl(E_Q(y)F_{x,y,v}\bigr)(g)
 =\sum_{\gamma\in Q(\Q)\backslash G(\Q)}
       [F_{x,y,v}(\gamma g)](1_{M_Q}),\qquad g\in G(\A),
\]
where $1_{M_Q}$ is the identity of the Levi subgroup $M_Q$.
By the preceding bound, this sum converges absolutely and uniformly
for $x$ in the chosen polydisc \cite[Section~2.2]{BernsteinLapid2024}.
For each $\gamma$, write $\gamma g=nmk$ with $n\in N_Q(\A)$,
$m=(m_b)_b\in M_Q(\A)$ and $k\in K$.  By the extension formula in
the proof of Lemma~\ref{lem:compact-realization},
\[
 [F_{x,y,v}(\gamma g)](1_{M_Q})
 =\delta_Q(m)^{1/2}\left(\prod_b|\det m_b|_\A^{c_b+y_b}\right)
       F_{x,y,v}(k)(m).
\]
The scalar factor is independent of $x$, and
$[x^\alpha]F_{x,y,v}(k)=f_{\alpha,v}(k)$ by definition.  Thus the
coefficient of this summand is $[(f_{\alpha,v})_y(\gamma g)](1_{M_Q})$
by the same extension formula.  Cauchy's formula allows us to take
the coefficient term by term, so
\[
 \begin{aligned}
 \bigl(E_Q(y)(f_{\alpha,v})_y\bigr)(g)
 &=\sum_{\gamma\in Q(\Q)\backslash G(\Q)}
        [(f_{\alpha,v})_y(\gamma g)](1_{M_Q})\\
 &=[x^\alpha]\bigl(E_Q(y)F_{x,y,v}\bigr)(g)
  =[x^\alpha]\bigl(H^\omega(x,y)v\bigr)(g).
 \end{aligned}
\]
The identity extends meromorphically in $y$.  Its right hand side is
holomorphic near $y=0$, so the spanning assertion gives holomorphy of
$E_Q(y)f_y$ for every $f\in\mathcal I_Q(0)^{K_f}$.  Moreover,
\begin{equation}\label{eq:block-coefficient-identity}
 \Phi\!\left(f_{\alpha,v}\otimes
            \frac{\partial_y^\beta|_0}{\beta!}\right)
       =[x^\alpha y^\beta]H^\omega(x,y)v.
\end{equation}
These coefficients span $\mathcal A_{\mathfrak c}^{K_f}$, so $\Phi$ is
surjective.  For $X\in\mathfrak g^1$, applying $\partial_y^\beta|_0$
to $X(E_Q(y)f_y)=E_Q(y)(\pi_y(X)f)_y$ gives the Leibniz rule
\eqref{eq:g-action-jets}.  The map $\Phi$ also commutes with
$K_\infty$, whose action is independent of $y$.
\end{proof}

Let $\Omega_i$ be the Harish--Chandra parameter of
$\tau_{i,\infty}|{\det}|^{t_i^0}$.
Let $\chi^0$ be the infinitesimal character of $M_Q$ with parameter
$\bigsqcup_{i\in B_b}\Omega_i$ on each block $b$.
For $\varphi\in\mathcal A_{\mathfrak c}$ and fixed $k\in K_\infty K_f$,
write
\[
 (\CT_Q\varphi)(mk)=\sum_\chi h_{\chi,k}(m),
\]
where $h_{\chi,k}$ is the generalized $\chi$-component in the Levi variable $m$.
Define the projection $\Pi_Q^0$ by
\begin{equation}
\label{eqn: projector-constant-term}
 (\Pi_Q^0\CT_Q\varphi)(mk):=h_{\chi^0,k}(m).
\end{equation}

Take $y_1,\ldots,y_{s-1}$ as coordinates on $Y$.  Write
$a=\operatorname{diag}(a_b I_{n_b})_{b=1}^{s}$ for a positive real block scalar,
trivial at the finite places, with $\prod_b a_b^{n_b}=1$.  Its height
coordinates and their inverse are
\[
 T_j=n_j\log(a_j/a_s)\quad(j<s),\qquad
 a_s(T)=\exp\!\left(-\frac1n\sum_{j<s}T_j\right),\qquad
 a_j(T)=a_s(T)e^{T_j/n_j}\quad(j<s).
\]
For $\xi\in\mathcal I_Q(0)^{K_f}\otimes\D_y$,
$m\in\prod_b\GL_{n_b}(\A)^1$ and $k\in K_\infty K_f$, define
\begin{equation}
\label{eqn; def-cal-P}
 \mathcal P_\xi(T;m,k)
 :=\left(\prod_b a_b(T)^{-n_b c_b}\right)
    (\Pi_Q^0\CT_Q\Phi(\xi))(a(T)mk).
\end{equation}

\begin{proposition}\label{lem:separated-cuspidal-parameters}
Under the hypothesis of Proposition~\ref{lem:block-coefficient-map},
the map $\Phi$ defined there is an isomorphism.
For $\xi=\sum_\beta F_\beta\otimes\partial_0^\beta\in \mathcal I_Q(0)^{K_f}\otimes\D_y$,
\begin{equation}\label{eq:separated-support-constant-term}
 \mathcal P_\xi(T;m,k)=\sum_\beta T^\beta[F_\beta(k)](m).
\end{equation}
Thus $[F_\beta(k)](m)=[T^\beta]\mathcal P_\xi(T;m,k)$.
\end{proposition}

\begin{proof}
By the proof of Proposition~\ref{lem:block-coefficient-map}, the sections
$f_{\alpha,v}$ span $\mathcal I_Q(0)^{K_f}$.  It suffices to take
$\xi=f_{\alpha,v}\otimes\partial_0^\beta$, where $v$ is a compact
section for $P^\omega$.  Put $\lambda=\left.\partial_y^\beta\right|_0[x^\alpha]$.
Multiplying \eqref{eq:block-coefficient-identity} by $\beta!$ gives
\[
 \Phi(\xi)=\lambda\bigl(H^\omega(x,y)v\bigr),
\]
By differentiation under the integral over
the compact quotient $N_Q(\Q)\backslash N_Q(\A)$, we obtain 
\[
 \CT_Q\Phi(\xi)=\lambda\bigl(\CT_Q(H^\omega(x,y)v)\bigr).
\]

Let $v_t\in I_{P^\omega}(\omega t)^{K_f}$ be the flat section through $v$,
and put $\omega\tau=\boxtimes_{j=1}^h\tau_{\omega^{-1}(j)}$.
By cuspidality, only $w\in{}_{Q}W_{P^\omega}$ with
$wM_{P^\omega}w^{-1}\subseteq M_Q$ occur in
\eqref{eq:full-normalized-constant-term}.  For such $w$, the standard
parabolic $Q_w\subseteq Q$ has Levi $M_{Q_w}=wM_{P^\omega}w^{-1}$.
For $m\in M_Q(\A)$ and $k\in K$, the $w$-summand of
$\CT_Q(H^\omega(x,y)v)$ is
\[
 I_w(t)(mk)=p(t)\,\delta_Q(m)^{-1/2}
   \sum_{\gamma\in Q_w(\Q)\backslash Q(\Q)}
   \bigl[(M(w,\omega\tau,\omega t)v_t)(\gamma mk)\bigr](1_{M_{Q_w}}).
\]
Here $M(w,\omega\tau,\omega t)v_t$ is a section for $Q_w$ with
exponent vector $w(\omega t)$, and the sum is continued meromorphically.
In $M_{P^\omega}$ the factor $\GL_{r_i}$ carrying $\tau_i$ occupies
the coordinate interval
\[
 J_i^\omega=
 \left[1+\sum_{\omega(j)<\omega(i)}r_j,\ \sum_{\omega(j)\le\omega(i)}r_j\right].
\]
Recall that $n_a=\sum_{j\in B_a}r_j$ is the size of the $a$th block
of $Q=P_{n_1,\ldots,n_s}$.  Viewing $w$ as a permutation of $\{1,\ldots,n\}$, define
$\sigma_w:\{1,\ldots,h\}\to\{1,\ldots,s\}$ by
\[
 \sigma_w(i)=b
 \quad\Longleftrightarrow\quad
 w(J_i^\omega)\subseteq
 \left[1+\sum_{a<b}n_a,\ \sum_{a\le b}n_a\right].
\]
The interval on the right consists of the coordinates of the $b$th
block of $M_Q$.  The inclusion $wM_{P^\omega}w^{-1}\subseteq M_Q$
ensures that $w(J_i^\omega)$ lies in exactly one such interval.
Since the Eisenstein operator on each block
intertwines the Levi action, $I_w(t)$ has infinitesimal character
$\chi_w(t)$ with Harish--Chandra parameter
\[
 \bigsqcup_{\sigma_w(i)=b}\bigl(\Omega_i+(t_i-t_i^0)\bigr)
 \qquad\text{on block }b.
\]
The $\Omega_i$ are pairwise disjoint by regularity.  Hence if
$\chi_w(t^0)=\chi_{w'}(t^0)$, then $\sigma_w=\sigma_{w'}$ and
$w'w^{-1}\in W_Q$.  Since $w,w'$ are minimal double-coset
representatives, $w=w'$.  For $w=1$ we have $\sigma_1^{-1}(b)=B_b$,
so $\chi_1(t^0)=\chi^0$.

Choose $Z\in Z(U(\mathfrak m_{Q,\C}))$ such that the values
$z_w(t^0):=\chi_w(t^0)(Z)$ are pairwise distinct.  Let
$z_w(t)=\chi_w(t)(Z)$.  By \eqref{eq:full-normalized-constant-term},
$\CT_Q(H^\omega(x,y)v)=\sum_u I_u(t)$, where $u$ runs over the
contributing Weyl elements.  Fix $w$ in this sum.  Since
$(Z-z_{w'}(t))I_u(t)=(z_u(t)-z_{w'}(t))I_u(t)$, we obtain
\[
 \begin{aligned}
 \left(\prod_{w'\ne w}\frac{Z-z_{w'}(t)}{z_w(t)-z_{w'}(t)}\right)
       \CT_Q\bigl(H^\omega(x,y)v\bigr)
 &=\sum_u\left(\prod_{w'\ne w}
       \frac{z_u(t)-z_{w'}(t)}{z_w(t)-z_{w'}(t)}\right) I_u(t)\\
 &=I_w(t).
 \end{aligned}
\]
 This identity holds
meromorphically, and the denominators are nonzero at $t^0$.  Thus
each $I_w(t)$ is holomorphic near $t^0$.  Fix $w$ and
$U\in Z(U(\mathfrak m_{Q,\C}))$.  Denote $I_{\mu,\nu}=[x^\mu y^\nu]I_w(t)$
and $c_{a,b}=[x^a y^b]\chi_w(t)(U)$, so $c_{0,0}=\chi_w(t^0)(U)$.
Since $U$ acts only in the Levi variable, comparing coefficients in
$U I_w(t)=\chi_w(t)(U)I_w(t)$ gives
\[
 (U-c_{0,0})I_{\mu,\nu}
 =\sum_{\substack{0\le a\le\mu,\;0\le b\le\nu\\(a,b)\ne(0,0)}}
     c_{a,b}\,I_{\mu-a,\nu-b}.
\]
Each application of $U-c_{0,0}$ lowers the total degree of the indices
by at least one and kills $I_{0,0}$, for which the sum is empty.
Thus $(U-c_{0,0})^{|\alpha|+|\beta|+1}$
kills $\lambda I_w=\beta!I_{\alpha,\beta}$.  This holds for every central
$U$, so $\lambda I_w$ has generalized infinitesimal character $\chi_w(t^0)$.
In $\CT_Q\Phi(\xi)=\sum_w\lambda I_w$ these characters are pairwise distinct,
and $\chi_1(t^0)=\chi^0$.  Hence \eqref{eqn: projector-constant-term} gives
\[
 (\Pi_Q^0\CT_Q\Phi(\xi))(mk)=h_{\chi^0,k}(m)=(\lambda I_1)(mk).
\]

For $w=1$, the inclusion $P^\omega\subseteq Q$ gives $Q_1=P^\omega$,
and $M(1,\omega\tau,\omega t)=\mathrm{id}$.  Apply
Lemma~\ref{lem:regularized-sum-through-Q} to the sum in $I_1(t)$
at $g=a(T)mk$.  We obtain
\[
 \begin{aligned}
 I_1(t)(a(T)mk)
 &=\delta_Q(a(T)m)^{-1/2}[F_{x,y,v}(a(T)mk)](1_{M_Q})\\
 &=\left(\prod_b|\det(a_b(T)m_b)|_\A^{c_b+y_b}\right)[F_{x,y,v}(k)](a(T)m).
 \end{aligned}
\]
On block $b$, the modulus character of $P_b$ is trivial on the scalar
$a_b(T)I_{n_b}$.  For $i\in B_b$, $\tau_i$ has trivial central character
and $|\det(a_b(T)I_{r_i})|_\A=a_b(T)^{r_i}$.  Thus for $u_b\in V_b$,
\[
 [H_b(x^{(b)})u_b](a_b(T)m_b)
 =a_b(T)^{\sum_{i\in B_b}r_i(t_i^0-c_b+x_i)}
       [H_b(x^{(b)})u_b](m_b).
\]
The exponent is zero since $\sum_{i\in B_b}r_it_i^0=n_bc_b$ and
$\sum_{i\in B_b}r_ix_i=0$.  Hence $[F_{x,y,v}(k)](a(T)m)=[F_{x,y,v}(k)](m)$.
Since $|\det m_b|_\A=1$, the formula for $I_1(t)$ becomes
\[
 I_1(t)(a(T)mk)
 =\left(\prod_b a_b(T)^{n_b(c_b+y_b)}\right)[F_{x,y,v}(k)](m).
\]
Now $\sum_b n_by_b=0$ and $T_j=n_j\log(a_j(T)/a_s(T))$ give
$\prod_b a_b(T)^{n_by_b}=e^{\sum_{j<s}y_jT_j}$.  By \eqref{eqn; def-cal-P}
and the definition of $f_{\alpha,v}$ in \eqref{eq:block-coefficient-section},
\[
 \begin{aligned}
 \mathcal P_\xi(T;m,k)
 &=\lambda\left(e^{\sum_{j<s}y_jT_j}[F_{x,y,v}(k)](m)\right)\\
 &=\left.\partial_y^\beta e^{\sum_{j<s}y_jT_j}\right|_{y=0}
       [f_{\alpha,v}(k)](m)
   =T^\beta[f_{\alpha,v}(k)](m).
 \end{aligned}
\]
Linearity gives \eqref{eq:separated-support-constant-term} for general
$\xi$.  If $\Phi(\xi)=0$, taking $[T^\beta]$ gives $[F_\beta(k)](m)=0$
for every $\beta,m,k$.  The block forms $F_\beta(k)$ are determined
by their restrictions to $\prod_b\GL_{n_b}(\A)^1$, so $F_\beta|_K=0$.
Hence $F_\beta=0$ for every $\beta$, and $\Phi$ is injective.
Proposition~\ref{lem:block-coefficient-map} gives surjectivity.
\end{proof}

For a holomorphic operator family $T(z)=\sum_\gamma T_\gamma z^\gamma$
in the compact realizations and $\delta_\alpha(f)=[z^\alpha]f$, define
\[
 T_{\D}(f\otimes\delta_\alpha)
 :=\sum_{\gamma\le\alpha}T_\gamma f\otimes\delta_{\alpha-\gamma}.
\]
For composable holomorphic families $B,T$, convolution of their
Taylor coefficients gives
\begin{equation}\label{eq:holomorphic-operator-composition}
 (BT)_{\D}=B_{\D}T_{\D},\qquad
 (T^{-1})_{\D}T_{\D}=\mathrm{id}
 \quad\text{if $T^{-1}$ is holomorphic}.
\end{equation}

For $I_P(t)$ as in \eqref{eqn: inducing-data-general-Levi}, let $V$
be the full-level compact realization of $I_P(t^0)$ and let $z$ be
coordinates centred at $t^0$. A \emph{full-level inducing jet} for
$P$ at $t^0$ is a finite sum
\begin{equation}\label{eq:full-level-inducing-jet}
 \eta=\sum_\alpha f_\alpha\otimes\delta_\alpha\in V\otimes\D_z,
 \qquad f_\alpha\in V,\qquad \delta_\alpha(g)=[z^\alpha]g .
\end{equation}

\begin{proposition}\label{lem:coefficient-kernel-constant-terms}
Use the notation of \eqref{eq:full-level-inducing-jet}, with affine-linear
coordinates $z$ and the inducing blocks in any order.  Suppose $t^0$ is real and the
Harish--Chandra parameter of $I_P(t^0)_\infty$ has distinct entries.
For a scalar holomorphic germ $a(z)$, put
\[
 H(z)=a(z)E_P(t(z)),\qquad
 T_{Q,\sigma}(z)=a(z)M(\sigma,\tau,t(z)),
\]
where $Q$ runs over the standard parabolics and $\sigma$ indexes the summands in
\eqref{eqn: constant-term-formula-special}.  If $H(z)f$ and all
$T_{Q,\sigma}(z)f$ are holomorphic near $0$ for every $f\in V$,
then the coefficient map $\Xi_H$ of Definition~\ref{def:coefficient-map} satisfies
\begin{equation}\label{eq:cuspidal-coefficient-kernel}
 \ker\Xi_H=\bigcap_{Q,\sigma}\ker (T_{Q,\sigma})_{\D}.
\end{equation}
If $t_1^0\ge\cdots\ge t_h^0$, we can take
$a(z)=p(z)=\prod_{e\in\mathcal H}u_e(z)$
by Proposition~\ref{prop:cuspidal-normal-crossings} and its proof.
\end{proposition}

\begin{proof}
Multiplying \eqref{eqn: constant-term-formula-special} by $a(z)$
gives $T_{Q,\sigma}(z)$ in the $\sigma$-summand of the cuspidal
constant term of $H(z)$. The locally uniform moderate-growth bounds
of \cite[Section~8.6, Corollary]{BernsteinLapid2024} permit finite
derivatives to commute with the cuspidal projections.

Fix $\eta\in\ker\Xi_H$ and a standard parabolic $Q$ occurring in
\eqref{eqn: constant-term-formula-special}. Regard $\sigma$ as a
permutation of the Levi blocks and put
\[
 \sigma\tau=\boxtimes_{b=1}^h\tau_{\sigma^{-1}(b)},\qquad
 \sigma t^0=(t_{\sigma^{-1}(b)}^0)_{b=1}^h .
\]
For $m\in M_Q(\A)^1$, $Y\in\mathfrak a_Q^G$ and $k\in K_\infty K_f$,
the normalized cuspidal constant term is
\[
 0=(\CT_Q\Xi_H(\eta))^{\mathrm{cusp}}(m\exp(Y)k)
   =\sum_\sigma e^{\langle\sigma t^0,Y\rangle}\Phi_\sigma(Y;m,k).
\]
As in \eqref{eq:constant-exponents}, $\Phi_\sigma$ is polynomial in
$Y$, since $\eta$ differentiates the exponential only finitely
many times. By \eqref{eq:normalized-intertwiner-notation} and $M=rN$,
the target of $T_{Q,\sigma}(z)$ is induced from
$\boxtimes_{b=1}^h\tau_{\sigma^{-1}(b)}|\det|^{t_{\sigma^{-1}(b)}(z)}$.
Put $F_z=T_{Q,\sigma}(z)f$. For $m=(m_b)_b\in M_Q(\A)^1$,
each $|\det m_b|_\A=1$ and $\delta_Q(m)=1$. Thus the extension
formula in the proof of Lemma~\ref{lem:compact-realization},  gives
\[
 e^{-\langle\rho_Q,Y\rangle}F_z(m\exp(Y)k)(1_{M_Q})
 =e^{\langle\sigma t(z),Y\rangle}F_z(k)(m).
\]
Let $\eta=\sum_\alpha v_\alpha\otimes\delta_\alpha$ and choose the
dual coordinates $X$ given by
$\langle\sigma(t(z)-t^0),Y\rangle=\langle z,X\rangle$.
After removing $e^{\langle\sigma t^0,Y\rangle}$, Definition~\ref{def:coefficient-map} give
\[
\begin{aligned}
 \Phi_\sigma(Y;m,k)
 &=\sum_\alpha [z^\alpha]\!
   \left(e^{-\langle\sigma t^0+\rho_Q,Y\rangle}
   [T_{Q,\sigma}(z)v_\alpha](m\exp(Y)k)(1_{M_Q})\right)\\
 &=\sum_\alpha [z^\alpha]\!
   \left(e^{\langle z,X\rangle}[T_{Q,\sigma}(z)v_\alpha](k)(m)\right)\\
 &=\sum_\beta \frac{X^\beta}{\beta!}\,f_\beta(k)(m),
\end{aligned}
\]
where $ f_\beta=\sum_{\alpha\ge\beta}
   [z^{\alpha-\beta}]\bigl(T_{Q,\sigma}(z)v_\alpha\bigr).$
Each $T_{Q,\sigma}(z)v_\alpha$ is holomorphic in the fixed target
compact realization. Hence $f_\beta$ belongs to that realization,
and $f_\beta(k)$ lies in the cuspidal space of $\sigma\tau$ by
Definition~\ref{def:compact-realization}\textup{(ii)}. Its restriction
to $M_Q(\A)^1$ is the coefficient of $X^\beta/\beta!$.

Fix a Weyl element $\sigma_0$ in \eqref{eqn: constant-term-formula-special} for $Q$. Using the Petersson inner product of
Definition~\ref{def:cuspidal-constant-term-component}, project onto the
$\sigma_0\tau$-isotypic component on $M_Q(\Q)\backslash M_Q(\A)^1$.
Since $\GL_r(\A)=\GL_r(\A)^1\cdot\R_{>0}I_r$ and the positive real
scalars act trivially on each $\tau_i$, its restriction to
$\GL_{r_i}(\A)^1$ determines $\tau_i$. By orthogonality, the projected identity is
\begin{equation}
\label{eqn: proj-sigma-tau-component}
 0=\sum_{\substack{\sigma\in W^{\supset Q}(P),\ P_\sigma=P\\
                   \sigma\tau\simeq\sigma_0\tau}}
       e^{\langle\sigma t^0,Y\rangle}\Phi_\sigma(Y;m,k).
\end{equation}
By \eqref{eqn: inducing-data-general-Levi},
$\sum_b r_{\sigma^{-1}(b)}t_{\sigma^{-1}(b)}^0=0$, so $\sigma t^0$
annihilates $\mathfrak a_G$ and defines a linear form on
$\mathfrak a_Q^G=\mathfrak a_Q/\mathfrak a_G$.
Use the coordinates on
$\mathfrak a_Q$ in which $H_Q(g)=(\log|\det g_j|_\A)_{j=1}^h$
for $g=(g_j)_j\in M_Q(\A)$, and let $Y_b$ be the class of the
$b$-th coordinate vector in $\mathfrak a_Q^G$.
Let $\sigma\ne\sigma_0$ be a Weyl element in
\eqref{eqn: constant-term-formula-special} with $\sigma\tau\simeq\sigma_0\tau$.
Choose a block $b$ with $i=\sigma^{-1}(b)\ne j=\sigma_0^{-1}(b)$.
On this block $\tau_i\simeq\tau_j$. Since the parameter of $I_P(t^0)_\infty$
has distinct entries, we have 
\[
 \langle\sigma t^0,Y_b\rangle=t_i^0
 \ne t_j^0=\langle\sigma_0t^0,Y_b\rangle .
\]
Put $Y_0=\sum_{b=1}^h t_{\sigma_0^{-1}(b)}^0Y_b$. Since $\sigma t^0$
and $\sigma_0t^0$ are permutations of the same real vector,
\[
 \langle\sigma_0t^0-\sigma t^0,Y_0\rangle
 =\sum_b t_{\sigma_0^{-1}(b)}^0
       \bigl(t_{\sigma_0^{-1}(b)}^0-t_{\sigma^{-1}(b)}^0\bigr)=\frac12\sum_b
       \bigl(t_{\sigma_0^{-1}(b)}^0-t_{\sigma^{-1}(b)}^0\bigr)^2>0.
\]
Fix $Y,m,k$. Substitute $Y+sY_0$ for $Y$ in the identity \eqref{eqn: proj-sigma-tau-component}
and divide by $e^{\langle\sigma_0t^0,Y+sY_0\rangle}$ to give
\[
 \Phi_{\sigma_0}(Y+sY_0;m,k)
 =-\sum_{\substack{\sigma\in W^{\supset Q}(P),\ P_\sigma=P\\
                    \sigma\tau\simeq\sigma_0\tau,\ \sigma\ne\sigma_0}}
 e^{\langle\sigma t^0-\sigma_0t^0,Y+sY_0\rangle}
 \Phi_\sigma(Y+sY_0;m,k).
\]
Each term on the right tends to zero as $s\to+\infty$, since the
exponential decays and $\Phi_\sigma(Y+sY_0;m,k)$ is polynomial in $s$.
Thus the polynomial $\Phi_{\sigma_0}(Y+sY_0;m,k)$ tends to zero
and is identically zero. Taking $s=0$ gives
$\Phi_{\sigma_0}(Y;m,k)=0$ for all $Y,m,k$.

Since $\Phi_\sigma(Y;m,k)=0$ for all $X,m,k$, every
$f_\beta(k)(m)=0$. Hence every $f_\beta=0$, so
$(T_{Q,\sigma})_{\D}\eta=\sum_\beta f_\beta\otimes\delta_\beta=0$.
This is the split-centre calculation in the proof of
Proposition~\ref{lem:separated-cuspidal-parameters}, applied to $M_Q$.
Conversely, if $(T_{Q,\sigma})_{\D}\eta=0$ for every $(Q,\sigma)$,
then every cuspidal projection of $\Xi_H(\eta)$ vanishes by
\eqref{eqn: constant-term-formula-special}, hence $\Xi_H(\eta)=0$
by \cite[Section~6.2, Proposition]{BernsteinLapid2024}.
\end{proof}

\begin{corollary}\label{cor:regularized-coefficient-image}
Use the data $(P,\tau,t^0,\mathcal H)$ of Proposition~\ref{prop:cuspidal-normal-crossings}
and put $p=\prod_{e\in\mathcal H}u_e$. Let $\mathfrak c$ be the
cuspidal support of $I_P(t^0)$ and $\chi$ the infinitesimal character
of $I_P(t^0)_\infty$ on $\mathfrak Z_G^1$. Then
\[
 \operatorname{im}\Xi_{pE_P}
 =\bigl(\operatorname{Fin}_{\ker\chi}\mathcal A(G)_{\mathfrak c}\bigr)^{K_f}.
\]
If $T(t):V'\to V$ is a holomorphic family of isomorphisms between
full-level compact realizations, with $V$ the realization of $I_P(t)^{K_f}$
and $T^{-1}$ holomorphic near $t^0$, then
$\operatorname{im}\Xi_{(pE_P)\circ T}=\operatorname{im}\Xi_{pE_P}$.
\end{corollary}

\begin{proof}
By \cite[Section~2, definition of $\mathcal A_{\{P\},\varphi(\pi)}$ by Taylor coefficients]{GrbacGrobner2024},
the Taylor coefficients of the regularized Eisenstein family at $t^0$
span the summand with support $\mathfrak c$.
Averaging the inducing sections over $K_f$ gives its $K_f$-invariants.
For these sections take $Q(s)=p$ by Proposition~\ref{prop:cuspidal-normal-crossings}.
For a $K_f$-fixed holomorphic section $f_t$, put $\Phi(t)=p(t)E_P(t)f_t$.
By Corollary~\ref{cor:residues-and-coefficients}, its Taylor coefficients
span the same space as the Laurent coefficients of $E_P(t)f_t$.
Write $\delta_\alpha=[(t-t^0)^\alpha]$ for coefficient extraction in the
independent coordinates $t_1,\ldots,t_{h-1}$. In the compact realization
write $f_t=\sum_\beta f_\beta(t-t^0)^\beta$. By Definition~\ref{def:coefficient-map},
\[
 \delta_\alpha\Phi
 =\sum_{\beta\le\alpha}\Xi_{pE_P}(f_\beta\otimes\delta_{\alpha-\beta}).
\]
Conversely, $\Xi_{pE_P}(f\otimes\delta_\alpha)=\delta_\alpha\Phi$ for
the flat section through $f$.

Let $\chi_t$ be the infinitesimal character of $I_P(t)_\infty$.
Equivariance gives $Z\Phi(t)=\chi_t(Z)\Phi(t)$ for $Z\in\mathfrak Z_G^1$.
For $Z\in\ker\chi$, the scalar $\chi_t(Z)$ is holomorphic by the
Harish--Chandra isomorphism and $\chi_{t^0}(Z)=0$. Thus for $|\alpha|=d$
and $Z_1,\ldots,Z_{d+1}\in\ker\chi$,
\[
 Z_1\cdots Z_{d+1}\,\delta_\alpha\Phi
 =\delta_\alpha\!\left(\Bigl(\prod_{j=1}^{d+1}\chi_t(Z_j)\Bigr)\Phi(t)\right)=0.
\]
Hence every coefficient belongs to $\operatorname{Fin}_{\ker\chi}\mathcal A(G)_{\mathfrak c}$.
Since the coefficients for $K_f$-fixed sections span $\mathcal A(G)_{\mathfrak c}^{K_f}$, we obtain
\[
 \operatorname{im}\Xi_{pE_P}
 =\mathcal A(G)_{\mathfrak c}^{K_f}
 =\bigl(\operatorname{Fin}_{\ker\chi}\mathcal A(G)_{\mathfrak c}\bigr)^{K_f}.
\]

Finally, put $z_i=t_i-t_i^0$ for $1\le i<h$. On the full-level inducing
jets of \eqref{eq:full-level-inducing-jet}, $T(t):V'\to V$ induces
$T_{\D}:V'\otimes\D_z\to V\otimes\D_z$, with inverse $(T^{-1})_{\D}$
by \eqref{eq:holomorphic-operator-composition}. Definition~\ref{def:coefficient-map} gives
\[
 \begin{aligned}
 \Xi_{(pE_P)\circ T}(\eta')&=\Xi_{pE_P}(T_{\D}\eta')
   &&(\eta'\in V'\otimes\D_z),\\
 \Xi_{pE_P}(\eta)&=\Xi_{(pE_P)\circ T}((T^{-1})_{\D}\eta)
   &&(\eta\in V\otimes\D_z).
 \end{aligned}
\]
These identities give both inclusions of the images.
\end{proof}

\section{The degree bound and proof of the main theorem}
\label{sec:numerical}

This section combines the degree calculation of Section~\ref{sec:packets}
with the separation of cuspidal parameters in
Propositions~\ref{lem:block-coefficient-map} and
\ref{lem:separated-cuspidal-parameters}.The only
results used from later sections are
Theorems~\ref{lem:one-copy-internal-bound} and
\ref{prop:two-cuspidal-bound}; both enter solely through
Lemma~\ref{lem:component-depth-at-most-two} below. 

Fix full-level cuspidal inducing data compatible with the infinitesimal
character of the trivial representation on $\mathfrak{sl}_n$. Unpacking
its Speh factors gives the partition
\begin{equation}\label{eqn: tiling-inf-char-7.1}
 U_n=\{0,\ldots,n-1\}=C\sqcup\bigsqcup_\nu B_\nu,
 \qquad B_\nu=\{b_{\nu,1}<\cdots<b_{\nu,k_\nu}\},\quad k_\nu\ge2.
\end{equation}
Here $C$ consists of the character weights, and $B_\nu-c_\nu$ is the
Harish--Chandra parameter of the underlying unitary cuspidal
representation $\pi_\nu$, where
\[
 c_\nu=\frac{b_{\nu,1}+b_{\nu,k_\nu}}2,\qquad
 h_\nu=\lfloor k_\nu/2\rfloor,\qquad
 b_\nu=\lfloor k_\nu^2/4\rfloor,\qquad
 w_\nu=w(\pi_\nu)=b_{\nu,k_\nu}-b_{\nu,1}.
\]
Symmetry gives $b_{\nu,i}+b_{\nu,k_\nu+1-i}=2c_\nu$ and
$c_\nu\in\tfrac12\Z$. If $k_\nu$ is odd, its middle weight is
$c_\nu\in B_\nu\subset\Z$, so $w_\nu$ is even. For $k_\nu=2$,
Remark~\ref{rem:rank-two-motivic-parity} instead gives $w_\nu$ odd and
$c_\nu\in\tfrac12+\Z$. Isomorphic cuspidal representations have the
same centered parameter, hence
\begin{equation}\label{eq:component-cuspidal-translations}
 \pi_\nu\simeq\pi_\mu\quad\Longrightarrow\quad
 B_\mu=B_\nu+(c_\mu-c_\nu).
\end{equation}

Let $K$ range over the connected components of
\[
 C\cup\{c_\nu\}_\nu
 \cup\bigcup_{a,a+1\in C}[a,a+1]_{\R}
 \cup\bigcup_{\substack{\pi_\nu\simeq\pi_\mu\\c_\mu=c_\nu+1}}
                  [c_\nu,c_\mu]_{\R}.
\]
For a component containing cuspidal factors put
\[
 m_K=|C\cap K|,\quad H_K=\sum_{c_\nu\in K}h_\nu,\quad
 \ell_K=\#\{\nu:c_\nu\in K\},\quad
 n_K=m_K+\sum_{c_\nu\in K}k_\nu.
\]
For $x\in U_n\setminus B_\nu$, define
\[
 d_\nu(x)=\min\bigl(\#\{b\in B_\nu:b<x\},
                         \#\{b\in B_\nu:b>x\}\bigr).
\]
This counts  intervals
$(b_{\nu,i},b_{\nu,k_\nu+1-i})$, $1\le i\le h_\nu$, containing $x$.
We call $(b_{\nu,i},b_{\nu,k_\nu+1-i}) (1\le i \le h_\nu)$
a  paired interval.

\begin{lemma}\label{lem:weighted-cuspidal-depth}
Every $c\in C\cap K$ satisfies
\[
 \sum_{c_\nu\in K}d_\nu(c)\ge\min(H_K,4).
\]
If any paired interval of a factor centered in $K$ misses $c$, the sum
is at least four.
\end{lemma}

\begin{proof}
Fix $c\in C\cap K$ and put $D=\sum_{c_\nu\in K}d_\nu(c)$.
If every paired interval contains $c$, then $D=H_K$.  Otherwise a
paired interval lies entirely on one side of $c$, since $c$ is a
character weight.  The reflection $x\mapsto n-1-x$ takes $(a,a')$ to
$(n-1-a',n-1-a)$ and $B_\nu$ to $n-1-B_\nu=B_\nu+(n-1-2c_\nu)$.
Thus it carries components to components.  For each $\nu$, the paired
intervals of $B_\nu$ containing $c$ correspond bijectively to those of
$n-1-B_\nu$ containing $n-1-c$, so both counts equal $d_\nu(c)$.
Reflecting if necessary, choose
\[
 b=\min\{a>c:a\in B_\nu,\ a<c_\nu,\ c_\nu\in K\},
 \qquad b\in B_{\nu_0},\qquad t=c_{\nu_0}>b .
\]
Then $c\le b-1$, $t\ge b+\tfrac12$ and $[c,t]_{\R}\subset K$.
By minimality of $b$ and $c\notin B_\nu$,
\begin{equation}\label{eqn: exclusion-ineq}
 c_\nu\in K,\quad a\in B_\nu,\quad a<c_\nu,\quad a<b
 \quad\Longrightarrow\quad a<c.
\end{equation}

Consider an interval $[s,s+1]_{\R}\subset K$ meeting $(b-1,b+1)$
with $c_\nu=s$, $c_\mu=s+1$ and $\pi_\nu\simeq\pi_\mu$.  Then
$b-2<s<b+1$.  By \eqref{eq:component-cuspidal-translations} and
the disjoint partition \eqref{eqn: tiling-inf-char-7.1},
$B_\mu=B_\nu+1$ and $B_\mu\cap B_\nu=\varnothing$.  Thus each
paired interval $(a,2s-a)$ of $B_\nu$ corresponds to
$(a+1,2s-a+1)$ of $B_\mu$.

\begin{itemize}
\item[(1)] If $s\in\Z$, then $s=b-1$ or $s=b$.  Of $B_\nu$ and $B_\mu$,
the weight set centred at $b$ is $B_\mu$ when $s=b-1$ and $B_\nu$
when $s=b$.  This set is disjoint from $B_{\nu_0}$ since $c_{\nu_0}=t>b$.
If its cardinality were odd, its middle weight would be $b\in B_{\nu_0}$.
In rank two the centre is strictly half-integral by
Remark~\ref{rem:rank-two-motivic-parity}.
Thus $k_\nu=k_\mu$ is even and at least four, so
$h_\nu=h_\mu\ge2$.  Since $a<s$ and both are integers,
\[
 a\le s-1\le b-1,\qquad 2s-a\ge s+1\ge b.
\]

\item[(2)] If $s\in\tfrac12+\Z$, then
$s\in\{b-\tfrac32,b-\tfrac12,b+\tfrac12\}$.  The positive
half-integer $s-a\neq \tfrac12$, since then
$2s-a=a+1$ would belong to both $B_\nu$ and $B_\mu$.  Hence
\begin{equation}\label{eqn: s-half-integer-ineq}
 s-a\ge\tfrac32,\qquad a\le s-\tfrac32\le b-1,\qquad
 2s-a\ge s+\tfrac32\ge b.
\end{equation}
\end{itemize}

In both cases, each paired interval $(a,2s-a)$ of $B_\nu$ has lower
endpoint $a\le b-1$, and the corresponding paired interval
$(a+1,2s-a+1)$ of $B_\mu$ has lower endpoint $a+1\le b$.
If $a+1=b$, then $a=b-1<c$ by \eqref{eqn: exclusion-ineq},
contrary to $c\le b-1$.  Thus both lower endpoints are $<b$, hence
$<c$ by \eqref{eqn: exclusion-ineq}.  Both upper endpoints are
$\ge b>c$.  It follows that $d_\nu(c)=h_\nu$ and $d_\mu(c)=h_\mu$.
In case \textup{(1)} this gives $D\ge h_\nu+h_\mu\ge4$.

It remains to consider the case where no such interval has integral
centres.  Since $b\notin C$, no interval joining consecutive elements
of $C$ meets $(b-1,b+1)$.  The only possible intervals meeting it are
therefore
\[
 [b-\tfrac32,b-\tfrac12]_{\R},\qquad
 [b-\tfrac12,b+\tfrac12]_{\R},\qquad
 [b+\tfrac12,b+\tfrac32]_{\R}.
\]
The intervals $(b-1,b-\tfrac12)$ and $(b-\tfrac12,b)$ lie in
$[c,t]_{\R}\subset K$.  Only the first and second intervals above,
respectively, cover them.  By the definition of $K$, choose indices with
\[
 [c_{\nu_1},c_{\nu_2}]_{\R}=[b-\tfrac32,b-\tfrac12]_{\R},
 \qquad
 [c_\mu,c_{\nu_3}]_{\R}=[b-\tfrac12,b+\tfrac12]_{\R},
\]
The centres $c_{\nu_1},c_{\nu_2},c_{\nu_3}$ are distinct, so the three
indices are distinct.  Case \textup{(2)} gives $d_{\nu_j}(c)=h_{\nu_j}\ge1$
for $1\le j\le3$.
If $t\ge b+1$, the third interval must also occur to cover $(b+\tfrac12,b+1)$.
Write it as $[c_\nu,c_{\nu_4}]_{\R}$ with $\pi_\nu\simeq\pi_{\nu_4}$.
Since $c_{\nu_4}=b+\tfrac32$, the index $\nu_4$ is distinct from the
first three.  Again $d_{\nu_4}(c)=h_{\nu_4}\ge1$, so $D\ge4$.

Otherwise $t=b+\tfrac12$, and $B_{\nu_0}$ contains the pair
$\{b,2t-b\}=\{b,b+1\}$.  The rank $k_{\nu_0}$ is even since $t\notin\Z$.
Rank two would give $w_{\nu_0}=1$, contrary to Lemma~\ref{lem:w11},
so $h_{\nu_0}\ge2$.  For any paired interval $(a,2t-a)$ of $B_{\nu_0}$
other than $(b,b+1)$, the lower endpoint $a$ is an integer with
$a<t=b+\tfrac12$ and $a\ne b$.  Hence $a\le b-1$, so
$a<c$ by \eqref{eqn: exclusion-ineq}, while $2t-a>t>c$.
Only $(b,b+1)$ misses $c$, giving $d_{\nu_0}(c)=h_{\nu_0}-1\ge1$.  Since
$d_{\nu_j}(c)=h_{\nu_j}$ for $1\le j\le3$, the index $\nu_0$ is
distinct from all three.  Hence
\[
 D\ge d_{\nu_0}(c)+\sum_{j=1}^3 d_{\nu_j}(c)\ge1+3=4.
\]
\end{proof}

For disjoint cuspidal weight sets $A,B$, define
\begin{equation}\label{eq:cuspidal-pair-counts}
\begin{aligned}
 p(A,B)&=\min\bigl(\#\{(a,b)\in A\times B:a<b\},
                      \#\{(a,b)\in A\times B:a>b\}\bigr),\\
 I(A,B)&=\sum_{a\in A}d_B(a)+\sum_{b\in B}d_A(b),\qquad
 \sigma(A,B)=p(A,B)-\tfrac12 I(A,B).
\end{aligned}
\end{equation}
Here $d_B(x)$ is defined by the same minimum as
$d_\nu(x)$.  Summing either term in that minimum gives
\begin{equation}\label{eqn: property-positivity-sigma}
 \sum_{a\in A}d_B(a)\le p(A,B),\qquad I(A,B)\le2p(A,B),
 \qquad\sigma(A,B)\ge0.
\end{equation}
Define
\begin{equation}\label{eq:weighted-depth-definition}
 S_\nu=\sum_{x\in U_n\setminus B_\nu}d_\nu(x)
 =\sum_{c\in C}d_\nu(c)+\sum_{\mu\ne\nu}\sum_{b\in B_\mu}d_\nu(b).
\end{equation}
Note that $B_\nu=\{b_1<\cdots<b_{k_\nu}\}$ gives
\begin{equation}\label{eq:weighted-depth-sum}
 S_\nu
 =\sum_{i=1}^{k_\nu-1}
       \min(i,k_\nu-i)(b_{i+1}-b_i-1)=\sum_{i=1}^{h_\nu}(b_{k_\nu+1-i}-b_i)-b_\nu
 \ge w(\pi_\nu)-k_\nu+1.
\end{equation}

Let $r_K$ count the maximal integer intervals of character weights in
$K$. For each isomorphism class of unitary cuspidal representations
occurring among the $\pi_\nu$, choose a representative $\pi$ and partition
$\{c_\nu\in K:\pi_\nu\simeq\pi\}$ into maximal consecutive sequences $\rho$.
Let $t_K$ count these sequences, and denote their lengths and ranks by $l_\rho,k_\rho$.
The weight union
\[
 \Omega_K=(C\cap K)\sqcup\bigsqcup_{c_\nu\in K}B_\nu
\]
Let $\bar c_K$ be its mean, let $\chi_K$ have
Harish--Chandra parameter $\Omega_K-\bar c_K$, and let $\mathfrak c_K$ be
the support represented by
\[
 \left(\boxtimes_{c\in C\cap K}|\cdot|^{c-\bar c_K}\right)
 \boxtimes
 \left(\boxtimes_{c_\nu\in K}\pi_\nu|\det|^{c_\nu-\bar c_K}\right).
\]
Denote
\[
 \mathcal A_K=\operatorname{Fin}_{\ker\chi_K}
                    \mathcal A(\GL_{n_K})_{\mathfrak c_K},\qquad
 G_K=\sum_{c\in C\cap K}\sum_{c_\nu\in K}d_\nu(c),\qquad
 P_K=\sum_{\substack{\nu<\mu\\c_\nu,c_\mu\in K}}p(B_\nu,B_\mu).
\]

For an inducing datum in the Franke filtration of $\mathcal A_K$, put
\begin{equation}\label{eqn: Franke-filtration-K}
 t=(R,\sigma_t,\lambda_t),\quad R=P_{n_1,\ldots,n_s},\quad
 \sigma_t=\boxtimes_{b=1}^sJ(l_b,\tau_b),\quad
 k_b=\operatorname{rank}(\tau_b),\quad n_b=k_bl_b.
\end{equation}
Let $z_b=(\lambda_t)_b$.  By the definition of $J(l_b,\tau_b)$ in
Subsection~\ref{subsec:speh-weights}, its cuspidal exponents in increasing
order are $-(l_b-1)/2+j$ for $0\le j<l_b$.  Twisting by $|\det|^{z_b}$
adds $z_b$ to each exponent, so the smallest is $z_b-(l_b-1)/2$.
Let $a_b$ be the center of the corresponding weight set before subtracting
$\bar c_K$.  By the definition of $\mathfrak c_K$, this same exponent is
$a_b-\bar c_K$.  Thus
\begin{equation}\label{eqn: relation-mean-shift-K}
 a_b-\bar c_K=z_b-\frac{l_b-1}{2},\qquad
 A_{b,j}=A_{b,0}+j\quad(0\le j<l_b),
\end{equation}
where $A_{b,j}$ is the shifted weight set for the cuspidal representation
$\tau_b|\det|^{a_b+j-\bar c_K}$ in block $b$.  It is obtained by adding
$a_b+j$ to each entry of the Harish--Chandra parameter of $\tau_{b,\infty}$.
The sets $A_{b,j}$
are the singletons in $C\cap K$ and the sets $B_\nu$ with $c_\nu\in K$,
each occurring exactly once.

\begin{lemma}\label{lem:component-block-count}
Every datum \eqref{eqn: Franke-filtration-K} has $s\ge r_K+t_K$.
\end{lemma}
\begin{proof}
By equality of cuspidal supports, the centers in a character Speh block
form a consecutive subset of $C\cap K$, hence lie in one maximal
character interval.  If $k_b\ge2$, comparing the cuspidal support of
$J(l_b,\tau_b)|\det|^{z_b}$ with $\mathfrak c_K$ gives, by
\eqref{eqn: relation-mean-shift-K}, an index $\nu_j$ for each $0\le j<l_b$ such that
\[
 \tau_b|\det|^{a_b+j-\bar c_K}
 \simeq \pi_{\nu_j}|\det|^{c_{\nu_j}-\bar c_K}.
\]
Since $\tau_b$ and $\pi_{\nu_j}$ are unitary, comparing the absolute values
of the central characters on positive real scalars gives $c_{\nu_j}=a_b+j$.
Removing the common twist gives $\pi_{\nu_j}\simeq\tau_b$ for every
$0\le j<l_b$.  By \eqref{eq:component-cuspidal-translations},
$B_{\nu_j}=B_{\nu_0}+(c_{\nu_j}-c_{\nu_0})=B_{\nu_0}+j$.
These shifted weight sets have respective centers
\[
 (c_{\nu_0},\ldots,c_{\nu_{l_b-1}})=(a_b,\ldots,a_b+l_b-1).
\]
These centers lie in a single maximal consecutive sequence $\rho$ in $K$
for the isomorphism class of $\tau_b$.  Since all characters and cuspidal factors are covered, every
one of the $r_K+t_K$ intervals or sequences requires a Levi block.
\end{proof}

Let
\[
 G_K=\sum_{c\in C\cap K}\sum_{\nu:\,c_\nu\in K}d_\nu(c),
 \qquad
 P_K=\sum_{\substack{\nu<\mu\\c_\nu,c_\mu\in K}}p(B_\nu,B_\mu).
\]

\begin{lemma}\label{lem:general-component-degree}
Let $K$ be a component containing cuspidal factors. 
For every $F$ with $\chi_{F^\vee}=\chi_K$ on $\mathfrak{sl}_{n_K}$,
\[
 H^{q_K}(\mathfrak{sl}_{n_K},\O_{n_K};
       \mathcal A_K^{K_{n_K,f}}\otimes F)\ne0
 \quad\Longrightarrow\quad d_{n_K}-q_K\ge J_K^0,
\]
where $J_K^0=G_K+\sum_{c_\nu\in K}b_\nu+P_K+r_K+t_K-1$.
\end{lemma}

\begin{proof}
We may take $F$ irreducible.  By Corollary~\ref{cor:franke-degree-reduction}
with $J=\ker\chi_K$, there is a datum $t$ as in \eqref{eqn: Franke-filtration-K}
for which $M(t)\otimes F$ has nonzero cohomology in degree $q_K$.  Choose a contributing
Kostant/K\"unneth summand with permutation $v$ and block degrees $q_b$.
Lemma~\ref{lem:split-centre-cohomology} removes the split-center degrees;
as in Proposition~\ref{prop:degree},
\begin{equation}\label{eqn: cohom-bound-K}
 q_K=\ell_W(v)+\sum_bq_b,\qquad
 d_{n_K}-q_K=\sum_b(d_{n_b}-q_b)
              +\dim\mathfrak n_R-\ell_W(v)+(s-1).
\end{equation}
To compute $\dim\mathfrak n_R-\ell_W(v)$ in \eqref{eqn: cohom-bound-K},
set $\gamma=\Lambda_F+\rho_{n_K}$.  For block $b$, let
$m_b=\sum_{b'<b}n_{b'}$ and denote by $h_{m_b+1}<\cdots<h_{m_b+n_b}$
 the entries of the Harish--Chandra parameter of
$J(l_b,\tau_{b,\infty})|\det|^{z_b}$.  By \eqref{eqn: relation-mean-shift-K},
these entries form the set $\bigsqcup_{j=0}^{l_b-1}(A_{b,j}-\bar c_K)$.
Equation~\eqref{eq:ordered-kostant-parameter-comparison} gives
$h_u-h_w=\gamma_{v(w)}-\gamma_{v(u)}$. Since $\gamma$ is strictly
decreasing, $h_u<h_w$ exactly when $v(u)<v(w)$.  As in the proof of
Proposition~\ref{prop:degree}, $\dim\mathfrak n_R$ counts all pairs of
positions $u,w$ in blocks $b<b'$, and $\ell_W(v)$ counts those with
$v(u)>v(w)$, since $v$ is increasing within each block.  
Thus, by the definition of $p$ in \eqref{eq:cuspidal-pair-counts},
\begin{equation}\label{eqn: dim-diff-blocks}
 \dim\mathfrak n_R-\ell_W(v)
 =\sum_{b<b'}\sum_{i,j}\#\{(x,y)\in A_{b,i}\times A_{b',j}:x<y\}
 \ge\sum_{b<b'}\sum_{i,j}p(A_{b,i},A_{b',j}).
\end{equation}
Here $0\le i<l_b$ and $0\le j<l_{b'}$.  Fix $c\in C\cap K$ and
$\nu$ with $c_\nu\in K$.  Both $\{c\}$ and $B_\nu$ occur exactly once
among the $A_{b,j}$.  They lie in different Levi blocks, since all
sets in block $b$ have cardinality $k_b$, while $|B_\nu|=k_\nu\ge2$.
Thus the sum over $b<b'$ in the lower bound counts this pair once.
In either block order its contribution is
\[
 p(\{c\},B_\nu)
 =\min\bigl(\#\{x\in B_\nu:x<c\},\#\{x\in B_\nu:x>c\}\bigr)
 =d_\nu(c).
\]
Summing over $c$ and $\nu$ gives $G_K$.  Cuspidal pairs in distinct
Levi blocks supply their terms in $P_K$.

Fix a block with $k_b\ge2$, and let $\lambda_1<\cdots<\lambda_{k_b}$
be the Harish--Chandra entries of $\tau_{b,\infty}$.  By
\eqref{eqn: relation-mean-shift-K},
\[
 A_{b,j}=\{a_b+\lambda_r+j:1\le r\le k_b\},\qquad 0\le j<l_b.
\]
For $r<k_b$, the gap $d=\lambda_{r+1}-\lambda_r$ is a positive integer,
since $a_b+\lambda_r$ and $a_b+\lambda_{r+1}$ belong to $U_n$.
If $d<l_b$, then $a_b+\lambda_{r+1}=a_b+\lambda_r+d$ belongs to both
$A_{b,0}$ and $A_{b,d}$, contrary to the disjoint partition
\eqref{eqn: tiling-inf-char-7.1}.  Hence $\lambda_{r+1}-\lambda_r\ge l_b$.
For $0\le i<j<l_b$ it follows that $\lambda_r+i>\lambda_{r'}+j$
exactly when $r>r'$.  Of the $k_b^2$ index pairs $(r,r')$ with
$1\le r,r'\le k_b$, exactly $\binom{k_b}{2}$ satisfy $r>r'$, giving
$p(A_{b,i},A_{b,j})=\binom{k_b}{2}$.
Proposition~\ref{prop:speh-cohomology-bound} gives
\[
 d_{n_b}-q_b\ge b(k_b,l_b)
 =l_b\left\lfloor\frac{k_b^2}{4}\right\rfloor
       +\binom{l_b}{2}\binom{k_b}{2}.
\]
Write $A_{b,j}=B_{\nu_j}$ for $0\le j<l_b$ as in the proof of
Lemma~\ref{lem:component-block-count}.  Since $k_{\nu_j}=k_b$, the
two terms in this bound are
\[
 l_b\left\lfloor\frac{k_b^2}{4}\right\rfloor
 =\sum_{j=0}^{l_b-1}b_{\nu_j},\qquad
 \binom{l_b}{2}\binom{k_b}{2}
 =\sum_{0\le i<j<l_b}p(B_{\nu_i},B_{\nu_j}).
\]
If $k_b=1$, the relative cochains vanish above degree $d_{n_b}$, so
$d_{n_b}-q_b\ge0$.  Here $A_{b,j}=\{a_b+j\}$, and for $i<j$ we have
$p(\{a_b+i\},\{a_b+j\})=\min(1,0)=0$.

Summing the block bounds counts $b_\nu$ once for every $c_\nu\in K$,
since each $B_\nu=A_{b, j}$ for exactly one block $b$ and at exactly one index $j$.
Since $p$ is symmetric, for $\nu<\mu$ with $c_\nu,c_\mu\in K$ the term $p(B_\nu,B_\mu)$
is counted once in the sum over $i<j$ if the two sets are in the same
block, and once in \eqref{eqn: dim-diff-blocks} if they are in distinct
blocks.  Together these terms give $P_K$.  The pairs $(\{c\},B_\nu)$
in \eqref{eqn: dim-diff-blocks} give $G_K$.  Substituting in
\eqref{eqn: cohom-bound-K} and using $s\ge r_K+t_K$ from
Lemma~\ref{lem:component-block-count}, we obtain
\[
 \begin{aligned}
 d_{n_K}-q_K
 &\ge \sum_{c_\nu\in K}b_\nu+P_K+G_K+s-1\\
 &\ge \sum_{c_\nu\in K}b_\nu+P_K+G_K+r_K+t_K-1=J_K^0.
 \end{aligned}
\]
\end{proof}

Define the component deficit
\begin{equation}\label{eq:component-deficit-definition}
\begin{aligned}
 \Phi_K={}&\frac{G_K}{2}-2m_K+2r_K
 +\sum_{c_\nu\in K}\left(b_\nu-2k_\nu+\frac{S_\nu}{2}\right)\\
 &+\sum_{\substack{\nu<\mu\\c_\nu,c_\mu\in K}}\sigma(B_\nu,B_\mu)
 -\sum_{\rho\subset K}(k_\rho-1).
\end{aligned}
\end{equation}

\begin{proposition}\label{lem:global-degree-from-components}
For each component containing cuspidal factors, suppose $J_K$ is a
lower bound for $d_{n_K}-q_K$ whenever
\[
 H^{q_K}(\mathfrak{sl}_{n_K},\O_{n_K};
     \mathcal A_K^{K_{n_K,f}}\otimes F)\ne0,
 \qquad \chi_{F^\vee}=\chi_K\text{ on }\mathfrak{sl}_{n_K}.
\]
Then every nonzero global class of degree $q$ in this support satisfies
\[
 d_n-q\ge2n-2+\sum_K(\Phi_K+J_K-J_K^0).
\]
In particular, $J_K\ge J_K^0-\Phi_K$ for every $K$ implies
$d_n-q\ge2n-2$.
\end{proposition}

\begin{proof}
As before Lemma~\ref{lem:weighted-cuspidal-depth}, take all connected
components of the subset of $\R$
\[
 C\cup\{c_\nu\}_\nu
 \cup\bigcup_{c,c+1\in C}[c,c+1]_{\R}
 \cup\bigcup_{\substack{\pi_\nu\simeq\pi_\mu\\c_\mu=c_\nu+1}}
       [c_\nu,c_\mu]_{\R}.
\]
Let $s$ be the total number of these connected components, including those
consisting only of characters. Order them from right to left. Thus $K$ precedes $K'$ exactly when
$\min K>\max K'$. For their total ranks $n_1,\ldots,n_s$, put
$Q=P_{n_1,\ldots,n_s}$. By Proposition~\ref{lem:separated-cuspidal-parameters},
\[
 \mathcal A_{\mathfrak c}^{K_f}\simeq\mathcal I_Q(0)^{K_f}\otimes\D_y.
\]
Fix a contributing $v\in W^Q$ and a nonzero Levi K\"unneth summand
with degrees $q_1,\ldots,q_s$. Consider $V^{M_Q}_{\mu_v}=\boxtimes_b E_{b,v}$ such that
\[
 q-\ell_W(v)=\sum_bq_b,\qquad
 H^{q_b}\!\left(\mathfrak{sl}_{n_b},\O_{n_b};
       \mathcal A_b^{K_{n_b,f}}\otimes E_{b,v}\right)\ne0
 \quad(1\le b\le s).
\]

For the $b$-th component $K$, we have
$ d_{n_b}-q_b\ge J_b$ with
\[
 J_b=\begin{cases}
 J_K, & K\text{ contains cuspidal factors},\\
 2n_b-2, & K\text{ consists only of characters}.
 \end{cases}
\]
In the first case the bound follows by assumption, since $\mathcal A_b=\mathcal A_K$ and
$\chi_{E_{b,v}^\vee}=\chi_K$ on $\mathfrak{sl}_{n_b}$.
In the second case let $m=n_b$. Then $\mathcal A_b=\mathcal A_{m,1}$
has trivial support, and the highest weight of $E_{b,v}^\vee$ restricts
to zero on $\mathfrak{sl}_m$. For either determinant parity
of $E_{b,v}$, Theorem~\ref{prop:constant} and duality~\eqref{eq:duality}
give $d_m-q_b\ge 2m-2$.
Using $d_n=\sum_b d_{n_b}+\dim\mathfrak n_Q+(s-1)$
(Proposition~\ref{prop:degree}) and $q=\ell_W(v)+\sum_b q_b$, we obtain
\begin{equation}\label{eq:separated-support-degree}
 d_n-q=\sum_{b=1}^s(d_{n_b}-q_b)+\dim\mathfrak n_Q-\ell_W(v)+s-1\ge\sum_{b=1}^sJ_b+\dim\mathfrak n_Q-\ell_W(v)+s-1.
\end{equation}

For symmetric weight sets $A,B$ with centres $a_0>b_0$, the map
\[
 (a,b)\longmapsto(2a_0-a,2b_0-b)
\]
injects the pairs $a<b$ into the pairs $a>b$. Hence
$\#\{(a,b)\in A\times B:a<b\}=p(A,B)$.
The analogous counts for character--cuspidal and
character--character pairs are $d_\nu(c)$ and $0$, respectively.
By the definition of $S_\nu$ in \eqref{eq:weighted-depth-definition},
the partition $U_n=C\sqcup\bigsqcup_\nu B_\nu$ gives
\[
\begin{aligned}
 \sum_\nu S_\nu
 &=\sum_\nu\sum_{c\in C}d_\nu(c)
   +\sum_{\nu<\mu}I(B_\nu,B_\mu),\\
 \dim\mathfrak n_Q-\ell_W(v)+\sum_K(G_K+P_K)
 &=\sum_\nu\sum_{c\in C}d_\nu(c)
   +\sum_{\nu<\mu}p(B_\nu,B_\mu).
\end{aligned}
\]
Using
\eqref{eq:cuspidal-pair-counts}, we obtain
\[
\begin{aligned}
 &\dim\mathfrak n_Q-\ell_W(v)+\sum_K(G_K+P_K)-\frac12\sum_\nu S_\nu\\
 &\quad=\frac12\sum_\nu\sum_{c\in C}d_\nu(c)
       +\sum_{\nu<\mu}\left(p(B_\nu,B_\mu)-\frac12 I(B_\nu,B_\mu)\right)\\
 &\quad=\frac12\sum_\nu\sum_{c\in C}d_\nu(c)
       +\sum_{\nu<\mu}\sigma(B_\nu,B_\mu).
\end{aligned}
\]
By \eqref{eqn: property-positivity-sigma}, 
\[
 \frac12\sum_K\sum_{\nu:\,c_\nu\in K}\sum_{c\in C\setminus K}d_\nu(c)
 +\sum_K\sum_{\substack{\nu<\mu\\c_\nu\in K,\,c_\mu\notin K}}
       \sigma(B_\nu,B_\mu)\ge0.
\]
On the other hand, by definition of $G_K$, 
\[
\sum_{\nu: c_\nu\in K}\sum_{c\in C}d_\nu(c)=\sum_{\nu:\,c_\nu\in K}\sum_{c\in C\setminus K}d_\nu(c)+G_K
\]

We obtain
\begin{equation}\label{eq:component-pair-bound}
 \dim\mathfrak n_Q-\ell_W(v)+\sum_KP_K\ge
 \frac12\sum_\nu S_\nu-\frac12\sum_KG_K
 +\sum_K\sum_{\substack{\nu<\mu\\c_\nu,c_\mu\in K}}
          \sigma(B_\nu,B_\mu).
\end{equation}

Let $r_C$ be the total number of maximal intervals of character weights
and $s_C$ the number of character-only components.
For a maximal sequence of consecutive centres $\rho=(\nu_1,\ldots,\nu_l)$
for $\pi$, put $B=B_{\nu_1}$. Let $\lambda_1<\cdots<\lambda_{k_\rho}$
be the Harish--Chandra parameter of $\pi$. Since
$\pi_{\nu_j}\simeq\pi$, subtracting $c_{\nu_j}$ from the weights in
$B_{\nu_j}$ gives these same entries by Subsection~\ref{subsec:speh-weights}.
Using $c_{\nu_j}=c_{\nu_1}+j-1$, we obtain
\[
 b_{\nu_j,i}=c_{\nu_j}+\lambda_i
            =b_{\nu_1,i}+j-1
 \qquad(1\le i\le k_\rho,\ 1\le j\le l).
\]
Thus $B_{\nu_j}=B+(j-1)$. Taking the union over $j$ gives
\[
 \bigcup_{j=1}^{l}B_{\nu_j}
 =\{b+j-1:b\in B,\ 1\le j\le l\}
 =\bigcup_{b\in B}\{b,b+1,\ldots,b+l-1\},\qquad |B|=k_\rho.
\]
Thus $U_n\setminus C$ is a union of at most
$\sum_\rho k_\rho$ integer intervals, which implies
\[
 r_C\le1+\sum_\rho k_\rho.
\]
Note that by definition, we have $ r_C=s_C+\sum_Kr_K$, hence
\begin{equation}\label{eq:character-component-count-bound}
 1-s_C\ge\sum_Kr_K-\sum_\rho k_\rho.
\end{equation}

The character-only components have total rank $n-\sum_Kn_K$, and
the total number of connected components
\[
 s=s_C+\#\{K:K\text{ is a connected component containing some }c_\nu\}.
\]
Since $J_b=2n_b-2$ on each character-only component, we have
\[
 \sum_{b=1}^s J_b
 =2\left(n-\sum_K n_K\right)-2s_C+\sum_K J_K.
\]
Hence 
\[
 \sum_{b=1}^s J_b+(s-1)-(2n-2)
 =1-s_C+\sum_K(J_K-2n_K+1).
\]
For each component $K$ containing a cuspidal factor, write
$J_K=(J_K-J_K^0)+J_K^0$. Substituting $J_K^0$ from
Lemma~\ref{lem:general-component-degree} and using the definition of $n_K$ gives
\[
 J_K-2n_K+1
 =(J_K-J_K^0)
   +\left(G_K+\sum_{\nu:\,c_\nu\in K}b_\nu+P_K+r_K+t_K-1\right)
 -2\left(m_K+\sum_{\nu:\,c_\nu\in K}k_\nu\right)+1.
\]
Collecting terms in \eqref{eq:separated-support-degree}, we obtain
\[
\begin{aligned}
 d_n-q-(2n-2)\ge{}&
 \dim\mathfrak n_Q-\ell_W(v)+\sum_KP_K+1-s_C\\
 &+\sum_K\left(
      G_K-2m_K+\sum_{\nu:\,c_\nu\in K}(b_\nu-2k_\nu)
      +r_K+t_K+J_K-J_K^0\right).
\end{aligned}
\]
Apply the lower bound \eqref{eq:component-pair-bound} and
\eqref{eq:character-component-count-bound}, 
grouping the sums by $K$ gives
\[
\begin{aligned}
 &d_n-q-(2n-2)\\
 &\quad\ge\sum_K\Biggl[
   \biggl(\frac12\sum_{\nu:\,c_\nu\in K}S_\nu-\frac{G_K}{2}
        +\sum_{\substack{\nu<\mu\\c_\nu,c_\mu\in K}}\sigma(B_\nu,B_\mu)\biggr)
   +\biggl(r_K-\sum_{\rho\subset K}k_\rho\biggr)\\
 &\qquad\qquad
   +G_K-2m_K+\sum_{\nu:\,c_\nu\in K}(b_\nu-2k_\nu)
   +r_K+t_K+J_K-J_K^0\Biggr].
\end{aligned}
\]
Since $t_K=\sum_{\rho\subset K}1$, we have
$t_K-\sum_{\rho\subset K}k_\rho=-\sum_{\rho\subset K}(k_\rho-1)$.
Collecting terms gives
\[
\begin{aligned}
 d_n-q-(2n-2)
 &\ge\sum_K\Biggl[\frac{G_K}{2}-2m_K
     +\sum_{\nu:\,c_\nu\in K}\left(b_\nu-2k_\nu+\frac{S_\nu}{2}\right)
     +\sum_{\substack{\nu<\mu\\c_\nu,c_\mu\in K}}\sigma(B_\nu,B_\mu)\\
 &\qquad\qquad
     +2r_K-\sum_{\rho\subset K}(k_\rho-1)+J_K-J_K^0\Biggr]\\
 &=\sum_K(\Phi_K+J_K-J_K^0),
\end{aligned}
\]
as required.
\end{proof}

\begin{lemma}\label{lem:component-arithmetic-margin}
Let $K$ be a component containing $\ell_K\ge1$ cuspidal factors.
With $\Phi_K$ as in Proposition~\ref{lem:global-degree-from-components},
\begin{equation}\label{eq:component-arithmetic-margin}
 \Phi_K\ge \frac{G_K}{2}-2m_K+2r_K+\ell_K.
\end{equation}
\end{lemma}

\begin{proof}
For a cuspidal factor $\pi_\nu$ of rank $k=k_\nu$,
\eqref{eq:weighted-depth-sum} gives $S_\nu\ge w(\pi_\nu)-k+1$.
If $k$ is odd, its integral centre makes $w(\pi_\nu)$ even.
The rank--weight bounds of Lemma~\ref{lem:w11} give
\[
 b_\nu-3k+\frac{S_\nu}{2}
 \ge\left\lfloor\frac{k^2}{4}\right\rfloor
     -\frac{7k}{2}+\frac{w(\pi_\nu)+1}{2}
 \ge
 \begin{cases}
 0,&k=2,\,4,\\
 3,&k=3,\\
 (k-6)(k-8)/4,&k\ge6\text{ even},\\
 (k-7)^2/4,&k\ge5\text{ odd}
 \end{cases}
 \ge0.
\]
Hence
\begin{equation}\label{eq:cuspidal-constant-margin}
 b_\nu-2k_\nu+\frac{S_\nu}{2}\ge k_\nu.
\end{equation}
For a maximal sequence $\rho$ of $l$ consecutive centres for
a cuspidal representation of rank $k$,
\[
 \sum_{\nu\in\rho}\left(b_\nu-2k+\frac{S_\nu}{2}\right)-(k-1)
 \ge (l-1)k+1\ge l.
\]
Summing over $\rho\subset K$ and using $\sigma(B_\nu,B_\mu)\ge0$ in
the formula for $\Phi_K$ in Proposition~\ref{lem:global-degree-from-components} gives
\eqref{eq:component-arithmetic-margin}.
\end{proof}

\begin{lemma}\label{lem:rank-two-component-budget}
Let $K$ be a component containing cuspidal factors, with
$m_K>0$ and $H_K\le3$.
\begin{enumerate}[label=(\roman*),ref=\roman*,leftmargin=2em]
\item\label{it:component-paired-intervals}
Every $c\in C\cap K$ lies in each
paired interval of every $B_\nu$ with $c_\nu\in K$, and
$G_K=m_KH_K$.
\item\label{it:component-rank-two-formula}
If all cuspidal ranks are two, put $W_K=\sum_{c_\nu\in K}w_\nu$ and
let $N_K$ count the unordered pairs of their weight intervals for which
one strictly contains the other. Then $\ell_K=H_K\le3$ and
\begin{equation}\label{eq:rank-two-component-budget}
 \Phi_K=\frac{G_K}{2}-2m_K+N_K
             +\frac{W_K-7\ell_K}{2}+2r_K-t_K.
\end{equation}
\item\label{it:component-even-ranks}
Every cuspidal factor $\nu$ with $c_\nu\in K$ has even rank: $k_\nu=2h_\nu$.
\end{enumerate}
\end{lemma}

\begin{proof}
For~\textup{(\ref{it:component-paired-intervals})}, fix
$c\in C\cap K$. Lemma~\ref{lem:weighted-cuspidal-depth}
gives
\[
 H_K\le\sum_\nu d_\nu(c)\le\sum_\nu h_\nu=H_K.
\]
Thus $d_\nu(c)=h_\nu$ for every $\nu$, and $G_K=mH_K$.
In particular, every internal paired interval contains $c$.

For~\textup{(\ref{it:component-rank-two-formula})}, suppose all cuspidal
factors have rank two. Then $h_\nu=1$ for each $\nu$, so $\ell_K=H_K\le3$.
Fix $\nu<\mu$ and write $B_\nu=\{a_1<a_2\}$ and $B_\mu=\{b_1<b_2\}$.
Since $m>0$, part~\textup{(\ref{it:component-paired-intervals})} gives
a common character weight in the two paired intervals. Their endpoints are
distinct, so after exchanging the two factors if necessary, the four
comparisons in \eqref{eq:cuspidal-pair-counts} give
\[
 \begin{array}{c|ccc}
 \text{Endpoint order} & \#\{(i,j):a_i<b_j\} & \#\{(i,j):a_i>b_j\}
     & p(B_\nu,B_\mu)\\ \hline
 a_1<b_1<b_2<a_2 & 2 & 2 & 2\\
 a_1<b_1<a_2<b_2 & 3 & 1 & 1
 \end{array}
\]
The endpoints lying inside the other interval are $b_1,b_2$ in the
first row and $b_1,a_2$ in the second. Thus $I(B_\nu,B_\mu)=2$
in both cases, and $\sigma(B_\nu,B_\mu)=p(B_\nu,B_\mu)-1$ equals
$1$ precisely when one interval contains the other and $0$ otherwise.
Summing over $\nu<\mu$ counts each such unordered pair once, giving
$\sum_{\nu<\mu}\sigma(B_\nu,B_\mu)=N_K$.
Since $b_\nu=1$ and $S_\nu=w(\pi_\nu)-1$, substitution in the formula for
$\Phi_K$ in Proposition~\ref{lem:global-degree-from-components} proves
\eqref{eq:rank-two-component-budget}.

For~\textup{(\ref{it:component-even-ranks})}, suppose $k_\eta=2h_\eta+1$ for some $c_\eta\in K$. Then
$c_\eta=b_{\eta,h_\eta+1}\in U_n\setminus C$.
Choose $c\in C\cap K$, which exists since by assumption.
The segment from $c_\eta$ to $c$ lies in $K$ and contains
consecutive integers $a,a+1$ with exactly one in $C$.
Thus $[a,a+1]_{\R}$ is not a character interval.

Fix $x\in\{a+\tfrac14,a+\tfrac34\}$. Then $x\in[a,a+1]_{\R}\subset K$.
If $x\in[d,d+1]_{\R}$ for some $d,d+1\in C$, then
\[
 d=\lfloor x\rfloor=a,\qquad d+1=a+1,\qquad \{a,a+1\}\subset C,
\]
contrary to the choice of $a$. The remaining union in the definition
of $K$ therefore gives an interval $[c_\nu,c_\mu]_{\R}$ containing $x$,
with $\pi_\nu\simeq\pi_\mu$ and $c_\mu=c_\nu+1$.  This interval lies in the same
connected component as $x$, and its endpoints lie in $\tfrac12\Z$. Hence
\[
 c_\nu<x<c_\mu=c_\nu+1,\qquad \pi_\nu\simeq\pi_\mu,\qquad
 2\lfloor k_\nu/2\rfloor\le H_K\le3.
\]
Thus $k_\nu=k_\mu\in\{2,3\}$. Rank three is excluded by $\#(\{a,a+1\}\cap C)=1$ and
\[
 k_\nu=3\ \Longrightarrow\
 (c_\nu,c_\mu)=(\lfloor x\rfloor,\lceil x\rceil)=(a,a+1)
 \in (U_n\setminus C)^2.
\]
In rank two, $c_\nu\in(x-1,x)\cap(\tfrac12+\Z)$, so
\[
 \begin{array}{c|cc}
 x & c_\nu & c_\mu=c_\nu+1\\ \hline
 a+\tfrac14 & a-\tfrac12 & a+\tfrac12\\
 a+\tfrac34 & a+\tfrac12 & a+\tfrac32
 \end{array}
\]
The three distinct centres in the table give three rank-two factors in $K$.
Together with the odd-rank factor $\eta$, they give
\[
 H_K=\sum_{\nu:\,c_\nu\in K} h_\nu\ge 1+1+1+h_\eta\ge4,
\]
contrary to $H_K\le3$. Hence $k_\nu=2h_\nu$ for every $c_\nu\in K$.
\end{proof}

\begin{lemma}\label{lem:component-depth-at-most-two}
Let $K$ be a component containing cuspidal factors, with
$m_K>0$ and $H_K\le2$. Let $J_K^0$ and $\Phi_K$ be as in
Lemma~\ref{lem:general-component-degree} and
Proposition~\ref{lem:global-degree-from-components}, respectively.
For every algebraic representation $F$ of $\GL_{n_K}$ whose dual
has infinitesimal character $\chi_K$ on $\mathfrak{sl}_{n_K}$,
and for either determinant parity,
\[
 H^q(\mathfrak{sl}_{n_K},\O_{n_K};
       \mathcal A_K^{K_{n_K,f}}\otimes F)\ne0
 \quad\Longrightarrow\quad d_{n_K}-q\ge J_K^0-\Phi_K.
\]
\end{lemma}

\begin{proof}
By Lemma~\ref{lem:rank-two-component-budget}\textup{(\ref{it:component-paired-intervals})}
and~\textup{(\ref{it:component-even-ranks})},
$G_K=m_KH_K$ and $k_\nu=2h_\nu$ for every $c_\nu\in K$.
Since $1\le\ell_K\le\sum_{\nu:\,c_\nu\in K}h_\nu=H_K\le2$,
there are one or two cuspidal factors.

\emph{One cuspidal factor of rank $2$ or $4$.}
Let $\nu$ be the unique cuspidal index with $c_\nu\in K$. Put
$h=h_\nu$, $k=k_\nu=2h\in\{2,4\}$, and
$B_\nu=\{\beta_1<\cdots<\beta_{2h}\}$, where $\beta_i=b_{\nu,i}$. There is no interval joining cuspidal
centres, so $C\cap K=I=\{a,\ldots,a+m_K-1\}$ and
$\min I\le c_\nu\le\max I$. By \eqref{eqn: tiling-inf-char-7.1},
$I\cap B_\nu=\varnothing$. Here
\[
 r_K=t_K=1,\qquad
 J_K^0=hm_K+h^2+1.
\]
Apply Theorem~\ref{lem:one-copy-internal-bound} with
$\pi=\pi_\nu$, $B=B_\nu$, $c=c_\nu$, and $m=|I|=m_K$.
Then $m+k=n_K$, $\bar c=\bar c_K$, $\chi=\chi_K$, and
$\mathcal A=\mathcal A_K^{K_{n_K,f}}$. Thus for $F$ as in the statement,
\[
 H^q(\mathfrak{sl}_{n_K},\O_{n_K};\mathcal A_K^{K_{n_K,f}}\otimes F)\ne0
 \quad\Longrightarrow\quad d_{n_K}-q\ge(h+2)m_K+h^2-3.
\]
By \eqref{eq:component-arithmetic-margin},
\[
 (h+2)m_K+h^2-3-(J_K^0-\Phi_K)
 =2m_K-4+\Phi_K
 \ge 2m_K-4+\frac{hm_K}{2}-2m_K+3
 =\frac{hm_K}{2}-1\ge0.
\]
For the pattern $(4)$, $h=2$, so the theorem gives $d_{n_K}-q\ge4m_K+1$
and the last inequality follows from $m_K\ge1$.
For $(2)$, $h=1$ and the strictly half-integral centre requires $m_K\ge2$.

\emph{Case two cuspidal factors.}
Both factors must have rank two. There is at most one interval joining the
two strictly half-integral centres. If it meets $[a,b]_{\R}$
for a maximal interval of character weights $[a,b]$, its unique
integer belongs to $[a,b]$. Thus it meets at most one such real
interval. By connectedness, $C\cap K=I$ is a single maximal
interval of character weights, so $r_K=1$. By
\eqref{eq:rank-two-component-budget} and
Lemma~\ref{lem:general-component-degree},
\[
 J_K^0=2m_K+3+N_K+t_K,\qquad \Phi_K=\frac{W_K}{2}-m_K-5+N_K-t_K.
\]
If $\Phi_K\ge0$, Lemma~\ref{lem:general-component-degree} gives
\[
 H^q(\mathfrak{sl}_{n_K},\O_{n_K};\mathcal A_K^{K_{n_K,f}}\otimes F)\ne0
 \quad\Longrightarrow\quad d_{n_K}-q\ge J_K^0\ge J_K^0-\Phi_K.
\]
Assume $\Phi_K<0$.
We verify the centre condition of Theorem~\ref{prop:two-cuspidal-bound}.
For $t_K=2$ there is no interval joining the centres. By connectedness,
$K=[\min I,\max I]_{\R}$, so both strictly half-integral centres lie in
$(\min I,\max I)$.
For $t_K=1$ the shifted weight sets are
\[
 \{b,b+w\},\qquad\{b+1,b+1+w\},\qquad w\ge11.
\]
The interval joining their centres contains the unique integer
$a=b+(w+1)/2$, and connectedness gives $a\in I$.
By Lemma~\ref{lem:rank-two-component-budget}\textup{(\ref{it:component-paired-intervals})},
\[
 I\subset (b,b+w)\cap(b+1,b+1+w)\cap\Z
 =(b+1,b+w)\cap\Z
 =\{b+2,\ldots,b+w-1\}.
\]
The centres are $a-\tfrac12$ and $a+\tfrac12$.
If either centre lies outside $[\min I,\max I]_{\R}$, then
$\min I=a$ or $\max I=a$. Since $a=b+(w+1)/2$, we have
\[
\begin{aligned}
 \min I=a &\quad\Longrightarrow\quad
 m_K=|I|\le (b+w-1)-a+1=\frac{w-1}{2},\\
 \max I=a &\quad\Longrightarrow\quad
 m_K=|I|\le a-(b+2)+1=\frac{w-1}{2}.
\end{aligned}
\]
Since $W_K=2w$ and $N_K=0$, this gives
$\Phi_K=w-m_K-6\ge(w-11)/2\ge0$, contrary to $\Phi_K<0$.
Thus for $t_K=1$ both centres also lie in $[\min I,\max I]_{\R}$.

Order the two factors so that $c_1\ge c_2$, and write their
shifted weight sets as $B_i=\{b_i,b_i+w_i\}$, with $c_i=b_i+w_i/2$.
By \eqref{eqn: tiling-inf-char-7.1}, $I,B_1,B_2$ are pairwise disjoint,
and the centre condition just proved gives $\min I\le c_2\le c_1\le\max I$.
Take $m=|I|=m_K$ in Theorem~\ref{prop:two-cuspidal-bound}. Then
$n_K=m_K+4$. The mean $\bar c_K$ of the shifted weights in
$I\sqcup B_1\sqcup B_2$ is
\[
 \bar c_K=\frac{\sum_{x\in I}x+\sum_{x\in B_1}x+\sum_{x\in B_2}x}{n_K}
         =\frac{\sum_{x\in I}x+2c_1+2c_2}{m_K+4}
         =\bar c
\]
by \eqref{eq:two-cuspidal-centering}. Hence
\eqref{eq:two-cuspidal-support}--\eqref{eq:two-cuspidal-automorphic-module}
give $\chi=\chi_K$ and $\mathcal A=\mathcal A_K^{K_{n_K,f}}$.
Thus for $F$ as in the statement, the theorem gives
\[
 H^q(\mathfrak{sl}_{n_K},\O_{n_K};\mathcal A_K^{K_{n_K,f}}\otimes F)\ne0
 \quad\Longrightarrow\quad d_{n_K}-q\ge4m_K-4.
\]
By Remark~\ref{rem:rank-two-motivic-parity}, both rank-two factors have
odd motivic weights. Let $w_i=w(\pi_i)=2a_i+1$ with $a_i\in\Z$
for $i=1,2$. Then
\[
 \frac{W_K}{2}=\frac{w_1+w_2}{2}
 =\frac{(2a_1+1)+(2a_2+1)}{2}=a_1+a_2+1\in\Z.
\]
Since $\Phi_K<0$, we have
\[
 \frac{W_K}{2}-m_K-5+N_K-t_K=\Phi_K<0
 \quad\Longrightarrow\quad m_K>\frac{W_K}{2}-5+N_K-t_K.
\]
Both sides of the strict inequality are integers, so
\[
 m_K\ge\left(\frac{W_K}{2}-5+N_K-t_K\right)+1
      =\frac{W_K}{2}-4+N_K-t_K.
\]
Using this bound, $W_K\ge22$ and $t_K\le2$, we find
\[
 (4m_K-4)-(J_K^0-\Phi_K)=W_K/2+m_K-12-2t_K\ge W_K-16+N_K-3t_K\ge0.
\]
Both internal estimates apply to $\mathcal A_K^{K_{n_K,f}}$  and to every $F$ in the statement.
\end{proof}

\begin{proposition}\label{lem:small-rank-two-component}
Let $K$ be a component containing cuspidal factors, with
$m_K>0$ and $H_K\le3$. Let $J_K^0$ and $\Phi_K$ be as in
Lemma~\ref{lem:general-component-degree} and
Proposition~\ref{lem:global-degree-from-components}, respectively.
For every algebraic representation $F$ of $\GL_{n_K}$ whose dual
has infinitesimal character $\chi_K$ on $\mathfrak{sl}_{n_K}$,
and for either determinant parity,
\[
 H^q(\mathfrak{sl}_{n_K},\O_{n_K};
       \mathcal A_K^{K_{n_K,f}}\otimes F)\ne0
 \quad\Longrightarrow\quad d_{n_K}-q\ge J_K^0-\Phi_K.
\]
\end{proposition}

\begin{proof}
For $H_K\le2$, the assertion follows from Lemma~\ref{lem:component-depth-at-most-two}.
It remains to consider $H_K=3$. By
Lemma~\ref{lem:rank-two-component-budget}\textup{(\ref{it:component-paired-intervals})}
and~\textup{(\ref{it:component-even-ranks})}, $G_K=3m_K$ and
$k_\nu=2h_\nu$ for every $c_\nu\in K$.
Fix $c\in C\cap K$. For $\nu\ne\mu$ with $c_\nu,c_\mu\in K$, every paired
interval from either $B_\nu$ or $B_\mu$ contains $c$.
Let $(a,b)$ and $(a',b')$ be one such pair of intervals.
Their four endpoints are distinct since $B_\nu\cap B_\mu=\varnothing$.
Exchanging the intervals if necessary, assume $a<a'$. Since both
contain $c$, there are two possible orders.
\[
\begin{array}{c|l}
 \text{Order} & \text{Endpoints in the other interval}\\ \hline
 a<a'<c<b<b' & a'\in(a,b),\quad b\in(a',b')\\
 a<a'<c<b'<b & a'\in(a,b),\quad b'\in(a,b)
\end{array}
\]
Since both ranks
are even, every weight in either set is an endpoint of exactly one
paired interval in that set. There are $h_\nu h_\mu$ choices of one paired
interval from $B_\nu$ and one from $B_\mu$. By the table, each choice
contributes exactly two to the sum in \eqref{eq:cuspidal-pair-counts}.
Hence $I(B_\nu,B_\mu)=2h_\nu h_\mu$.
Using \eqref{eq:weighted-depth-definition} and \eqref{eq:cuspidal-pair-counts}
with $d_\nu(x)\ge0$, we obtain
\[
\begin{aligned}
 \sum_{\nu:\,c_\nu\in K} S_\nu
 &=\sum_{\nu:\,c_\nu\in K}\left(\sum_{c\in C}d_\nu(c)
        +\sum_{\mu\ne\nu}\sum_{b\in B_\mu}d_\nu(b)\right)\\
 &\ge\sum_{\nu:\,c_\nu\in K}\sum_{c\in C\cap K}d_\nu(c)
   +\sum_{\substack{\nu<\mu\\c_\nu,c_\mu\in K}}
      \left(\sum_{b\in B_\mu}d_\nu(b)+\sum_{a\in B_\nu}d_\mu(a)\right)\\
 &=G_K+\sum_{\substack{\nu<\mu\\c_\nu,c_\mu\in K}}I(B_\nu,B_\mu)\\
 &=3m_K+\sum_{\substack{\nu<\mu\\c_\nu,c_\mu\in K}}2h_\nu h_\mu
  =3m_K+9-\sum_{\nu:\,c_\nu\in K}h_\nu^2.
\end{aligned}
\]
By assumption 
\[
H_K=\sum_{\nu:\,c_\nu\in K}h_\nu=3.
\]
Using $b_\nu=h_\nu^2$ and
$\sum_{\rho\subset K}(k_\rho-1)\le\sum_{\nu:\,c_\nu\in K}(2h_\nu-1)$ in
the formula for $\Phi_K$ in Proposition~\ref{lem:global-degree-from-components}, we obtain
\[
\begin{aligned}
 \Phi_K
 &\ge-\frac{m_K}{2}+\sum_{\nu:\,c_\nu\in K}(h_\nu^2-4h_\nu)
 +\frac{3m_K+9-\sum_{\nu:\,c_\nu\in K} h_\nu^2}{2}
 +2r_K-\sum_{\nu:\,c_\nu\in K}(2h_\nu-1)\\
 &=m_K-\frac{27}{2}+2r_K+\frac12\sum_{\nu:\,c_\nu\in K}(h_\nu^2+2).
\end{aligned}
\]
With $G_K=3m_K$,  \eqref{eq:component-arithmetic-margin} gives
\[
 2\Phi_K\ge -m_K+4r_K+2\ell_K.
\]
Adding this to the preceding bound for $\Phi_K$,
we obtain, since $\ell_K=\sum_{\nu:\,c_\nu\in K}1$,
\[
\begin{aligned}
 3\Phi_K&\ge 6r_K-\frac{27}{2}
       +\frac12\sum_{\nu:\,c_\nu\in K}(h_\nu^2+2)+2\ell_K,\\
 6\Phi_K&\ge 12r_K-27+\sum_{\nu:\,c_\nu\in K}(h_\nu^2+2)+4\ell_K\\
        &=12r_K-27+\sum_{\nu:\,c_\nu\in K}(h_\nu^2+6).
\end{aligned}
\]
Since
$h^2+6-5h=(h-2)(h-3)\ge0$ for positive integers $h$,
\[
 6\Phi_K\ge 12r_K-27+5H_K=12(r_K-1)\ge0.
\]
Thus $J_K^0\ge J_K^0-\Phi_K$, and Lemma~\ref{lem:general-component-degree}
gives the proposition.
\end{proof}

\begin{theorem}\label{prop:general-mixed-support-vanishing}
Every full-level cuspidal support $\mathfrak c$ compatible with $J_0$
satisfies
\[
 H^q(Y_n;\C)_{\mathfrak c}=0\qquad(q>d_n-2n+2).
\]
\end{theorem}

\begin{proof}
The trivial support is fixed by contragredience, so
Theorem~\ref{prop:constant} and duality~\eqref{eq:duality} give
the assertion for both determinant parities.
Fix a nontrivial support. For each component $K$ containing
cuspidal factors, \eqref{eq:component-arithmetic-margin} gives
\[
 \Phi_K\ge
 \begin{cases}
 \ell_K,&m_K=0,\\
 2r_K+\ell_K,&m_K>0,\ H_K\ge4,
 \end{cases}
\]
since in the second case Lemma~\ref{lem:weighted-cuspidal-depth} gives
\[
 \sum_{\nu:\,c_\nu\in K}d_\nu(c)\ge\min(H_K,4)=4
 \qquad(c\in C\cap K).
\]
Summing over the $m_K=|C\cap K|$ character weights gives
\[
 G_K=\sum_{c\in C\cap K}\sum_{\nu:\,c_\nu\in K}d_\nu(c)
     \ge\sum_{c\in C\cap K}4=4m_K.
\]
For every $q_K\ge0$ and every irreducible algebraic representation $F$
of $\GL_{n_K}$ with $\chi_{F^\vee}=\chi_K$ on $\mathfrak{sl}_{n_K}$,
Lemma~\ref{lem:general-component-degree} gives
\[
 H^{q_K}(\mathfrak{sl}_{n_K},\O_{n_K};\mathcal A_K^{K_{n_K,f}}\otimes F)\ne0
 \quad\Longrightarrow\quad d_{n_K}-q_K\ge J_K^0.
\]
Thus $J_K=J_K^0$ satisfies the hypothesis of
Proposition~\ref{lem:global-degree-from-components}. With this choice,
the two bounds above give
\[
 \Phi_K+J_K-J_K^0=\Phi_K\ge
 \begin{cases}
  \ell_K,&m_K=0,\\
  2r_K+\ell_K,&m_K>0,\ H_K\ge4
 \end{cases}
 \ \ge0.
\]
If $m_K>0$ and $H_K\le3$, Proposition~\ref{lem:small-rank-two-component}
gives the internal bound $J_K=J_K^0-\Phi_K$. Thus if
$H^q(Y_n;\C)_{\mathfrak c}\ne0$, Proposition~\ref{lem:global-degree-from-components} gives
\[
 d_n-q\ge2n-2+\sum_K(\Phi_K+J_K-J_K^0)\ge2n-2.
\]
Hence $q\le d_n-2n+2$, so
\[
 H^q(Y_n;\C)_{\mathfrak c}=0\qquad(q>d_n-2n+2).
\]
\end{proof}

\begin{corollary}[Theorem~\ref{thm:main}]\label{cor:main-theorem}
For every $n\ge2$ and every $0\le i<n-1$,
\[
 H^{N_n-i}(\Gamma_n;\Q)=0.
\]
\end{corollary}

\begin{proof}
Fix $n\ge2$ and $0\le i<n-1$, and put $q=N_n-i$.  Since
\[
 q-(d_n-2n+2)=n-1-i>0,
\]
Theorem~\ref{prop:general-mixed-support-vanishing} gives
$H^q(Y_n;\C)_{\mathfrak c}=0$ for every full-level cuspidal support
compatible with $J_0$.  Franke's comparison
(Theorem~\ref{thm:franke-comparison}) and
\eqref{eq:shapiro-decomposition} identify $H^q(\Gamma_n;\C)$ with
the direct sum of these summands.  Hence
\[
 H^{N_n-i}(\Gamma_n;\C)=0.
\]
 Consequently
\[
 H^{N_n-i}(\Gamma_n;\Q)\otimes_\Q\C
 \cong H^{N_n-i}(\Gamma_n;\C)=0.
\]
Faithfulness of scalar extension proves the theorem.
\end{proof}

\section{Residue quotients and kernels of coefficient maps}
\label{sec:residue-quotients}

In an ordered list of inducing factors, a \emph{cut after $i$} splits
the block containing the consecutive factors $\tau_i,\tau_{i+1}$ into
\begin{equation}\label{eq:cut-after-index}
 (\tau_a,\ldots,\tau_b)\ \longmapsto\
 \bigl((\tau_a,\ldots,\tau_i),\,(\tau_{i+1},\ldots,\tau_b)\bigr)
 \qquad(a\le i<b).
\end{equation}
For a set $D$ of cuts in an interval of characters, we identify the
inducing module $\mathcal M_D$ and bound its cohomological degree.
For two copies of a cuspidal representation $\pi$, we determine the
kernel of the coefficient map onto $\mathcal A_{2,\pi}^{K_{2k,f}}$.
Section~\ref{sec:internal-cohomology} uses both calculations.
All coefficient modules retain finite derivatives in every spectral
variable unless a kernel condition is explicitly imposed.

\subsection{Cuspidal supports of Speh representations}
\label{subsec:speh-support-modules}

We consider the full automorphic module with the
cuspidal support of a Speh representation.

For a regular algebraic unitary cuspidal representation $\pi$ of
$\GL_k(\A)$, with $k\ge1$, and an integer $\ell\ge1$, put $N=k\ell$.
Let $\mathfrak c_{\ell,\pi}$ be the cuspidal support with Levi
$\GL_k^\ell$ and inducing representation
\begin{equation}
\label{eqn: cuspidal-support-fix-cusp}
 \mathop{\boxtimes}_{i=1}^{\ell}
              \pi|{\det}|^{(\ell+1-2i)/2}.
\end{equation}
Let $\chi_{\ell,\pi}:\mathfrak Z_{\GL_N}^1\to\C$ be the
infinitesimal character of $J(\ell,\pi)_\infty$ on $\mathfrak{sl}_N$.
Using the notation of Subsections~\ref{subsec:automorphic-comparison}
and~\ref{subsec:topological-trivial-datum}, define
\[
 \mathcal A_{\ell,\pi}
 :=\operatorname{Fin}_{\ker\chi_{\ell,\pi}}
       \mathcal A(\GL_N)_{\mathfrak c_{\ell,\pi}}.
\]
The first index counts the inducing factors and the second specifies
the cuspidal representation, as in $J(\ell,\pi)$.
This module includes the Eisenstein coefficients and their finite
derivatives in the relative spectral variables $x_i$, which vary the
exponent of the $i$-th factor from $(\ell+1-2i)/2$ to
$(\ell+1-2i)/2+x_i$, with $\sum_{i=1}^\ell x_i=0$.
Full-level invariants are written $\mathcal A_{\ell,\pi}^{K_{N,f}}$.

For $\ell=1$, the support is represented by $(\GL_k,\pi)$.
The Eisenstein operator for $P=\GL_k$ is the identity, and the relation
gives $x_1=0$. Thus cuspidal multiplicity
one for $\GL_k$ gives
\[
 \mathcal A_{1,\pi}\simeq\pi.
\]

For the trivial character $\pi=1$ of $\GL_1(\A)$ and $\ell=m$,
the support is ${\rm triv}$ from
Subsection~\ref{subsec:topological-trivial-datum}.  Thus
\[
 \mathcal A_{m,1}
   =\operatorname{Fin}_{J_0}\mathcal A(\GL_m)_{\rm triv},
 \qquad \mathcal A_{1,1}=\C.
\]
We use $\mathcal A_{\ell,\pi}$ to emphasize the $\ell$ consecutive twists
of $\pi$ in \eqref{eqn: cuspidal-support-fix-cusp}. Thus $\mathcal A_{m,1}$
denotes $\mathcal A(\GL_m)_{\rm triv}$, and $\mathcal A_{2,\pi}$ denotes
the module on $\GL_{2k}$ with cuspidal support represented by
$\pi|{\det}|^{1/2}\boxtimes\pi|{\det}|^{-1/2}$ on $\GL_k\times\GL_k$.

\subsection{Residues for characters and repeated cuspidal factors}
\label{subsec:character-residue-image}
\label{subsec:repeated-cuspidal-residues}

For $m\ge2$, put $B_m=P_{1,\ldots,1}$ and
$\rho_m=(\rho_{m,1},\ldots,\rho_{m,m})$, where
\begin{equation}
\label{eqn: notation-section-5.2}
 \rho_{m,i}=\frac{m+1-2i}{2},\qquad
 \sum_{i=1}^m x_i=0,\qquad
 u_i=x_i-x_{i+1}\quad(1\le i<m).
\end{equation}
The $u_i$ are coordinates on $\sum_i x_i=0$.
Let $V_m$ be the full-level compact realization of
$\mathop{\times}_{i=1}^m|\cdot|^{\rho_{m,i}+x_i}$, as in
Definition~\ref{def:compact-realization}\textup{(iii)};
it is independent of $x$ by Lemma~\ref{lem:compact-realization}.
By Proposition~\ref{prop:cuspidal-normal-crossings},
\begin{equation}\label{eq:character-regularized-family}
 H_m(x)=\left(\prod_{i=1}^{m-1}u_i\right)
                 E_{B_m}^{\GL_m}(\rho_m+x)
\end{equation}
is holomorphic near $x=0$.
Let $\D_{u_1,\ldots,u_{m-1}}$ denote the space of finite linear
combinations of derivative evaluations at $u=0$, and set
$\mathcal T_m=V_m\otimes\D_{u_1,\ldots,u_{m-1}}$.
With $\mathcal A_{m,1}$ as in
Subsection~\ref{subsec:speh-support-modules}, write the coefficient map
of $H_m$ in the $u$-coordinates (Definition~\ref{def:coefficient-map}) as
\[
 \Xi_m:=\Xi_{H_m}:\mathcal T_m\longrightarrow
                  \mathcal A_{m,1}^{K_{m,f}},\qquad
 \Xi_m(f\otimes\delta)=\delta\bigl(H_m(\,\cdot\,)f\bigr).
\]
Its image is $\mathcal A_{m,1}^{K_{m,f}}$.
For $S\subseteq\{1,\ldots,m-1\}$, define
\begin{equation}\label{eq:character-residue-module}
 \mathcal R_{m,S}
   :=\Xi_m\left(\bigcap_{i\in S}\ker(u_i:\mathcal T_m\to\mathcal T_m)\right).
\end{equation}
Proposition~\ref{prop:cuspidal-normal-crossings} gives
\[
 \operatorname{Res}_S\!\left(E_{B_m}^{\GL_m}(\rho_m+x)f\right)
 =\frac{\left.H_m(x)f\right|_{u_i=0,\ i\in S}}
        {\prod_{i\notin S}u_i}.
\]
Thus $\mathcal R_{m,S}$ is spanned by the Laurent coefficients of
these residues as $f$ varies in $V_m$, or equivalently by the
Taylor coefficients of $\left.H_m(x)f\right|_{u_i=0,\ i\in S}$
in all remaining variables
(Corollary~\ref{cor:residues-and-coefficients}).
Denote also $\mathcal R_{m,a}=\mathcal R_{m,\{a\}}$.
For $m=1$, set $\mathcal R_{1,\varnothing}=\C$.

For $m\ge2$, put $t_j=\rho_{m,j}+x_j$.  For $1\le i<m$, on the
hyperplane $u_i=0$ we have $t_i-t_{i+1}=1$.
Set $c_i=(t_i+t_{i+1})/2$.  The quotient
\eqref{eq:global-residue-quotient}, twisted by $|{\det}|^{c_i}$, is
\[
 |\cdot|^{t_i}\times|\cdot|^{t_{i+1}}
 \xrightarrow{\ R_{1,1}\otimes|{\det}|^{c_i}\ }
 |{\det}|^{c_i}.
\]
Let $V_m^{(i)}$ be the full-level compact realization of
\[
 \left(\mathop{\times}_{j<i}|\cdot|^{t_j}\right)
 \times |{\det}|^{c_i}\times
 \left(\mathop{\times}_{j>i+1}|\cdot|^{t_j}\right),
\]
where the middle factor is a character of $\GL_2(\A)$.
Inducing the displayed quotient with the remaining characters gives
\begin{equation}
\label{eqn: def-R_i}
R_i:V_m\twoheadrightarrow V_m^{(i)}.
\end{equation}
By Lemma~\ref{lem:adjacent-character-quotient}, this map in the compact
realizations is independent of the parameters on $u_i=0$.
Let $L_i=\ker R_i$.
For $D\subseteq\{1,\ldots,m-1\}$, set
\begin{equation}
\label{eqn: def-k^D}
 \mathcal K^D
   =\sum_{\substack{1\le i<m\\i\notin D}}L_i\otimes
          \C\operatorname{ev}_{u_i}\otimes\D_{(u_j)_{j\ne i}},
 \qquad
 \mathcal M_D=\mathcal T_m/\mathcal K^D.
\end{equation}
Here $\operatorname{ev}_{u_i}(a)=a(0)$ for a germ in $u_i$.
Since $R_i$ intertwines for every parameter on $u_i=0$,
each summand of $\mathcal K^D$ is stable under the action
\eqref{eq:g-action-jets}, including all derivatives in the other variables.
Franke's kernel and coefficient-spanning theorems
\cite[Section~5]{Franke2008} give
\begin{equation}\label{eq:character-coefficient-kernel}
 \ker\Xi_m=\mathcal K^\varnothing,\qquad
 \mathcal M_\varnothing
    \xrightarrow[\ \Xi_m\ ]{\ \sim\ }
       \mathcal A_{m,1}^{K_{m,f}}.
\end{equation}
For $i\in\{1,\ldots,m-1\}\setminus D$, we have
$u_i\operatorname{ev}_{u_i}=0$, while multiplication by $u_i$
preserves every summand indexed by $j\ne i$.
Thus $u_i\mathcal K^D\subseteq\mathcal K^{D\cup\{i\}}$, and the map
for the additional cut after $i$
\begin{equation}\label{eq:character-cut-map}
 P_i^D:\mathcal M_D\longrightarrow\mathcal M_{D\cup\{i\}},\qquad
 [\eta]_D\longmapsto[u_i\eta]_{D\cup\{i\}}
\end{equation}
is well defined, where $[\eta]_D$ denotes the class of
$\eta\in\mathcal T_m$ modulo $\mathcal K^D$.
These maps commute and are equivariant, since scalar multiplication
on distributions commutes with \eqref{eq:g-action-jets}.

For $\alpha\in\Z_{\ge0}^{m-1}$, define
\[
 \delta_\alpha(g):=[u^\alpha]g
   =\frac{1}{\alpha!}\partial_u^\alpha g(0),
 \qquad \alpha!=\prod_{i=1}^{m-1}\alpha_i!,
\]
for a holomorphic germ $g$ at $u=0$.
Then $u_i\delta_\alpha=\delta_{\alpha-e_i}$, where $e_i$
is the $i$th standard basis vector of $\Z^{m-1}$, and a distribution with a
negative index is zero.  Let $(\mathcal M_D)_\alpha$ denote the
image of $V_m\otimes\C\delta_\alpha$ under the quotient map
$\mathcal T_m\to\mathcal M_D$.

\begin{lemma}\label{lem:character-cut-grading}
For every $D\subseteq\{1,\ldots,m-1\}$,
\[
 \mathcal M_D=\bigoplus_{\alpha\in\Z_{\ge0}^{m-1}}(\mathcal M_D)_\alpha,\qquad
 (\mathcal M_D)_\alpha\cong
 V_m\Big/\sum_{\substack{i\notin D\\\alpha_i=0}}L_i .
\]
For $i\notin D$, $P_i^D$ is zero on this summand if $\alpha_i=0$;
if $\alpha_i>0$, it maps the summand isomorphically onto
$(\mathcal M_{D\cup\{i\}})_{\alpha-\mathbf e_i}$.
\end{lemma}

\begin{proof}
In \eqref{eqn: def-k^D}, the term indexed by $j\notin D$
contributes $L_j\otimes\C\delta_\alpha$ in degree $\alpha$
precisely when $\alpha_j=0$. For $\alpha_i>0$ we have
\[
 \{j\notin D:\alpha_j=0\}
 =\{j\notin D\cup\{i\}:(\alpha-\mathbf e_i)_j=0\}.
\]
The equality $u_i\delta_\alpha=\delta_{\alpha-\mathbf e_i}$ gives
the claimed isomorphism.
\end{proof}

\begin{proposition}\label{prop:character-block-cut-degree}
For $m\ge1$, $D\subseteq\{1,\ldots,m-1\}$ and either
determinant parity $\varepsilon$, we have
\[
 H^q(\mathfrak{sl}_m,\O_m;\mathcal M_D\otimes E_\varepsilon)\ne0
 \quad\Longrightarrow\quad q\le d_m-2m+2+|D|.
\]
Here $\mathcal M_\varnothing=\mathcal A_{m,1}^{K_{m,f}}$, also for $m=1$.
\end{proposition}

\begin{proof}
For $m=1$, $\mathcal M_\varnothing=\C$ and $d_1=0$.
For $m\ge2$, write $D=\{d_1<\cdots<d_t\}$ and set
\[
 d_0=0,\quad d_{t+1}=m,\qquad
 m_\nu=d_\nu-d_{\nu-1},\qquad
 c_\nu=\frac{m-d_{\nu-1}-d_\nu}{2}.
\]
Cut the ordered indices $1,\ldots,m$ between $i$ and $i+1$
for each $i\in D$.  This gives $t+1$ consecutive blocks
$\{d_{\nu-1}+1,\ldots,d_\nu\}$, of lengths $m_\nu$,
for $1\le\nu\le t+1$.
For $1\le i<m$, the indices $i$ and $i+1$ lie in the same
block exactly when $i\notin D$.
Let $R_D=P_{m_1,\ldots,m_{t+1}}$ and
$\mathfrak n_D=\Lie(N_{R_D})_\C$.
For each block, define
\[
 y_\nu=\frac1{m_\nu}\sum_{a=1}^{m_\nu}x_{d_{\nu-1}+a},
 \qquad z_{\nu,a}=x_{d_{\nu-1}+a}-y_\nu.
\]
Then $\sum_a z_{\nu,a}=0$, $\sum_\nu m_\nu y_\nu=0$, and by \eqref{eqn: notation-section-5.2},
\[
 \rho_{m,d_{\nu-1}+a}+x_{d_{\nu-1}+a}
 =\rho_{m_\nu,a}+z_{\nu,a}+c_\nu+y_\nu.
\]
For $1\le a<m_\nu$, by \eqref{eqn: notation-section-5.2}
$u_{d_{\nu-1}+a}=z_{\nu,a}-z_{\nu,a+1}$.
Reindexing the defining sum for $\mathcal K^D$ of \eqref{eqn: def-k^D} gives
\[
 \mathcal K^D
 =\sum_{\nu=1}^{t+1}\sum_{a=1}^{m_\nu-1}
    L_{d_{\nu-1}+a}\otimes
    \C\operatorname{ev}_{u_{d_{\nu-1}+a}}\otimes
    \D_{(u_j)_{j\ne d_{\nu-1}+a}}.
\]
For $1\le a<m_\nu$, let
$R_a^{(\nu)}:V_{m_\nu}\twoheadrightarrow V_{m_\nu}^{(a)}$
be induced, within $\GL_{m_\nu}$, from
\eqref{eq:global-residue-quotient} in
Lemma~\ref{lem:adjacent-character-quotient}, twisted by
$|{\det}|^{(m_\nu-2a)/2+z_{\nu,a}}$:
\[
 |\cdot|^{\rho_{m_\nu,a}+z_{\nu,a}}
 \times|\cdot|^{\rho_{m_\nu,a+1}+z_{\nu,a+1}}
 \twoheadrightarrow
 |{\det}|^{(m_\nu-2a)/2+z_{\nu,a}}
 \qquad(z_{\nu,a}=z_{\nu,a+1}).
\]
Inducing this map with the other blocks and twisting block $\mu$
by $|{\det}|^{c_\mu+y_\mu}$ gives $R_{d_{\nu-1}+a}$.

Let $\D_{z_\nu}$ denote the finite-order distributions at zero on
$\sum_a z_{\nu,a}=0$, and let $\D_y$ denote those at zero on
$\sum_\nu m_\nu y_\nu=0$.
For $\Xi_{m_\nu}$ defined with the variables $z_\nu$,
Franke's kernel formula
\cite[Theorem~5.2, Equation~(5.18)]{Franke2008} reads
\[
 \ker\Xi_{m_\nu}
 =\sum_{a=1}^{m_\nu-1}\ker R_a^{(\nu)}\otimes
   \ker\bigl((z_{\nu,a}-z_{\nu,a+1}):
                   \D_{z_\nu}\to\D_{z_\nu}\bigr).
\]
For $m_\nu=1$, the coefficient map is $\operatorname{id}_\C$
and its kernel is zero.
The coordinate change $x\leftrightarrow(z,y)$ identifies
\[
 \D_{u_1,\ldots,u_{m-1}}
 \cong\left(\bigotimes_{\nu=1}^{t+1}\D_{z_\nu}\right)\otimes\D_y.
\]
Put $W_\nu=V_{m_\nu}\otimes\D_{z_\nu}$, the domain of
$\Xi_{m_\nu}$, with $W_\nu=\C$ when $m_\nu=1$.
Using the real compact induction $\mathcal I_{R_D,K}$ of
Definition~\ref{def:compact-realization}\textup{(i)}, induction in stages
identifies the underlying $\O_m$-modules as
\[
 \mathcal T_m\cong
 \mathcal I_{R_D,K}\bigl(W_1\otimes\cdots\otimes W_{t+1}\bigr)
 \otimes\D_y.
\]
Under this identification, the formula for $\ker\Xi_{m_\nu}$ and
$u_{d_{\nu-1}+a}=z_{\nu,a}-z_{\nu,a+1}$ give
\[
 \mathcal K^D\cong
 \mathcal I_{R_D,K}\left(
   \sum_{\nu=1}^{t+1}
   W_1\otimes\cdots\otimes\ker\Xi_{m_\nu}
      \otimes\cdots\otimes W_{t+1}\right)\otimes\D_y.
\]
Passing to $W_\nu/\ker\Xi_{m_\nu}$ for each $\nu$ and
using exactness of normalized induction gives
\[
 \begin{aligned}
 \mathcal I_D(y)
   &=\mathop{\times}_{\nu=1}^{t+1}
       \bigl(\mathcal A_{m_\nu,1}|{\det}|^{c_\nu+y_\nu}\bigr),\\
 \mathcal M_D&\cong\mathcal I_D(0)^{K_{m,f}}\otimes\D_y.
 \end{aligned}
\]
The action on $\mathcal I_D(0)^{K_{m,f}}\otimes\D_y$ is
\eqref{eq:g-action-jets} for the family $\mathcal I_D(y)$.

The Harish--Chandra parameter of
$\mathcal A_{m_\nu,1}|{\det}|^{c_\nu}$ is
$\{\rho_{m_\nu,a}+c_\nu:1\le a\le m_\nu\}$.
The infinitesimal-character comparison in the proof of
Proposition~\ref{prop:degree} therefore forces the Kostant representative
\[
 v(d_{\nu-1}+a)=m-d_\nu+a
 \qquad(1\le\nu\le t+1,\ 1\le a\le m_\nu).
\]
It preserves the order inside each block and reverses the order of
the blocks.  Hence
\[
 \ell_W(v)=\sum_{\mu<\nu}m_\mu m_\nu=\dim\mathfrak n_D.
\]
Moreover,
\[
 (v^{-1}\rho_m-\rho_m)_{d_{\nu-1}+a}
 =d_{\nu-1}+d_\nu-m=-2c_\nu.
\]
Thus Kostant's coefficient on each block has trivial
$\mathfrak{sl}_{m_\nu}$-action and one of the two determinant
parities on $\O_{m_\nu}$.
Theorem~\ref{prop:constant}, applied to
$\mathcal A_{m_\nu,1}^{K_{m_\nu,f}}$ with this coefficient,
gives $q_\nu\le d_{m_\nu}-2m_\nu+2$ for nonzero cohomology.
The coordinates of $\rho_{R_D}$ on block $\nu$ all equal
$c_\nu$.  Thus the constant part of the relative split-centre
action is $c_\nu+c_\nu-2c_\nu=0$.
Lemma~\ref{lem:split-centre-cohomology} gives
\[
 H^a((\mathfrak a_{R_D}^{\GL_m})_\C;\D_y)
 \cong
 \begin{cases}
 \C,&a=0,\\
 0,&a>0.
 \end{cases}
\]
Hence
\[
 \begin{aligned}
 H^q(\mathfrak{sl}_m,\O_m;\mathcal M_D\otimes E_\varepsilon)
       \ne0
 &\Longrightarrow
 q\le\sum_\nu(d_{m_\nu}-2m_\nu+2)+\dim\mathfrak n_D\\
 &=d_m-2m+2+t.
 \end{aligned}
\]
We used
$d_m=\sum_\nu d_{m_\nu}+\dim\mathfrak n_D+t$.
\end{proof}

Let $\pi$ be a regular algebraic unitary cuspidal representation of
$\GL_k(\A)$, with $k\ge2$, and assume that the Harish--Chandra
parameter of $J(2,\pi)_\infty$ has distinct entries.
Let $V$ be the full-level compact
realization of Definition~\ref{def:compact-realization}\textup{(iii)} for
\[
 I(v)=\pi|{\det}|^{(1+v)/2}\times
                         \pi|{\det}|^{-(1+v)/2}.
\]
For $f\in V$, let $f_v\in I(v)^{K_{2k,f}}$ be the section
with compact value $f$.  Put $P=P_{k,k}$ and define
\[
 (E(v)f)(g)=
 \sum_{\gamma\in P(\Q)\backslash\GL_{2k}(\Q)}
       [f_v(\gamma g)](1_{M_P}),
\]
initially in the region of absolute convergence and then by
meromorphic continuation, as in \eqref{eq:speh-eisenstein-family}.
In the notation \eqref{eq:normalized-intertwiner-notation}, put
\[
 \begin{aligned}
 N(v)&=N\left((12),\pi\boxtimes\pi,
          \left(\tfrac{1+v}{2},-\tfrac{1+v}{2}\right)\right),\\
 N(v)&:I(v)\longrightarrow
       \pi|{\det}|^{-(1+v)/2}\times\pi|{\det}|^{(1+v)/2}.
 \end{aligned}
\]
Lemma~\ref{lem:speh-residue-constant-term} identifies
$\operatorname{im}N(0)$ with $J(2,\pi)$.
On $K_{2k,f}$-fixed vectors, define $R_\pi$ by the factorization
\begin{equation}\label{eq:two-cuspidal-residue-quotient}
 \begin{gathered}
 V\xrightarrow{\ R_\pi\ }J(2,\pi)^{K_{2k,f}}
 \overset{\iota}{\hookrightarrow}
 \bigl(\pi|{\det}|^{-1/2}\times\pi|{\det}|^{1/2}\bigr)^{K_{2k,f}},\\
 N(0)f_0=\iota(R_\pi f).
 \end{gathered}
\end{equation}
Here $R_\pi$ is surjective.  In the automorphic realization
of $J(2,\pi)$ by residues, the same lemma gives
\begin{equation}\label{eq:two-cuspidal-residue-identity}
 \operatorname{Res}_{v=0}E(v)f=c_\pi R_\pi f,
 \qquad c_\pi\ne0.
\end{equation}
The two exponents have sum zero, and all derivatives below are
taken in the single variable $v$.
The space of finite-order distributions at $v=0$ is
\[
 \D_v=\bigoplus_{j\ge0}\C\delta_j,\qquad
 \delta_j(h)=\frac{h^{(j)}(0)}{j!},\qquad
 \delta_0=\operatorname{ev}_v,
\]
where $h$ is a holomorphic germ at zero.

\begin{lemma}\label{lem:two-cuspidal-coefficient-kernel}
Let $V$ be the common compact realization of $I(v)^{K_{2k,f}}$,
\[
 V=\mathcal I_{P_{k,k},K}(\pi_\infty\boxtimes\pi_\infty)
       \otimes_\C \bigl(V_{\pi_f}^{K_{k,f}}\bigr)^{\otimes 2}.
\]
Here the first factor is the space of $\O_{2k}$-finite functions
in \eqref{eq:compact-picture}, and $V_{\pi_f}^{K_{k,f}}$ is the
space of $K_{k,f}$-fixed vectors in $\pi_f$.
The family $H_\pi(v)=vE(v)$ is holomorphic near zero.
Define its coefficient map (Definition~\ref{def:coefficient-map}) as
\[
 \Xi_\pi:=\Xi_{H_\pi}:V\otimes\D_v
       \longrightarrow\mathcal A_{2,\pi}^{K_{2k,f}},
 \qquad
 \Xi_\pi(f\otimes\delta_j)
   =\left.\frac1{j!}\partial_v^j\bigl(vE(v)f\bigr)\right|_{v=0}.
\]
Then the sequence
\[
 0\longrightarrow(\ker R_\pi)\otimes\operatorname{ev}_v
 \longrightarrow V\otimes\D_v
 \xrightarrow{\ \Xi_\pi\ }
       \mathcal A_{2,\pi}^{K_{2k,f}}\longrightarrow0
\]
is exact as a sequence of $(\mathfrak{sl}_{2k},\O_{2k})$-modules,
where $V\otimes\D_v$ carries the action \eqref{eq:g-action-jets}.
\end{lemma}

\begin{proof}
Consider the intertwining operator  $M(v)=r(v)N(v)$, where
\[
 r(v)=\epsilon(v)^{-1}
 \frac{\Lambda(1+v,\pi\times\pi^\vee)}
      {\Lambda(2+v,\pi\times\pi^\vee)},
\]
where $\epsilon(v)$ is holomorphic and nonvanishing.
The numerator of $r(v)$ has a simple pole at $v=0$, and the denominator
is holomorphic and nonzero there.  Hence
$a(v):=vr(v)$ is holomorphic with $a(0)\ne0$.
The normalized operator $N(v)$ is holomorphic near zero
\cite[Sections~I.10--I.11 and II.2--II.3]{MoeglinWaldspurger1989}.

For $s\in\R$, let $a_s\in M_P(\A)$ have real component
$\operatorname{diag}(e^{s/k}1_k,e^{-s/k}1_k)$ and finite
component $1$.  Take $m=(m_1,m_2)\in M_P(\A)$ with
$|{\det} m_1|_\A=|{\det} m_2|_\A=1$, and
$\kappa\in\O_{2k}K_{2k,f}$.
Cuspidality of $\pi$ leaves only the identity and the
block transposition in \eqref{eq:full-normalized-constant-term}, so
\begin{equation}\label{eq:two-cuspidal-cleared-constant-term}
 (\CT_P H_\pi(v)f)(ma_s\kappa)
 =v e^{(1+v)s}[f(\kappa)](m)+a(v)e^{-(1+v)s}[(N(v)f)(\kappa)](m).
\end{equation}
Here $\CT_P$ is normalized.  If $\omega_{\pi,\infty}$
denotes the central character of $\pi_\infty$, then
\[
 (\pi\boxtimes\pi)(a_s)
 =\omega_{\pi,\infty}(e^{s/k})
  \omega_{\pi,\infty}(e^{-s/k})\operatorname{id}
 =\operatorname{id}.
\]
Both terms are holomorphic in $v$.  All other standard cuspidal
constant terms vanish by
Proposition~\ref{prop:bernstein-lapid-constant-term}.
The pole-detection criterion
\cite[Sections~I.3.4 and I.4.10]{MoeglinWaldspurger1995}
therefore gives holomorphy of $H_\pi$.

Let $\eta=\sum_{j=0}^d f_j\otimes\delta_j$ satisfy
$\Xi_\pi(\eta)=0$.  Taking Taylor coefficients in
\eqref{eq:two-cuspidal-cleared-constant-term} gives
\[
 0=e^s\sum_{j=1}^d\frac{s^{j-1}}{(j-1)!}[f_j(\kappa)](m)+e^{-s}\sum_{j=0}^d[v^j]\!\left(
       e^{-vs}a(v)[(N(v)f_j)(\kappa)](m)\right).
\]
After multiplication by $e^{-s}$, the second term tends to
zero as $s\to+\infty$.  Therefore we must have 
\[
\sum_{j=1}^d\frac{s^{j-1}}{(j-1)!}[f_j(\kappa)](m)=0.
\]
This implies $[f_j(\kappa)](m)=0$ for every $j\ge1$,
$m$, and $\kappa$.  The fixed central characters determine
the values of these cusp forms from their restrictions to
$|{\det} m_1|_\A=|{\det} m_2|_\A=1$, so
$f_1=\cdots=f_d=0$.
The residue identity \eqref{eq:two-cuspidal-residue-identity} now gives
\[
 \Xi_\pi(f_0\otimes\delta_0)=H_\pi(0)f_0=c_\pi R_\pi f_0.
\]
Since $c_\pi\ne0$, this vanishes exactly when
$f_0\in\ker R_\pi$.  Consequently,
\[
 \ker\Xi_\pi=(\ker R_\pi)\otimes\C\delta_0.
\]

Let $F(v)=\sum_{\ell\ge0}v^\ell F_\ell$, with $F_\ell\in V$,
be a holomorphic $V$-valued germ at $v=0$.
The coefficient-spanning result
\cite[Section~2]{GrbacGrobner2024}, applied to the support
$\pi|{\det}|^{1/2}\boxtimes\pi|{\det}|^{-1/2}$ and followed by
averaging over $K_{2k,f}$, says that
$\mathcal A_{2,\pi}^{K_{2k,f}}$ is spanned by Laurent
coefficients of $E(v)F(v)$ at zero as $F$ varies.
Since $E(v)$ has at most a simple pole,
\[
 [v^n]\bigl(E(v)F(v)\bigr)
 =\sum_{\ell=0}^{n+1}
       \Xi_\pi(F_\ell\otimes\delta_{n+1-\ell})
 \qquad(n\ge-1).
\]
Thus $\Xi_\pi$ is surjective.

Finally, let $\rho_v$ denote the induced action on $I(v)$,
expressed in the compact realization $V$.
For $X\in\mathfrak{sl}_{2k}$, Eisenstein equivariance gives
\[
 X\Xi_\pi(f\otimes\delta_j)
 =\sum_{\ell=0}^j
   \Xi_\pi\left(([v^\ell]\rho_v(X)f)
                    \otimes\delta_{j-\ell}\right).
\]
This is the action \eqref{eq:g-action-jets} in the basis
$\delta_j$.  The action of $\O_{2k}$ on $V$ is independent of $v$,
so $\Xi_\pi$ is also $\O_{2k}$-equivariant.
This proves the asserted exact sequence of modules.
\end{proof}

\subsection{Translation of a relative parameter}

\begin{lemma}\label{lem:finite-translation}
Let $U$ be a $(\mathfrak l,K)$-module, and let $N$ be a
locally nilpotent endomorphism commuting with this action.  The
action on each distribution space is trivial.  For a variable
$y$, write $\delta_j^y=\partial_y^j|_0/j!$, and similarly
for $\zeta$.  The map
\begin{equation}\label{eq:finite-translation-map}
\begin{aligned}
 \Theta_N:U\otimes\D_\zeta&\longrightarrow U\otimes\D_y,\\
 \Theta_N(f\otimes\delta_j^\zeta)
   &=\sum_{k\ge0}\binom{j+k}{k}N^kf\otimes\delta_{j+k}^y
\end{aligned}
\end{equation}
is an isomorphism of $(\mathfrak l,K)$-modules.  Its inverse
interchanges $y,\zeta$ and replaces $N$ by $-N$, and
\begin{equation}\label{eq:finite-translation-conjugation}
 \Theta_N^{-1}(y-N)\Theta_N=\zeta,\qquad
 \Theta_N^{-1}y\Theta_N=\zeta+N.
\end{equation}
Consequently, for every $r\ge0$,
\begin{equation}\label{eq:finite-translation-kernel}
 \ker(y-N)^r
 =\Theta_N\!\left(U\otimes\ker(\zeta^r:\D_\zeta\to\D_\zeta)\right),
 \qquad
 \ker\zeta^r=\bigoplus_{j=0}^{r-1}\C\delta_j^\zeta,
\end{equation}
where the sum is zero for $r=0$.
Multiplication by $y-N$ is surjective, with right inverse
\[
 S_N(f\otimes\delta_j^y)
   =\sum_{k\ge0}N^kf\otimes\delta_{j+k+1}^y,\qquad
 (y-N)\circ S_N=\mathrm{id}_{U\otimes\D_y}.
\]
For every $N$-invariant subspace $U_0\subseteq U$,
$\Theta_N$ maps $U_0\otimes\D_\zeta$ onto
$U_0\otimes\D_y$, and $S_N$ preserves
$U_0\otimes\D_y$.
\end{lemma}

\begin{proof}
Every displayed sum is finite because $N$ is locally nilpotent.
In the composite of $\Theta_N$ with the map obtained by replacing
$N$ by $-N$, the coefficient of
$N^nf\otimes\delta_{j+n}$ is
\[
 \sum_{k=0}^n(-1)^k
       \binom{j+k}{k}\binom{j+n}{n-k}
 =\binom{j+n}{n}\sum_{k=0}^n(-1)^k\binom nk
 =\begin{cases}1,&n=0,\\0,&n>0.\end{cases}
\]
The same calculation applies in the opposite order.  This proves
the inverse formula.

Since $y\delta_j^y=\delta_{j-1}^y$, with
$\delta_{-1}^y=0$, for $j\ge1$ we have
\[
\begin{aligned}
 (y-N)\Theta_N(f\otimes\delta_j^\zeta)
 &=f\otimes\delta_{j-1}^y
   +\sum_{k\ge1}
       \left(\binom{j+k}{k}-\binom{j+k-1}{k-1}\right)
         N^kf\otimes\delta_{j+k-1}^y\\
 &=\sum_{k\ge0}\binom{j+k-1}{k}
         N^kf\otimes\delta_{j+k-1}^y
 =\Theta_N(f\otimes\delta_{j-1}^\zeta).
\end{aligned}
\]
For $j=0$, the consecutive terms
$N^{k+1}f\otimes\delta_k^y$ cancel, giving zero.
Thus $(y-N)\Theta_N=\Theta_N\zeta$.
Since $N\Theta_N=\Theta_NN$, this proves
\eqref{eq:finite-translation-conjugation}.
Taking the kernel of the $r$-th power gives
\eqref{eq:finite-translation-kernel}, since multiplication by
$\zeta$ lowers the distribution index by one.

The right-inverse formula follows from
\begin{equation}\label{eq:finite-translation-right-inverse}
 (y-N)\sum_{k\ge0}N^kf\otimes\delta_{j+k+1}^y
 =\sum_{k\ge0}N^kf\otimes\delta_{j+k}^y
   -\sum_{k\ge0}N^{k+1}f\otimes\delta_{j+k+1}^y
 =f\otimes\delta_j^y.
\end{equation}
Each map uses only powers of $N$ on its coefficient.  It therefore
commutes with the $(\mathfrak l,K)$-action and preserves every
stated $N$-invariant subspace; the inverse formula gives the
asserted equality for $\Theta_N$.
\end{proof}

\begin{remark}
Suppose $U\ne0$. Since $N$ is locally nilpotent, choose
$0\ne f\in\ker N$. For $\eta=f\otimes\delta_0^y\ne0$, we have
\[
 (y-N)\eta=f\otimes y\delta_0^y-Nf\otimes\delta_0^y=0,
 \qquad S_N((y-N)\eta)=0\ne\eta .
\]
Thus $y-N$ has a nonzero kernel and admits no left inverse.
\end{remark}

\section{Cohomology of mixed cuspidal supports}
\label{sec:internal-cohomology}

\subsection{Intertwining operators and cuts in the character interval}
\label{subsec:mixed-character-interval}

Let $I$ be an interval of consecutive character weights. We treat a
cuspidal factor of even rank and a pair of cuspidal factors of rank two.
For a half-integral cuspidal centre $c$ with $c\pm\tfrac12\in I$,
we compute the factor contributed by a central zero from the order of
the cuspidal factor and the characters at $c\pm\tfrac12$ in the constant term.
In both cases the cut maps define a filtration of the Eisenstein
coefficient module. We give the construction at the end of this subsection.

\begin{lemma}\label{lem:character-cuspidal-intertwiner}
Let $\pi$ be a full-level regular algebraic unitary cuspidal
representation of $\PGL_k(\A)$, where $k=2h$, and write
$\pi_\infty\simeq\mathop{\times}_{i=1}^hD_{d_i}$, with $d_i>0$.
Suppose $t^0,s^0\in\R$ satisfy
\[
 |t^0-s^0|<\frac{d_i}{2}\qquad(1\le i\le h).
\]
Then the normalized intertwining operator
\[
 N(t,s):\pi|{\det}|^s\times|\cdot|^t
       \longrightarrow|\cdot|^t\times\pi|{\det}|^s
\]
and its inverse are holomorphic near $(t^0,s^0)$ on
$K_{k+1,f}$-fixed sections in the compact realizations of
Definition~\ref{def:compact-realization}.
\end{lemma}

\begin{proof}
At the real place, the inequality makes
$D_{d_i}|{\det}|^{s^0}\times|\cdot|^{t^0}$ irreducible
\cite[Section~I.7, Lemma, p.~624]{MoeglinWaldspurger1989}.
Thus the normalized intertwining operators
\[
 N_i(t,s):D_{d_i}|{\det}|^s\times|\cdot|^t
       \longrightarrow|\cdot|^t\times D_{d_i}|{\det}|^s
\]
and their inverses are holomorphic near $(t^0,s^0)$
\cite[Section~I.1, Equation~(1), and Section~I.5(ii)]{MoeglinWaldspurger1989}.
Induce each $N_i(t,s)$ with the identity maps on the remaining
factors and apply them in the order $i=h,h-1,\ldots,1$.
Their composition is $N(t,s)_\infty$
\cite[Section~I.2, Equation~(5)]{MoeglinWaldspurger1989}, so it and its inverse
are holomorphic.
At each finite prime $p$, let $f_p^\circ$ and $f_p^{\prime\circ}$
be the normalized spherical vectors in the source and target of
$N(t,s)_p$, respectively. By \cite[Theorem~3.1, p.~397]{Casselman1980},
\[
 N(t,s)_p f_p^\circ=f_p^{\prime\circ},\qquad
 N(t,s)_f:=\bigotimes_{p<\infty}'N(t,s)_p:
 \bigotimes_{p<\infty}f_p^\circ\longmapsto
 \bigotimes_{p<\infty}f_p^{\prime\circ}.
\]
Thus $N(t,s)_f$ has scalar $1$ in these bases of the two
$K_{k+1,f}$-fixed lines. The assertion follows from the real-place calculation.
\end{proof}

\medskip

Let us introduce a filtration defined by the cut maps. We start with the
general setup. Let $X$ be a finite set whose elements we call cuts.
To each $Y\subseteq X$, we assign a $(\mathfrak g,K)$-module $M_Y$
with a grading as a vector space
\[
M_Y=\bigoplus_{\alpha\in\mathbb Z_{\ge0}^{X}}(M_Y)_\alpha.
\]
For any $j\in X\setminus Y$, let $P_j^Y:M_Y\to M_{Y\cup\{j\}}$
be a morphism of $(\mathfrak g,K)$-modules which we call a cut map.
We assume
\begin{itemize}
\item For distinct $i,j\in X\setminus Y$,
\[
 P_i^{Y\cup\{j\}}\circ P_j^Y=P_j^{Y\cup\{i\}}\circ P_i^Y.
\]
\item For every $\alpha\in\mathbb Z_{\ge0}^{X}$ and $j\in X\setminus Y$,
\[
 \begin{aligned}
 P_j^Y((M_Y)_\alpha)&=0 &&(\alpha_j=0),\\
 P_j^Y\big|_{(M_Y)_\alpha}\colon (M_Y)_\alpha
 &\xrightarrow{\ \sim\ }(M_{Y\cup\{j\}})_{\alpha-\mathbf e_j}
 &&(\alpha_j>0),
 \end{aligned}
\]
where $\mathbf e_j$ is the $j$-th standard basis vector of $\mathbb Z^X$.
\end{itemize}
Write $M=M_\varnothing$. We call the data
$(M_Y,\{P_j^Y\}_{j\in X\setminus Y})_{Y\subseteq X}$ satisfying
these two assumptions a cut system for $(X,M)$.

For an ordering $j_1,\ldots,j_r$ of $Y$, put $Y_0=\varnothing$
and $Y_k=\{j_1,\ldots,j_k\}$ for $1\le k\le r$. The composite
\[
 P_Y=P_{j_r}^{Y_{r-1}}\circ\cdots\circ P_{j_1}^{Y_0}\colon M\longrightarrow M_Y
\]
is independent of the ordering by the first assumption, since adjacent
cuts may be interchanged. Set $P_\varnothing=\mathrm{id}_M$.
For each target degree $\beta$, the second assumption gives an isomorphism
from $(M_Y)_{\beta+\mathbf e_j}$ onto $(M_{Y\cup\{j\}})_\beta$.
Thus every cut $P_j^Y$ is surjective, and so is each composite $P_Y$.

Given a cut system as above, define for each $C\subseteq X$
\[
D_C=\bigcap_{j\in C}\ker P_j^\varnothing,\qquad
 q_C=P_{X\setminus C}|_{D_C},\qquad E_C=q_C(D_C).
\]
We introduce the cut filtration on $M$ by setting
\[
 F_sM=\sum_{|C|=s}D_C\qquad(0\le s\le |X|).
\]
Since $D_{C\cup\{j\}}\subseteq D_C$, the filtration is decreasing.
Here $D_\varnothing=F_0M=M$, and we set $F_sM=0$ for $s>|X|$.
Write $\operatorname{gr}_s=F_s/F_{s+1}$.

\begin{lemma}\label{lem:cut-filtration-graded}
For $0\le s\le |X|$, the map with components $P_{X\setminus C}$
induces an isomorphism
\[
 \operatorname{gr}_s M \xrightarrow{\ \sim\ }\bigoplus_{|C|=s}E_C,\qquad
 m+F_{s+1}M\longmapsto (P_{X\setminus C}m)_{|C|=s}
 \quad (m\in F_sM).
\]
\end{lemma}

\begin{proof}
For $Z(\alpha)=\{j\in X:\alpha_j=0\}$, the cut rule gives
\[
 D_C=\bigoplus_{C\subseteq Z(\alpha)} M_\alpha,\qquad
 F_sM=\bigoplus_{|Z(\alpha)|\ge s}M_\alpha,
 \qquad M_\alpha=(M_\varnothing)_\alpha .
\]
For $|C|=s$, $P_{X\setminus C}$ kills every summand of $F_sM$
except those with $Z(\alpha)=C$. On these it is injective, since
each cut is an isomorphism onto its target summand and distinct
$\alpha$ give distinct target degrees. Their images sum to $E_C$.
In particular
\[
 \ker q_C=D_C\cap F_{s+1}M=\sum_{j\in X\setminus C}D_{C\cup\{j\}} .
\]
Thus the map on $F_sM$ has kernel $F_{s+1}M$ and image
$\bigoplus_{|C|=s}E_C$.
\end{proof}

Let $B\subseteq M$ be a $(\mathfrak g,K)$-submodule and $A=M/B$,
with $F_sB=B\cap F_sM$ and $F_sA=(F_sM+B)/B$.
For each $C\subseteq X$, let $\alpha_C:M\to M$ be an equivariant endomorphism
satisfying $\alpha_C(\ker P_Y)\subseteq\ker P_Y$ for $Y=X\setminus C$.
Since $P_Y$ is surjective and $\alpha_C$ preserves its kernel, $\alpha_C$
induces an endomorphism of $M_Y$ given by
\[
 \alpha_C\colon M_Y\longrightarrow M_Y,\qquad \alpha_C(P_Yx)=P_Y(\alpha_Cx)\quad(x\in M).
\]
\begin{proposition}\label{prop:cut-filtration}
Suppose in addition that for every $C\subseteq X$,
\begin{align}
 \ker(\alpha_C|_{D_C})&\subseteq B,\qquad
 P_{X\setminus C}(\alpha_C x)=0\quad(x\in B\cap F_{|C|}M),
 \label{eq:cut-kernel-tests}\\
 \alpha_CD_C&=D_C,\qquad \alpha_C(\ker q_C)=\ker q_C .
 \label{eq:cut-scalar-surjectivity}
\end{align}
Then for $0\le s\le |X|$, under the isomorphism of
Lemma~\ref{lem:cut-filtration-graded},
\[
 \operatorname{gr}_s B=\bigoplus_{|C|=s}\ker(\alpha_C|_{E_C}),\qquad
 \operatorname{gr}_s A\xrightarrow{\ \sim\ }\bigoplus_{|C|=s}E_C .
\]
The component of the last map indexed by $C$ sends the class of
$x\in F_sM$ to $\alpha_CP_{X\setminus C}x$.
\end{proposition}

\begin{proof}
For $|C|=s$, \eqref{eq:cut-scalar-surjectivity} gives
$\alpha_C E_C=q_C(\alpha_C D_C)=E_C$. If $x\in B\cap F_sM$, then
$\alpha_CP_{X\setminus C}x=P_{X\setminus C}(\alpha_Cx)=0$ by
\eqref{eq:cut-kernel-tests}. This puts the image of $\operatorname{gr}_s B$
in the indicated direct sum of kernels.
Conversely, fix $\xi\in E_C$ with $\alpha_C\xi=0$, and lift it to
$d\in D_C$. Then $\alpha_Cd\in\ker q_C$. By
\eqref{eq:cut-scalar-surjectivity}, choose $d'\in\ker q_C$ with
$\alpha_Cd'=\alpha_Cd$. The first test gives $d-d'\in B\cap D_C$.
Its $C$-component is $\xi$, and all its other components vanish,
since their complementary cuts include an index in $C$.
This proves the formula for $\operatorname{gr}_s B$.
Taking the quotient by
$\operatorname{gr}_s B$ and multiplying by $\alpha_C$ on each
component gives the asserted isomorphism for $A$.
\end{proof}

Let $I\subset\Z$ consist of $m\ge1$ consecutive integers. Arrange the
inducing characters  $\chi_j=|\cdot|^{t_j}$ in decreasing order of
their  exponents, so that $t_j^0-t_{j+1}^0=1$. Put
\[
 u_j=t_j-t_{j+1}-1\quad(1\le j<m).
\]
For $Y\subseteq\{1,\ldots,m-1\}$, let $\mathcal M_Y$ be the module
of \eqref{eqn: def-k^D}; thus $\mathcal M_\varnothing=\mathcal A_{m,1}^{K_{m,f}}$ by
\eqref{eq:character-coefficient-kernel}. For $j\in\{1,\ldots,m-1\}\setminus Y$,
the cut between $\chi_j$ and $\chi_{j+1}$ uses the map $P_j^Y$
of \eqref{eq:character-cut-map} with $D=Y$ and $i=j$.

\begin{lemma}\label{lem:character-cut-system}
Let $X=\{1,\ldots,m-1\}$. The modules $\mathcal M_Y$ for $Y\subseteq X$
and the maps $P_j^Y$ for $j\in X\setminus Y$ form a cut system
for $(X,\mathcal A_{m,1}^{K_{m,f}})$.
\end{lemma}

\begin{proof}
Each $P_j^Y$ is equivariant since multiplication by $u_j$ commutes
with the action \eqref{eq:g-action-jets}. For distinct $i,j\in X\setminus Y$
and $\eta\in\mathcal T_m$, \eqref{eq:character-cut-map} gives
\[
 P_i^{Y\cup\{j\}}P_j^Y[\eta]_Y
 =[u_iu_j\eta]_{Y\cup\{i,j\}}
 =[u_ju_i\eta]_{Y\cup\{i,j\}}
 =P_j^{Y\cup\{i\}}P_i^Y[\eta]_Y .
\]
For $f\in V_m$ and $\delta_\alpha(g)=[u^\alpha]g$,
Lemma~\ref{lem:character-cut-grading} gives
\[
 P_j^Y([f\otimes\delta_\alpha]_Y)=
 \begin{cases}
  0,&\alpha_j=0,\\
  [f\otimes\delta_{\alpha-\mathbf e_j}]_{Y\cup\{j\}},&\alpha_j>0,
 \end{cases}
\]
where $\mathbf e_j$ is the $j$-th standard basis vector of $\Z^{m-1}$.
For $\alpha_j>0$, Lemma~\ref{lem:character-cut-grading} shows that
$P_j^Y$ maps $(\mathcal M_Y)_\alpha$ isomorphically onto
$(\mathcal M_{Y\cup\{j\}})_{\alpha-\mathbf e_j}$.
\end{proof}

For $C\subseteq\{1,\ldots,m-1\}\setminus Y$ put
$D_C^Y=\bigcap_{i\in C}\ker P_i^Y$, with $D_\varnothing^Y=\mathcal M_Y$.

\begin{proposition}\label{lem:character-cut-modules}
For $j\notin C\cup Y$ there is an exact sequence
\[
 0\longrightarrow D_{C\cup\{j\}}^Y\longrightarrow D_C^Y
 \xrightarrow{P_j^Y}D_C^{Y\cup\{j\}}\longrightarrow0.
\]
For $i\in C$ we have $u_iD_C^Y=0$. For $\varepsilon\in\{1,\det\}$,
\begin{equation}\label{eq:character-cut-degree}
 H^p(\mathfrak{sl}_m,\O_m;D_C^Y\otimes\varepsilon)\ne0
 \quad\Longrightarrow\quad p\le d_m-2m+2+|Y|+2|C|.
\end{equation}
For $Y=\varnothing$ we have $D_C^\varnothing=\mathcal R_{m,C}$.
Taking $C=\{j\}$ gives $\ker P_j^\varnothing=\mathcal R_{m,j}$.
\end{proposition}

\begin{proof}
By the divided-power calculation of Lemma~\ref{lem:character-cut-grading},
$P_j^Y$ is zero on the summand of multidegree $\alpha$ if
$\alpha_j=0$, and otherwise maps that summand isomorphically onto the
one of degree $\alpha-\mathbf e_j$. Thus $D_C^Y$ is the sum of the
summands with $\alpha_i=0$ for $i\in C$. This proves exactness and
$u_iD_C^Y=0$; for $Y=\varnothing$ the same summands give
$D_C^\varnothing=\mathcal R_{m,C}$ by \eqref{eq:character-residue-module}.
For $C=\varnothing$, \eqref{eq:character-cut-degree} is the bound for
$\mathcal M_Y$ of Proposition~\ref{prop:character-block-cut-degree}
($d_1=0$ when $m=1$). For $j\notin C\cup Y$, the long exact
sequence bounds the degree on $D_{C\cup\{j\}}^Y$ by
\[
 d_m-2m+2+\max\bigl\{|Y|+2|C|,\ 1+|Y|+1+2|C|\bigr\}
 =d_m-2m+2+|Y|+2|C\cup\{j\}|.
\]
\end{proof}

\subsection{One cuspidal factor}
\label{subsec:central-zero-merger}
Let $\pi$ be a full-level regular algebraic unitary cuspidal
representation of $\PGL_k(\A)$, where $k=2h$. Let
$B=\{\beta_1<\cdots<\beta_k\}$ be its shifted integral weight set
as in \eqref{eq:packet} with $\ell=1$. Let $I\subset\Z$ be an interval
of $m\ge1$ consecutive character weights with $I\cap B=\varnothing$. Put
\[
 c=\frac{\beta_i+\beta_{k+1-i}}2\quad(1\le i\le h),\qquad
 \bar c=\frac{\sum_{x\in I}x+\sum_{\beta\in B}\beta}{m+k}.
\]
Let $\mathfrak c$ be represented by
\[
 \left(\boxtimes_{x\in I}|\cdot|^{x-\bar c}\right)
     \boxtimes\pi|{\det}|^{c-\bar c}.
\]
Let $\chi$ have parameter $(I\sqcup B)-\bar c$ on
$\mathfrak{sl}_{m+k}$, and put
\[
\mathcal A=(\operatorname{Fin}_{\ker\chi}
          \mathcal A(\GL_{m+k})_{\mathfrak c})^{K_f}.
          \]

\begin{theorem}\label{lem:one-copy-internal-bound}
For one cuspidal factor of rank $k=2h$, assume $\min I\le c\le\max I$.
For an algebraic representation $F$ whose dual has the 
infinitesimal character $\chi$ on $\mathfrak{sl}_{m+k}$,
\[
 H^q(\mathfrak{sl}_{m+k},\O_{m+k};\mathcal A\otimes F)=0
 \qquad\bigl(d_{m+k}-q<(h+2)m+h^2-3\bigr).
\]
\end{theorem}
The proof will be given at the end of this subsection.

Under the hypotheses of Theorem~\ref{lem:one-copy-internal-bound},
the interval $[\min I,\max I]_{\R}$ contains $c=(\beta_h+\beta_{h+1})/2$
and contains neither $\beta_h$ nor $\beta_{h+1}$, since its integer
points form $I$. Hence $I\subset(\beta_h,\beta_{h+1})$.
Put $Q=P_{m,k}$,
$t_j^0=\max I+1-j-\bar c$, $s_\pi=c-\bar c$, and set
\[
 s_I=\frac1m\sum_jt_j^0,\quad y_I=\frac1m\sum_j(t_j-t_j^0),\quad
 x_j=t_j-t_j^0-y_I,\quad y_\pi=t_\pi-s_\pi,\quad y=y_\pi-y_I.
\]
Thus $y$ is a coordinate on $my_I+ky_\pi=0$. For normalized induction
from $Q$, put
\begin{equation}\label{eq:one-cusp-induced-module}
\begin{aligned}
 \mathcal I_Q(y)&:=
 \bigl(\mathcal A_{m,1}|\det|^{s_I+y_I}\bigr)
 \times\bigl(\pi|\det|^{s_\pi+y_\pi}\bigr),\\
 \mathcal M&:=\mathcal I_Q(0)^{K_f}\otimes\D_y.
\end{aligned}
\end{equation}
The action on $\mathcal M$ is the differentiated action
\eqref{eq:g-action-jets}. We have fixed a finite spherical vector
to identify $\pi^{K_{k,f}}$ with $\pi_\infty$. Put
$r=\operatorname{ord}_{s=1/2}\Lambda(s,\pi)$.

When $c\in\tfrac12+\Z$ and $r>0$, let
$a=\#\{x\in I:x>c\}$. The characters at $c+1/2,c-1/2$ are
$\chi_a,\chi_{a+1}$, respectively. Put
\[
 u=u_a,\qquad w=t_\pi-t_{a+1}-\tfrac12=y-x_{a+1},\qquad
 t_a-t_\pi-\tfrac12=u-w,\qquad p=\prod_{j=1}^{m-1}u_j.
\]
Let $E_P$ and $E_{P'}$ be the Eisenstein families with the following
standard parabolics and ordered inducing factors. Write $1^d$ for
$d$ blocks of size one and $\chi_j=|\cdot|^{t_j}$. Then
\[
\begin{aligned}
 P  &= P_{1^m,k},       &\quad&(\chi_1,\ldots,\chi_m,\pi|\det|^{t_\pi}),\\
 P' &= P_{1^a,k,1^{m-a}}, &&(\chi_1,\ldots,\chi_a,\pi|\det|^{t_\pi},
                              \chi_{a+1},\ldots,\chi_m).
\end{aligned}
\]
Since $t_a^0=s_\pi+\tfrac12$ and $t_{a+1}^0=s_\pi-\tfrac12$,
the factors for $P'$ are in decreasing order of their exponents at the base point.
For $P$, the inducing data of \eqref{eqn: inducing-data-general-Levi} are
$\tau=\mathbf1^{\boxtimes m}\boxtimes\pi$ and $t=(t_1,\ldots,t_m,t_\pi)$.

For a full-level inducing jet $\eta$ for $P$ as in \eqref{eq:full-level-inducing-jet}, with
$z=(u_1,\ldots,u_{m-1},y)$, write $[\eta]_\varnothing$
for its image in $\mathcal M$ under the map induced by
$\Xi_m\otimes\mathrm{id}_{\pi_\infty}\otimes\mathrm{id}_{\D_y}$.
The subscript $\varnothing$ denotes the empty set of character cuts.

\begin{lemma}\label{lem:one-cusp-coefficient-map}
\textup{(i)} If $c\in\Z$ or $r=0$, the coefficient map is an isomorphism
\[
 \mathcal M\xrightarrow{\ \sim\ }\mathcal A,\qquad
 f\otimes\delta\longmapsto\delta\bigl(E_Q(y)f_y\bigr),
\]
where $f_y$ is the flat section through $f$.

\textup{(ii)} Suppose $c\in\tfrac12+\Z$ and $r>0$. Put $H=w^rpE_P$.
For each $f$ in the full-level compact realization of $I_P(t^0)$,
let $f_t$ be the flat section through $f$. The family
\[
 H(z)f=
 \bigl(t_\pi-t_{a+1}-\tfrac12\bigr)^r
 \Bigl(\prod_{j=1}^{m-1}(t_j-t_{j+1}-1)\Bigr)E_P(t)f_t
\]
is holomorphic in $z=(u_1,\ldots,u_{m-1},y)$ near $z=0$, which
corresponds to $t^0=(t_1^0,\ldots,t_m^0,s_\pi)$ on
$\sum_{j=1}^m t_j+kt_\pi=0$. The coefficient map $\Xi_H$ induces
\[
 \mathcal M\twoheadrightarrow\mathcal A,\qquad
 [\eta]_\varnothing\longmapsto\Xi_H(\eta).
\]
For $\sigma\in S_{m+1}$, $\pi$ crosses $\chi_j$ in the
$\sigma$-summand of \eqref{eqn: constant-term-formula-special} exactly when
$\sigma(m+1)<\sigma(j)$, by \eqref{eq:normalized-intertwiner-notation}.
Write $s_\sigma$ for the entry in the second row below. Then $w^r$
times the product of the normalizing scalars for these crossings is
$s_\sigma$ times a function holomorphic and nonzero at $t^0$.
\begin{equation}\label{eq:one-cusp-exterior-scalars}
\begin{array}{c|cccc}
 \text{position of }\pi\text{ relative to }\chi_a,\chi_{a+1}&
 \text{before both}&\chi_a,\pi,\chi_{a+1}&
 \chi_{a+1},\pi,\chi_a&\text{after both}\\ \hline
 s_\sigma & (u-w)^r&1&(u-w)^rw^r&w^r .
\end{array}
\end{equation}
\end{lemma}

\begin{proof}
Write $d_i=\beta_{k+1-i}-\beta_i$ for $1\le i\le h$. Since
$I\subset(\beta_h,\beta_{h+1})$,
\[
 |t_j^0-s_\pi|=|\max I+1-j-c|<d_i/2
 \qquad(1\le j\le m,\ 1\le i\le h).
\]
Lemma~\ref{lem:character-cuspidal-intertwiner}
gives holomorphic normalized intertwiners and inverses for all
character--$\pi$ exchanges.

\textup{(i)} At $t=t^0$, the characters have exponents
$t_1^0>\cdots>t_m^0$ and $\pi$ has exponent $s_\pi$.
Start with the order
\[
 \bigl((\chi_j)_{t_j^0\ge s_\pi},\ \pi|\det|^{t_\pi},\
       (\chi_j)_{t_j^0<s_\pi}\bigr),
\]
with $j$ increasing in each group. Moving $\pi|\det|^{t_\pi}$ to
the last position crosses precisely the $\chi_j$ with $t_j^0<s_\pi$.
For each such $j$, the rank-one normalization in the proof of
Proposition~\ref{prop:cuspidal-normal-crossings} gives
\[
 r\bigl((12),\pi\boxtimes\mathbf1,(t_\pi,t_j)\bigr)
 =e_j(t)\frac{\Lambda(t_\pi-t_j,\pi)}{\Lambda(1+t_\pi-t_j,\pi)},
 \qquad e_j\in\mathcal O_{t^0}^{\times}.
\]
At $t^0$ the two arguments are $s_\pi-t_j^0$ and $1+s_\pi-t_j^0$.
Both lie in $\tfrac12\Z$, and both are integers if $c\in\Z$.
Lemma~\ref{lem:cuspidal-L-nonvanishing} gives
$\Lambda(s,\pi)\ne0$ for $s\in\tfrac12\Z\setminus\{\tfrac12\}$.
If $r=0$, then
$\Lambda(\tfrac12,\pi)\ne0$ as well. Thus each scalar and its reciprocal
are holomorphic at $t^0$. Each standard intertwiner is the above
scalar times the normalized intertwiner of
Lemma~\ref{lem:character-cuspidal-intertwiner}. Hence their composition
and its inverse are holomorphic.

Group the $m$ characters in the first block and $\pi$ in the second,
so $Q=P_{m,k}$. Since $k>1$ and $t_i^0-t_j^0=j-i$ for the
characters, the index pairs in $\mathcal H$ of \eqref{eq:cuspidal-polar-pairs}
correspond precisely to
$(\chi_j,\chi_{j+1})$, $1\le j<m$, all in the first block.
The exponents after subtracting the two block means are
\[
 t_j^0-s_I=\frac{m+1-2j}{2}\quad(1\le j\le m),
 \qquad t_\pi^0-s_\pi=0.
\]
Thus $\boxtimes_{j=1}^m|\cdot|^{t_j^0-s_I}$ is the inducing
representation in \eqref{eqn: cuspidal-support-fix-cusp} for
$(\ell,\pi)=(m,\mathbf1)$. By the definition in
Subsection~\ref{subsec:speh-support-modules}, the corresponding module is
$\mathcal A_{m,1}=\operatorname{Fin}_{J_0}\mathcal A(\GL_m)_{\rm triv}$.
The second block has inducing datum $(\GL_k,\pi)$, for which the same
subsection gives $\mathcal A_{1,\pi}\simeq\pi$. By
\eqref{eq:one-cusp-induced-module}, the source in
\eqref{eq:separated-support-induction} is $\mathcal M$.
Proposition~\ref{lem:block-coefficient-map} gives
the surjection $\Phi:\mathcal M\twoheadrightarrow\mathcal A$.
Take the $Q$-constant term $\CT_Q\Phi(\xi)$ of the automorphic form
$\Phi(\xi)$. By \eqref{eqn: projector-constant-term}, $\Pi_Q^0$
selects its generalized infinitesimal-character component with parameters
$I-\bar c$ on $\GL_m$ and $B-\bar c$ on $\GL_k$.
If $\Phi(\xi)=0$, then $\Pi_Q^0\CT_Q\Phi(\xi)=0$, hence
$\mathcal P_\xi=0$ by definition \eqref{eqn; def-cal-P}.
Use the coordinate $y_I=-ky/(m+k)$ as in
Proposition~\ref{lem:separated-cuspidal-parameters}, and write
$\xi=\sum_{d=0}^D F_d\otimes\partial_{y_I}^d|_0$.
For $g\in M_Q(\A)^1$ and $\kappa\in K_\infty K_f$,
\eqref{eq:separated-support-constant-term} gives
\[
 0=\mathcal P_\xi(T;g,\kappa)
   =\sum_{d=0}^D T^d[F_d(\kappa)](g)\qquad(T\in\R).
\]
Taking the coefficient of $T^d$ gives $[F_d(\kappa)](g)=0$
for every $g,\kappa$. Thus each $F_d=0$, so $\xi=0$ and $\ker\Phi=0$.

For $m=1$, write $I=\{x\}$. The hypothesis $\min I\le c\le\max I$
gives $x=c$. Since $t_\pi=s_\pi+y_\pi$, at the base point we have
\[
 t_1^0=x-\bar c=c-\bar c=s_\pi=t_\pi^0.
\]
The chosen order is already $(\chi_1,\pi|\det|^{t_\pi})$, so the
permutation is the identity. Since $t_1^0-t_\pi^0=0$,
$\mathcal H=\varnothing$ by \eqref{eq:cuspidal-polar-pairs}.
The same coefficient map therefore gives the isomorphism for $Q=P_{1,k}$.

\textup{(ii)} Write $\mathcal O=\mathcal O_{t^0}$ for the local ring
in the proof of Proposition~\ref{prop:cuspidal-normal-crossings}, here on
$\sum_{j=1}^m t_j+kt_\pi=0$ at $t^0=(t_1^0,\ldots,t_m^0,s_\pi)$.
In the coordinates $z=(u_1,\ldots,u_{m-1},y)$, this is the ring
$R_z=\C\{u_1,\ldots,u_{m-1},y\}$ of germs at $z=0$ used before
Definition~\ref{def:coefficient-map}.
In passing from $P'$ to $P$, $\pi$ crosses precisely
$\chi_{a+1},\ldots,\chi_m$. For $j>a$ we have
$s_\pi-t_j^0=j-a-\tfrac12$. Hence the rank-one normalization in the
proof of Proposition~\ref{prop:cuspidal-normal-crossings} and
Lemma~\ref{lem:cuspidal-L-nonvanishing} give
\[
 \prod_{j=a+1}^{m}
 \frac{\Lambda(t_\pi-t_j,\pi)}{\Lambda(1+t_\pi-t_j,\pi)}
 \ \in\ w^r\mathcal O^\times ,
 \qquad t_\pi-t_{a+1}=\tfrac12+w .
\]
The standard intertwiner from $P'$ to $P$ is therefore $w^rbU$,
where $b\in\mathcal O^\times$ and $U$ is the normalized intertwiner.
The functional equation gives $E_{P'}=E_P\circ(w^rbU)$ and
\[
b^{-1}pE_{P'}U^{-1}=w^rpE_P=H.
\]
For $P'$ the polar pairs are $(\chi_j,\chi_{j+1})$, $1\le j<m$,
with equations $u_j=0$. Thus $pE_{P'}$ is holomorphic by
Proposition~\ref{prop:cuspidal-normal-crossings}, and its Taylor
coefficients span $\mathcal A$ \cite[Section~2]{GrbacGrobner2024}.
For $a<j\le m$, let $N_j$ be the normalized exchange of
$\pi|\det|^{t_\pi}$ with $\chi_j$, induced with the identity on the
other factors. The composition law
\cite[Section~I.1, Equation~(1), p.~607]{MoeglinWaldspurger1989} gives
\[
 U=N_m\cdots N_{a+1},\qquad
 (bU)^{-1}=b^{-1}N_{a+1}^{-1}\cdots N_m^{-1}.
\]
Since $|t_j^0-s_\pi|<d_i/2$ for every $i$,
Lemma~\ref{lem:character-cuspidal-intertwiner} gives holomorphy of
$N_j$ and $N_j^{-1}$ near $t^0$. Thus $b\in\mathcal O^\times$
gives holomorphy of $bU$ and $(bU)^{-1}$. Hence $H$ is holomorphic and
$\operatorname{im}\Xi_H=\operatorname{im}\Xi_{pE_{P'}}=\mathcal A$.

By the same rank-one normalization, the scalar in
\eqref{eq:one-cusp-exterior-scalars} belongs to
\[
 w^r\prod_{\substack{1\le j\le m\\\sigma(m+1)<\sigma(j)}}
 \frac{\Lambda(t_j-t_\pi,\pi^\vee)}{\Lambda(1+t_j-t_\pi,\pi^\vee)}
 \,\mathcal O^\times .
\]
The two exceptional factors are
\[
\begin{aligned}
 j=a:\quad&
 \frac{\Lambda(\tfrac12+u-w,\pi^\vee)}
      {\Lambda(\tfrac32+u-w,\pi^\vee)}\in (u-w)^r\mathcal O^\times,\\
 j=a+1:\quad&
 \frac{\Lambda(-\tfrac12-w,\pi^\vee)}
      {\Lambda(\tfrac12-w,\pi^\vee)}\in w^{-r}\mathcal O^\times .
\end{aligned}
\]
All other factors for these exchanges are units by
Lemma~\ref{lem:cuspidal-L-nonvanishing}. Multiplying $w^r$ by the
factor for each character crossed gives \eqref{eq:one-cusp-exterior-scalars}.

The map $\eta\mapsto[\eta]_\varnothing$ is surjective by \eqref{eq:character-coefficient-kernel}
and exactness of compact induction. For $\sigma\in S_{m+1}$, the block of
rank $k$ is in position $b=\sigma(m+1)$. The standard parabolic $Q'$
in \eqref{eqn: constant-term-formula-special} and the position $\kappa(j)$
of $\chi_j$ after deleting this block are given by
\begin{equation}\label{eq:one-cusp-output-order}
 Q'=P_{1^{b-1},k,1^{m+1-b}},\qquad
 \kappa(j)=
 \begin{cases}
  \sigma(j),&\sigma(j)<b,\\
  \sigma(j)-1,&\sigma(j)>b
 \end{cases}
 \qquad(1\le j\le m).
\end{equation}
Thus $\kappa\in S_m$, and the characters occur in the order
$(\chi_{\kappa^{-1}(1)},\ldots,\chi_{\kappa^{-1}(m)})$.
By \eqref{eq:cuspidal-coefficient-kernel} on $\GL_m$,
$(T_{P_{1^m},\kappa})_{\D}$ factors through $\Xi_m$, hence induces an
operator $\rho_\kappa$ on $\mathcal M$. Let $B_\sigma$ be the
product of the normalized character--$\pi$ exchanges multiplied by
their normalizing scalar and $w^r$. By \eqref{eq:one-cusp-exterior-scalars},
$B_\sigma$ is holomorphic. Thus \eqref{eq:holomorphic-operator-composition}
and Proposition~\ref{lem:coefficient-kernel-constant-terms} for $H$ give
\[
 \begin{aligned}
 (T_{Q',\sigma})_{\D}\eta&=(B_\sigma)_{\D}\bigl(\rho_\kappa([\eta]_\varnothing)\bigr),\\
 \{\eta:[\eta]_\varnothing=0\}&\subseteq\bigcap_{Q',\sigma}\ker(T_{Q',\sigma})_{\D}
 =\ker\Xi_H .
 \end{aligned}
\]
Since $\operatorname{im}\Xi_H=\mathcal A$, the stated map
$[\eta]_\varnothing\mapsto\Xi_H(\eta)$ is well defined and surjective.
\end{proof}

As in Subsection~\ref{subsec:character-residue-image}, $B_n=P_{1^n}$
denotes the upper triangular Borel subgroup of $\GL_n$.

\begin{lemma}\label{lem:character-cut-common-kernel}
Fix $Y\subseteq\{1,\ldots,m-1\}$ and $j\in\{1,\ldots,m-1\}\setminus Y$.
Let $d,d'$ be consecutive in $Y\cup\{0,m\}$ with $d<j<d'$.
For a permutation $\kappa$ of $d+1,\ldots,d'$, let $\rho_\kappa$
be the map on $\mathcal M_Y$ induced by $(T_{B_{d'-d},\kappa})_{\D}$
of Proposition~\ref{lem:coefficient-kernel-constant-terms} on the block
$(\chi_{d+1},\ldots,\chi_{d'})$, acting by the identity on the other blocks.
Then the cut map $P_j^Y$ of \eqref{eq:character-cut-map} satisfies
\begin{equation}\label{eq:character-cut-common-kernel}
 \ker P_j^Y=\bigcap_{\kappa(j)<\kappa(j+1)}\ker\rho_\kappa .
\end{equation}
The ordered inducing characters in the target are
$(\chi_{\kappa^{-1}(i)})_{i=d+1}^{d'}$, as in
\eqref{eq:normalized-intertwiner-notation}. The condition
$\kappa(j)<\kappa(j+1)$ is equivalent to
$(j,j+1)\notin\operatorname{Inv}(\kappa)$.
\end{lemma}

\begin{proof}
By Lemma~\ref{lem:character-cut-grading}, every class in $\ker P_j^Y$
has a representative $\eta$ supported in degrees $\alpha_j=0$, so
$u_j\eta=0$. Suppose $\kappa(j)<\kappa(j+1)$. In the proof of
Proposition~\ref{prop:cuspidal-normal-crossings}, only the factors $u_i$
with $\kappa(i)>\kappa(i+1)$ are needed to clear the poles of the
standard intertwiner. Hence $T_{B_{d'-d},\kappa}/u_j$ is holomorphic.
Apply the first identity of \eqref{eq:holomorphic-operator-composition}
with $B=T_{B_{d'-d},\kappa}/u_j$ and $T=u_j\,\mathrm{id}$.
Since $T_{\D}\eta=u_j\eta$, we obtain
\[
 (T_{B_{d'-d},\kappa})_{\D}\eta
 =B_{\D}T_{\D}\eta
 =(T_{B_{d'-d},\kappa}/u_j)_{\D}(u_j\eta)=0.
\]
Thus $\rho_\kappa[\eta]_Y=0$.

Conversely, split the block after $j$. Let $\alpha$ permute
$d+1,\ldots,j$ and $\beta$ permute $j+1,\ldots,d'$.
Their concatenation $\alpha\oplus\beta$ places $\chi_j$ before
$\chi_{j+1}$. Since
\[
 \prod_{d<i<d'}u_i
 =u_j\Bigl(\prod_{d<i<j}u_i\Bigr)\Bigl(\prod_{j<i<d'}u_i\Bigr),
\]
the definition $P_j^Y[\eta]_Y=[u_j\eta]_{Y\cup\{j\}}$ gives
\[
 \bigl(\rho_{\alpha\oplus\beta}[\eta]_Y\bigr)(k)
 =(\mathrm{id}\otimes\rho_\alpha\otimes\rho_\beta\otimes\mathrm{id})
   \bigl((P_j^Y[\eta]_Y)(k)\bigr),\qquad k\in\O_m .
\]
For the block $(\chi_{b+1},\ldots,\chi_{b+\ell})$, take
$(b,\ell)=(d,j-d)$ on the left and $(b,\ell)=(j,d'-j)$ on the right.
Put $\bar t=\ell^{-1}\sum_{i=1}^{\ell}t_{b+i}$ and $s_i=t_{b+i}-\bar t$.
After removing the common twist $|\det|^{\bar t}$, apply
\eqref{eq:cuspidal-coefficient-kernel} with
\[
 G=\GL_\ell,\qquad P=B_\ell,\qquad \tau=\mathbf1^{\boxtimes\ell},
 \qquad I_{B_\ell}(s)=\mathop{\times}_{i=1}^{\ell}|\cdot|^{s_i}.
\]
At the base point $s_i^0=(\ell+1-2i)/2=\rho_{\ell,i}$, as in
\eqref{eqn: notation-section-5.2}. Each $\gamma\in S_\ell$
permutes the diagonal entries, so $M_Q=\gamma M_{B_\ell}\gamma^{-1}=M_{B_\ell}$.
Since $Q$ is standard, $Q=B_\ell$, and hence
$\ker\Xi_\ell=\bigcap_{\gamma\in S_\ell}\ker(T_{B_\ell,\gamma})_{\D}$.
Since $\Xi_\ell$ is surjective by \eqref{eq:character-coefficient-kernel},
this gives the injection
\[
 \iota_\ell:\mathcal A_{\ell,1}^{K_{\ell,f}}\hookrightarrow
 \bigoplus_{\gamma\in S_\ell}\operatorname{im}\rho_\gamma,\qquad
 \Xi_\ell(\zeta)\longmapsto
 \bigl((T_{B_\ell,\gamma})_{\D}\zeta\bigr)_{\gamma\in S_\ell}.
\]
Thus $\iota_\ell(v)=(\rho_\gamma v)_{\gamma\in S_\ell}$, with
$\iota_1=\mathrm{id}_\C$. On the original block, $\rho_{\alpha\oplus\beta}$
is the $(\alpha\oplus\beta)$-component of $\iota_{d'-d}$.
For $k\in\O_m$,  denote $v=(P_j^Y[\eta]_Y)(k)$. It follows
\[
\begin{aligned}
 \bigl[(\mathrm{id}\otimes\iota_{j-d}\otimes\iota_{d'-j}\otimes\mathrm{id})v
      \bigr]_{\alpha,\beta}
 &=(\mathrm{id}\otimes\rho_\alpha\otimes\rho_\beta\otimes\mathrm{id})v\\
 &=(\rho_{\alpha\oplus\beta}[\eta]_Y)(k).
\end{aligned}
\]
Since tensoring over $\C$ preserves injectivity, we obtain
\[
 \rho_{\alpha\oplus\beta}[\eta]_Y=0\ \text{for all }\alpha,\beta
 \quad\Longleftrightarrow\quad
 (P_j^Y[\eta]_Y)(k)=0\ \text{for all }k\in\O_m
 \quad\Longleftrightarrow\quad P_j^Y[\eta]_Y=0.
\]
This gives $\bigcap_{\alpha,\beta}\ker\rho_{\alpha\oplus\beta}=\ker P_j^Y$
and proves the reverse inclusion.
\end{proof}

Put $\mathcal B=\ker(\mathcal M\twoheadrightarrow\mathcal A)$ for the map of
Lemma~\ref{lem:one-cusp-coefficient-map}. To filter $\mathcal A$ by
Proposition~\ref{prop:cut-filtration}, we restrict the character cut system of
Lemma~\ref{lem:character-cut-system} to $X=\varnothing$ if $c\in\Z$ or
$r=0$, and to $X=\{a\}$ if $c\in\tfrac12+\Z$ and $r>0$.
We then take normalized induction from $Q=P_{m,k}$. For $Y\subseteq X$,
the modules of the cut system for $(X,\mathcal M)$ are
\[
 M_Y=\bigl((\mathcal M_Y|\det|^{s_I})\times
           (\pi_\infty|\det|^{s_\pi})\bigr)\otimes\D_y,\qquad
 M_\varnothing=\mathcal M .
\]
The action on jets is obtained by differentiating the family with exponents
$s_I+y_I$ and $s_\pi+y_\pi$ as in \eqref{eq:one-cusp-induced-module}.
For $X=\{a\}$, denote by $P_a$ the character cut $P_a^\varnothing
:\mathcal M_\varnothing\to\mathcal M_{\{a\}}$. Its induced map
$\mathcal M\to M_{\{a\}}$, also denoted by $P_a$, is given by
\[
 (P_af)(k_0)=(P_a\otimes\mathrm{id}_{\pi_\infty}
              \otimes\mathrm{id}_{\D_y})f(k_0)\qquad(k_0\in\O_{m+k}).
\]

\begin{lemma}\label{lem:one-cusp-inducing-quotient}
If $c\in\tfrac12+\Z$ and $r>0$, then
\[
 \mathcal B=\ker P_a\cap\ker(w^r:\mathcal M\to\mathcal M),
 \qquad w=y-x_{a+1}.
\]
\end{lemma}

\begin{proof}
For $f\in\mathcal M$, choose a full-level inducing jet $\eta$ as in
\eqref{eq:full-level-inducing-jet} for $P=P_{1^m,k}$ with $f=[\eta]_\varnothing$.
By Lemma~\ref{lem:one-cusp-coefficient-map}\textup{(ii)}, the map
$\mathcal M\twoheadrightarrow\mathcal A$ sends $f$ to $\Xi_H(\eta)$,
where $H=w^rpE_P$. For each $\sigma\in S_{m+1}$, take $Q'$ and
$\kappa$ as in \eqref{eq:one-cusp-output-order}. The factorization in
the proof of that lemma gives
\[
 (T_{Q',\sigma})_{\D}\eta
 =(B_\sigma)_{\D}\,\rho_\kappa(f)
 =(B_\sigma/s_\sigma)_{\D}\,\rho_\kappa(s_\sigma f).
\]
The last equality uses $s_\sigma f=[s_\sigma\eta]_\varnothing$ and
\eqref{eq:holomorphic-operator-composition} for the scalar germ $s_\sigma$.
The factor $B_\sigma/s_\sigma$ is a unit times the product of the normalized
exchanges of $\pi|\det|^{t_\pi}$ with the characters it crosses. Its inverse
is holomorphic by Lemma~\ref{lem:character-cuspidal-intertwiner}, so its
action on jets is invertible by \eqref{eq:holomorphic-operator-composition}.
Thus the definition of $\mathcal B$ and \eqref{eq:cuspidal-coefficient-kernel} give
\[
 \begin{aligned}
 f\in\mathcal B
 &\iff \Xi_H(\eta)=0\\
 &\iff (T_{Q',\sigma})_{\D}\eta=0
       \quad\text{for every }\sigma\in S_{m+1}\\
 &\iff \rho_\kappa(s_\sigma f)=0
       \quad\text{for every }\sigma\in S_{m+1}.
 \end{aligned}
\]
On compact section values, the injection $\iota_m$ in the proof of
Lemma~\ref{lem:character-cut-common-kernel} and \eqref{eq:character-cut-common-kernel}
with $Y=\varnothing$, $j=a$, give
\[
 \bigcap_{\kappa\in S_m}\ker\rho_\kappa=\{0\},\qquad
 \bigcap_{\kappa(a)<\kappa(a+1)}\ker\rho_\kappa=\ker P_a .
\]
Hence \eqref{eq:one-cusp-exterior-scalars} gives
\[
 f\in\mathcal B \quad\iff\quad (u_a-w)^rf=0,\quad P_af=0,\quad w^rf=0 .
\]
In fact, putting $\pi$ first or last gives $(u_a-w)^r$ or $w^r$ for every
$\kappa$. The entry $1$ occurs precisely for $\kappa(a)<\kappa(a+1)$,
by inserting $\pi$ between $\chi_a$ and $\chi_{a+1}$. The remaining
entry $(u_a-w)^rw^r$ kills $\ker w^r$.
On $\ker P_a$, $u_a=0$ by Proposition~\ref{lem:character-cut-modules},
so $(u_a-w)^r=(-1)^rw^r$, as required.
\end{proof}

\begin{proof}[Proof of Theorem~\ref{lem:one-copy-internal-bound}]
By additivity we may take $F$ irreducible. Use the cut system defined
before Lemma~\ref{lem:one-cusp-inducing-quotient}, with $\alpha_\varnothing=1$.

If $X=\varnothing$, then $c\in\Z$ or $r=0$. By
Lemma~\ref{lem:one-cusp-coefficient-map}\textup{(i)}, the coefficient map
$\mathcal M\to\mathcal A$ is an isomorphism, so $\mathcal B=0$.
The only subset is $\varnothing$, hence $F_0\mathcal M=\mathcal M$
and $F_1\mathcal M=0$. Taking images gives
$\operatorname{gr}_0\mathcal A=\mathcal A\simeq\mathcal M$.

If $X=\{a\}$, take $\alpha_C$ in Proposition~\ref{prop:cut-filtration}
for $C=\{a\}$ to be multiplication by $w^r$ on $\mathcal M$,
where $w=y-x_{a+1}$.
Group the multidegrees $\alpha$ of Lemma~\ref{lem:character-cut-system}
by $\alpha_a$. On representatives $P_a$ multiplies the character jet
by $u_a$ and leaves the $y$-distribution unchanged. Since $x_{a+1}$
is linear in the $u_j$, $x_{a+1}^{|\alpha|+1}\delta_\alpha=0$,
where $|\alpha|=\sum_j\alpha_j$. Taking $N$ to be multiplication by
$x_{a+1}$ in Lemma~\ref{lem:finite-translation} gives the right inverse
$S_{x_{a+1}}$ of $w=y-x_{a+1}$,
\[
 S_{x_{a+1}}(f\otimes\delta_j^y)
 =\sum_{k\ge0}x_{a+1}^k f\otimes\delta_{j+k+1}^y
 \qquad(f\in\mathcal I_Q(0)^{K_f},\ j\ge0).
\]
Since $u_a$ commutes with $x_{a+1}$, $P_a$ commutes with
$S_{x_{a+1}}$ term by term. Thus $S_{x_{a+1}}$ preserves $\ker P_a$,
and $w^r S_{x_{a+1}}^r=\mathrm{id}$ gives $w^r(\ker P_a)=\ker P_a$.
By Lemma~\ref{lem:one-cusp-inducing-quotient},
$\mathcal B=\ker P_a\cap\ker w^r$. Proposition~\ref{prop:cut-filtration}
now gives $F_1\mathcal A=\ker P_a/\mathcal B$, $F_2\mathcal A=0$, and
\[
 P_a: \operatorname{gr}_0\mathcal A\xrightarrow{\ \sim\ }M_{\{a\}},
 \qquad
 w^r: \operatorname{gr}_1\mathcal A\xrightarrow{\ \sim\ }\ker P_a .
\]

The modules $M_\varnothing=\mathcal M$ and $M_{\{a\}}$ are induced
from $Q=P_{m,k}$ by the definition of $M_Y$ before
Lemma~\ref{lem:one-cusp-inducing-quotient}. On the character module,
Proposition~\ref{lem:character-cut-modules} gives
\[
 \ker\bigl(P_a:\mathcal A_{m,1}^{K_{m,f}}\to\mathcal M_{\{a\}}\bigr)
 =\mathcal R_{m,a}.
\]
For $f\in\mathcal M$, the formula for $(P_af)(k_0)$
therefore gives $P_af=0$ exactly when
$f(k_0)\in\mathcal R_{m,a}\otimes\pi_\infty\otimes\D_y$
for every $k_0\in\O_{m+k}$. Hence
\[
 \ker\bigl(P_a:\mathcal M\to M_{\{a\}}\bigr)=
 \bigl((\mathcal R_{m,a}|\det|^{s_I})\times
       (\pi_\infty|\det|^{s_\pi})\bigr)\otimes\D_y ,
\]
where $s_I=\frac1m\sum_{j=1}^m t_j^0$ is the mean character exponent.
The jet action is obtained by differentiating the family in
\eqref{eq:one-cusp-induced-module} with exponents $s_I+y_I$ and $s_\pi+y_\pi$.

Suppose $H^q(\mathfrak{sl}_{m+k},\O_{m+k};\mathcal A\otimes F)\ne0$.
By the long exact sequence, some $\operatorname{gr}_s\mathcal A\otimes F$
has nonzero cohomology in degree $q$. For a contributing summand the relative
split centre has cohomology only in degree zero by
Lemma~\ref{lem:split-centre-cohomology}, since $y\delta_{j+1}^y=\delta_j^y$.
By Shapiro's lemma, \eqref{eq:kostant-decomposition} and
\eqref{eq:levi-kunneth-centre},
\[
 q=\ell_W(v)+q_1+q_2,
\]
where $v\in W^Q$ and $q_1,q_2$ are the cohomological degrees on
the $\GL_m$- and $\GL_k$-blocks, respectively.

By Lemma~\ref{lem:finite-derivative-vanishing}, the ordered block parameters
are $t_{m+1-j}^0=\min I+j-1-\bar c$ and
$\beta_i-c+s_\pi=\beta_i-\bar c$. Since $v$ is increasing within
each block, \eqref{eq:ordered-kostant-parameter-comparison} gives
\[
 \ell_W(v)=\#\{(j,i):1\le j\le m,\ 1\le i\le k,\ v(j)>v(m+i)\}
 =\sum_{x\in I}\#\{i:\beta_i<x\}=hm,
\]
as $\beta_h<x<\beta_{h+1}$ for every $x\in I$.

Since $I$ is consecutive, the coefficient on $\GL_m$ is a determinant
power by Lemma~\ref{lem:consecutive-kostant-character}. Its restriction
to $\O_m$ is $1$ or $\det$. Proposition~\ref{lem:character-cut-modules}
and \eqref{eq:character-cut-degree} therefore give
\[
\begin{array}{c|c|c|c}
 X&s&\text{character factor of }\operatorname{gr}_s\mathcal A
     &\text{upper bound for }q_1\\ \hline
 \varnothing&0&D_\varnothing^\varnothing=\mathcal A_{m,1}^{K_{m,f}}&d_m-2m+2\\
 \{a\}&0&D_\varnothing^{\{a\}}=\mathcal M_{\{a\}}&d_m-2m+3\\
 \{a\}&1&D_{\{a\}}^\varnothing=\mathcal R_{m,a}&d_m-2m+4
\end{array}
\]
On the cuspidal block, Proposition~\ref{prop:speh-cohomology-bound} with
$\ell=1$ gives $q_2\le d_k-h^2$. Thus
$q\le hm+d_m-2m+2+2|X|+d_k-h^2$. Using
$d_{m+k}=d_m+d_k+2hm+1$ and $|X|\le1$, we obtain
\[
 d_{m+k}-q\ge(h+2)m+h^2-1-2|X|\ge(h+2)m+h^2-3. \qedhere
\]
\end{proof}

\subsection{Two cuspidal factors}\label{subsec:two-cuspidal-distinct-centres}

Recall from Subsection~\ref{subsec:mixed-character-interval} that
$I\subset\Z$ consists of the $m\ge1$ consecutive weights of the
inducing $\GL_1$-characters.
Let $\pi_1,\pi_2$ be full-level unitary cuspidal representations of
$\PGL_2(\A)$, with $\pi_{i,\infty}=D_{w_i}$. Let $B_i=\{b_i,b_i+w_i\}$
be their shifted weight pairs, with $b_i\in\Z$ and centres $c_i=b_i+w_i/2$.
Assume that $I,B_1,B_2$ are pairwise disjoint and put
\begin{equation}\label{eq:two-cuspidal-centering}
 \bar c=\frac{\sum_{x\in I}x+2c_1+2c_2}{m+4}.
\end{equation}
Let $\mathfrak c$ be the cuspidal support represented by
\begin{equation}\label{eq:two-cuspidal-support}
 \left(\boxtimes_{x\in I}|\cdot|^{x-\bar c}\right)
 \boxtimes\pi_1|{\det}|^{c_1-\bar c}
 \boxtimes\pi_2|{\det}|^{c_2-\bar c}.
\end{equation}
Let $\chi:\mathfrak Z_{\GL_{m+4}}^1\to\C$ have Harish--Chandra parameter
\begin{equation}\label{eq:two-copy-infinitesimal-parameter}
                   (I\sqcup B_1\sqcup B_2)-\bar c,
\end{equation}
and denote
\begin{equation}\label{eq:two-cuspidal-automorphic-module}
 \mathcal A=\bigl(\operatorname{Fin}_{\ker\chi}
                  \mathcal A(\GL_{m+4})_{\mathfrak c}\bigr)^{K_f}.
\end{equation}

\begin{theorem}\label{prop:two-cuspidal-bound}
For the two cuspidal factors above, assume $\min I\le c_2\le c_1\le\max I$.
Define $\bar c,\mathfrak c,\chi,\mathcal A$ by
\eqref{eq:two-cuspidal-centering}--\eqref{eq:two-cuspidal-automorphic-module}.
For an algebraic representation $F$ whose dual has infinitesimal
character $\chi$ on $\mathfrak{sl}_{m+4}$,
\[
 H^q(\mathfrak{sl}_{m+4},\O_{m+4};\mathcal A\otimes F)=0
                    \qquad(d_{m+4}-q<4m-4).
\]
\end{theorem}
The proof will be given at the end of this subsection.

We keep the hypotheses of Theorem~\ref{prop:two-cuspidal-bound} throughout
this subsection. For $i=1,2$, the real interval $[\min I,\max I]_{\R}$
contains $c_i$ and contains neither $b_i$ nor $b_i+w_i$, since its integer
points form $I$ and $I\cap B_i=\varnothing$. As $b_i<c_i<b_i+w_i$,
we have $I\subset(b_1,b_1+w_1)\cap(b_2,b_2+w_2)$.
Take $(t_1,\ldots,t_m,s_1,s_2)\in\C^{m+2}$ on the
hyperplane $\sum_{j=1}^m t_j+2s_1+2s_2=0$ near the base point
\[
 (t_1^0,\ldots,t_m^0,s_1^0,s_2^0)
 =(\max I-\bar c,\ldots,\min I-\bar c,\ c_1-\bar c,\ c_2-\bar c).
\]
The inducing characters are $\chi_j=|\cdot|^{t_j}$. For $i=1,2$,
$s_i$ is the complex exponent in the determinant twist
$\Pi_i=\pi_i|\det|^{s_i}$ of $\pi_i$. At $t_j=t_j^0$ and
$s_i=s_i^0=c_i-\bar c$, these factors give \eqref{eq:two-cuspidal-support}.
The trivial central characters imply $w_i$ odd. Thus $c_i\in\tfrac12+\Z$
and both $c_i\pm\tfrac12$ belong to $I$.
If $\pi_1\simeq\pi_2$, then $w_1=w_2$, so disjointness of
$B_i=\{c_i-w_i/2,c_i+w_i/2\}$ gives $c_1>c_2$.
For each $Q$, we will consider the induction of the following modules
\begin{equation}\label{eq:two-cuspidal-cases}
\begin{array}{c|c|c|c}
 &Q&M_Q&\text{full-level inducing module}\\ \hline
 \textup{(i)}\quad \pi_1\not\simeq\pi_2\ \text{or}\ c_1-c_2\ge2
  &P_{m,2,2}&\GL_m\times\GL_2\times\GL_2
  &\mathcal A_{m,1}^{K_{m,f}}\boxtimes\pi_{1,\infty}\boxtimes\pi_{2,\infty}\\
 \textup{(ii)}\quad \pi_1\simeq\pi_2=\pi,\ c_1-c_2=1
  &P_{m,4}&\GL_m\times\GL_4
  &\mathcal A_{m,1}^{K_{m,f}}\boxtimes\mathcal A_{2,\pi}^{K_{4,f}} .
\end{array}
\end{equation}
In Subsection~\ref{subsec:central-zero-merger}, $\mathcal M$ denotes the
source module in \eqref{eq:separated-support-induction} for $Q=P_{m,k}$
with one cuspidal factor $\pi$. For the two cuspidal factors
$\pi_1,\pi_2$, we take $Q$ from \eqref{eq:two-cuspidal-cases} and let
\[
 \mathcal M=\mathcal I_Q(0)^{K_{m+4,f}}\otimes\D_y,
\]
where $y_b=a_b-a_b^0$ is the change in the mean exponent on block $b$
from its value at the base point. The action is given by \eqref{eq:g-action-jets}.

Denote
\[
 \begin{gathered}
 e_i=\#\{x\in I:x>c_i\},\qquad r_i=\operatorname{ord}_{s=1/2}\Lambda(s,\pi_i),
 \qquad u_j=t_j-t_{j+1}-1,\\
 g_{i,+}=t_{e_i}-s_i-\tfrac12,\qquad
 g_{i,-}=s_i-t_{e_i+1}-\tfrac12.
 \end{gathered}
\]
Then $ g_{i,+}+g_{i,-}=u_{e_i}.$

Let $P=P_{1^m,2,2}\subset\GL_{m+4}$ be the standard parabolic with Levi
$M_P=\GL_1^m\times\GL_2\times\GL_2$ and ordered inducing factors
$(\chi_1,\ldots,\chi_m,\Pi_1,\Pi_2)$. In the notation of
Proposition~\ref{lem:coefficient-kernel-constant-terms}, take
\[
 H=aE_P,\qquad a=g_{1,-}^{r_1}g_{2,-}^{r_2}p,\qquad
 p=\begin{cases}
 \prod_{j=1}^{m-1}u_j,&\textup{(i)},\\
 v\prod_{j=1}^{m-1}u_j,&\textup{(ii)},\quad v=s_1-s_2-1.
 \end{cases}
\]
Here $\prod_{j=1}^{m-1}u_j$ is the regularizing factor of $H_m$ on
$\GL_m$ in Subsection~\ref{subsec:character-residue-image}. In case
\textup{(ii)}, the additional factor $v$ is the regularizing factor of
$H_\pi(v)=vE(v)$ on $\GL_4$ in Lemma~\ref{lem:two-cuspidal-coefficient-kernel}.

We use the notation $[\eta]_\varnothing$ introduced before
Lemma~\ref{lem:one-cusp-coefficient-map} for one cuspidal factor.
For the two cuspidal factors in \eqref{eq:two-cuspidal-cases}, put
$u=(u_1,\ldots,u_{m-1})$ and take $z=(u,y)$ in case \textup{(i)}
and $z=(u,v,y)$ in case \textup{(ii)}.
Let $K=K_\infty K_f$, let $n_1,\ldots,n_s$ be the block sizes of $Q$,
and put $K_b=K\cap\GL_{n_b}(\A)$. For a full-level compact section
$f$ of $I_P(t^0)$ and $k\in K$, its \emph{compact section value}
$f^\sharp(k)$ is the function on $K_1\times\cdots\times K_s$ given by
\[
 f^\sharp(k)(k_1,\ldots,k_s)
   =f\bigl(\operatorname{diag}(k_1,\ldots,k_s)k\bigr)
 \qquad(k_b\in K_b).
\]
This is $v^\sharp(k)$ of Remark~\ref{rem: Summation through Q} with $v=f$.
For the two parabolics in \eqref{eq:two-cuspidal-cases} we have
\[
 f^\sharp(k)\in
 \begin{cases}
 V_m\otimes\pi_{1,\infty}\otimes\pi_{2,\infty},&Q=P_{m,2,2},\\
 V_m\otimes V,&Q=P_{m,4},
 \end{cases}
\]
where $V$ is the full-level compact realization on $\GL_4$ from
Lemma~\ref{lem:two-cuspidal-coefficient-kernel}.
On the full-level inducing jets of \eqref{eq:full-level-inducing-jet}, define
\begin{subequations}
\begin{equation}\label{eq:two-cusp-inducing-jet-quotient}
 I_P(t^0)^{K_{m+4,f}}\otimes\D_z\longrightarrow\mathcal M,\qquad
 \eta\longmapsto[\eta]_\varnothing
\end{equation}
by linear extension of
\begin{equation}\label{eq:two-cusp-block-quotient}
\begin{aligned}
 \textup{(i)}\quad [f\otimes\lambda\otimes\delta]_\varnothing(k)
 &=\bigl(\Xi_m(\,\cdot\,\otimes\lambda)\otimes
   \mathrm{id}_{\pi_{1,\infty}}\otimes\mathrm{id}_{\pi_{2,\infty}}\bigr)
   f^\sharp(k)\otimes\delta,\\
 \textup{(ii)}\quad [f\otimes\lambda\otimes\mu\otimes\delta]_\varnothing(k)
 &=\bigl(\Xi_m(\,\cdot\,\otimes\lambda)\otimes
   \Xi_\pi(\,\cdot\,\otimes\mu)\bigr)f^\sharp(k)\otimes\delta .
\end{aligned}
\end{equation}
\end{subequations}
Here $\lambda\in\D_u$, $\mu\in\D_v$, and $\delta\in\D_y$.
The map \eqref{eq:two-cusp-inducing-jet-quotient} is the block coefficient surjection
in the proof of Proposition~\ref{lem:block-coefficient-map} tensored with $\mathrm{id}_{\D_y}$.

\begin{lemma}\label{lem:two-cusp-coefficient-map}
The family $H$ is holomorphic at the base point. Its coefficient map
$\Xi_H$ descends to $\mathcal M$ through the quotient map
\eqref{eq:two-cusp-inducing-jet-quotient}, giving the surjection
\[
 \mathcal M\twoheadrightarrow\mathcal A,\qquad
 [\eta]_\varnothing\longmapsto\Xi_H(\eta).
\]
Put
\[
 \mathcal B=\{[\eta]_\varnothing\in\mathcal M:\Xi_H(\eta)=0\}.
\]
Then we have the exact sequence
\begin{equation}\label{eq:two-cuspidal-kernel-exact}
 0\longrightarrow\mathcal B\longrightarrow\mathcal M
 \longrightarrow\mathcal A\longrightarrow0 .
\end{equation}
\end{lemma}

\begin{proof}
\textup{(1)} We show that $\operatorname{im}\Xi_H=\mathcal A$.
Since $c_1\ge c_2$, we have $e_1\le e_2$. Let
\[
 P'=P_{1^{e_1},\,2,\,1^{e_2-e_1},\,2,\,1^{m-e_2}}\subset\GL_{m+4}.
\]
In the order of the inducing factors of $P'$, the exponents are
\[
 \bigl((t_j^0)_{1\le j\le e_1},\,s_1^0,\,
       (t_j^0)_{e_1<j\le e_2},\,s_2^0,\,(t_j^0)_{e_2<j\le m}\bigr).
\]
They are weakly decreasing since
$t_{e_i}^0-s_i^0=s_i^0-t_{e_i+1}^0=\tfrac12$ and $s_1^0\ge s_2^0$.
For $P'$, the set $\mathcal H$ in \eqref{eq:cuspidal-polar-pairs}  consists of the pairs of positions occupied by
\[
 (\chi_j,\chi_{j+1})\quad(1\le j<m),\qquad
 (\Pi_1,\Pi_2)\quad\text{in case \textup{(ii)} only}.
\]
In fact, $t_j^0-t_k^0=k-j$, and $\pi_1\simeq\pi_2$ with $s_1^0-s_2^0=1$
exactly in case \textup{(ii)}. Thus $p=\prod_{e\in\mathcal H}u_e$.

In passing from $P'$ to $P$,
$\Pi_i$ crosses precisely the $\chi_j$ with $j>e_i$. Since
$s_i^0-t_j^0=j-e_i-\tfrac12$, the rank-one normalization in the proof of
Proposition~\ref{prop:cuspidal-normal-crossings} and Lemma~\ref{lem:cuspidal-L-nonvanishing}
give
\[
 \prod_{j>e_i}\frac{\Lambda(s_i-t_j,\pi_i)}{\Lambda(1+s_i-t_j,\pi_i)}
 \in g_{i,-}^{r_i}\mathcal O_{t^0}^\times,\qquad
 s_i-t_{e_i+1}=\tfrac12+g_{i,-}.
\]
Denote by $N$ the normalized intertwining operator from $P'$ to $P$.
By Lemma~\ref{lem:character-cuspidal-intertwiner}, $N,N^{-1}$ are holomorphic
since $|t_j^0-s_i^0|<w_i/2$. The functional equation gives
\[
 H=b^{-1}pE_{P'}N^{-1},\qquad b\in\mathcal O_{t^0}^\times .
\]
Thus
$pE_{P'}$ and $H$ are holomorphic by Proposition~\ref{prop:cuspidal-normal-crossings}.
At $t^0$, the cuspidal support and the Harish--Chandra parameter for $P'$
are \eqref{eq:two-cuspidal-support} and \eqref{eq:two-copy-infinitesimal-parameter}.
Corollary~\ref{cor:regularized-coefficient-image}
applied to $P'$ with $T=(bN)^{-1}$ gives
\[
 \operatorname{im}\Xi_H=\operatorname{im}\Xi_{pE_{P'}}
 =\bigl(\operatorname{Fin}_{\ker\chi}\mathcal A(\GL_{m+4})_{\mathfrak c}\bigr)^{K_f}
 =\mathcal A
\]
by \eqref{eq:two-cuspidal-automorphic-module}.

\textup{(2)} It remains to show that $\Xi_H$ factors through $\mathcal M$.
Suppose $[\eta]_\varnothing=0$.
Recall from \eqref{eq:two-cuspidal-cases} that $P=P_{1^m,2,2}\subseteq Q$,
where $Q=P_{m,2,2}$ in case \textup{(i)} and $Q=P_{m,4}$ in case \textup{(ii)}.
For the $w$-summand of \eqref{eqn: constant-term-formula-special} with output parabolic $Q'$,
the condition $P_w=P$ gives
\[
 w^{-1}M_{Q'}w=M_P\subseteq M_Q .
\]
Thus $Q'$ has blocks of size two at positions $w(m+1),w(m+2)$
and blocks of size one at all other positions.
Let $D$ index the inducing factors in a block of $Q$. They occupy the
positions $w(D)$ in $Q'$.
Take $\tau|_D:D\to w(D)$ increasing and put $\kappa=\tau^{-1}w$.
Then $\kappa(D)=D$. The exchanges within each $D$ occur in $\kappa$
and those between distinct blocks occur in $\tau$, so
\begin{equation}\label{eq:two-cusp-intertwiner-factorization}
 M(w)=M(\tau)M(\kappa),\qquad
 T_{Q',w}=\bigl(g_{1,-}^{r_1}g_{2,-}^{r_2}M(\tau)\bigr)\,pM(\kappa).
\end{equation}

In the inducing order $(\chi_1,\ldots,\chi_m,\Pi_1,\Pi_2)$,
$\chi_j$ is in position $j$ and $\Pi_i$ is in position $m+i>j$.
By \eqref{eq:normalized-intertwiner-notation}, their positions in the
$w$-summand are $w(j)$ and $w(m+i)$, respectively. The pair contributes
an exchange exactly when $\Pi_i$ now precedes $\chi_j$, that is,
when $w(m+i)<w(j)$.
By the calculation of \eqref{eq:one-cusp-exterior-scalars}, the scalar
for this exchange is, up to a holomorphic unit,
$\Lambda(t_j-s_i,\pi_i^\vee)/\Lambda(1+t_j-s_i,\pi_i^\vee)$.
Since $1+t_j^0-s_i^0=\tfrac32+e_i-j$, the denominator is a unit for
$j\ne e_i+1$ by Lemma~\ref{lem:cuspidal-L-nonvanishing}. For $j=e_i+1$,
\[
 \Lambda(1+t_{e_i+1}-s_i,\pi_i^\vee)
 =\Lambda(\tfrac12-g_{i,-},\pi_i^\vee)
 \in g_{i,-}^{r_i}\mathcal O_{t^0}^{\times}.
\]
The numerators are holomorphic, so for $i=1,2$,
\[
 g_{i,-}^{r_i}\prod_{\substack{1\le j\le m\\w(m+i)<w(j)}}
 \frac{\Lambda(t_j-s_i,\pi_i^\vee)}{\Lambda(1+t_j-s_i,\pi_i^\vee)}
 \in\mathcal O_{t^0}.
\]
The normalized intertwining operators for these exchanges are holomorphic by
Lemma~\ref{lem:character-cuspidal-intertwiner}, since $|t_j^0-s_i^0|<w_i/2$.
In case \textup{(i)} of \eqref{eq:two-cuspidal-cases}, $s_1^0-s_2^0=c_1-c_2\ge0$.
By the proof of
Proposition~\ref{prop:cuspidal-normal-crossings}, the normalized operator
exchanging $\Pi_1,\Pi_2$ is holomorphic. A pole of its scalar requires both
$\pi_1\simeq\pi_2$ and $c_1-c_2\in\{0,1\}$, a combination excluded by case \textup{(i)}.
In case \textup{(ii)} of \eqref{eq:two-cuspidal-cases}, $Q=P_{m,4}$
and $D=\{m+1,m+2\}$ indexes the two factors $\Pi_1,\Pi_2$ of the
$\GL_4$-block. Since $w=\tau\kappa$ and $\tau|_D$ is increasing,
\[
 w(m+1)>w(m+2)\quad\Longleftrightarrow\quad
 \kappa(m+1)>\kappa(m+2).
\]
The regularizing factor of $H_m$ on $\GL_m$ is $\prod_{j=1}^{m-1}u_j$
by \eqref{eq:character-regularized-family}. That of $H_\pi(v)=vE(v)$
on $\GL_4$ is $v=s_1-s_2-1$ by Lemma~\ref{lem:two-cuspidal-coefficient-kernel}.
Their product gives $p=v\prod_{j=1}^{m-1}u_j$.
Since $s_1-s_2=1+v$, the induction on the $\GL_4$-block is
\[
 \Pi_1\times\Pi_2
 \simeq
 \underbrace{\bigl(\pi|\det|^{(1+v)/2}\times
                   \pi|\det|^{-(1+v)/2}\bigr)}_{I(v)}
 \otimes|\det|^{(s_1+s_2)/2}.
\]
Here $I(v)$ is the family introduced before Lemma~\ref{lem:two-cuspidal-coefficient-kernel}.
The determinant twist on $\GL_4$ factors out of the standard intertwining
integral, so the operator exchanging $\Pi_1,\Pi_2$ on compact sections
is $M(v)$ from the proof of that lemma. Thus the $\GL_4$-factor of $pM(\kappa)$ is
\[
 \begin{cases}
 v\,\mathrm{id},&(\kappa(m+1),\kappa(m+2))=(m+1,m+2),\\
 vM(v)=\bigl(vr(v)\bigr)N(v),&(\kappa(m+1),\kappa(m+2))=(m+2,m+1).
 \end{cases}
\]
In the second case the normalizing scalar is
\[
 r(v)=\epsilon(v)^{-1}
 \frac{\Lambda(1+v,\pi\times\pi^\vee)}{\Lambda(2+v,\pi\times\pi^\vee)}.
\]
Its numerator has a simple pole at $v=0$, and its denominator and
$\epsilon(v)$ are units there. Hence $vr(v)$ is the holomorphic unit
$a(v)$ multiplying $N(v)$ in \eqref{eq:two-cuspidal-cleared-constant-term}.
Thus $g_{1,-}^{r_1}g_{2,-}^{r_2}M(\tau)$ is holomorphic.

Put $\kappa_m=\kappa|_{\{1,\ldots,m\}}$. On $\GL_m$ we have
$T_{B_m,\kappa_m}=(\prod_{j=1}^{m-1}u_j)M(\kappa_m)$ by
\eqref{eq:character-regularized-family}. In case \textup{(i)}, $\kappa$
acts by the identity on each $\pi_{i,\infty}$. In case \textup{(ii)},
the two possibilities on $\GL_4$ are $v\,\mathrm{id}$ and $vM(v)$
by \eqref{eq:two-cuspidal-cleared-constant-term}.
By \eqref{eq:cuspidal-coefficient-kernel} of Proposition~\ref{lem:coefficient-kernel-constant-terms}
applied to $H_m$, $(T_{B_m,\kappa_m})_{\D}$ factors through $\Xi_m$.
The same formula for $H_\pi$ in case \textup{(ii)} shows that both
$(v\,\mathrm{id})_{\D}$ and $(vM(v))_{\D}$ factor through $\Xi_\pi$.
These operators depend only on exponent differences within each block.
Taking their tensor product with $\mathrm{id}_{\D_y}$ on compact section values
as in \eqref{eq:two-cusp-block-quotient} therefore gives
\[
 \ker\bigl(\eta\mapsto[\eta]_\varnothing\bigr)\subseteq
 \ker\bigl((pM(\kappa))_{\D}\bigr).
\]
Hence we obtain a well-defined map on $\mathcal M$ by
\begin{equation}\label{eq:two-cusp-kappa-on-jets}
 \rho_\kappa([\eta]_\varnothing):=(pM(\kappa))_{\D}\eta .
\end{equation}
Since $[\eta]_\varnothing=0$, \eqref{eq:holomorphic-operator-composition} gives
\[
 (T_{Q',w})_{\D}\eta
 =(g_{1,-}^{r_1}g_{2,-}^{r_2}M(\tau))_{\D}\,
 \rho_\kappa([\eta]_\varnothing)=0 .
\]
By \eqref{eq:cuspidal-coefficient-kernel} for $H$, $\Xi_H(\eta)=0$.
Thus $[\eta]_\varnothing\mapsto\Xi_H(\eta)$ defines a map
$\mathcal M\to\mathcal A$ with kernel $\mathcal B$.
It is surjective by \textup{(1)}, giving \eqref{eq:two-cuspidal-kernel-exact}.
\end{proof}

Recall that $e_i=\#\{x\in I:x>c_i\}=\max I-c_i+\tfrac12$ for $i=1,2$.
In the decreasing order on $I$, positions $e_i,e_i+1$ have weights
$c_i+\tfrac12,c_i-\tfrac12$, respectively. Let $X_I$ index the cuts
of \eqref{eq:cut-after-index} between $\chi_{e_i}$ and $\chi_{e_i+1}$
for $i=1,2$. In case \textup{(ii)}, put
$z=m+1$. In the inducing order $(\chi_1,\ldots,\chi_m,\Pi_1,\Pi_2)$,
the cut after position $z$ separates $\Pi_1$ from $\Pi_2$ in the $\GL_4$-block. Then
\[
 X_I=\{e_1,e_2\},\qquad
 X=\begin{cases}
 X_I,&\textup{(i)},\\
 X_I\sqcup\{z\},&\textup{(ii)}.
 \end{cases}
\]
Thus $X_I=\{e_1\}$ if $e_1=e_2$. On $\mathcal M$, let
$P_{e_i}^\varnothing$ be induced by the character cut map
\eqref{eq:character-cut-map}, also when $r_i=0$. In case \textup{(ii)},
let $P_z^\varnothing$ be induced by the map
\begin{equation}\label{eq:two-cusp-cut-map}
 P_v:\mathcal A_{2,\pi}^{K_{4,f}}\longrightarrow V\otimes\D_v,\qquad
 [\eta]\longmapsto v\eta.
\end{equation}
Here $V$ is the full-level compact realization of $I(v)$ on $\GL_4$,
and $\mathcal A_{2,\pi}^{K_{4,f}}$ is the quotient of $V\otimes\D_v$
by $(\ker R_\pi)\otimes\operatorname{ev}_v$ in
Lemma~\ref{lem:two-cuspidal-coefficient-kernel} with $k=2$.
For $C\subseteq X$, let
\begin{equation}
\label{eqn: def-D-C}
 D_C=\bigcap_{j\in C}\ker P_j^\varnothing\subseteq\mathcal M,\qquad D_\varnothing=\mathcal M.
\end{equation}

For $Y\subseteq X$, put $D=Y\cap X_I$ and
$\mathcal M_D=(V_m\otimes\D_u)/\mathcal K^D$ as in \eqref{eqn: def-k^D}.
Recall that the character exponents at the base point are
$t_j^0=\max I+1-j-\bar c$ for $1\le j\le m$. Their mean is
\[
 a_1^0=\frac1m\sum_{j=1}^m t_j^0=\frac{\max I+\min I}{2}-\bar c.
\]

\begin{subequations}\label{eq:two-cusp-cut-modules}
For $Q=P_{m,2,2}$, take $(a_2^0,a_3^0)=(s_1^0,s_2^0)$ and define
\begin{equation}\label{eq:two-cusp-cut-modules-m22}
 \mathcal M_{Q,Y}:=
 \bigl((\mathcal M_D|\det|^{a_1^0})
 \times(\pi_{1,\infty}|\det|^{a_2^0})
 \times(\pi_{2,\infty}|\det|^{a_3^0})\bigr)\otimes\D_y .
\end{equation}
For $Q=P_{m,4}$, take $a_2^0=(s_1^0+s_2^0)/2$ and define
\begin{equation}\label{eq:two-cusp-cut-modules-m4}
 \mathcal M_{Q,Y}:=
 \begin{cases}
 \bigl((\mathcal M_D|\det|^{a_1^0})\times
    (\mathcal A_{2,\pi}^{K_{4,f}}|\det|^{a_2^0})\bigr)\otimes\D_y,
       &z\notin Y,\\[2pt]
 \bigl((\mathcal M_D|\det|^{a_1^0})\times
    ((V\otimes\D_v)|\det|^{a_2^0})\bigr)\otimes\D_y,
       &z\in Y .
 \end{cases}
\end{equation}
\end{subequations}
Here $\D_y$ consists of finite-order distributions at zero on
$my_1+2y_2+2y_3=0$ or $my_1+4y_2=0$, respectively.
In particular $\mathcal M_{Q,\varnothing}=\mathcal M$.

Let $\eta$ be a full-level inducing jet for $P=P_{1^m,2,2}$, and write
\[
 \eta=\sum_\alpha f_\alpha\otimes\delta_\alpha
 \in I_P(t^0)^{K_{m+4,f}}\otimes
 \begin{cases}
 \D_u\otimes\D_y,&Q=P_{m,2,2},\\
 \D_u\otimes\D_v\otimes\D_y,&Q=P_{m,4}.
 \end{cases}
\]
Its image under \eqref{eq:two-cusp-inducing-jet-quotient} is
$[\eta]_\varnothing\in\mathcal M$. For each $\alpha$, the values
$f_\alpha^\sharp(k)(1,\ldots,1)=f_\alpha(k)$
for $k\in K=K_\infty K_f$ determine $f_\alpha$ by
Lemma~\ref{lem:compact-realization}. As in \eqref{eq:two-cusp-block-quotient},
we define the image $[\eta]_Y\in\mathcal M_{Q,Y}$ of the same inducing jet
by giving the map on pure tensors with first factor $f$.
For $h\in V_m$ and $\lambda\in\D_u$, denote
$[h\otimes\lambda]_D=h\otimes\lambda+\mathcal K^D$ in $\mathcal M_D$.

\textup{(i)} For $Q=P_{m,2,2}$, if $f^\sharp(k)=h\otimes f_1\otimes f_2$
with $f_i\in\pi_{i,\infty}$, put
\[
 [f\otimes\lambda\otimes\delta]_Y(k)
   =[h\otimes\lambda]_D\otimes f_1\otimes f_2\otimes\delta
   \qquad(\delta\in\D_y).
\]
\textup{(ii)} For $Q=P_{m,4}$, take $\mu\in\D_v$ and $\delta\in\D_y$.
If $f^\sharp(k)=h\otimes g$ with $g\in V$, put
\[
 [f\otimes\lambda\otimes\mu\otimes\delta]_Y(k)=
 \begin{cases}
 [h\otimes\lambda]_D\otimes\Xi_\pi(g\otimes\mu)\otimes\delta,&z\notin Y,\\
 [h\otimes\lambda]_D\otimes(g\otimes\mu)\otimes\delta,&z\in Y.
 \end{cases}
\]
Extend both formulas by linearity. Write $r=\operatorname{diag}(r_1,\ldots,r_s)$
in $M_Q(\R)\cap K_\infty=\prod_{b=1}^s\O_{n_b}$. By the definition of
$f^\sharp(k)$,
\[
 f^\sharp(rk)(k_1,\ldots,k_s)=f^\sharp(k)(k_1r_1,\ldots,k_sr_s).
\]
Thus $r_b$ acts by right translation on the $b$-th factor and fixes
the distribution factors. On $\mathcal M_D$ this gives
$r_1[h\otimes\lambda]_D=[(r_1h)\otimes\lambda]_D$, since
$\mathcal K^D$ in \eqref{eqn: def-k^D} is $\O_m$-stable.
In case \textup{(i)}, $r_2,r_3$ act on $f_1,f_2$ by
$\pi_{1,\infty}(r_2),\pi_{2,\infty}(r_3)$, respectively.
In case \textup{(ii)}, $r_2\in\O_4$ acts on $g$, and
$\Xi_\pi((r_2g)\otimes\mu)=r_2\Xi_\pi(g\otimes\mu)$ by
Lemma~\ref{lem:two-cuspidal-coefficient-kernel}. In both cases we obtain
\[
 [\eta]_Y(rk)=r\cdot[\eta]_Y(k)
 \qquad\bigl(r\in M_Q(\R)\cap K_\infty\bigr).
\]
These pointwise maps preserve right $K_f$-invariance and $K_\infty$-finiteness.
By the Iwasawa decomposition
$\GL_{m+4}(\R)=Q(\R)K_\infty$, these values determine a unique
section in $\mathcal M_{Q,Y}$.
For $Y=\varnothing$,
\eqref{eq:character-coefficient-kernel} identifies
$[h\otimes\lambda]_\varnothing$ with $\Xi_m(h\otimes\lambda)$, giving
\eqref{eq:two-cusp-block-quotient}.

For $j\in X\setminus Y$, the additional cut
$P_j^Y:\mathcal M_{Q,Y}\to\mathcal M_{Q,Y\cup\{j\}}$ is given by
\begin{subequations}
\begin{align}
 P_j^Y[\eta]_Y&=[u_j\eta]_{Y\cup\{j\}} &&(j\in X_I\setminus Y),
 \label{eq:two-cusp-character-cut-on-jets}\\
 P_z^Y[\eta]_Y&=[v\eta]_{Y\cup\{z\}} &&(\text{case \textup{(ii)}},\ z\notin Y).
 \label{eq:two-cusp-z-cut-on-jets}
\end{align}
\end{subequations}
Let $p_j^Y$ be the corresponding map on the inducing Levi module, so that
\begin{equation}\label{eq:two-cusp-cut-on-levi-values}
 (P_j^Y\xi)(k)=p_j^Y\bigl(\xi(k)\bigr)
 \qquad(\xi\in\mathcal M_{Q,Y},\ k\in K_\infty).
\end{equation}
The multiplications commute. Define $P_Y:\mathcal M\to\mathcal M_{Q,Y}$
as their composite and $P_\varnothing=\mathrm{id}$. For $C\subseteq X$
put $E_C=P_{X\setminus C}(D_C)$, with $D_C$ as in \eqref{eqn: def-D-C}.

Fix $C\subseteq X$.
For $w=\tau\kappa$ as in \eqref{eq:two-cusp-intertwiner-factorization},
$\kappa(D)=D$ for each set $D$ indexing the inducing factors in a
block of $Q$, and $\tau|_D$ is increasing. In case \textup{(i)},
$\kappa$ fixes the two cuspidal positions since each forms a separate block.
In case \textup{(ii)}, $\Pi_1,\Pi_2$ occupy positions $z=m+1,z+1=m+2$
in the same block. Consider those $\kappa$ for which
\begin{equation}\label{eq:two-cusp-internal-orders}
\begin{aligned}
 \kappa(\{1,\ldots,e\})&=\{1,\ldots,e\} &&(e\in X_I\setminus C),\\
 \kappa(e+1)&<\kappa(e) &&(e\in C\cap X_I),\\
 (\kappa(z),\kappa(z+1))&=(z,z+1) &&(\text{case \textup{(ii)}},\ z\notin C),\\
 (\kappa(z),\kappa(z+1))&=(z+1,z) &&(\text{case \textup{(ii)}},\ z\in C).
\end{aligned}
\end{equation}
Since $\tau(z)<\tau(z+1)$ in case \textup{(ii)}, the cuspidal order under
$w=\tau\kappa$ is
$(\Pi_1,\Pi_2)$ if $z\notin C$ and $(\Pi_2,\Pi_1)$ if $z\in C$.

\begin{lemma}\label{lem:two-cusp-complementary-cut-detection}
Let $Y=X\setminus C$ and $x\in\mathcal M$ with $P_Yx\in E_C$.
Then $P_Yx=0$ if and only if $\rho_\kappa(x)=0$ for every permutation
$\kappa$ preserving each block of $Q$ and satisfying
\eqref{eq:two-cusp-internal-orders}.
\end{lemma}

\begin{proof}
\textup{(1)} Put $D=Y\cap X_I$ and $\xi=P_Yx$. For $x=[\eta]_\varnothing$,
iterating \eqref{eq:two-cusp-character-cut-on-jets} gives
\[
 P_Dx=\Bigl[\Bigl(\prod_{j\in D}u_j\Bigr)\eta\Bigr]_D
       \in\mathcal M_{Q,D}.
\]
Thus $\xi=P_Dx$ if $Y=D$. In case \textup{(ii)} with $z\in Y$,
$\xi=P_z^D(P_Dx)$.
Let $J$ run over the intervals obtained by cutting after each $j\in D$.
The first condition of \eqref{eq:two-cusp-internal-orders} is
\[
 \kappa(\{1,\ldots,j\})=\{1,\ldots,j\}\quad(j\in D)
 \quad\Longleftrightarrow\quad
 \kappa|_{\{1,\ldots,m\}}=\bigoplus_J\gamma_J,\qquad \gamma_J:J\xrightarrow{\sim}J .
\]
On the $\GL_m$-factor, $P_D$ is the
composite $\mathcal M_\varnothing\to\mathcal M_D$ of the cut maps
\eqref{eq:character-cut-map} at $j\in D$, with
$\mathcal M_\varnothing\simeq\mathcal A_{m,1}^{K_{m,f}}$ by
\eqref{eq:character-coefficient-kernel}.
For the character coordinates $x_j$ of \eqref{eqn: notation-section-5.2}, put
$\bar x_J=|J|^{-1}\sum_{j\in J}x_j$.
Let $\D_{\bar x}$ be the space of finite-order distributions at the origin
of $\sum_J|J|\bar x_J=0$, as in the proof of
Proposition~\ref{prop:character-block-cut-degree}.
Take compact section values in the realization induced in stages through the
intervals $J$.
For $f\in\mathcal A_{m,1}^{K_{m,f}}$, iterating the formula for
$\rho_{\alpha\oplus\beta}$ in the proof of Lemma~\ref{lem:character-cut-common-kernel} gives
\[
 (\rho_{\oplus_J\gamma_J}f)(k)
 =\Bigl(\Bigl(\bigotimes_J\rho_{\gamma_J}\Bigr)\otimes\mathrm{id}_{\D_{\bar x}}\Bigr)
    \bigl((P_Df)(k)\bigr),\qquad k\in\O_m .
\]
On each term of $(P_Df)(k)$, this gives
\[
 \Bigl(\bigotimes_J f_J\Bigr)\otimes\delta
 \longmapsto \Bigl(\bigotimes_J\rho_{\gamma_J}f_J\Bigr)\otimes\delta,
 \qquad f_J\in\mathcal A_{|J|,1}^{K_{|J|,f}},\quad\delta\in\D_{\bar x}.
\]
In $\mathcal M_{Q,D}$, tensor this map with the identity on the cuspidal
factors and on $\D_y$.
By \eqref{eq:cuspidal-coefficient-kernel}, as in the proof of
Lemma~\ref{lem:character-cut-common-kernel},
$\iota_{|J|}:f_J\mapsto(\rho_{\gamma_J}f_J)_{\gamma_J}$ is injective. Thus
\[
 \bigcap_{(\gamma_J)_J}\ker\Bigl(\bigotimes_J\rho_{\gamma_J}\Bigr)
 =\ker\Bigl(\bigotimes_J\iota_{|J|}\Bigr)=0 .
\]
Since $\xi\in E_C=P_Y(D_C)$, choose $x_C\in D_C$ with $\xi=P_Yx_C$.
For every $j\in C$, we may apply the cut at $j$ first. Since
$P_j^\varnothing x_C=0$, we have
\[
 P_j^Y\xi=P_{Y\cup\{j\}}x_C=0 .
\]
For $j\in C\cap X_I$, let $J_0$ be the interval containing $j,j+1$.
On compact section values, $P_j^Y$ acts by $P_j^D$ on $\mathcal M_D$
and by the identity on the other factors, by \eqref{eq:two-cusp-character-cut-on-jets}.
If $\gamma_{J_0}(j)<\gamma_{J_0}(j+1)$, then
$\ker P_j^D\subseteq\ker\rho_{\gamma_{J_0}}$ by
\eqref{eq:character-cut-common-kernel}. Applying $\rho_{\gamma_J}$ on
the remaining intervals gives
\[
 \Bigl(\Bigl(\bigotimes_J\rho_{\gamma_J}\Bigr)\otimes\mathrm{id}\Bigr)\xi(k)=0
 \qquad(k\in K_\infty K_f).
\]

\textup{(2)} To obtain $\rho_\kappa(x)$ from $\xi=P_Yx$,
recall that the $\GL_m$-factor of $\rho_\kappa$ in
\eqref{eq:two-cusp-kappa-on-jets} acts on $x$ by $\rho_{\kappa_m}$,
where $\kappa_m=\bigoplus_J\gamma_J$ and $\gamma_J=\kappa|_J$.
Part \textup{(1)} factors this action through $P_D$. Since
$\xi=P_Dx$ for $Y=D$ and $\xi=P_z^D(P_Dx)$ for $z\in Y$,
we obtain, on compact values at $k\in K_\infty K_f$,
\[
 \Bigl(\Bigl(\bigotimes_J\rho_{\gamma_J}\Bigr)\otimes\mathrm{id}\Bigr)\xi(k)
 =
 \begin{cases}
  (\rho_{\kappa_m}\otimes\mathrm{id})\,x(k),&Y=D,\\
  (\rho_{\kappa_m}\otimes P_v\otimes\mathrm{id}_{\D_y})\,x(k),&z\in Y .
 \end{cases}
\]
On the cuspidal factors, \eqref{eq:two-cusp-kappa-on-jets} gives
$\mathrm{id}_{\pi_{1,\infty}}\otimes\mathrm{id}_{\pi_{2,\infty}}$ in case \textup{(i)}.
In case \textup{(ii)} the two operators on inducing jets are
$(v\,\mathrm{id}_V)_\D$ and $(T_{P_{2,2},(12)})_\D$ by
\eqref{eq:two-cuspidal-cleared-constant-term}.
For $\zeta\in V\otimes\D_v$, denote by $[\zeta]$  its class in
\[
 (V\otimes\D_v)/\ker\Xi_\pi \simeq \mathcal A_{2,\pi}^{K_{4,f}},
 \qquad \ker\Xi_\pi=(\ker R_\pi)\otimes\C\operatorname{ev}_v,
\]
as in Lemma~\ref{lem:two-cuspidal-coefficient-kernel}.

If $z\in Y$, then $\kappa$ fixes $z,z+1$, so
$T_{P_{2,2},1}(v)=v\,\mathrm{id}_V$ by \eqref{eq:two-cuspidal-cleared-constant-term}.
By \eqref{eq:two-cusp-cut-map}, $P_z^D$ sends each $\GL_4$-coefficient
$[\zeta]$ of $P_Dx$ to $v\zeta$. In passing from $\xi=P_z^D(P_Dx)$
to $\rho_\kappa(x)$ in \eqref{eq:two-cusp-kappa-on-jets}, the map on
the $\GL_4$-factor of $\xi$ is
\[
 \operatorname{im}P_v\longrightarrow V\otimes\D_v,\qquad
 P_v([\zeta])=v\zeta\longmapsto(T_{P_{2,2},1})_\D\zeta=v\zeta .
\]
Thus whenever $\kappa$ fixes the cuspidal positions,
$\rho_\kappa(x)(k)=((\bigotimes_J\rho_{\gamma_J})\otimes\mathrm{id})\xi(k)$.

If $z\in C$, then $Y=D$ and $\kappa$ exchanges $\Pi_1,\Pi_2$ in
positions $z,z+1$. On the two $\GL_2$-blocks of $P_{2,2}\subset\GL_4$,
this permutation is $(12)$, with $T_{P_{2,2},(12)}(v)=vM(v)$
by \eqref{eq:two-cuspidal-cleared-constant-term}.
The cut map $P_Y$ acts by the identity on the $\GL_4$-factor.
Fix $k\in K_\infty K_f$ and choose a decomposition
$\xi(k)=\sum_i b_i\otimes[\zeta_i]$, with
$[\zeta_i]\in\mathcal A_{2,\pi}^{K_{4,f}}$ and the $b_i$ linearly independent in
$(\bigotimes_J\mathcal A_{|J|,1}^{K_{|J|,f}})\otimes\D_{\bar x}\otimes\D_y$.
By \textup{(1)} and \eqref{eq:two-cusp-z-cut-on-jets},
\[
 0=(P_z^Y\xi)(k)=\sum_i b_i\otimes P_v([\zeta_i]).
\]
Linear independence gives $P_v([\zeta_i])=0$ for every $i$.
The $\GL_4$-factor of $\rho_\kappa(x)$ in \eqref{eq:two-cusp-kappa-on-jets}
is therefore obtained by applying to each $[\zeta_i]$ the map
\[
 \ker P_v\longrightarrow\operatorname{im}\bigl((T_{P_{2,2},(12)})_\D\bigr),\qquad
 [\zeta]\longmapsto(T_{P_{2,2},(12)})_\D\zeta=(vM(v))_\D\zeta .
\]
Hence in this case
\[
 \rho_\kappa(x)(k)
 =\sum_i\Bigl(\Bigl(\bigotimes_J\rho_{\gamma_J}\Bigr)
       \otimes\mathrm{id}_{\D_{\bar x}\otimes\D_y}\Bigr)b_i
       \otimes(T_{P_{2,2},(12)})_\D\zeta_i .
\]

\textup{(3)} If $\xi=0$, the formulas in \textup{(1)} and \textup{(2)}
give $\rho_\kappa(x)=0$ for every $\kappa$ in \eqref{eq:two-cusp-internal-orders}.
Suppose $\xi\ne0$ and choose $k\in K_\infty K_f$ with $\xi(k)\ne0$.
Since $(\bigotimes_J\iota_{|J|})\otimes\mathrm{id}$ is injective,
there is a tuple $(\gamma_J)_J$ for which
\[
 t=\Bigl(\Bigl(\bigotimes_J\rho_{\gamma_J}\Bigr)\otimes\mathrm{id}\Bigr)
       \xi(k)\ne0 .
\]
By \textup{(1)}, this tuple has $\gamma_J(j+1)<\gamma_J(j)$ for every
$j\in C\cap X_I$ in $J$. Set $\kappa|_J=\gamma_J$ and take the
cuspidal order specified in \eqref{eq:two-cusp-internal-orders}.
If $\kappa$ fixes the cuspidal positions, then $\rho_\kappa(x)(k)=t$.
If it exchanges them, choose a decomposition $t=\sum_i b_i\otimes[\zeta_i]$ with
$[\zeta_i]\in\ker P_v$ and $b_i$ linearly independent in the tensor
product of the other factors. Since $t\ne0$, choose $i$ with
$[\zeta_i]\ne0$. We have $v\zeta_i=P_v([\zeta_i])=0$, and
\eqref{eq:cuspidal-coefficient-kernel} for $H_\pi$ gives
\[
 \ker\Xi_\pi
 =\{\zeta\in V\otimes\D_v:v\zeta=0,\ (T_{P_{2,2},(12)})_\D\zeta=0\}.
\]
Since $\zeta_i\notin\ker\Xi_\pi$, it follows that
$(T_{P_{2,2},(12)})_\D\zeta_i\ne0$. By linear independence of the $b_i$,
\[
 \rho_\kappa(x)(k)=\sum_i b_i\otimes(T_{P_{2,2},(12)})_\D\zeta_i\ne0 .
\]
This proves the converse.
\end{proof}

In the inducing order for $P=P_{1^m,2,2}$, the $j$-th factor is
$\chi_j$ for $1\le j\le m$, $\Pi_1$ for $j=m+1$, and $\Pi_2$ for $j=m+2$.
Given $w\in S_{m+2}$, let $Q'$ have blocks of size two at positions
$w(m+1),w(m+2)$ and blocks of size one elsewhere.
Let $\widetilde w\in S_{m+4}$ map the coordinates of each block $j$ of $P$
in increasing order onto those of block $w(j)$ of $Q'$.
This is the Weyl element indexing the corresponding summand of
\eqref{eqn: constant-term-formula-special} in the proof of
Lemma~\ref{lem:two-cusp-coefficient-map}, where
$\widetilde w^{-1}M_{Q'}\widetilde w=M_P$.
As in that proof, we denote both $\widetilde w$ and its block permutation by $w$.
The $j$-th inducing factor occupies the $w(j)$-th block
of $M_{Q'}$ by \eqref{eq:normalized-intertwiner-notation}.
For inducing factors $A,B$ in positions $a,b$ of $P$, respectively,
write $A\prec_w B$ if $w(a)<w(b)$. In particular,
\[
 \Pi_i\prec_w\chi_j\quad\Longleftrightarrow\quad w(m+i)<w(j)
 \qquad (i=1,2,\ 1\le j\le m).
\]

For $w=\tau\kappa$ and a jet $\eta$, put $x=[\eta]_\varnothing$.
The factors $g_{i,\pm}$ below act on quotient classes by
$g_{i,\pm}[\eta]_Y=[g_{i,\pm}\eta]_Y$. Indeed, for any holomorphic
scalar germ $f$ at $t^0$ we have
\[
 f(\operatorname{ev}_{u_j}\otimes\delta)
 =\operatorname{ev}_{u_j}\otimes(f|_{u_j=0}\,\delta),
\]
and likewise for $v$. Thus $f$ preserves the relations in
\eqref{eqn: def-k^D} and Lemma~\ref{lem:two-cuspidal-coefficient-kernel}.
By \eqref{eq:one-cusp-exterior-scalars}, the scalar for crossing
$\Pi_i$ with $\chi_j$ belongs to
\[
 \mathcal O_{t^0}^{\times}\cdot
 \begin{cases}
  g_{i,+}^{r_i},&j=e_i,\\
  g_{i,-}^{-r_i},&j=e_i+1,\\
  1,&j\notin\{e_i,e_i+1\}.
 \end{cases}
\]
Crossing $\chi_{e_i+1}$ cancels the prefactor $g_{i,-}^{r_i}$.
Thus \eqref{eq:two-cusp-intertwiner-factorization} and
\eqref{eq:two-cusp-kappa-on-jets} give
\begin{equation}\label{eq:two-cusp-shuffle-jets}
 (T_{Q',w})_\D\eta=(B_w)_\D\rho_\kappa\left(
 \prod_{i:\,\Pi_i\prec_w\chi_{e_i}}g_{i,+}^{r_i}
 \prod_{i:\,\chi_{e_i+1}\prec_w\Pi_i}g_{i,-}^{r_i}\,x\right).
\end{equation}
Here
\[
 B_w:=
 \frac{g_{1,-}^{r_1}g_{2,-}^{r_2}}
 {\displaystyle
  \prod_{\substack{1\le i\le2\\ \Pi_i\prec_w\chi_{e_i}}}g_{i,+}^{r_i}
  \prod_{\substack{1\le i\le2\\ \chi_{e_i+1}\prec_w\Pi_i}}g_{i,-}^{r_i}}
 \,M(\tau).
\]

\begin{lemma}\label{lem:two-cusp-shuffle-orders}
Fix $C\subseteq X$.
In case \textup{(i)} of \eqref{eq:two-cuspidal-cases}, set
\[
 \alpha_C=\prod_{i:\,e_i\in C}g_{i,+}^{r_i}.
\]
In case \textup{(ii)}, write $r=r_1=r_2$. Here $e_2=e_1+1$.
The following table gives $\alpha_C$, together with an order used in the
proof. The labels $1,2,3,4,5$ stand for
$\chi_{e_1},\Pi_1,\chi_{e_2},\Pi_2,\chi_{e_2+1}$, respectively.
The output order records only the relative order of these five factors.
\[
\begin{array}{ccc|ccc}
 C&\alpha_C&\text{output order}&C&\alpha_C&\text{output order}\\ \hline
 \varnothing&1&(1,2,3,4,5)&\{z\}&g_{2,+}^r&(1,4,2,3,5)\\
 \{e_1\}&g_{1,+}^r&(2,3,1,4,5)&\{e_1,z\}&g_{1,-}^r&(3,4,1,2,5)\\
 \{e_2\}&g_{2,-}^r&(1,2,5,3,4)&\{e_2,z\}&g_{2,+}^r&(1,4,2,5,3)\\
 \{e_1,e_2\}&g_{1,+}^r g_{2,-}^r&(2,5,3,1,4)&\{e_1,e_2,z\}&g_{1,+}^r g_{2,+}^r&(4,2,5,3,1)
\end{array}
\]

\textup{(1)} If $x=[\eta]_\varnothing\in D_C$ and $\alpha_Cx=0$, then
$(T_{Q',w})_\D\eta=0$ for every $w$.

\textup{(2)} For each $\kappa$ in
Lemma~\ref{lem:two-cusp-complementary-cut-detection} there is a shuffle
$\tau$ such that $w=\tau\kappa$ satisfies, for every $\eta$,
\begin{equation}\label{eq:two-cusp-detecting-shuffle}
 (T_{Q',w})_\D\eta=(B_w)_\D\rho_\kappa(\alpha_C[\eta]_\varnothing),
 \qquad B_w^{-1}\ \text{holomorphic at }t^0 .
\end{equation}
\end{lemma}

\begin{proof}
\textup{(1)} Let $x\in D_C$ and suppose $(T_{Q',w})_\D\eta\ne0$.
For each pair $\Pi_i\prec_w\chi_j$, the corresponding factor of $B_w$
in \eqref{eq:two-cusp-shuffle-jets} acts by the normalized intertwining operator
\[
 N(t_j,s_i)^{-1}:\chi_j\times\Pi_i\longrightarrow\Pi_i\times\chi_j
\]
of Lemma~\ref{lem:character-cuspidal-intertwiner} with $\pi=\pi_i$,
and by the identity on the other inducing factors.
The remaining factors are a unit and, if $\tau$ exchanges $\Pi_1$ and $\Pi_2$ in case \textup{(i)},
their standard intertwiner. These are holomorphic by the proof of
Lemma~\ref{lem:two-cusp-coefficient-map}.
For $e_i\in C$, the cut formula \eqref{eq:two-cusp-character-cut-on-jets}
gives $P_{e_i}^\varnothing(fx)=fP_{e_i}^\varnothing x=0$
for every scalar germ $f$. Since $\tau$ preserves the character order,
\eqref{eq:character-cut-common-kernel} and Proposition~\ref{lem:character-cut-modules} give
\[
 \chi_{e_i+1}\prec_w\chi_{e_i},\qquad
 (g_{i,+}+g_{i,-})x=u_{e_i}x=0 .
\]
Thus $\Pi_i\prec_w\chi_{e_i}$ or $\chi_{e_i+1}\prec_w\Pi_i$.
For each such $i$, the product in \eqref{eq:two-cusp-shuffle-jets}
contains $g_{i,+}^{r_i}$ or $g_{i,-}^{r_i}$, and
$g_{i,-}^{r_i}x=(-1)^{r_i}g_{i,+}^{r_i}x$.
Choose one such factor for each $i$ with $e_i\in C$, also when
$e_1=e_2$. Since the scalar factors commute, we may replace each chosen
$g_{i,-}^{r_i}$ on $x$ by $(-1)^{r_i}g_{i,+}^{r_i}$.
Thus the product in \eqref{eq:two-cusp-shuffle-jets} sends $x$ to
\[
 h\Bigl(\prod_{i:\,e_i\in C}g_{i,+}^{r_i}\Bigr)x,
\]
where $h\in\mathcal O_{t^0}$ is the product of the unchosen factors
and the signs from these replacements. In case \textup{(i)} the above
product of $g_{i,+}^{r_i}$ over $i$ with $e_i\in C$ is $\alpha_C$.

In case \textup{(ii)} with $z\notin C$ or $C=X$, the table gives
\[
 \Bigl(\prod_{i:\,e_i\in C}g_{i,+}^{r}\Bigr)x
 =
 \begin{cases}
  \alpha_Cx,& C=\varnothing,\ \{e_1\},\ X,\\
  (-1)^r \alpha_Cx,& C=\{e_2\},\ \{e_1,e_2\}.
 \end{cases}
\]
Indeed, in the second line $e_2\in C$, so
$g_{2,+}^{r}x=(-1)^r g_{2,-}^{r}x$. Absorbing this sign in $h$
gives $h \alpha_Cx$.

In case \textup{(ii)} with $z\in C$ and $C\ne X$, suppose
$\Pi_1\prec_w\Pi_2$, so $w(z)<w(z+1)$. Since $w=\tau\kappa$,
$\kappa$ preserves $\{z,z+1\}$ and $\tau(z)<\tau(z+1)$,
$\kappa$ fixes both positions. The $\GL_4$-factor of
$\rho_\kappa$ is therefore induced by $T_{P_{2,2},1}(v)=v\,\mathrm{id}_V$
in \eqref{eq:two-cuspidal-cleared-constant-term}. On a coefficient
$[\zeta]\in\mathcal A_{2,\pi}^{K_{4,f}}$ this gives
\[
 (T_{P_{2,2},1})_\D\zeta=v\zeta=P_v([\zeta])
\]
by \eqref{eq:two-cusp-cut-map}. In \eqref{eq:two-cusp-shuffle-jets}, let
\[
 f=\prod_{\substack{1\le i\le 2\\ \Pi_i\prec_w\chi_{e_i}}}g_{i,+}^{r}
   \prod_{\substack{1\le i\le 2\\ \chi_{e_i+1}\prec_w\Pi_i}}g_{i,-}^{r}.
\]
Since $x\in D_C$ and $z\in C$,
$P_z^\varnothing(fx)=fP_z^\varnothing x=0$. Thus
$\rho_\kappa(fx)=0$, contrary to $(T_{Q',w})_\D\eta\ne0$
by \eqref{eq:two-cusp-shuffle-jets}. Hence
\[
 \Pi_2\prec_w\Pi_1,\qquad\text{so}\qquad
 \chi_{e_2}\prec_w\Pi_1\ \text{ or }\ \Pi_2\prec_w\chi_{e_2}.
\]
Since $e_2=e_1+1$, these inequalities supply $g_{1,-}^r$ or
$g_{2,+}^r$, respectively. Now
\[
 \begin{gathered}
 (g_{1,-}+g_{2,+})x=vx=0,\\
 g_{1,-}^rx=(-1)^rg_{2,+}^rx,
 \end{gathered}
 \qquad
 \alpha_C=
 \begin{cases}
 g_{1,-}^r,&C=\{e_1,z\},\\
 g_{2,+}^r,&C=\{z\},\,\{e_2,z\}.
 \end{cases}
\]
Here too the product sends $x$ to $h \alpha_Cx$. Thus $\alpha_Cx=0$
forces $(T_{Q',w})_\D\eta=0$ in every case.

\textup{(2)} In case \textup{(i)}, let $a_i-1,b_i$ be consecutive
in $\{0,m\}\cup(X_I\setminus C)$ with $a_i\le e_i+1\le b_i$.
By \eqref{eq:two-cusp-internal-orders},
\[
 \kappa^{-1}(\{a_i,\ldots,b_i\})=\{a_i,\ldots,b_i\}.
\]
In the character order $(\chi_{\kappa^{-1}(1)},\ldots,\chi_{\kappa^{-1}(m)})$,
insert $\Pi_i$ before $\chi_{\kappa^{-1}(a_i)}$, in the
order $\Pi_1,\Pi_2$ if $a_1=a_2$. Since $e_1\le e_2$ gives
$a_1\le a_2$, we have
\[
 \Pi_1\prec_w\Pi_2,\qquad \Pi_i\prec_w\chi_{e_i+1},\qquad
 \Pi_i\prec_w\chi_{e_i}\ \Longleftrightarrow\ a_i\le e_i\
 \Longleftrightarrow\ e_i\in C .
\]
The product in \eqref{eq:two-cusp-shuffle-jets} is
$\prod_{i:e_i\in C}g_{i,+}^{r_i}=\alpha_C$, also when $e_1=e_2$.

In case \textup{(ii)}, start with the full character list
$(\chi_{\kappa^{-1}(1)},\ldots,\chi_{\kappa^{-1}(m)})$.
The characters labelled $1,3,5$ occur in positions
$\kappa(e_1),\kappa(e_2),\kappa(e_2+1)$, respectively.
For $e\in X_I\setminus C$, \eqref{eq:two-cusp-internal-orders} places
all $\chi_j$ with $j\le e$ before all $\chi_j$ with $j>e$.
For $e\in X_I\cap C$, it places $\chi_{e+1}$ before $\chi_e$.
Since $e_2=e_1+1$, this gives the first row below.
\[
 \begin{array}{c|cccc}
 C\cap X_I & \varnothing & \{e_1\} & \{e_2\} & \{e_1,e_2\}\\ \hline
 \text{order of }1,3,5 &(1,3,5)&(3,1,5)&(1,5,3)&(5,3,1)\\
 \{j:j+1\prec_w j\},\ z\notin C &\varnothing&\{1\}&\{4\}&\{1,4\}\\
 \{j:j+1\prec_w j\},\ z\in C     &\{3\}&\{2\}&\{3\}&\{1,3\}
 \end{array}
\]
Start with  the character order $\Delta=(\chi_{\kappa^{-1}(1)},\ldots,\chi_{\kappa^{-1}(m)})$. 
Consider the output order for $C$ in the statement from left to right.
Insert each of $\Pi_1,\Pi_2$ in the full list immediately before the next
character of $\Delta$ labelled $1,3$ or $5$ in the output order. If there is no
such character after it, append it to $\Delta$. For example, for
$C=\{z\}$ the output order is $(1,4,2,3,5)$, so we insert
$\Pi_2,\Pi_1$ immediately before $\chi_{e_2}$.
Let $w(j)$ be the position of the $j$-th inducing factor in the resulting
list. Deleting $\Pi_1,\Pi_2$ recovers the full character list $\Delta$.

Then $\Pi_1\prec_w\Pi_2$ if $z\notin C$ and $\Pi_2\prec_w\Pi_1$ if
$z\in C$, as required by \eqref{eq:two-cusp-internal-orders}.
The indices $j=1,2,3,4$ with $j+1\prec_w j$ contribute
$g_{1,+}^r,g_{1,-}^r,g_{2,+}^r,g_{2,-}^r$, respectively, to the scalar
multiplying $x$ in \eqref{eq:two-cusp-shuffle-jets}.
For $j\notin\{e_1,e_2,e_2+1\}$, the normalizing scalar for crossing
$\chi_j$ with either $\Pi_i$ is a unit by \eqref{eq:one-cusp-exterior-scalars}.
The last two rows therefore give $\alpha_C$ in all eight cases.
In both cases $\tau=w\kappa^{-1}$ exchanges only characters with
cuspidal factors. Hence $B_w^{-1}$ is holomorphic by
Lemma~\ref{lem:character-cuspidal-intertwiner}. Equation~\eqref{eq:two-cusp-shuffle-jets} gives
\eqref{eq:two-cusp-detecting-shuffle}.
\end{proof}

\begin{lemma}\label{lem:two-cusp-jet-nilpotence}
For $P=P_{1^m,2,2}$, let $\eta$ be a full-level inducing jet in the
source of \eqref{eq:two-cusp-inducing-jet-quotient}. Let $T$ be a linear
form in $u_1,\ldots,u_{m-1}$ in case \textup{(i)}, or in
$u_1,\ldots,u_{m-1},v$ in case \textup{(ii)}.
If $\eta$ has total degree at most $q$ in these variables, then $T^{q+1}\eta=0$.
\end{lemma}

\begin{proof}
In case \textup{(i)}, choose linear coordinates $(t,y_1)$ on
$my_1+2y_2+2y_3=0$. In case \textup{(ii)}, choose a linear coordinate
$t$ on $my_1+4y_2=0$. In these coordinates,
\[
 \eta\in I_P(t^0)^{K_{m+4,f}}\otimes
 \begin{cases}
 \D_u\otimes\D_{y_1}\otimes\D_t,&\textup{(i)},\\
 \D_u\otimes\D_v\otimes\D_t,&\textup{(ii)}.
 \end{cases}
\]
Thus $\eta$ is a finite sum of terms
\[
 f\otimes\lambda\otimes\delta_d^t=
 \begin{cases}
 f\otimes\delta_\beta\otimes\delta_c^{y_1}\otimes\delta_d^t,&\textup{(i)},\\
 f\otimes\delta_\beta\otimes\delta_a^v\otimes\delta_d^t,&\textup{(ii)},
 \end{cases}
\]
where $f\in I_P(t^0)^{K_{m+4,f}}$ and $\beta\in\Z_{\ge0}^{m-1}$.
The distribution factor extracts the Taylor coefficient of $u^\beta y_1^c t^d$
in case \textup{(i)} and of $u^\beta v^a t^d$ in case \textup{(ii)}.
Thus the order hypothesis gives $|\beta|=\sum_p\beta_p\le q$ in case
\textup{(i)} and $|\beta|+a\le q$ in case \textup{(ii)}.
Multiplication by $u_p$ acts on $\delta_\beta$ by
\[
 (u_p\delta_\beta)(h)=[u^\beta](u_ph)=
 \begin{cases}
  \delta_{\beta-\mathbf e_p}(h),&\beta_p>0,\\
  0,&\beta_p=0.
 \end{cases}
\]
In case \textup{(ii)}, likewise $v\delta_a^v=\delta_{a-1}^v$, with
$\delta_{-1}^v=0$.
Since $T$ is a linear combination of the $u_p$ and, in case \textup{(ii)},
of $v$, each multiplication by $T$ lowers this total degree by one
in every surviving term. Hence $T^{q+1}\lambda=0$ for each term, so $T^{q+1}\eta=0$.
\end{proof}

\begin{lemma}\label{lem:two-cusp-shuffle-summands}
\textup{(1)} Take $X=X_I$ in case \textup{(i)} and $X=X_I\sqcup\{z\}$
in case \textup{(ii)}. The modules $\mathcal M_{Q,Y}$, $Y\subseteq X$,
and the cuts $P_j^Y$ of \eqref{eq:two-cusp-character-cut-on-jets}
and \eqref{eq:two-cusp-z-cut-on-jets} form a cut system for $(X,\mathcal M)$.

\textup{(2)} For $C\subseteq X$, put $Y=X\setminus C$ and take $\alpha_C$
as in Lemma~\ref{lem:two-cusp-shuffle-orders}. Multiplication by $\alpha_C$
is equivariant on each module in \textup{(1)} and commutes with the cuts.
For $\mathcal B=\ker(\mathcal M\to\mathcal A)$, it satisfies the kernel
conditions \eqref{eq:cut-kernel-tests} and the surjectivity conditions
\eqref{eq:cut-scalar-surjectivity} of Proposition~\ref{prop:cut-filtration}.
\end{lemma}

\begin{proof}
\textup{(1)} Fix $Y\subseteq X$ and put $D=Y\cap X_I$.
The multidegree $\alpha\in\mathbb Z_{\ge0}^X$ records the exponents
$\alpha_j=\beta_j$ for $j\in X_I$ in $\delta_\beta(g)=[u^\beta]g(u)$.
In case \textup{(ii)} it also records $\alpha_z=d$ in $\delta_d(g)=[v^d]g(v)$.
The character factor of multidegree $\alpha$ is
\[
 \bigoplus_{\substack{\beta\in\mathbb Z_{\ge0}^{m-1}\\
               \beta_{e_1}=\alpha_{e_1},\ \beta_{e_2}=\alpha_{e_2}}}
       (\mathcal M_D)_\beta .
\]
The remaining coordinates of $\beta$ are arbitrary.
For $j\in X_I\setminus Y$, the cut multiplies the character jet by $u_j$.
Thus for $h\in V_m$ it sends
\[
 [h\otimes\delta_\beta]_D\longmapsto
 \begin{cases}
  0,&\alpha_j=0,\\
  [h\otimes\delta_{\beta-\mathbf e_j}]_{D\cup\{j\}},&\alpha_j>0 .
 \end{cases}
\]
For $\alpha_j>0$, the target indices $\gamma$ satisfy
$\gamma_j=\alpha_j-1$ and $\gamma_i=\alpha_i$ for $i\in X_I\setminus\{j\}$.
By Lemma~\ref{lem:character-cut-system}, the inverse on the $\gamma$-summand
sends $[h\otimes\delta_\gamma]_{D\cup\{j\}}$ to $[h\otimes\delta_{\gamma+\mathbf e_j}]_D$.

In case \textup{(ii)} with $z\notin Y$, the cut on the $\GL_4$-block
is given by \eqref{eq:two-cusp-cut-map},
\[
 P_v([f\otimes\delta_d])=
 \begin{cases}
  0,&d=0,\\
  f\otimes\delta_{d-1},&d>0,
 \end{cases}
 \qquad f\in V,\quad d=\alpha_z.
\]
The only relations before this cut are $(\ker R_\pi)\otimes\C\delta_0$
by Lemma~\ref{lem:two-cuspidal-coefficient-kernel}, so the second line is an
isomorphism for every $d>0$.

On representatives the cuts multiply by $u_j$ or $v$, so they commute
and are equivariant by \eqref{eq:g-action-jets}. Since $M_Q(\R)\cap K_\infty$
fixes the distribution factors, the grading passes to compact sections,
with the full factor $\D_y$ in each multidegree.

\textup{(2)} Multiplication by the scalar germ $\alpha_C$ on jets, $(\alpha_C\delta)(h)=\delta(\alpha_Ch)$,
descends to the quotient classes as in \eqref{eq:two-cusp-shuffle-jets}
and is equivariant by \eqref{eq:g-action-jets}. For $A\subseteq X$,
\eqref{eq:two-cusp-character-cut-on-jets} and \eqref{eq:two-cusp-z-cut-on-jets} give
\[
\begin{aligned}
 P_j^A(\alpha_C[\eta]_A)&=[u_j\alpha_C\eta]_{A\cup\{j\}}
 =[\alpha_Cu_j\eta]_{A\cup\{j\}}=\alpha_CP_j^A[\eta]_A
 &&(j\in X_I\setminus A),\\
 P_z^A(\alpha_C[\eta]_A)&=[v \alpha_C\eta]_{A\cup\{z\}}
 =[\alpha_C v\eta]_{A\cup\{z\}}=\alpha_CP_z^A[\eta]_A
 &&(\text{case \textup{(ii)}},\ z\notin A).
\end{aligned}
\]
Composing the cuts gives $P_Y(\alpha_Cx)=\alpha_CP_Yx$, so $\alpha_C$ preserves $\ker P_Y$.

Recall that $t_j$ and $s_i$ are the exponents in
$\chi_j=|\cdot|^{t_j}$ for $1\le j\le m$ and
$\Pi_i=\pi_i|\det|^{s_i}$ for $i=1,2$. Their values at the base point
in \eqref{eq:two-cuspidal-support} are
\[
 t_j^0=\max I+1-j-\bar c,\qquad s_i^0=c_i-\bar c .
\]
Denote $\Delta t_j=t_j-t_j^0,\ \Delta s_i=s_i-s_i^0$.
For each block $b$ of $Q$, $y_b$ is the change in its mean exponent.
Thus $y_1=m^{-1}\sum_{j=1}^m\Delta t_j$.
In case \textup{(i)}, $y_{i+1}=\Delta s_i$, and in case \textup{(ii)},
$y_2=(\Delta s_1+\Delta s_2)/2$.
For each linear factor $l=g_{i,+}$ or $g_{i,-}$ of $\alpha_C$, let $b$
index the block containing $\Pi_i$. Then $l=t+T$, where
\[
 \begin{array}{c|c|l}
 l&t&T\\ \hline
 g_{i,+}&y_1-y_b&(\Delta t_{e_i}-y_1)-(\Delta s_i-y_b)\\
 g_{i,-}&y_b-y_1&(\Delta s_i-y_b)-(\Delta t_{e_i+1}-y_1).
 \end{array}
\]
Here $\Delta t_j-y_1=m^{-1}\sum_p(\Delta t_j-\Delta t_p)$ is linear
in the $u_p$, and $\Delta s_i-y_b$ is $0$ in case \textup{(i)}
and $(-1)^{i+1}v/2$ in case \textup{(ii)}.
In case \textup{(i)}, $b=i+1$. Use $(t,y_1)$
as coordinates on $my_1+2y_2+2y_3=0$, with $t=y_1-y_b$ for
$g_{i,+}$ and $t=y_b-y_1$ for $g_{i,-}$.
In case \textup{(ii)}, $b=2$.
Since $my_1+4y_2=0$, $y_1-y_2=(m+4)y_1/4$, so $t$ alone is a
coordinate on this line.

As in the proof of Lemma~\ref{lem:two-cusp-jet-nilpotence}, each full-level inducing jet
is a sum of terms $f\otimes\lambda\otimes\delta_d^t$ with
$T^{q+1}\lambda=0$, where $q$ is the total degree of $\lambda$ in the $u_p$
and, in case \textup{(ii)}, in $v$ as well.
For $A\subseteq X$, put $D=A\cap X_I$. On compact section values,
the character factor is taken modulo the subspace
\[
 \mathcal K^D=
 \sum_{\substack{1\le j<m\\j\notin D}}
 (\ker R_j)\otimes\C\operatorname{ev}_{u_j}\otimes\D_{(u_p)_{p\ne j}}
\]
of \eqref{eqn: def-k^D}. In case \textup{(ii)} with $z\notin A$, the
$\GL_4$-factor is also taken modulo
$\ker\Xi_\pi=(\ker R_\pi)\otimes\C\operatorname{ev}_v$ by
Lemma~\ref{lem:two-cuspidal-coefficient-kernel}. Each of these subspaces
is tensored with the full distribution space $\D_t$.
For a distribution $\mu$ in the remaining variables other than $t$,
\[
 (-T)^k(\operatorname{ev}_{u_j}\otimes\mu)
 =\operatorname{ev}_{u_j}\otimes(-T|_{u_j=0})^k\mu ,
\]
and the same identity holds with $v$ in place of $u_j$.
Thus $(-T)^k$ preserves the evaluation factor and fixes the section
coefficient in $\ker R_j$ or $\ker R_\pi$. Replacing $\delta_d^t$
by $\delta_{d+k+1}^t$ also preserves each relation subspace, since
every $t$-order is allowed. Hence the right inverse of
Lemma~\ref{lem:finite-translation} descends termwise to $\mathcal M_{Q,A}$, with
\[
 S_{-T}[f\otimes\lambda\otimes\delta_d^t]_A
 =\sum_{k=0}^q[f\otimes(-T)^k\lambda\otimes\delta_{d+k+1}^t]_A,
 \qquad lS_{-T}=1 .
\]
Since $u_j,v$ commute with $T$ and the cut maps act by the identity on $\D_t$,
$P_j^A S_{-T}=S_{-T}P_j^A$ termwise for $j\in X\setminus A$.
For $W=D_C$ or $W=\ker q_C=D_C\cap\ker P_Y$,
we obtain $W=lS_{-T}W\subseteq lW\subseteq W$. Taking the product gives $\alpha_CW=W$.

Let $x=[\eta]_\varnothing$. By Lemma~\ref{lem:two-cusp-shuffle-orders}\textup{(1)}
and \eqref{eq:cuspidal-coefficient-kernel},
\[
 x\in D_C,\quad \alpha_Cx=0
 \quad\Longrightarrow\quad (T_{Q',w})_\D\eta=0\ \ (\forall w)
 \quad\Longrightarrow\quad x\in\mathcal B .
\]
Now take $x=\sum_{|C'|=|C|}x_{C'}\in\mathcal B$, with $x_{C'}\in D_{C'}$.
For $C'\ne C$, a cut in $C'\setminus C\subseteq Y$ kills $x_{C'}$.
Hence $P_Y(\alpha_Cx)=\alpha_Cq_Cx_C\in E_C$. For each $\kappa$ in
Lemma~\ref{lem:two-cusp-complementary-cut-detection}, choose $w$ by
Lemma~\ref{lem:two-cusp-shuffle-orders}\textup{(2)}. Then
\[
 \rho_\kappa(\alpha_Cx)=(B_w^{-1})_\D(T_{Q',w})_\D\eta=0 .
\]
Lemma~\ref{lem:two-cusp-complementary-cut-detection} gives $P_Y(\alpha_Cx)=0$.
\end{proof}

\begin{lemma}\label{lem:two-cusp-induced-graded-pieces}
For $C\subseteq X$, put $C_I=C\cap X_I$ and $Y_I=X_I\setminus C_I$,
where $X_I=\{e_1,e_2\}$. With the block means $a_b^0$ in
\eqref{eq:two-cusp-cut-modules}, the module $E_C=P_{X\setminus C}(D_C)$ is given by
\[
 E_C\simeq
 \begin{cases}
 \bigl((D_{C_I}^{Y_I}|\det|^{a_1^0})\times
       (\pi_{1,\infty}|\det|^{a_2^0})\times
       (\pi_{2,\infty}|\det|^{a_3^0})\bigr)\otimes\D_y,
       &Q=P_{m,2,2},\\
 \bigl((D_{C_I}^{Y_I}|\det|^{a_1^0})\times
       ((\ker P_v)|\det|^{a_2^0})\bigr)\otimes\D_y,
       &Q=P_{m,4},\ z\in C,\\
 \bigl((D_{C_I}^{Y_I}|\det|^{a_1^0})\times
       ((\operatorname{im}P_v)|\det|^{a_2^0})\bigr)\otimes\D_y,
       &Q=P_{m,4},\ z\notin C.
 \end{cases}
\]
The action on jets is obtained by varying the block exponents as $a_b^0+y_b$
in \eqref{eq:g-action-jets}.
\end{lemma}

\begin{proof}
For $Y\subseteq X$, put $D=Y\cap X_I$. Consider module $U_Y$
on $M_Q(\R)$ 
\[
 U_Y=
 \begin{cases}
 \bigl((\mathcal M_D|\det|^{a_1^0})\boxtimes
       (\pi_{1,\infty}|\det|^{a_2^0})\boxtimes
       (\pi_{2,\infty}|\det|^{a_3^0})\bigr)\otimes\D_y,
       &Q=P_{m,2,2},\\
 \bigl((\mathcal M_D|\det|^{a_1^0})\boxtimes
       (\mathcal A_{2,\pi}^{K_{4,f}}|\det|^{a_2^0})\bigr)\otimes\D_y,
       &Q=P_{m,4},\ z\notin Y,\\
 \bigl((\mathcal M_D|\det|^{a_1^0})\boxtimes
       ((V\otimes\D_v)|\det|^{a_2^0})\bigr)\otimes\D_y,
       &Q=P_{m,4},\ z\in Y .
 \end{cases}
\]
Here $\D_y$ consists of finite-order distributions at the origin of
$my_1+2y_2+2y_3=0$ or $my_1+4y_2=0$, respectively. 
By definition of  $\mathcal M_{Q,Y}$ (cf. \eqref{eq:two-cusp-cut-modules-m22} and
\eqref{eq:two-cusp-cut-modules-m4}) , we have the following identification in the compact
realization \eqref{eq:compact-picture}, 
\[
 \mathcal M_{Q,Y}\xrightarrow{\ \sim\ }\operatorname{Ind}_Q^{\GL_{m+4}(\R)}U_Y,
 \qquad f\otimes\delta\longmapsto
       \bigl(k\longmapsto f(k)\otimes\delta\bigr)
 \qquad(k\in K_\infty,\ \delta\in\D_y).
\]
Recall the Levi map $p_j^Y:U_Y\to U_{Y\cup\{j\}}$ of
\eqref{eq:two-cusp-cut-on-levi-values}, for $j\in X\setminus Y$.
If $j\in X_I$, its map on the $\GL_m$-factor is
$P_j^D([\eta]_D)=[u_j\eta]_{D\cup\{j\}}$ by \eqref{eq:character-cut-map}.
If $j=z$ in case \textup{(ii)}, its map on the $\GL_4$-factor is
$P_v([\zeta])=v\zeta$ by \eqref{eq:two-cusp-cut-map}.
In each case it is the identity on the other tensor factors, including $\D_y$.
Thus $P_j^Y=\operatorname{Ind}_Q(p_j^Y)$. Let
\[
 L=\ker\!\left(U_\varnothing
       \xrightarrow{(p_j^\varnothing)_{j\in C}}
       \bigoplus_{j\in C}U_{\{j\}}\right).
\]
Inducing this map gives $f\mapsto(P_j^\varnothing f)_{j\in C}$,
whose kernel is $D_C$ by \eqref{eqn: def-D-C}. Thus exactness gives
$D_C=\operatorname{Ind}_Q L$.
Let $p_{X\setminus C}:U_\varnothing\to U_{X\setminus C}$ be the composite
of the remaining cuts on the Levi. Its restriction to $L$ factors as
\[
 L\twoheadrightarrow p_{X\setminus C}(L)\hookrightarrow U_{X\setminus C}.
\]
Induction preserves the surjection and the injection, so the image of $P_{X\setminus C}|_{D_C}$
is $E_C=\operatorname{Ind}_Q\bigl(p_{X\setminus C}(L)\bigr)$.
On the character block, $P_{Y_I}(D_{C_I}^\varnothing)=D_{C_I}^{Y_I}$ by
the exact sequence of Proposition~\ref{lem:character-cut-modules}.
In case \textup{(ii)}, the cut on the $\GL_4$-block therefore gives
$\ker P_v$ if $z\in C$ and $\operatorname{im}P_v$ if $z\notin C$.
The cuts leave $\D_y$ unchanged, giving the stated formulas.
\end{proof}

\begin{proof}[Proof of Theorem~\ref{prop:two-cuspidal-bound}]
We may take $F$ irreducible. Denote
$H^q(U)=H^q(\mathfrak{sl}_{m+4},\O_{m+4};U\otimes F)$.

\emph{The filtration of $\mathcal A$.}
Lemma~\ref{lem:two-cusp-shuffle-summands}\textup{(1)} gives the cut system
$(\mathcal M_{Q,Y},\{P_j^Y\}_{j\in X\setminus Y})_{Y\subseteq X}$
for $(X,\mathcal M)$. Part~\textup{(2)} verifies \eqref{eq:cut-kernel-tests}
and \eqref{eq:cut-scalar-surjectivity} for $\mathcal B=\ker(\mathcal M\to\mathcal A)$.
Recall that $F_s\mathcal A$ is the image of
$F_s\mathcal M=\sum_{|C|=s}D_C$ under the quotient map
$\mathcal M\twoheadrightarrow\mathcal A$ in \eqref{eq:two-cuspidal-kernel-exact},
so $F_s\mathcal A=(F_s\mathcal M+\mathcal B)/\mathcal B$.
Proposition~\ref{prop:cut-filtration} therefore gives
\[
 F_s\mathcal A/F_{s+1}\mathcal A\simeq\bigoplus_{|C|=s}E_C .
\]
The long exact sequences of this finite filtration give
\[
 H^q(\mathcal A)\ne0\quad\Longrightarrow\quad H^q(E_C)\ne0
 \quad\text{for some }C\subseteq X .
\]

\emph{The degrees on the Levi blocks.}
Fix $C\subseteq X$ with $H^q(E_C)\ne0$. By
Lemma~\ref{lem:two-cusp-induced-graded-pieces}, $E_C$ is induced from $Q$
with character factor $D_{C_I}^{Y_I}$. The relative split centre contributes
only in degree zero by Lemma~\ref{lem:split-centre-cohomology} applied to $\D_y$.

By Shapiro's lemma (Lemma~\ref{lem:parabolic-cohomology}), the Kostant
decomposition \eqref{eq:kostant-decomposition} and the Levi K\"unneth
formula \eqref{eq:levi-kunneth-centre}, $H^q(E_C)\ne0$ requires a
nonzero tensor product of block cohomology groups for some $v_Q\in W^Q$.
Consider the decomposition $V_{\mu_{v_Q}}^{M_Q}=\boxtimes_b E_{b,v_Q}$
over the Levi blocks, and let $q_b$ be the degree of the cohomology group
on block $b$ in the chosen nonzero tensor product. Then
\[
 q=\ell_W(v_Q)+\sum_b q_b,\qquad
 H^{q_1}(\mathfrak{sl}_m,\O_m;D_{C_I}^{Y_I}\otimes E_{1,v_Q})\ne0 .
\]
The character block has parameter $I-\bar c$ by \eqref{eq:two-cuspidal-support}
and the jet filtration in the proof of Lemma~\ref{lem:finite-derivative-vanishing}.
Since $I$ is consecutive, $E_{1,v_Q}$ is a power of the determinant by
Lemma~\ref{lem:consecutive-kostant-character}. Its restriction to $\O_m$
is $1$ or $\det$. Hence \eqref{eq:character-cut-degree} gives
\[
 q_1\le d_m-2m+2+|Y_I|+2|C_I|\le d_m-2m+6,
\]
where $ |Y_I|+|C_I|=|X_I|=|\{e_1,e_2\}|\le2 .$

\emph{The bound for $q$.}
Recall that the shifted weight pair of $\pi_i$ is
\[
 B_i=\{b_i,b_i+w_i\}=\{c_i-w_i/2,c_i+w_i/2\}\qquad(i=1,2).
\]
The ordered block parameters, up to the common shift $-\bar c$, are
$(I,B_1,B_2)$ in case \textup{(i)} and $(I,B_1\sqcup B_2)$ in case
\textup{(ii)}, with each block in increasing order.
For each $x\in I$, $b_i<x<b_i+w_i$ gives one inversion with $B_i$.
The inversion count from \eqref{eq:ordered-kostant-parameter-comparison} is therefore
\[
 \ell_W(v_Q)=
 \begin{cases}
  2m+\#\{(b,b')\in B_1\times B_2:b>b'\}\le 2m+3,&\textup{(i)},\\
  2m,&\textup{(ii)} .
 \end{cases}
\]
Indeed $b_1<\min I<b_2+w_2$, so at most three of the four pairs in
$B_1\times B_2$ are inversions.

In case \textup{(i)}, the cuspidal factors are $\pi_{1,\infty}$ and
$\pi_{2,\infty}$ by \eqref{eq:two-cuspidal-cases}, so
Proposition~\ref{prop:speh-cohomology-bound} gives
$q_2,q_3\le d_2-b(2,1)=1$. In case \textup{(ii)} with $z\in C$,
the $\GL_4$-factor is $\ker P_v\simeq J(2,\pi)^{K_{4,f}}$ by
Lemma~\ref{lem:two-cuspidal-coefficient-kernel} and \eqref{eq:two-cusp-cut-map},
so the same proposition gives $q_2\le d_4-b(2,2)=9-3=6$.
In case \textup{(ii)} with $z\notin C$, apply Shapiro and Kostant to $P_{2,2}$ on
$\operatorname{im}P_v=V\otimes\D_v$. Write $b=b_2$, $w=w_2=w_1$.
Since $c_1-c_2=1$, the ordered blocks are
$(b+1,b+1+w)$ and $(b,b+w)$. Their four entries have ranks
$(2,4,1,3)$, giving three inversions. By Lemma~\ref{lem:split-centre-cohomology},
the split-centre cohomology with coefficients in $\D_v$ vanishes in positive degrees, so
$q_2\le 3+1+1+0=5$.

Combining the block degrees with $\ell_W(v_Q)$ gives
\begin{equation}\label{eq:two-cuspidal-cut-degrees}
 H^q(\mathcal A)\ne0\quad\Longrightarrow\quad
 q\le\begin{cases}
 (2m+3)+(d_m-2m+6)+1+1=d_m+11,&\textup{(i)},\\
 2m+(d_m-2m+6)+\max\{6,5\}=d_m+12,&\textup{(ii)} .
 \end{cases}
\end{equation}
Subtracting from $d_{m+4}=d_m+4m+10$ gives
\[
 d_{m+4}-q\ge
 \begin{cases}4m-1,&\textup{(i)},\\4m-2,&\textup{(ii)}.\end{cases}
\]
Both bounds imply $H^q(\mathcal A)=0$ for $d_{m+4}-q<4m-4$.
\end{proof}

\subsection{An example with two cuspidal factors in rank fourteen}
\label{subsec:two-cuspidal-example}

We illustrate the proof of Theorem~\ref{prop:two-cuspidal-bound} with the two
distinct level-one cuspidal representations $\pi_1=\Delta_{17}$ and
$\pi_2=\Delta_{15}$ recalled in Lemma~\ref{lem:w11}. Put
\[
 B_1=\{-1,16\},\qquad B_2=\{0,15\},\qquad I=\{3,\ldots,12\}.
\]
Here $m=10$, the sets $I,B_1,B_2$ are pairwise disjoint, and
$3<c_1=c_2=\bar c=15/2<12$. Thus the cuspidal support
$\mathfrak c$ in \eqref{eq:two-cuspidal-support} is represented by
\[
 \left(\boxtimes_{x=3}^{12}|\cdot|^{x-15/2}\right)
             \boxtimes\pi_1\boxtimes\pi_2 .
\]
Let $\chi:\mathfrak Z_{\GL_{14}}^1\to\C$ have parameter
$(I\sqcup B_1\sqcup B_2)-15/2$. Then \eqref{eq:two-cuspidal-automorphic-module} gives
\[
 \mathcal A=\left\{\phi\in\mathcal A(\GL_{14})_{\mathfrak c}^{K_{14,f}}:
         (\ker\chi)^N\phi=0\text{ for some }N\ge1\right\}.
\]
Take $F$ to be the irreducible algebraic representation of $\GL_{14}$
of highest weight
$\Lambda_F=(a+2,a+2,\underbrace{a,\ldots,a}_{10},a-2,a-2)$, $a\in\Z$.
Since $\rho_{14,j}=15/2-j$, the trace-zero parameter of $F^\vee$ is
\[
 \{a-\Lambda_{F,j}-\rho_{14,j}:1\le j\le14\}
 =\{-1,0,3,\ldots,12,15,16\}-15/2=(I\sqcup B_1\sqcup B_2)-\bar c .
\]
Thus $F^\vee$ has infinitesimal character $\chi$ on $\mathfrak{sl}_{14}$.

\emph{The filtration of $\mathcal A$.}
Since $\pi_1\not\simeq\pi_2$, we are in case \textup{(i)} of
\eqref{eq:two-cuspidal-cases} with $Q=P_{10,2,2}\subset\GL_{14}$.
Order the characters by decreasing weight, $12,11,\ldots,3$. Then
\[
 e_1=e_2=\#\{x\in I:x>15/2\}=5,\qquad X=X_I=\{5\}.
\]
Thus both cuspidal factors give the cut between $\chi_5$ and $\chi_6$,
of weights $8$ and $7$. The ten character exponents at the base point
are $9/2,7/2,\ldots,-9/2$, with mean zero, and $c_i-\bar c=0$ for $i=1,2$.
For inducing exponents $(t_1,\ldots,t_{10},s_1,s_2)$, the three block
means are $y_1=\frac1{10}\sum_{j=1}^{10}t_j$, $y_2=s_1$ and $y_3=s_2$.
Using normalized induction from $Q=P_{10,2,2}$, we have
\[
\begin{aligned}
 Y&=\{(y_1,y_2,y_3)\in\C^3:10y_1+2y_2+2y_3=0\},\\
 \mathcal I_Q(y)&=(\mathcal A_{10,1}|\det|^{y_1})
           \times(\pi_1|\det|^{y_2})\times(\pi_2|\det|^{y_3}),\\
 \mathcal M&=\mathcal I_Q(0)^{K_{14,f}}\otimes\D_y .
\end{aligned}
\]
Here $\mathcal A_{10,1}$ is the $\GL_{10}$-module of
Subsection~\ref{subsec:speh-support-modules} for the ten character exponents above,
and $\D_y$ consists of the finite-order distributions at $0\in Y$.
For $P=P_{1^{10},2,2}$, let $E_P(t)$ be the Eisenstein family induced
from $(|\cdot|^{t_1},\ldots,|\cdot|^{t_{10}},
\pi_1|\det|^{s_1},\pi_2|\det|^{s_2})$.
The central vanishing orders $r_i=\operatorname{ord}_{s=1/2}\Lambda(s,\pi_i)$
are $r_1=1$ and $r_2=0$.

With $u_j=t_j-t_{j+1}-1$, the family in Lemma~\ref{lem:two-cusp-coefficient-map} is
\[
 H(u,y)=(s_1-t_6-\tfrac12)\Bigl(\prod_{j=1}^{9}u_j\Bigr)E_P(t).
\]
At $t^0=(9/2,7/2,\ldots,-9/2,0,0)$, we have $(u,y)=(0,0)$.
For a full-level inducing jet of \eqref{eq:full-level-inducing-jet},
\eqref{eq:two-cusp-inducing-jet-quotient} and Lemma~\ref{lem:two-cusp-coefficient-map}
give
\[
\begin{array}{ccccc}
 I_P(t^0)^{K_{14,f}}\otimes\D_u\otimes\D_y &\twoheadrightarrow&
 \mathcal M &\twoheadrightarrow& \mathcal A\\
 \eta &\longmapsto& [\eta]_\varnothing &\longmapsto& \Xi_H(\eta).
\end{array}
\]
Here $\D_u$ consists of the finite-order distributions at $0\in\C^9$.
For the two subsets of $X=\{5\}$, the definition
$E_C=P_{X\setminus C}(D_C)$ gives
\[
\begin{aligned}
 C=\varnothing:\quad &D_\varnothing=\mathcal M, &
 E_\varnothing&=P_5^\varnothing(\mathcal M),\\
 C=\{5\}:\quad &D_{\{5\}}=\ker P_5^\varnothing, &
 E_{\{5\}}&=P_\varnothing(D_{\{5\}})=\ker P_5^\varnothing .
\end{aligned}
\]
Denote $\mathcal B=\ker(\mathcal M\twoheadrightarrow\mathcal A)$.
Since $e_1=e_2=5$, Lemma~\ref{lem:two-cusp-shuffle-orders} gives
\[
 \alpha_\varnothing=1,\qquad
 \alpha_{\{5\}}=g_{1,+}^{r_1}g_{2,+}^{r_2}=t_5-s_1-\tfrac12 .
\]
Lemma~\ref{lem:two-cusp-shuffle-summands}\textup{(2)} verifies
\eqref{eq:cut-kernel-tests} and \eqref{eq:cut-scalar-surjectivity} for $\mathcal B$.
For $C=\varnothing$, the second condition in \eqref{eq:cut-kernel-tests} gives
$P_5^\varnothing\mathcal B=0$, hence $\mathcal B\subseteq D_{\{5\}}$.
For $C=\{5\}=X$, we have $X\setminus C=\varnothing$, so
$P_{X\setminus C}=P_\varnothing=\mathrm{id}_{\mathcal M}$. The two conditions
in \eqref{eq:cut-kernel-tests} now give
$\mathcal B=\ker(\alpha_{\{5\}}|_{D_{\{5\}}})$, and
\eqref{eq:cut-scalar-surjectivity} gives $\alpha_{\{5\}}D_{\{5\}}=D_{\{5\}}$.
Thus $F_0\mathcal A=\mathcal A$, $F_1\mathcal A=D_{\{5\}}/\mathcal B$
and $F_2\mathcal A=0$. Writing $\bar x$ for the image of $x$ in $\mathcal A$,
the two isomorphisms of Proposition~\ref{prop:cut-filtration} are
\[
 \begin{aligned}
 \mathcal A/F_1\mathcal A&\xrightarrow{\ \sim\ }E_\varnothing,
 &\bar x+F_1\mathcal A&\longmapsto P_5^\varnothing x &&(x\in\mathcal M),\\
 F_1\mathcal A&\xrightarrow{\ \sim\ }E_{\{5\}},
 &\bar x&\longmapsto (t_5-s_1-\tfrac12)x &&(x\in D_{\{5\}}).
 \end{aligned}
\]

For $U\in\{\mathcal A,E_\varnothing,E_{\{5\}}\}$, write
$H^q(U)=H^q(\mathfrak{sl}_{14},\O_{14};U\otimes F)$.
The filtration and its long exact cohomology sequence therefore give
\[
 \begin{gathered}
 0\longrightarrow E_{\{5\}}\longrightarrow\mathcal A
   \longrightarrow E_\varnothing\longrightarrow0,\\
 H^q(\mathcal A)\ne0
 \quad\Longrightarrow\quad H^q(E_C)\ne0
 \quad\text{for some }C\subseteq\{5\}.
 \end{gathered}
\]

\emph{The inducing modules for $E_\varnothing$ and $E_{\{5\}}$.}
At $t^0$, the mean on the $\GL_{10}$-block is
\[
 y_1=\frac1{10}\sum_{j=1}^{10}t_j^0
    =\frac1{10}\sum_{j=1}^{10}\left(\frac{11}{2}-j\right)=0 .
\]
On the two $\GL_2$-blocks, $y_2=s_1^0=c_1-\bar c=0$ and
$y_3=s_2^0=c_2-\bar c=0$, since $c_1=c_2=\bar c=15/2$. By
Lemma~\ref{lem:two-cusp-induced-graded-pieces}, normalized induction from
$P_{10,2,2}$ gives
\[
 \begin{aligned}
 E_\varnothing&\simeq
   (\mathcal M_{\{5\}}\times\pi_{1,\infty}\times\pi_{2,\infty})\otimes\D_y,\\
 E_{\{5\}}&\simeq
   (\mathcal R_{10,5}\times\pi_{1,\infty}\times\pi_{2,\infty})\otimes\D_y .
 \end{aligned}
\]

By the proof of Proposition~\ref{prop:character-block-cut-degree}, normalized
induction from $P_{5,5}$ to $\GL_{10}$ gives
\[
 \begin{aligned}
 \mathcal I_{\{5\}}(z)
   &=\bigl(\mathcal A_{5,1}|\det|^{5/2+z}\bigr)
      \times\bigl(\mathcal A_{5,1}|\det|^{-5/2-z}\bigr),\\
 \mathcal M_{\{5\}}
   &\simeq \mathcal I_{\{5\}}(0)^{K_{10,f}}\otimes\D_z .
 \end{aligned}
\]
Here $z$ is the change in the first block mean and $\D_z$ consists of
finite-order distributions at $z=0$.
 Proposition~\ref{lem:character-cut-modules} gives
\[
 \mathcal R_{10,5}=\ker\bigl(P_5^\varnothing:
 \mathcal A_{10,1}^{K_{10,f}}\longrightarrow\mathcal M_{\{5\}}\bigr).
\]

\emph{The degrees on the Levi blocks of $Q$.}
Up to the common shift $-15/2$, the three blocks have ordered weights
$(3,\ldots,12;\,-1,16;\,0,15)$. Ranking them in $I\sqcup B_1\sqcup B_2$
gives the contributing $v_Q\in W^{Q}$ by
\eqref{eq:ordered-kostant-parameter-comparison} with
\[
 (v_Q(j))_{j=1}^{14}=(3,\ldots,12;\,1,14;\,2,13),\qquad
 \ell_W(v_Q)=10+10+2=22 .
\]
In fact, each $x\in I$ exceeds one element of each $B_i$, and the
inversions in $B_1\times B_2$ are $(16,0)$ and $(16,15)$.
For $1\le j\le10$, $v_Q(j)=j+2$ and the $(j+2)$-nd entry of the
highest weight of $F$ is $a$.
Thus the highest weight on the character block in \eqref{eq:kostant-decomposition}
has $j$-th entry $a+\rho_{14,j+2}-\rho_{14,j}=a-2$.

Fix $C$ with $H^q(E_C)\ne0$ and put $Y_I=\{5\}\setminus C$.
By Shapiro's lemma (Lemma~\ref{lem:parabolic-cohomology}),
\eqref{eq:levi-kunneth-centre} and Lemma~\ref{lem:split-centre-cohomology},
\[
 q=22+q_1+q_2+q_3,\qquad
 H^{q_1}(\mathfrak{sl}_{10},\O_{10};
       D_C^{Y_I}\otimes{\det}_{10}^{\,a-2})\ne0 .
\]
On $\O_{10}$, $\det_{10}^{\,a-2}=\det_{10}^{\,a}$.
Since $d_{10}=54$, \eqref{eq:character-cut-degree} gives
\[
\begin{array}{c|c|c|c}
 C & Y_I & D_C^{Y_I} & q_1\le d_{10}-20+2+|Y_I|+2|C| \\ \hline
 \varnothing & \{5\} & \mathcal M_{\{5\}} & 54-20+2+1=37\\
 \{5\} & \varnothing & \mathcal R_{10,5} & 54-20+2+2=38
\end{array}
\]
On the two $\GL_2$-factors, $\pi_{1,\infty}$ and $\pi_{2,\infty}$
give $q_2,q_3\le d_2-b(2,1)=1$ by Proposition~\ref{prop:speh-cohomology-bound}.
Consequently
\[
 q\le
 \begin{cases}
 22+37+1+1=61,& C=\varnothing,\\
 22+38+1+1=62,& C=\{5\}.
 \end{cases}
\]
Since $d_{14}=104$, this gives $d_{14}-q\ge104-62=42$ whenever
$H^q(\mathcal A)\ne0$. In particular,
\[
 H^q(\mathfrak{sl}_{14},\O_{14};\mathcal A\otimes F)=0
 \qquad(d_{14}-q<4m-4=36),
\]
as in Theorem~\ref{prop:two-cuspidal-bound}.

\bibliographystyle{unsrt}
\bibliography{Bibliography}

\iffalse

\fi

\end{document}